\documentclass[10pt]{article}
\usepackage{Handlebodies}

\title{\Large\bfseries\scshape stable moduli spaces of odd-dimensional manifold triads}
\author{\textsc{joão lobo fernandes}}
\date{}

\begin{document}


\maketitle


\begin{abstract}
    \noindent
    We establish a homotopy-theoretic description of the homology of stable moduli spaces of $(2n+1)$-dimensional manifold triads $(N, \hb N, \vb N)$ with fixed $\vb N$, whenever $n \geq 3$ and $(N, \hb N)$ is $1$-connected. Stabilization is performed by taking boundary connected sum with $S^n \times D^{n+1}$. This is an analog of earlier work of Galatius and Randal-Williams for even-dimensional manifolds with fixed boundary, and it extends a previous result by Botvinnik and Perlmutter. As a byproduct, we obtain an analog for odd-dimensional triads of Kreck's stable diffeomorphism classification of even-dimensional manifolds.
\end{abstract}

\tableofcontents

\section{Introduction and statement of results.}

Let $M$ be a compact, smooth manifold of dimension $d$ with boundary $\partial M$, and let $\Diff_{\partial}(M)$ denote the topological group of diffeomorphisms of $M$ that restrict to the identity on $\partial M$, equipped with the $C^\infty$-topology. The classifying space $\BDiff_{\partial}(M)$ classifies smooth bundles with fiber $M$, together with a trivialization of the subbundle with fiber $\partial M.$  Over the past two decades, significant progress has been made in understanding the homotopy type of $\BDiff_{\partial}(M)$. One of the major catalysts for this progress was the work of Galatius and Randal-Williams \cite{GRWStableModuli,GRWI,GRWII}, resulting in a higher-dimensional generalization of Madsen--Weiss' resolution of the \textit{Mumford conjecture} \cite{madsenweiss} to all even-dimensional manifolds. This program has found a wide and significant range of applications, such as \cite{kupersfinite,BerglundMadsen,krannich2021diffeomorphismsdiscssecondweiss,dalian,Bustamante2021FinitenessPO,evendisc,solidtori}. As the work of Galatius and Randal-Williams applies specifically to even-dimensional manifolds, many of these applications are consequently restricted to even dimensions as well. The purpose of this work is to provide an analog to their program in the setting of odd-dimensional manifold triads.   

\subsubsection*{The work of Galatius and Randal-Williams.} We believe that our result is best explained in analogy to the one of Galatius and Randal-Williams, which we now recall: Assume $d=2n$, and that $M$ is connected and has non-empty boundary. Let $W_{g,1}$ be the manifold obtained by taking a $g$-fold boundary connected sum of $W_{1,1}\coloneqq (S^n\times S^n)\backslash \text{int}(D^{2n})$. By choosing a $(2n-1)$-disc in $\partial M,$ we define the boundary connected sum $M\mathbin{\natural} W_{g,1}$ along this disc. We define the space $\BDiff_\partial(M_\infty)$ as the homotopy colimit over the stabilization maps $\BDiff_\partial(M\mathbin{\natural} W_{g,1})\to \BDiff_\partial(M\mathbin{\natural} W_{g+1,1})$, induced by boundary connected sum with $W_{1,1}.$ A parameterized version of the Pontrjagin--Thom construction yields a map
\begin{equation}\label{alpha}
    \BDiff_\partial(M_\infty)\to \left(\Omega^{\infty}\MTh_M\right)_{\Aut_{\partial M}(\theta_M)}
\end{equation}to a purely homotopy-theoretic object, which we now describe: Consider a \textit{Moore--Postnikov $n$-factorization} $\tau_M=\theta_M\circ l_M:M\to B\to \BO(2n)$ of a tangent classifier of $M$. In other words, the map $l_M$ is an $n$-connected cofibration and the map $\theta_M$ is an $n$-coconnected fibration. The spectrum $\MTh_M$ denotes the \textit{Thom spectrum} of the formal inverse of $\theta_M$. The latter admits an action by the group-like monoid $\Aut_{\partial M}(\theta_M)$ of self weak equivalences of $B$ over $\BO(2n)$ and under $\partial M.$ The target of the map \eqref{alpha} is the homotopy orbits of $\Omega^\infty\MTh_M$ by this action. The main result of \cite{GRWII} is that the map \eqref{alpha} is \textit{acyclic} onto the path component it hits, that is, it induces an isomorphism on homology groups with any local coefficient system pulled back from the target. The right-hand side of \eqref{alpha} is amenable to standard methods in homotopy theory. We refer to \cite{usersguide} for some explicit computations of the homology of the left-hand side via this result.

\subsubsection*{An odd-dimensional analog.}
We continue by explaining the setting of our work, in analogy to the even-dimensional case: Our main objects of study are compact smooth manifold triads $(\tri,\hb \tri,\vb \tri)$ of dimension $d=2n+1,$ that is, compact smooth $d$-dimensional manifolds $\tri$ with boundary, along with a decomposition $\partial \tri=\hb \tri\cup \vb \tri$ of their boundary into the union of two codimension $0$ compact submanifolds $\hb \tri$ and $\vb \tri$ such that $\partial(\hb \tri)=\hb \tri\cap \vb \tri=\partial(\vb \tri). $ We assume that $\hvb \tri\coloneqq \hb \tri\cap \vb \tri$ is non-empty and that $N$ is connected. The role of the space $\BDiff_\partial(M)$ above is now played by the classifying space $\BDiff_{\vb}(\tri)$ of the topological group of diffeomorphisms of $\tri$ which are the identity on $\vb \tri$ but are allowed to move $\hb \tri$.
\\\\
In this context, the role of $W_{g,1}$ above is now played by the triad $(V_g,W_{g,1}, D^{2n})$ where $V_g$ is the manifold obtained by $g$-fold boundary connected sum of $S^n\times D^{n+1}$ and $D^{2n}$ is a codimension $0$ disc in $\partial V_g$: The boundary of $V_g$ is the $g$-fold connected sum of $S^n\times S^n$, so we are using the fact that $W_{g,1}$ is diffeomorphic to $\sharp_g (S^n\times S^n)$ after removing an open $2n$-disc. Similarly, by picking a $2n$-dimensional half-disc in $\vb  \tri$ whose equatorial $(2n-1)$-disc lies in $\hb N$, boundary connected sum yields the triad $(\tri\mathbin{\natural} V_g,\hb \tri\mathbin{\natural} W_{g,1}, \vb \tri).$ We define the space $\BDiff_{\vb}(\tri_\infty)$ to be the homotopy colimit over the maps $\BDiff_{\vb}(\tri\mathbin{\natural} V_g)\to \BDiff_{\vb}(\tri\mathbin{\natural} V_{g+1})$ induced by boundary connected sum with $(V_1,W_{1,1},D^{2n})$. This is the analog of the left hand side in \eqref{alpha}. We proceed by defining the analog of right hand side in our context.
\\\\
We say that a factorization $(X,X')\to (Y,Y')\to (Z,Z')$ of a map of pairs $(X,X')\to (Z,Z')$ is a \textit{Moore--Postnikov $k$-factorization of pairs} if $X'\to Y'\to Z'$ is a Moore--Postnikov $k$-factorization and $X\cup_{X'} Y'\to Y\to Z$ a Moore--Postnikov $(k+1)$-factorization, both in the absolute sense. Such a factorization always exists (see \cref{mp exist and are unique} below), so we may fix a Moore--Postnikov $n$-factorization of pairs $\Theta_\tri\circ l_\tri:(\tri,\hb \tri)\to (B,B^\partial)\to (\BO(2n+1),\BO(2n))$ of a compatible pair of tangent classifiers $\tau_\tri: (\tri,\hb \tri)\to (\BO(2n+1),\BO(2n))$. We denote by $\MTH_\tri$ the cofiber of the canonical map of Thom spectra $\Sigma^{-1}\MTh_{\hb \tri}\to \MTh_\tri$ induced by the map of pairs $\Theta_N=(\theta_N,\theta_{\hb N}).$ The group-like monoid $\Aut_{\vb \tri}(\Theta_\tri)$ of self equivalences of the pair $(B,B^\partial)$  under $(\vb \tri,\hvb \tri)$ and over $(\BO(2n+1),\BO(2n))$ acts on this spectrum. A parameterized version of the Pontrjagin--Thom construction for pairs yields a map
\begin{equation}\label{alpha b}
    \BDiff_{\vb}(\tri_\infty)\to \left(\Omega^{\infty}\MTH_\tri\right)_{\Aut_{\vb \tri}(\Theta_\tri)}.
\end{equation}We now state our main result, which can be seen as an analog for odd-dimensional triads to the one of Galatius and Randal-Williams described above.

\begin{mainteo}\label{main no tang}
    Let $(\tri,\hb \tri,\vb \tri)$ be a compact smooth $(2n+1)$-dimensional manifold triad where $N$ is connected, $(\tri,\hb \tri)$ is $1$-connected, $\hvb \tri\neq \emptyset$, and $n\geq 3$. For a Moore--Postnikov $n$-factorization $\tau_\tri=\Theta_\tri\circ l_\tri$, the map
    \[\BDiff_{\vb}(\tri_\infty)\to \left(\Omega^{\infty}\MTH_\tri\right)_{\Aut_{\vb \tri}(\Theta_\tri)}\]is acyclic onto the path component it hits.
\end{mainteo}
\noindent
\begin{rmkk} We highlight the following points:
\begin{enumerate}[itemsep=4pt, label=(\textit{\roman*})]
    \item Specialized to $n\geq 4$ and to the triad $(D^{2n+1},D^{2n}_+,D^{2n}_-)$, \cref{main no tang} recovers a previous result by Botvinnik and Perlmutter \cite[Thm.\ A]{BP}, where $D^{2n}_\pm$ are the upper and lower hemispheres of $\partial D^{2n+1}=S^{2n}$ (see \cref{examples section}).  Our work can be seen as a generalization of this result to a more general class of triads (and incidentally to $n\geq 3$--although the extension of \cite[Thm.\ A]{BP} to dimension $7$ was also previously obtained by Krannich and Kupers). However, our methods are much closer in spirit to the proof of the analogous result by Galatius and Randal-Williams. In addition, we also consider a more general notion of tangential structures for triads (see \ref{where i talk about theta structures for the first time} below) than the one considered in \cite{BP}, and thus also recover the more general result Theorem A* of loc.cit.  
    \item The right-hand side of \eqref{alpha b} simplifies considerably in many cases. For example, if $\tri$ is an $h$-cobordism, that is, if the inclusions $\hb \tri\hookrightarrow \tri \hookleftarrow \vb \tri$ are equivalences, then the right-hand side is equivalent to $\Omega^\infty\Sigma^\infty_+B$ for $B$ fitting in a Moore--Postnikov $n$-factorization $\tri\to B\to \BO(2n+1)$ (see \cref{examples section}). In particular, the rational cohomology ring of any of its path components is isomorphic to the free graded commutative $\QQ$-algebra $\QQ[\tilde{\H}^*(B;\QQ)],$ so the same holds for $\BDiff_{\vb}(N_\infty).$ 
    \item\label{where i talk about theta structures for the first time} Galatius and Randal-Williams' result above is a corollary of a more general result on the moduli space $\BDiff_\partial^\theta(M)$ of $\theta$-structures on $M$: For a map $\theta:B\to \BO(2n)$, a $\theta$-structure on $M$ is a fiberwise isomorphism $TM\to \theta^*\gamma_{2n}$, where $\gamma_{2n}$ is the universal $2n$-vector bundle. In this spirit, we also prove a more general result, stated as \cref{main} in \cref{general}, that identifies the homology of stabilization of the moduli space $\BDiff_{\vb}^{\Theta}(\tri)$ of $\Theta$-structures on $\tri:$ For a map of pairs $\Theta=(\theta,\theta^\partial):(B,B^\partial)\to (\BO(2n+1),\BO(2n))$, a $\Theta$-structure on $(\tri,\hb \tri,\vb \tri)$ is a pair of compatible $\theta$- and $\theta^\partial$-structures on $\tri$ and $\hb \tri$, respectively (which includes compatible orientations, spin structures and framings, for example). \cref{main no tang} follows by applying this result to $\Theta=\id.$ See \cref{examples section} for other examples of $\Theta$-structures and their corresponding moduli spaces.
    \item Our result is not only analogous to the one of Galatius--Randal-Williams, but also compatible with it in the sense that there exists a commutative square
    \[\begin{tikzcd}
        \BDiff_{\vb}(\tri_\infty)\arrow[d]\arrow[r, "\eqref{alpha b}"] & \left(\Omega^{\infty}\MTH_\tri\right)_{\Aut_{\vb \tri}(\Theta_\tri)}\arrow[d] \\ 
        \BDiff_\partial(\hb \tri_\infty)\arrow[r, "\eqref{alpha}"] & \left(\Omega^{\infty}\MTh_{\hb \tri}\right)_{\Aut_{\hvb \tri}(\theta_{\hb \tri})}
    \end{tikzcd}\]where the left vertical map is given by restricting diffeomorphisms and the right vertical map is induced by the canonical map $\MTH_\tri\to \MTh_{\hb \tri}$ and the restriction of self-equivalences of $(B,B^\partial)$ to $B^\partial.$ 
    
\end{enumerate}\end{rmkk}

\subsubsection*{Stable moduli spaces of even-dimensional nullbordisms.}

Our main result \cref{main no tang} will be deduced from a stronger, although more technical, group completion-type result on moduli spaces of nullbordisms, which we now move to explain. A comparable deduction happens in the work of Galatius and Randal-Williams \cite{GRWStableModuli,GRWII}, so we explain this case first. One of the most useful objects exploited in loc.cit. is the \textit{cobordism category $\Cob_\theta$ of $\theta$-manifolds} for a map of spaces $\theta:B\to \BO(d).$ This is a category internal to topological spaces: Roughly speaking, the objects of $\Cob_\theta$ are $(d-1)$-dimensional closed smooth manifolds $P$ together with a vector bundle map $l_P:TP\oplus \varepsilon^1\to \theta^*\gamma_d$ which is a fiberwise isomorphism, where $\gamma_d$ is the universal $d$-dimensional vector bundle over $\BO(d)$. A morphism $W:P\leadsto Q$ in $\Cob_\theta$ is a $d$-dimensional compact smooth manifold $W$ such that $\partial W=P\sqcup Q$, \textit{i.e.} a \textit{cobordism}, together with a vector bundle map $l_W:TW\to \theta^*\gamma_d$, extending $l_P$ and $l_Q$ on its boundary (see \cite[Defn. 2.6]{GRWStableModuli} for more details). We highlight the following facts about $\Cob_{\theta}:$
\begin{enumerate}[topsep=3pt, itemsep=3pt, label=\protect\circled{\arabic*}]
    \item\label{bdiff is inside cob} Given a compact smooth $d$-dimensional manifold $M$ as above, the space of interest $  \smash{\BDiff_\partial(M)}$ is a path component of the morphism space $\Cob(\emptyset, \partial M)$, where $\Cob\coloneqq \Cob_{\id}$ for $\id:\BO(d)\to\BO(d)$ the identity map. More generally, for any map $\theta:B\to \BO(d)$, the space $\BDiff^\theta_\partial(M)$ as in \ref{where i talk about theta structures for the first time} above is a union of path components of $\Cob_\theta(\emptyset,\partial M)$.
    \item\label{cob is understood} For any map $\theta:B\to \BO(d)$, Galatius, Madsen, Tillmann, and Weiss \cite{GMTW} proved that the classifying space $\B\Cob_\theta $ is equivalent to $\Omega^\infty\MTc$, where $\MTc$ is the Thom spectrum of the formal inverse $-\theta$.  
    \item\label{stabilization intro} When $d=2n$, there exists a \textit{distinguished} morphism $H_P:P\leadsto P$ for every non-empty connected object $P\in \Cob_\theta$ such that the post-composition map $H_P\circ (-):\Cob_\theta(\emptyset, P)\to \Cob_\theta(\emptyset, P)$ of $H_P$ extends the stabilization map $\BDiff_\partial(M\natural W_{g,1})\to \BDiff_\partial(M\natural W_{g+1,1})$ defined above along the inclusion from \ref{bdiff is inside cob}, in case $P=\partial M$ and $\theta=\id.$ The underlying cobordism of $H_P$ is the connected sum $(P\times [0,1])\sharp W_1$, where $W_1\coloneqq S^n\times S^n=\partial V_1.$ Moreover, the space $\BDiff_\partial(M_\infty)$ is a path component of the space
    \[\Cob(\emptyset, P)[(H_{P})^{-1}]\coloneqq \hocolim \left(\Cob(\emptyset, P) \overset{H_{P}\circ (-)}{\rightarrow} \Cob(\emptyset, P) \to \cdots \right)\]given by "inverting" the action of this morphism, when $P=\partial M$.
\end{enumerate}In \cite{GRWII}, the authors reduce the proof that the map \eqref{alpha} is acyclic, using standard homotopy-theoretic methods together with \ref{bdiff is inside cob}--\ref{stabilization intro}, to the following statement about the spaces $\Cob_\theta(\emptyset,P)[(H_P)^{-1}]$ for a general map $\theta:B\to \BO(2n)$: The canonical map \begin{equation}\label{null grw intro}
    \Cob_{\theta,n}(\emptyset,P)[(H_P)^{-1}]\to \Omega_{[\emptyset,P]}\B\Cob_\theta
\end{equation}taking a morphism to the path represented by it, is acyclic. Here, the source is the subspace of $\Cob_{\theta}(\emptyset,P)$ of those morphisms $W:\emptyset\leadsto P$ such that $l_W:W\to B$ is $n$-connected. Most of the work in \cite{GRWII} resides in the proof of the latter statement, for which the authors use a mixture of parameterized surgery (as developed in \cite{GRWStableModuli}), a variant of the classical group completion theorem, and geometric/surgery-theoretic arguments.

\subsubsection*{Interlude: Stable diffeomorphism classification.}

In contrast to the acyclicity of the map \eqref{alpha}, the statement that \eqref{null grw intro} is acyclic already has content on $\H_0$: It recovers a classical result of Kreck \cite{Kreck1999SurgeryAD} on the classification of even-dimensional manifolds up to stable diffeomorphism, which we briefly recall. Let $P$ be a closed $(2n-1)$-dimensional manifold, we say that two compact connected $2n$-dimensional manifolds $M_0$ and $M_1$ with boundary, together with identifications of $\partial M_i$ with $P$ for $i=0,1$, are \textit{stable diffeomorphic relative to $P$} if there exists an integer $g\geq 0$ such that $M_0\sharp W_{g}$ and $ M_1\sharp W_{g}$ are diffeomorphic relative to $P$. In \cite[Thm. 2]{Kreck1999SurgeryAD}, Kreck established a necessary and sufficient condition for two such manifolds to be stably diffeomorphic, which is purely bordism-theoretic. It turns out that the acyclicity of \eqref{null grw intro} on $\H_0$ implies this result (see \cref{grw and kreck} for more details). Therefore, one can view the latter as a \textit{family version} of Kreck's result. Below we will explain that a similar perspective applies in our setting.

\subsubsection*{Stable moduli spaces of odd-dimensional triad nullbordisms.}

We move now to explaining the aforementioned stronger but more technical, group completion type result in our context. We also consider a category internal to topological spaces, namely the \textit{cobordism category $\Cobbt$ of $\Theta$-manifolds with boundary} for a map of pairs of spaces $\Theta=(\theta,\theta^\partial):(B,B^\partial)\to (\BO(d),\BO(d-1)):$ The objects of this category are $(d-1)$-dimensional compact smooth manifolds with boundary $(P,\partial P)$ together with a map of pairs of vector bundles $l_P:(TP\oplus \varepsilon^1,T\partial P\oplus\varepsilon^1)\to (\theta^*\gamma_d,\theta^{\partial}\gamma_{d-1})$ (see \ref{where i talk about theta structures for the first time}). A morphism $W:P\leadsto Q$ is a \textit{triad cobordism} between $P$ and $Q$, that is, a $d$-dimensional compact manifold triad $(W,\hb W,\vb W)$ as above where $\vb W=P\sqcup Q$, along with a map of pairs of vector bundles $l_W:(TW,T\hb W)\to (\theta^*\gamma_d,(\theta^\partial)^*\gamma_{d-1})$ compatible with $l_P$ and $l_Q$ (see \cref{cobordism categories definition} below for more details). We highlight now facts about this category, which are analogous to \ref{bdiff is inside cob}--\ref{stabilization intro} from above:

\begin{enumerate}[topsep=3pt, itemsep=3pt, label=\protect\circled{\arabic*$^{\scriptscriptstyle \partial}$}]
    \item\label{bdiff triads is in} Given a compact smooth $d$-dimensional manifold triad $(N,\hb N,\vb N)$ as above, then $\smash{\BDiff_{\vb}(N)}$ is a path component of the morphism space $\smash{\Cobb(\emptyset, \vb N)}.$ Again, the moduli space $\smash{\BDiff_{\vb}^\Theta(N)}$ defined in \ref{where i talk about theta structures for the first time} is a union of path components of $\Cobbt(\emptyset, \vb N)$ for any map of pairs $\Theta:(B,B^\partial)\to (\BO(d),\BO(d-1)).$
    \item Genauer \cite{genauer} showed that $\B\Cobbt$ is equivalent to $\Omega^\infty\MTH$, where $\MTH$ is the cofiber of the canonical map of Thom spectra $\Sigma^{-1} \MThb\to \MTc$ induced by $\Theta=(\theta,\theta^\partial)$.
    \item When $d=2n+1,$ there exists a \textit{distinguished} morphism $\smash{H_P:P\leadsto P}$ for every object $P\in \smash{\Cobbt}$ with non-empty connected boundary such that the left action map $\smash{H_P\circ (-):\Cobbt(\emptyset, P)\to \Cobbt(\emptyset, P)}$ extends the stabilization map $\BDiff_{\vb}(N\natural V_g)\to \BDiff_{\vb}(N\natural V_{g+1})$ defined above along the inclusion in \smash{\ref{bdiff triads is in}}, in case $P=\vb N$ and $\Theta=\id.$ The underlying cobordism of $H_P$ is the boundary connected sum $(P\times [0,1])\natural V_1$, performed away from $P\times \{0,1\}.$ Moreover, the space $\BDiff_{\vb}(N_\infty)$ is a path component of
    \[\Cobb(\emptyset, P)[(H_{P})^{-1}]\coloneqq \hocolim \left(\Cobb(\emptyset, P) \overset{H_{P}\circ (-)}{\rightarrow} \Cobb(\emptyset, P) \to \cdots \right)\]given by "inverting" the action of this morphism, when $P=\vb N$.
\end{enumerate}
In much of the same spirit as in the even-dimensional case, we reduce the proof of our main result (\cref{main no tang}) to the following theorem. For $k\geq -2$ an integer, we say that a map of pairs of spaces $f:(X,X')\to (Y,Y')$ is \textit{strongly k-connected} if $f|_{X'}:X'\to Y'$ is $k$-connected and the induced map $X\cup_{X'} Y'\to Y$ is $(k+1)$-connected, where the source denotes the homotopy pushout.

\begin{mainteo}[\cref{final with null}]\label{final with null in intro}
    Fix $n\geq 3$ and a map of pairs $\Theta:(B,B^\partial)\to (\BO(2n+1),\BO(2n))$. If the pair $(B,B^\partial)$ is $1$-connected and $B^\partial$ is path-connected, then for any $P\in \Cobbt$ such that $\partial P\neq \emptyset$ and $\Cob_{\Theta,n}^\partial(\emptyset,P)\neq \emptyset$, the canonical map
    \[\Cob_{\Theta,n}^\partial(\emptyset,P)[(H_P)^{-1}]\to \Omega_{[\emptyset,P]}\B\Cobbt\]is acyclic. Here, $\Cob_{\Theta,n}^\partial(\emptyset, P)$ denotes the subspace of $\Cobbt(\emptyset, P)$ of those morphisms $W:\emptyset\leadsto P$ such that $l_W:(W,\hb W)\to (B,B^\partial)$ is strongly $n$-connected.
\end{mainteo}

This result is the core of the present work. Its proof follows the overall strategy in the even-dimensional case of Galatius and Randal-Williams \cite{GRWStableModuli,GRWII}, but requires a number of generalizations (of varying difficulty) of geometric/surgery-theoretic techniques for even-dimensional manifolds with boundary (some classical and some developed in loc.cit.) to the context of odd-dimensional manifold triads.

\subsubsection*{Stable diffeomorphism classification of triads.}

As in the even-dimensional case, \cref{final with null in intro} already has content on $\H_0$: It yields an analog to Kreck's result \cite{Kreck1999SurgeryAD} for odd-dimensional triads. Since this result has not been (to our knowledge) previously established and might be of independent interest, we explain it now in full detail: Let $n\geq 3$ be an integer, $P$ be a compact $2n$-manifold with boundary, and $(N_i,\hb N_i,\vb N_i)$ be $(2n+1)$-dimensional manifold triads for $i=0,1$, together with an identification of $\vb N_i$ with $P$. Assume also that $\hb N_i$ is connected for $i=0,1.$ We say that:
\begin{enumerate}[itemsep=3pt, label=(\textit{\alph*})]
    \item\label{stably diffeo defn} The triads $N_0$ and $N_1$ are \textit{stably diffeomorphic} if there exists an integer $g\geq 0$ such that the triads $(N_0\natural V_g,\hb N_0\sharp W_{g},\vb N_0)$ and $(N_1\natural V_g,\hb N_1\sharp W_{g},\vb N_1)$ are diffeomorphic as triads relative to $(P,\partial P)$ (where the boundary connected sum is performed away from $\vb N_i$).
    \item\label{same normal type} The triads $N_0$ and $N_1$ have \textit{the same stable normal $n$-type} if there exist a map of pairs $\Theta^{\perp}:(B,B^\partial)\to (\BO,\BO)$ and strongly $n$-connected maps of pairs $l_i:(N_i,\hb N_i)\to (B,B^\partial)$ for $i=0,1$ making the following diagram commute up to homotopy
    \[\begin{tikzcd}
    	& {(N_0,\hb N_0)} \\
    	{(P,\partial P)} & {(B,B^\partial)} & {(\BO,\BO)} \\
    	& {(N_1,\hb N_1)}
    	\arrow["{l_1}", dashed, from=1-2, to=2-2]
    	\arrow["{\nu_0}", bend left=20, from=1-2, to=2-3]
    	\arrow[bend left=20, hook, from=2-1, to=1-2]
    	\arrow[bend right=20, hook, from=2-1, to=3-2]
    	\arrow["{\Theta^\perp}", from=2-2, to=2-3]
    	\arrow["{l_0}"', dashed, from=3-2, to=2-2]
    	\arrow["{\nu_1}"', bend right=20, from=3-2, to=2-3]
    \end{tikzcd}\]where $\nu_i$ is the classifying map for the stable normal bundle $(\nu_{N_i},\nu_{\hb N_i})$ of the pair $(N_i,\hb N_i).$
    \item\label{bordant normal smoothings} Given a map of pairs $\Theta^{\perp}:(B,B^\partial)\to (\BO,\BO)$, the triads $N_0$ and $N_1$ \textit{admit bordant $\Theta^\perp$-smoothings} if there exist choices of strongly $n$-connected maps $\smash{l_i:(N_i,\hb N_i)\to (B,B^\partial)}$ making the diagram above commute, such that the bordism class of manifolds with boundary \[ [N_0\cup_{P} (-N_1)]\in \Omega_{2n+1}^{\Theta^\perp}\]vanishes, where $\Omega_{k}^{\Theta^\perp}$ is the relative bordism group of $k$-dimensional compact manifolds with boundary $(M,\partial M)$ together with a lift $l_M:(M,\partial M)\to (B,B^\partial)$ of its stable normal bundle along $\Theta^\perp$. The bordism relation is triad bordism (in the sense of \cref{relative cobordism} but $\partial_0=\emptyset$) together with lifts of the stable normal bundle of the bordism compatible with the given lifts in its vertical boundary. Here, $N_0\cup_P(-N_1)$ denotes the manifold obtained by the gluing of $N_0$ and $N_1$ along the identifications with $P$, together with the map $l_0\cup_P (-l_1)$ to $(B,B^\partial)$, where $-l_1$ is the restriction of the map $l_1\oplus \varepsilon^1:(N_1\times [0,1],\hb N_1\times [0,1])\to (B,B^\partial)$ to $N_1\times \{1\},$ using an inwards-pointing trivialization of the normal bundle of $N_1\times \{1\}$ in $N_1\times [0,1]$ (see \cite[17]{Stong2015-dc}). 
\end{enumerate}The classical Pontrjagin-Thom construction for relative bordisms yields an isomorphism $\Omega_k^{\Theta^\perp}\cong \pi_k(\MM\Theta^\perp)$, where $\MM\Theta^\perp$ is the cofiber of the map of Thom spectra $\MM(\theta^\perp)^\partial\to \MM\theta^\perp$ induced by $\Theta^\perp=(\theta^\perp,(\theta^\perp)^\partial)$ (see e.g. \cite[25]{Stong2015-dc} or \cite[Thm. 3.1.5]{laures}). The following result is deduced from \cref{final with null in intro} below in \cref{proof of kreck section}.

\begin{maincor}[Stable diffeomorphism classification]\label{stable class thm}
    Let $n\geq 3$ be an integer, $P$ be a compact $2n$-manifold with boundary, and $(N_i,\hb N_i,\vb N_i)$ be $(2n+1)$-dimensional manifold triads for $i=0,1$, together with an identification of $\vb N_i$ with $P$. Assume:
    \begin{enumerate}[label=(\textit{\roman*})]
        \item The manifold $\hb N_i$ is non-empty and connected for $i=0,1.$
        \item The pair $(N_i,\hb N_i)$ is $1$-connected for $i=0,1.$
    \end{enumerate}Then the triads $N_0$ and $N_1$ are stably diffeomorphic relative to $P$ if and only if 
    \begin{enumerate}
        \item their relative Euler characteristics agree, i.e. $\chi(N_0,\hb N_0)=\chi(N_1,\hb N_1)$, and
        \item they have the same stable normal $n$-type $\Theta^\perp$ and admit bordant $\Theta^\perp$-smoothings.
    \end{enumerate}
\end{maincor}
\noindent\textbf{On the assumptions.} We conclude this introduction by commenting on the assumptions of our results:
\begin{enumerate}[topsep=2pt,label=(\textit{\roman*})]
    \item It would be interesting to see whether \cref{final with null in intro}, and consequently \cref{main no tang,stable class thm}, also holds in the excluded dimensions $3$ and $5$. One starting point could be to see whether \cref{stable class thm} holds in these dimensions, possibly via more direct methods.
    \item We believe that the $1$-connectivity conditions in our results are necessary. For example, it seems plausible that one can find an explicit counterexample of \cref{stable class thm} when the $1$-connectivity condition is dropped.
\end{enumerate}

\subsubsection*{Structure of the paper.}
\cref{preliminaries} introduces various categories of pairs that are used throughout, while introducing the notion of \textit{strong connectivity} of maps of pairs of spaces. We recall the work of Borodzik, Némethi, and Ranicki on handle decompositions for triads \cite{BNR} and prove a \textit{geometrical connectivity} statement in \cref{strong geometrical connectivity} relating strong connectivity with existence of certain handle decompositions. In \cref{cob cat section}, we recall the definition of the \textit{cobordism category of manifolds with boundary} from \cite{genauer,BP,steimle} and variants thereof. The main result of this section is \cref{surgery on morphisms}, where we prove a generalization of results of \cite{BP} for this larger class of tangential structures, which is an analog of \cite[Thm. 3.1]{GRWStableModuli}. In \cref{sta sta section} and \cref{proof of closure section} resides the core of the proof of \cref{final with null in intro,main}, where we prove a "stable stability" statement, which is analogous to \cite[Thm. 2.15]{GRWII}. The main geometrical input for this proof is a modification result for embeddings of discs into manifolds with boundary via surgery. We deduce \cref{final with null in intro} and a generalization thereof from the stable stability phenomena, in \cref{group comp section}. In \cref{general}, we deduce \cref{main} from \cref{final with null in intro}, by exploiting the notion of Moore--Postnikov factorizations of pairs. We finish by presenting further simplifications of \cref{main} in the context of certain examples in \cref{examples section}.

\subsubsection*{Acknowlegments. }

First and foremost, I would like to thank my Ph.D. advisor Manuel Krannich for suggesting this problem and providing invaluable support and insight throughout the whole project. I would also like to express my gratitude to Oscar Randal-Williams for his interest and expertise which led to many enlightening discussions, as well as for his hospitality at the University of Cambridge. I am also grateful to Jan Steinebrunner, Mark Powell, Robin Stoll, and Samuel Muñoz-Echániz for helpful discussions and comments in earlier drafts. Finally, I acknowledge an enormous intellectual debt to Søren Galatius and Oscar Randal-Williams. Without their work on stable moduli spaces of manifolds, this work would not exist.

I was supported by the Deutsche Forschungsgemeinschaft (DFG, German Research Foundation) through the Research Training Group DFG 281869850 (RTG 2229).

\section{Pairs of spaces and manifolds.}\label{preliminaries}

This section introduces and studies properties of pairs of spaces, manifolds, and vector bundles. We start by introducing the notion of \textit{strong connectivity} of maps of pairs of spaces and prove closure properties of this notion. We introduce manifold pairs and an appropriate notion of \textit{handle decompositions} on them and prove that this notion behaves well with the notion of strong connectivity.

\subsection{Strong connectivity.}\label{strong connectivity section}

    Throughout this paper, by \textit{space} we mean a compactly generated space in the sense of \cite[Defn. 1.1]{strickland}. We consider the category of spaces $\Top$ as a model category with the Quillen model structure (see \cite[Defn. 7.10.6]{Hirschhorn2003-du}), which turns it into a self-enriched cartesian closed model category (see \cite[Notation 7.10.2]{Hirschhorn2003-du} and \cite[Prop. 2.12]{strickland}\footnote{We point out a typo in this statement. The conclusion should read "the category of compactly generated spaces is cartesian closed". In fact, the statement about compactly generated weak Hausdorff spaces is given in \cite[Prop. 2.24]{strickland}.}). In this model category, we denote the homotopy pushout of a diagram $X\leftarrow Y\to Z$ as $X\cup_YZ$. When the map $Y\to X$ is a Hurewicz cofibration, the homotopy pushout is equivalent to the strict pushout (this follows from \cite[Thm. A.7]{Dugger2004}, as homotopy colimits in the Quillen and the Strom model structure agree). Thus, we can model $X\cup_Y Z$ using the double mapping cylinder construction (see \cite[Defn. 3.6.3]{Munson_Volić_2015}). When $Y\to X$ is a Hurewicz cofibration, we implicitly denote the strict pushout by $X\cup_YZ.$ In particular, we do so in the case where $Y\to X$ is the inclusion of a compact smooth submanifold of a smooth manifold. We denote by $X\times_YZ$ the homotopy pullback of a diagram $X\to Y\leftarrow Z$. Similarly, we often model this space using the path space construction (see \cite[Defn. 3.2.4]{Munson_Volić_2015}).

    By a \textit{pair of spaces}, we mean a map of spaces $A'\to A.$ If the map is implicit, we denote it by $(A,A').$ A map of pairs $f=(\alpha,\alpha'):(A,A')\to (B,B')$ is the data of two maps of spaces $\alpha:A\to B$ and $\alpha':A'\to B' $ such that the square
    \[\begin{tikzcd}
        A'\arrow[d]\arrow[r, "\alpha'"] & B' \arrow[d] \\
        A\arrow[r, "\alpha"] & B  
    \end{tikzcd}\]commutes. Given a map of pairs $f=(\alpha,\alpha'):(A,A')\to (B,B')$, we often denote the individual maps $\alpha$ and $\alpha'$ by $f|_{A}$ and by $f|_{A'}$, respectively. Given two pairs $(A,A')$ and $(B,B'),$ we denote the mapping space of pairs by $\Map((A,A'),(B,B'))$ seen as subspace of $\Map(A',B')\times \Map(A,B)$. Given a map $f':A'\to B',$ we denote the (strict) fiber of the projection $\Map((A,A'),(B,B'))\to \Map(A',B')$ at $f'$ by $\Map_{f'}(A,B)$ or $\Map_{A'}(A,B)$ if $f'$ is implicit.
\subsubsection{Connectivity for pairs.}

We introduce a convenient notion of connectivity in the category of pairs of spaces, which will be heavily used throughout the entire paper. Recall that for $k\geq -1$, a map of spaces is $k$-connected if all its homotopy fibers are $(k-1)$-connected. In particular, any map is $(-1)$-connected.

\begin{defn}\label{strong conn defn}
    Let $k\geq -1$. A map of pairs of spaces $f:(A,A')\to (B,B')$ is \textit{strongly $k$-connected} if $f|_{A'}$ is $k$-connected and the map $A\cup_{A'}B'\to B$ is $(k+1)$-connected.
\end{defn}

\begin{rmk}\label{cocartesian definition}
    We point out that this notion implies that the individual maps $f|_{A'}$ and $f|_{A}$ are $k$-connected, the first by definition and the second by the following argument. The second condition in \cref{strong conn defn} is, by definition, the condition that the square
    \[\begin{tikzcd}
        A' \arrow[d]\arrow[r] & B'\arrow[d] \\
        A \arrow[r] & B
    \end{tikzcd}\]is homotopy $(k+1)$-cocartesian, in the sense of \cite[Defn. 3.7.1]{Munson_Volić_2015}. We conclude that $f|_A$ is $k$-connected by applying \cite[Prop. 3.7.13.2.(a)]{Munson_Volić_2015} to this square, using the fact that $f|_{A'}$ is $k$-connected.
\end{rmk}

\begin{lemma}[2 out of 3]\label{closure properties strongly}
    Let $k\geq -1$. Let $f:(A,A')\to (B,B')$ and $g:(B,B')\to (C,C')$ be maps of pairs, then the following holds:
    \begin{enumerate}
        \item If $f$ and $g$ are strongly $k$-connected, then $g\circ f$ is strongly $k$-connected.
        
        \item If $f$ is strongly $(k-1)$-connected and $g\circ f$ are strongly $k$-connected, then $g$ is strongly $k$-connected.
        \item If $g$ is strongly $k$-connected and $g\circ f$ are strongly $(k-1)$-connected, then $f$ is strongly $(k-1)$-connected, provided $k\geq 2.$
    \end{enumerate} 
\end{lemma}
\begin{proof}
    Properties $1$ and $2$ follow from classical properties of connectivity of maps of spaces (see \cite[Prop. 2.6.15]{Munson_Volić_2015}) and \cite[Prop. 3.7.26]{Munson_Volić_2015}. For property $3$, the same properties of connectivity of maps of space imply that the map $f|_{A'}$ is $(k-1)$-connected. It remains to prove that $A\cup_{A'} B'\to B$ is $k$-connected. Consider the following homotopy pushout square of spaces
    \[\begin{tikzcd}
        A\cup_{A'}B'\arrow[d]\arrow[r] & B\arrow[d] \\
        A\cup_{A'}C'\arrow[r] & B\cup_{B'} C' 
    \end{tikzcd}\]given by identifying the bottom left corner with $A\cup_{A'}B'\cup_{B'}C'.$ By hypothesis, the map $A\cup_{A'} C'\to C$ is $k$-connected. Since the map ${B}\cup_{B'} {C'}\to C$ is $(k+1)$-connected, it follows that $A\cup_{A'} C'\to {B}\cup_{B'} {C'}$ is $k$-connected. On the other hand, $A\cup_{A'} B'\to A\cup_{A'} C'$ is $k$-connected, since $g|_{B'}:B'\to C' $ is $k$-connected. Since $k\geq 2,$ the former induces an isomorphism on fundamental groupoids and thus, the map $A\cup_{A'} B'\to B$ is $k$-connected by \cite[Prop. A]{WallGeometrical}.    
\end{proof}

Given a commutative square of maps of pairs, we call it a \textit{homotopy pushout square} if both the square consisting of the maps between the targets of the pairs and the square of the maps between the sources of the pairs are homotopy pushout squares of spaces. Analogously, a \textit{homotopy pullback square} is a commutative square that restricts to a homotopy pullback square of spaces both in the targets and the sources.   

\begin{lemma}[Pushouts]\label{strong conn and pushouts}
    Let $k\geq -1$. Assume that the square
    \begin{equation*}
        \begin{tikzcd}
        (A,A')\arrow[r,"\alpha"]\arrow[d, "f"'] &(C,C')\arrow[d,"g"]\\
        (B,B')\arrow[r] &(D,D').
    \end{tikzcd}
    \end{equation*}is a homotopy pushout square of pairs. If the map $f$ is strongly $k$-connected, then so is $g$. Conversely, if $g$ is strongly $k$-connected, then so is $f$, provided both $\alpha|_A$ and $\alpha|_{A'}$ induce isomorphisms on fundamental groupoids.
\end{lemma}
\begin{proof}
    We start by observing that, if $f|_{A'}$ is $k$-connected, then so is $g|_{C'}$, as $k$-connected maps are closed under taking homotopy pushouts. Conversely, if $g|_{C'}$ is $k$-connected and $\alpha|_{A'}$ induces an isomprhism of fundamental groupoids, then $f|_{A'}$ is $k$-connected by \cite[Prop. A]{WallGeometrical}. We proceed by showing that if $f$ is strongly $k$-connected, than $C\cup_{C'} D'\to D$ is $(k+1)$-connected, and thus $g$ is strongly $k$-connected, and its converse under the additional assumptions stated in the result. To do so, we prove now the subclaim that the following left square in
    \[\begin{tikzcd}\label{the square we want to prove}
        A\cup_{A'} B'\arrow[d]\arrow[r] & C\cup_{C'} D'\arrow[d] \\
        B\arrow[r] & D
    \end{tikzcd} \quad \begin{tikzcd}
        A\arrow[d]\arrow[r] & C\arrow[d] \\
        A\cup_{A'} B'\arrow[r] & C\cup_{C'} D'
    \end{tikzcd}\]is a homotopy pushout square. By the pasting lemma for homotopy pushout squares \cite[Prop. 3.7.26]{Munson_Volić_2015}, it suffices to prove that the right square above is a homotopy pushout square. Consider the following diagrams 
    \[\begin{tikzcd}
    A'\arrow[r]\arrow[d] & A\arrow[r] \arrow[d] & C\arrow[d] \\
    B'\arrow[r] & A\cup_{A'} B' \arrow[r] & C\cup_{C'} D'
    \end{tikzcd} \quad \begin{tikzcd}
    A'\arrow[r]\arrow[d] & C'\arrow[r] \arrow[d] & C\arrow[d] \\
    B'\arrow[r] & D' \arrow[r] & C\cup_{C'} D'
    \end{tikzcd}. \]Since the left square in the left diagram is a homotopy pushout by definition, it suffices to prove that the pasted square is a homotopy pushout. For that, we observe that it can be factored as in the right diagram above as the pasting of two homotopy pushout squares. This finishes the proof of the subclaim. We conclude that, if $A\cup_{A'}B'\to B $ is $(k+1)$-connected, then so is $C\cup_{C'}D'\to D,$ thus establishing that $g$ is strongly $k$-connected. 
    
    For the converse claim (under the additional assumption that $\alpha|_A$ and $\alpha|_A'$ induce fundamental groupoid isomorphisms), we have that $A\cup_{A'} B'\to C\cup_{C'} D'$ is a fundamental groupoid isomorphism, since $\alpha|_A$ also is. By \cite[Prop A]{WallGeometrical}, the map $A\cup_{A'} B'\to B$ is $(k+1)$-connected, as $C\cup_{C'}D'\to D$ also is. Hence, $f$ is strongly $k$-connected.
\end{proof}

For $k\geq -1$, we say that a pair of spaces $(A,A')$ is $k$-connected if the map of spaces $A'\to A$ is $k$-connected.

\begin{lemma}\label{1 conn: old same as new}
     Let $f:(A,A')\to (B,B')$ be a map of pairs such that $f|_{A'}$ is $k$-connected. 
     \begin{enumerate}[label=(\roman*)]
         \item  If $(A,A')$ is $0$-connected and $f$ is strongly $k$-connected, then the map $f_*:\pi_i(A,A',x)\to \pi_i(B,B',f(x))$ is an isomorphism for $i<k$ and a surjection for $i=k$ for all basepoints $x\in A'.$
         \item If $(A,A')$ is $1$-connected, then $f$ is strongly $k$-connected if and only if the map $f_*:\pi_i(A,A',x)\to \pi_i(B,B',f(x))$ is an isomorphism for $i<k+1$ and a surjection for $i=k+1$ for all basepoints $x\in A'.$
     \end{enumerate}
\end{lemma}
\begin{proof}
    We start by proving the first claim. Note that the condition on $f_*$ is equivalent to the square
    \[\begin{tikzcd}
        A' \arrow[d]\arrow[r] & B'\arrow[d] \\
        A \arrow[r] & B
    \end{tikzcd}\]being homotopy $(k-1)$-cartesian (in the sense of \cite[Defn. 3.3.1]{Munson_Volić_2015}), while the assumption is equivalent to it being homotopy $(k+1)$-cocartesian. We know that $A'\to B'$ is $k$-connected and $A'\to A$ is $0$-connected. Thus, the Blakers-Massey theorem implies that homotopy $(k+1)$-cocartesianness of the square above implies homotopy $(k-1)$-cartesianness (see \cite[Thm. 4.2.3] {Munson_Volić_2015}). This finishes the proof of the claim. The proof of the second claim follows similarly by observing that, since $(A,A')$ is $1$-connected, homotopy $(k+1)$-cocartesianness of the square above is equivalent to homotopy $k$-cartesianness by the Blakers-Massey theorem and its dual (see \cite[Thm. 4.2.3/4] {Munson_Volić_2015}).
\end{proof}

\subsubsection{Triad homotopy groups.}

It will be convenient to have an analog of the relative homotopy groups in the context of maps of pairs of spaces, which we start by recalling. When studying maps of spaces $f:X\to Y$, it is convenient to define the \textit{relative homotopy groups} $\pi_k(Y,X,x)$ as $\pi_0$ of the mapping space of maps $(D^k,\partial D^k)\to (Y,X)$ in the category of pairs of pointed spaces $\smash{\Top_*^{[1]}}$, where $X$ is pointed at $x$ and $Y$ at $f(x)$. These agree with the homotopy groups at $x$ of the homotopy fiber of $f$ at $f(x)$ in one degree lower. We now consider an analogous concept for maps of pairs and record some facts about them. Let $D^k_+\subset \RR^k$ be the subset of those tuples $x=(x_1,\cdots, x_k)$ such that $||x||\leq 1$ and $x_1\geq 0$. Let $\partial_0 D^k_+$ be the subset of $D^k_+$ of those points where $x_1=0$ and $\partial_1 D^k_+$ of those points with norm $1$. Let $\partial_{01}D^k$ be the intersection $\partial_0 D^k_+\cap \partial_1D^k_+.$

\begin{defn}[Triad homotopy groups]\label{triad homotopy}
    Let $k\geq 2$ and $f:(A,A')\to (B,B')$ a map of pairs and $a\in A' $. The $k$-th \textit{triad homotopy group} $\pi_k(B,B',A,a)$ is the set of path components of the space of squares
    \[\begin{tikzcd}
        (\partial_0 D^k_+,\partial_{01} D^k_+) \arrow[r, "\alpha"]\arrow[d] & (A,A')\arrow[d, "f"] \\
        (D^k_+,\partial_{1} D^k_+) \arrow[r, "\beta"] & (B,B')
    \end{tikzcd}\]seen as the mapping space in the category of functors from $[1]\times [1]$ to $\Top_*$ from the left vertical map of pairs to the rightmost vertical map of pairs. 
\end{defn}

\begin{war}
    We note that, although $A'$ is not in the notation, these groups depend heavily on $A'$. We believe that adding $A'$ in the notation would make it too cumbersome, so we choose to drop it.
\end{war}

See \cite{Blakers1949} for this definition when $A,B'\subset B$ and $A'=A\cap B$. As in this special case, when $k=2$ this is just a set with no extra structure. When $k\geq 3$, one can define a group structure given by stacking and for $k\geq 4$, this group structure is abelian. We can map this group (or set) to the $(k-2)$-nd homotopy group of the total homotopy fiber $\tofib_a(f)\coloneqq \hofib_a(A'\to B'\times_B A)$, modelled using the path-space construction, in the following way. By unwrapping the definition, we see that $\tofib_a(f)$ is the space of tuples $(x,\gamma_0,\gamma_1,\Gamma)\in A'\times (B')^I\times A^I\times B^{I\times I}$ such that $\gamma_0$ is a path in $B'$ from $f(x)$ to $f(a)$, $\gamma_1$ is a path in $A$ from $\iota_A(x)$ to $\iota_A(a)$, where $\iota_A:A'\to A$, and $\Gamma$ is a homotopy from $\iota_{B'}\circ \gamma_0$ to $f\circ \gamma_1,$ where $\iota_B:B'\to B$ (see also \cite[Prop. 5.5.8]{Munson_Volić_2015}). We consider this space as based in $a$, that is, at the tuple $(a,\const_{f(a)},\const_{\iota_A(a)},\const_{f(\iota(a))})$. We can define a map
\(\pi_k(B,B',A,a)\to \pi_{k-2}(\tofib_a(f),a)\) taking the class $[\beta,\alpha]$ of a diagram as above to the tuple $(\alpha|_{\partial_{01}D^k_+},\beta|_{\partial_1D^k_+},\alpha,\beta) $ by looking at $\alpha$ as a map from $\smash{\frac{\partial_{01}D^k_+\times I}{\partial_{01}D^k_+\times \{1\}}}$ and all the other maps analogously. It is not difficult to see that this map is a bijection for $k=2$ and a group isomorphism for $k\geq 3$. 

These triad homotopy groups are related to relative homotopy groups in the following way. We have a sequence of maps
\[\smash{\cdots \to \pi_k(A,A',a)\overset{f_*}{\to} \pi_k(B,B',f(a))\overset{p}{\to} \pi_k(B,B',A,a)\overset{\delta}{\to} \pi_{k-1}(A,A',a)\to \cdots} \]where $p$ takes the class of a map $(D^k,\partial D^k,*)\to (B,B',f(a))$ to the diagram where $\alpha$ is constant, and $\delta$ sends $[\beta,\alpha]$ to $[\alpha].$ This sequence is compatible with the identification of the relative and triad homotopy groups with the homotopy groups of the fibers of $\iota_A$ and $\iota_B$ and of the total homotopy fibers. We conclude that this sequence is exact. We finish with the following lemma, which will be useful in constructing maps of pairs.

\begin{lemma}\label{strong implies lifting}
    Let $k\geq 2$ and $i\leq  k$ be integers. Assume $(A,A')$ is $0$-connected and $f:(A,A')\to (B,B')$ is a  strongly $k$-connected map of pairs. Then every commutative square
    \[\begin{tikzcd}
        (\partial_0 D^i_+,\partial_{01} D^i_+) \arrow[r]\arrow[d] & (A,A')\arrow[d, "f"] \\
        (D^i_+,\partial_{1} D^i_+) \arrow[r]\arrow[ru, dashed, "g"] & (B,B')
    \end{tikzcd}\]admits a filler $g$ making the upper triangle commute strictly and the bottom one up to homotopy of pairs. If $(A,A')$ is $1$-connected, then the same conclusion holds for $i\leq k+1$.
\end{lemma}
\begin{proof}
    This follows by combining \cref{1 conn: old same as new} with the long exact sequence from above to deduce that $\pi_i(B,B',A)$ vanishes for every basepoint. From a nullhomotopy of the homotopy class of this square, we can construct such a filler $g$. We leave this check to the reader. 
\end{proof}

\subsection{Pairs of vector bundles.}

We consider now an appropriate notion of a vector bundle over a pair of spaces. The following is essentially a recollection of \cite[Section 4]{steimle}. For the entirety of this paper, a \textit{bundle map} between vector bundles will always be assumed to be fiberwise injective.

\begin{defn}\label{pairs of vector bundles}
    A \textit{pair of vector bundles} $(\xi,\xi')$ is a map of pairs $(\pi,\pi'):(\xi,\xi')\to (X,X')$ such that both $\pi$ and $\pi'$ are vector bundles and the map $\xi'\to \xi$ is a bundle map. Given two vector bundle pairs $\pi_X:(\xi,\xi')\to (X,X')$ and $\pi_Y:(\eta,\eta')\to (Y,Y')$, a \textit{map of pairs of vector bundles} $f:(\xi,\xi')\to (\eta,\eta')$ is the data of a commutative square of pairs of spaces
    \[\begin{tikzcd}
        (\xi,\xi')\arrow[r, "{(f,f')}"]\arrow[d, "\pi_X"] & (\eta,\eta')\arrow[d, "\pi_Y"]\\
        (X,X')\arrow[r] & (Y,Y')
    \end{tikzcd}\]such that $f$ and $f'$ are bundle maps.
\end{defn}

\begin{defn}\label{collared pairs}
    Let $(\xi,\xi')$ be a pair of vector bundles. A \textit{collar} on $(\xi,\xi')$ is a bundle map $c:\xi'\oplus \varepsilon^1\to \xi$ extending $\xi'\to \xi$ which is a fiberwise isomorphism, where $\varepsilon^1$ denotes the trivial vector bundle of rank $1$. We call the pair $((\xi,\xi'),c)$ a \textit{collared pair of vector bundles}. A map $f:(\xi,\xi')\to (\eta,\eta')$ between collared pairs of vector bundles is \textit{collared} if the resulting square commutes:
    \[\begin{tikzcd}
        \xi'\oplus \varepsilon^1\arrow[r, "f'\oplus \id"]\arrow[d] & \eta'\oplus \varepsilon^1\arrow[d]\\
        \xi\arrow[r,"f"] & \eta
    \end{tikzcd}.\]We denote the space of collared bundle maps by $\Bunc(\xi,\eta)$, topologized with the compact-open topology. Given a collared pair of vector bundles $(\xi,\xi')$ over $(X,X')$ and $\iota:(A,A')\to (X,X'),$ we have a collared pair of vector bundles $((\iota|_A)^*\xi,(\iota|_{A'})^*\xi')$ over $(A,A')$ defined by pulling back the individual bundles and the collar. Given a collared map $f|_A:((\iota|_A)^*\xi,(\iota|_{A'})^*\xi')\to (\eta,\eta')$ of collared vector bundle pairs, let $\Bunca(\xi,\eta)$ be the subspace of collared bundle maps extending $f|_A.$
\end{defn}

\begin{rmk}\label{bun col is pullback}
    Unravelling the definition, one sees that we have a strict pullback square of spaces
    \[\begin{tikzcd}
        \Bunca(\xi,\eta) \arrow[d] \arrow[r] & \Bun_A(\xi,\eta) \arrow[d] \\
        \Bun_{A'}(\xi',\eta') \arrow[r] & \Bun_{A'}(\xi'\oplus \varepsilon^1,\eta)
    \end{tikzcd}\]where $\Bun_A(\xi,\eta)$ denotes the space of bundle maps of vector bundles $\xi\to\eta$ extending a map $f|_A:\xi|_A\to \eta$. The left vertical map takes a collared bundle map $(f,f')$ and maps it to $f'$. The right vertical map takes a bundle map $f:\xi\to \eta$ to the composite $\smash{f|_{X'}=f\circ c:\eta'\oplus \varepsilon^1\to \eta|_{X'}\to \eta}$.
\end{rmk}

    Throughout the entirety of this paper, we model the classifying space $\BO(d)$ of the $d$-dimensional orthogonal group as the colimit of the Grassmanian manifolds $\Gr_d(d+k)$ of $d$-planes along the standard inclusions $\Gr_d(d+k)\to \Gr_d(d+k+1)$. Let $\gamma_d$ be the universal vector bundle of dimension $d$ over $\BO(d)$ given by the colimit over $k$ of the space of pairs $(V,v)$ of $V\in \Gr_d(d+k)$ a $d$-plane in $\RR^{d+k}$ and $v\in V$. There is a pair of vector bundles $\gamma_d:(\gamma_d,\gamma_{d-1})\to (\BO(d),\BO(d-1))$, whose underlying map of vector bundles is the usual stabilization map. This vector bundle pair admits a preferred collar given by the canonical map $\gamma_{d-1}\oplus \varepsilon^1\to \gamma_d$. Given a map $\Theta=(\theta,\theta'):(B,B')\to (\BO(d),\BO(d-1))$ of pairs, then the vector bundle $\Theta^*\gamma_d:(\theta^*\gamma_d,\theta'^*\gamma_{d-1})\to (B,B')$ admits a collar given by pulling back the collar along $\Theta.$ Given a collared vector bundle pair $(\xi,\xi')$ of dimension $l<d$, we have a preferred collar for $(\xi,\xi')\oplus \varepsilon^{d-l}\coloneqq(\xi\oplus \varepsilon^{d-l},\xi'\oplus \varepsilon^{d-l})$ induced by the collar $\xi'\oplus\varepsilon^1\to \xi.$

    \begin{rmk}[Universality of $\gamma$]\label{gamma is universal}
        Let $((\xi,\xi'),c)$ be a collared $d$-dimensional vector bundle pair over $(X,X')$, $(A,A')\to (X,X')$ be a cofibration and $f|_A:(\xi|_A,\xi'|_{A'})\to (\gamma_d,\gamma_{d-1})$ a collared bundle map, then the space $\Bunca(\xi,\gamma_d)$ is weakly contractible. This follows by first observing that the square in \cref{bun col is pullback} is a homotopy pullback square, since the right-hand vertical map is a fibration of spaces (see the proof of \cref{bundle map fibration}). If we take $(\eta,\eta')=(\gamma_d,\gamma_{d-1})$, then all of the spaces except the one in the initial vertex are weakly contractible, by universality of $\gamma_d$ and $\gamma_{d-1}.$ We conclude that the space in the initial vertex, that is $\Bunca(\xi,\gamma_d)$, is weakly contractible.
    \end{rmk}

\subsection{Pairs of manifolds and triads.}\label{section of pairs of manifolds} Throughout this paper, we follow the convention that \textit{manifolds} are smooth and allowed to have boundary, without further mention. A $d$-dimensional \textit{manifold pair} is a pair $(M,K)$ where $M$ is a $d$-dimensional compact manifold and $K$ is a codimension $0$ compact submanifold of $\partial M.$ The map of pairs $(TM,TK)\to (M,K)$ given by the differential of the inclusion is a pair of vector bundles. An inwards-pointing vector field on $K$ induces a collar on this vector bundle pair. We implicity carry such an inwards-pointing vector field and induced collar whenever we have a manifold pair. A map of manifold pairs $f:(M,K)\to (M',K')$ is a map of pairs of spaces where both maps are smooth, $f^{-1}(K')=K$, and $Df:(TM,TK)\to (TM',TK')$ is a collared bundle map. A related notion that we will use is that of a triad. A $d$-dimensional \textit{manifold triad} is a triple $(W,\hb W,\vb W)$ where $(W,\hb W)$ and $(W,\vb W)$ are manifold pairs such that $\hb W\cup \vb W=\partial W$ and $\partial (\hb W)=\hb W\cap \vb W=\partial (\vb W).$ We denote $\hb W\cap \vb W$ by $\hvb W$. When dealing with these objects, we will often only remember the manifold pairs $(W,\hb W)$ and $(\vb W,\hvb W)$, and the inclusion map between them. The following spaces of smooth maps will appear often throughout:
\begin{enumerate}[label=(\alph*)]
    \item Given two manifold pairs $(M,K)$ and $(M',K')$, the space $\Emb((M,K),(M',K'))$ (or $\Emb(M,M')$ when the subspaces are implicit) denotes the space of smooth embeddings of pairs, that is, smooth maps of pairs which are topological embeddings and whose differentials are fiberwise injective. For a subpair $(N,L)\subseteq (M,K) $ and an embedding $e:(N,L)\to (M',K')$, we denote by $\Emb_N((M,K),(M',K'))$ the space of smooth embeddings which extend $e$ on $(N,L)$.
    \item Given a manifold triad $(W,\hb W,\vb W),$ we denote by $\Diff_{\vb}(W)$ the space of diffeomorphisms of $W$ fixing $\vb W$ pointwise, whose differential is a collared bundle map. This space is a topological group where multiplication is given by composition of diffeomorphism and inverses are given by taking the inverse diffeomorphism.
\end{enumerate}Both of these spaces (and any subspace of the space of smooth maps) are equipped with the Whitney $C^\infty$-topology. By definition, the diffeomorphism group $\Diff_{\vb}(W)$ acts on $\Buncp(TW,\eta)$ for any vector bundle pair.

\subsubsection{Handle decompositions for pairs of manifolds.}

In this section, we introduce handle decompositions for cobordisms between triads (as defined below, see also \cref{pictures of handle}). To do so, we are inspired by the work of Borodzik, Némethi, and Ranicki in \cite{BNR}. Recall that a \textit{4-ad} is a compact $d$-dimensional manifold $W$ together with a decomposition $\partial W=\partial_0 W\cup \partial_1 W\cup \partial_2 W$ into codimension $0$ compact submanifolds such that $\partial(\partial_i W)=\smash{\bigcup_{i\neq j} \partial_i W\cap \partial_{\!j} W}$ and that $\partial(\partial_i W\cap \partial_{\!j} W)=\smash{\bigcap_{k} \partial_k W}.$ We use the notation $\partial_{ij} W\coloneqq \partial_{i}W\cap \partial_{\!j} W$ and $\partial_{012} W\coloneqq \partial_0 W\cap \partial_1 W\cap \partial_2 W.$

\begin{defn}\label{relative cobordism}
    Let $(M,\hb M,\vb M)$ and $(N,\hb N,\vb N)$ be two $(d-1)$-dimensional manifold triads. A $d$-dimensional manifold $4$-ad $(W,\partial_0 W,\partial_1 W,\partial_2 W)$ is called a \textit{triad cobordism from $M$ to $N$} if there is an equality of triads $(\partial_2 W,\partial_{02} W,\partial_{12} W)=(M,\hb M,\vb M)\sqcup (N,\hb N,\vb N)$ and there is a diffeomorphism of pairs $(\partial_1 W,\partial_{01} W)\cong (\vb M\times [0,1],\hvb M\times [0,1])$ relative to $(\vb M,\hvb M)$, seen as a subspace of the second factor given by the inclusion at $\{0\}$. In this case, we denote $\hb W\coloneqq \partial_0 W$ and $\vb W=\partial_2 W $. We denote a triad cobordism $W$ from $M$ to $N$ by $W:M\leadsto N.$  Observe that the definition implies that $\hb W$ is a cobordism with trivial boundary (or a \textit{r$\partial$}-cobordism in the sense of \cite{WallGeometrical}) from $\hb M$ to $\hb N$ (see \cref{relative cobordism picture}).
\end{defn}

\begin{figure}
    \centering
    \tikzset{every picture/.style={line width=0.75pt}} 

\begin{tikzpicture}[x=0.75pt,y=0.75pt,yscale=-1,xscale=1]

\draw [color={rgb, 255:red, 43; green, 0; blue, 255 }  ,draw opacity=1 ]   (230.29,140.67) .. controls (230.39,114.64) and (233.61,102.58) .. (237.5,101.35) .. controls (243.59,99.42) and (251.32,123.97) .. (251.32,162.95) ;
\draw [color={rgb, 255:red, 144; green, 19; blue, 254 }  ,draw opacity=1 ]   (230.29,140.67) -- (251.32,162.95) ;
\draw [color={rgb, 255:red, 208; green, 2; blue, 27 }  ,draw opacity=1 ]   (251.32,162.95) -- (391.55,161.19) ;
\draw    (237.5,101.35) .. controls (267.02,83.76) and (347.27,119.56) .. (376.79,101.97) ;
\draw [color={rgb, 255:red, 43; green, 0; blue, 255 }  ,draw opacity=1 ]   (376.79,101.97) .. controls (390.44,92.88) and (393.39,145.36) .. (391.55,161.19) ;
\draw [color={rgb, 255:red, 208; green, 2; blue, 27 }  ,draw opacity=1 ] [dash pattern={on 4.5pt off 4.5pt}]  (230.29,140.67) -- (372.36,140.67) ;
\draw [color={rgb, 255:red, 43; green, 0; blue, 255 }  ,draw opacity=1 ] [dash pattern={on 4.5pt off 4.5pt}]  (372.36,140.67) .. controls (373.83,123.96) and (372.36,111.06) .. (376.79,101.97) ;
\draw [color={rgb, 255:red, 144; green, 19; blue, 254 }  ,draw opacity=1 ] [dash pattern={on 4.5pt off 4.5pt}]  (372.36,140.67) -- (391.55,161.19) ;
\draw [color={rgb, 255:red, 43; green, 0; blue, 255 }  ,draw opacity=1 ]   (228.08,62.1) .. controls (228.17,36.07) and (231.4,24.01) .. (235.29,22.78) .. controls (241.38,20.85) and (249.11,45.4) .. (249.11,84.38) ;
\draw [color={rgb, 255:red, 208; green, 2; blue, 27 }  ,draw opacity=1 ]   (249.11,84.38) -- (389.33,82.62) ;
\draw    (235.29,22.78) .. controls (264.81,5.19) and (345.05,40.99) .. (374.57,23.4) ;
\draw [color={rgb, 255:red, 43; green, 0; blue, 255 }  ,draw opacity=1 ]   (374.57,23.4) .. controls (388.23,14.31) and (391.18,66.79) .. (389.33,82.62) ;
\draw [color={rgb, 255:red, 208; green, 2; blue, 27 }  ,draw opacity=1 ] [dash pattern={on 4.5pt off 4.5pt}]  (228.08,62.1) -- (370.14,62.1) ;
\draw [color={rgb, 255:red, 43; green, 0; blue, 255 }  ,draw opacity=1 ] [dash pattern={on 4.5pt off 4.5pt}]  (370.14,62.1) .. controls (371.62,45.39) and (370.14,32.49) .. (374.57,23.4) ;
\draw [color={rgb, 255:red, 144; green, 19; blue, 254 }  ,draw opacity=1 ][fill={rgb, 255:red, 208; green, 2; blue, 27 }  ,fill opacity=1 ]   (233.61,186.11) -- (254.65,208.4) ;
\draw [color={rgb, 255:red, 208; green, 2; blue, 27 }  ,draw opacity=1 ][fill={rgb, 255:red, 208; green, 2; blue, 27 }  ,fill opacity=1 ]   (254.65,208.4) -- (394.87,206.64) ;
\draw [color={rgb, 255:red, 208; green, 2; blue, 27 }  ,draw opacity=1 ][fill={rgb, 255:red, 208; green, 2; blue, 27 }  ,fill opacity=1 ]   (233.61,186.11) -- (375.68,186.11) ;
\draw [color={rgb, 255:red, 144; green, 19; blue, 254 }  ,draw opacity=1 ][fill={rgb, 255:red, 208; green, 2; blue, 27 }  ,fill opacity=1 ]   (375.68,186.11) -- (394.87,206.64) ;
\draw [color={rgb, 255:red, 43; green, 0; blue, 255 }  ,draw opacity=1 ]   (195.37,142.43) .. controls (195.47,116.4) and (198.69,104.34) .. (202.58,103.11) .. controls (208.67,101.18) and (216.4,125.73) .. (216.4,164.71) ;
\draw [color={rgb, 255:red, 144; green, 19; blue, 254 }  ,draw opacity=1 ]   (195.37,142.43) -- (216.4,164.71) ;
\draw [color={rgb, 255:red, 43; green, 0; blue, 255 }  ,draw opacity=1 ]   (193.52,62.1) .. controls (193.62,36.07) and (196.84,24.01) .. (200.73,22.78) .. controls (206.82,20.85) and (214.56,45.4) .. (214.56,84.38) ;
\draw [color={rgb, 255:red, 144; green, 19; blue, 254 }  ,draw opacity=1 ]   (195.37,185.82) -- (216.4,208.1) ;
\draw [color={rgb, 255:red, 43; green, 0; blue, 255 }  ,draw opacity=1 ]   (415.88,102.85) .. controls (429.53,93.76) and (432.48,146.24) .. (430.64,162.07) ;
\draw [color={rgb, 255:red, 43; green, 0; blue, 255 }  ,draw opacity=1 ]   (411.45,141.55) .. controls (412.92,124.84) and (411.45,111.94) .. (415.88,102.85) ;
\draw [color={rgb, 255:red, 144; green, 19; blue, 254 }  ,draw opacity=1 ]   (411.45,141.55) -- (430.64,162.07) ;
\draw [color={rgb, 255:red, 43; green, 0; blue, 255 }  ,draw opacity=1 ]   (411.08,23.11) .. controls (424.73,14.02) and (427.68,66.5) .. (425.84,82.33) ;
\draw [color={rgb, 255:red, 43; green, 0; blue, 255 }  ,draw opacity=1 ]   (406.65,61.81) .. controls (408.13,45.1) and (406.65,32.2) .. (411.08,23.11) ;
\draw [color={rgb, 255:red, 144; green, 19; blue, 254 }  ,draw opacity=1 ]   (414.4,186.99) -- (433.59,207.52) ;
\draw  [draw opacity=0][fill={rgb, 255:red, 233; green, 89; blue, 89 }  ,fill opacity=0.4 ] (394.87,206.64) -- (254.65,208.4) -- (233.61,186.11) -- (375.68,186.11) -- cycle ;
\draw  [draw opacity=0][fill={rgb, 255:red, 233; green, 89; blue, 89 }  ,fill opacity=0.4 ] (391.55,161.19) -- (251.32,162.95) -- (230.29,140.67) -- (372.36,140.67) -- cycle ;
\draw  [draw opacity=0][fill={rgb, 255:red, 43; green, 0; blue, 255 }  ,fill opacity=0.18 ] (216.4,164.71) -- (195.37,142.43) -- (196.66,120.88) -- (197.58,114.29) -- (199.43,107.98) -- (200.72,104.76) -- (202.58,103.11) -- (205.52,104.02) -- (209.39,110.77) -- (213.26,126.16) -- (215.85,140.97) -- cycle ;
\draw  [draw opacity=0][fill={rgb, 255:red, 43; green, 0; blue, 255 }  ,fill opacity=0.18 ] (251.32,162.95) -- (230.29,140.67) -- (231.58,119.12) -- (232.51,112.53) -- (234.35,106.22) -- (235.64,103) -- (237.5,101.35) -- (240.44,102.27) -- (244.31,109.01) -- (248.19,124.4) -- (250.77,139.21) -- cycle ;
\draw  [draw opacity=0][fill={rgb, 255:red, 43; green, 0; blue, 255 }  ,fill opacity=0.18 ] (430.64,162.07) -- (411.45,141.55) -- (412.74,120) -- (413.17,113.28) -- (414.65,106.83) -- (415.88,102.85) -- (418.52,101.55) -- (422.39,103.31) -- (425.9,109.46) -- (429.34,125.28) -- (431.07,139.66) -- cycle ;
\draw  [draw opacity=0][fill={rgb, 255:red, 43; green, 0; blue, 255 }  ,fill opacity=0.18 ] (391.55,161.19) -- (372.36,140.67) -- (373.65,119.12) -- (374.08,112.4) -- (375.56,105.95) -- (376.79,101.97) -- (379.43,100.67) -- (383.31,102.43) -- (386.81,108.59) -- (390.26,124.4) -- (391.98,138.78) -- cycle ;
\draw  [draw opacity=0][fill={rgb, 255:red, 245; green, 166; blue, 35 }  ,fill opacity=0.44 ] (389.95,73.16) -- (389.33,82.62) -- (249.11,84.38) -- (248.99,74.19) -- (247.33,50.62) -- (244.01,34.79) -- (239.95,25.26) -- (236.2,22.16) -- (244.13,19.52) -- (255.57,17.91) -- (272.36,18.65) -- (288.41,20.55) -- (316.08,24.8) -- (336.75,27.44) -- (355.75,28.61) -- (367.01,26.85) -- (374.57,23.4) -- (378.26,22.31) -- (382.14,25.39) -- (385.83,35.06) -- (388.59,50.31) -- cycle ;
\draw  [draw opacity=0][fill={rgb, 255:red, 245; green, 166; blue, 35 }  ,fill opacity=0.44 ] (392.16,151.73) -- (391.55,161.19) -- (251.32,162.95) -- (251.2,152.76) -- (249.54,129.19) -- (246.22,113.36) -- (242.16,103.83) -- (238.41,100.73) -- (246.34,98.1) -- (257.78,96.48) -- (274.57,97.22) -- (290.62,99.12) -- (318.3,103.37) -- (338.96,106.01) -- (357.97,107.18) -- (369.22,105.43) -- (376.79,101.97) -- (380.48,100.88) -- (384.35,103.96) -- (388.04,113.63) -- (390.81,128.88) -- cycle ;
\draw  [draw opacity=0][fill={rgb, 255:red, 245; green, 166; blue, 35 }  ,fill opacity=0.44 ] (371.07,45.62) -- (370.14,62.1) -- (228.08,62.1) -- (228.26,51.19) -- (229,42.54) -- (230.48,32.86) -- (233.06,25.24) -- (236.2,22.16) -- (244.13,19.52) -- (255.57,17.91) -- (272.36,18.65) -- (288.41,20.55) -- (316.08,24.8) -- (336.75,27.44) -- (355.75,28.61) -- (367.01,26.85) -- (374.57,23.4) -- (378.26,22.31) -- (374.57,23.4) -- (373.28,26.27) -- (371.8,33.89) -- cycle ;
\draw  [draw opacity=0][fill={rgb, 255:red, 245; green, 166; blue, 35 }  ,fill opacity=0.44 ] (373.34,124.12) -- (372.42,140.61) -- (230.35,140.61) -- (230.54,129.69) -- (231.28,121.05) -- (232.75,111.37) -- (235.33,103.75) -- (238.47,100.67) -- (246.4,98.03) -- (257.84,96.42) -- (274.63,97.15) -- (290.68,99.06) -- (318.3,104.16) -- (338.59,106.79) -- (351.32,107.82) -- (369.28,105.36) -- (376.85,101.91) -- (380.54,100.82) -- (376.85,101.91) -- (375.56,104.77) -- (374.08,112.4) -- cycle ;

\draw (170.8,122.4) node [anchor=north west][inner sep=0.75pt]  [font=\footnotesize]  {$\textcolor[rgb]{0.17,0,1}{M}$};
\draw (161.04,36.85) node [anchor=north west][inner sep=0.75pt]  [font=\footnotesize]  {$\textcolor[rgb]{0.17,0,1}{\partial ^{h} M}$};
\draw (168.41,195.45) node [anchor=north west][inner sep=0.75pt]  [font=\footnotesize]  {$\textcolor[rgb]{0.17,0,1}{\partial }\textcolor[rgb]{0.17,0,1}{^{v}}\textcolor[rgb]{0.17,0,1}{M}$};
\draw (439.87,113.57) node [anchor=north west][inner sep=0.75pt]  [font=\footnotesize,color={rgb, 255:red, 43; green, 0; blue, 255 }  ,opacity=1 ]  {$N$};
\draw (434.61,40.51) node [anchor=north west][inner sep=0.75pt]  [font=\footnotesize]  {$\textcolor[rgb]{0.17,0,1}{\partial }\textcolor[rgb]{0.17,0,1}{^{h} N}$};
\draw (431.98,182.11) node [anchor=north west][inner sep=0.75pt]  [font=\footnotesize]  {$\textcolor[rgb]{0.17,0,1}{\partial }\textcolor[rgb]{0.17,0,1}{^{v} N}$};
\draw (344.15,9.06) node [anchor=north west][inner sep=0.75pt]  [font=\footnotesize]  {$\textcolor[rgb]{0.74,0.5,0.11}{\partial ^{h} W}$};
\draw (312.28,116.09) node [anchor=north west][inner sep=0.75pt]  [font=\footnotesize]  {$\textcolor[rgb]{0.55,0.34,0.16}{W}$};
\draw (252.65,191.95) node [anchor=north west][inner sep=0.75pt]  [font=\footnotesize]  {$\textcolor[rgb]{0.82,0.01,0.11}{\partial _{1} W\cong \partial ^{v} M\times [ 0,1]}$};

\end{tikzpicture}

    \caption{This is a triad cobordism $(W, \hb W,\partial_1 W, \vb W)$ from $(M,\hb M,\vb M)$ to $(N,\hb N,\vb N).$}
    \label{relative cobordism picture}
\end{figure}

We say that two triad cobordisms $W:M\leadsto N$ and $W':M\leadsto N$ are \textit{diffeomorphic} if they are diffeomorphic as $4$-ads where the diffeomorphism sends $M$ in $W$ to $M$ in $W'$ (and thus $N$ to $N$), not necessarily via the identity. We call the subspaces $\hb W$ and $\vb W$ of a triad $W$ by \textit{horizontal} and \textit{vertical boundary}, respectively. 

\begin{rmk}
    This notion can be seen as a generalization of a cobordism between manifolds with boundary with trivialized boundary (or $r\partial$-cobordisms in the sense of \cite{WallGeometrical}) by taking $\hb M=\hb N=\emptyset.$ Moreover, given a triad cobordism $W$, then $\hb W$ is a $r\partial$-cobordism.
\end{rmk}

For $k\geq 0$, recall the definition of $D^k_+$ from \cref{triad homotopy}. We see this as a triad $(D^k_+,\partial_0D^k_+,\partial_1D^k_+)$ given by $\partial_0D^k$ to be the subset of $\partial D^k_+$ where the first coordinate is $0$. Denote by $\partial_{01}D^k_+=\partial_0D^k_+\cap \partial_1D^k_+.$ For $k=0,$ this is the triad $(\{0\},\emptyset,\emptyset)$.

\begin{defn}[Handles]\label{defn of handles}
    Let $d\geq k$ be non-negative integers and let $W:M\leadsto N$ be a $d$-dimensional triad cobordism.
    \begin{enumerate}[label=(\textit{\roman*})]
        \item Let $f:S^{k-1}\times D^{d-k}\to N$ be an embedding disjoint from $\partial N$. Define the $4$-ad $(W',\hb W',\partial_1 W',\vb W')$ given by
        \begin{align*}
            W'&\coloneqq W\cup_{f} D^k\times D^{d-k},
            &\hb W'\coloneqq\hb W, \\  
            \vb W'&\coloneqq\vb W\backslash f(S^{k-1}\times D^{d-k})^\circ\cup D^k\times S^{d-k-1}, &\partial_1 W'\coloneqq\partial_1 W.
        \end{align*}There is a unique triad structure on $N'\coloneqq N\backslash f(S^{k-1}\times D^{d-k})^\circ\cup D^k\times S^{d-k-1}$ making $W'$ into a triad cobordism from $M$ to $N'$. We say $W'$ is obtained by \textit{attaching an interior $k$-handle} to $W$ along $f.$ The subspace $(D^k\times \{0\},\emptyset)\subset (W',\hb W')$ is called the \textit{core} and $(\{0\}\times D^{d-k},\emptyset)$ the \textit{cocore} of this handle.
        \item Let $f:(S^{k-1}\times D^{d-k}_+,S^{k-1}\times \partial_0D^{d-k})\to (N,\hb N)$ be an embedding disjoint from $\vb N$. Define the $4$-ad $(W',\hb W',\partial_1 W',\vb W')$ given by
        \begin{align*}
            W'&\coloneqq W\cup_{f} D^k\times D^{d-k}_+,
            &\hb W'\coloneqq\hb W\cup_f D^k\times \partial_0D^{d-k}_+, \\  
            \vb W'&\coloneqq\vb W\backslash f(S^{k-1}\times D^{d-k})^\circ\cup_f D^k\times \partial_1D^{d-k}_+, &\partial_1 W'\coloneqq\partial_1 W.
        \end{align*}There is a unique triad structure on $N'\coloneqq N\backslash f(S^{k-1}\times D^{d-k}_+)^\circ\cup_f D^k\times \partial_1D^{d-k}_+$ making $W'$ into a triad cobordism from $M$ to $N'$. We say $W'$ is obtained by \textit{attaching a right $k$-handle} to $W$ along $f.$ The subspace $D^k\times \{0\}\subset \hb W'$ is called the \textit{core} and $(\{0\}\times D^{d-k}_+,\{0\}\times \partial_0D^{d-k}_+)\subset (W',\hb W')$ the \textit{cocore} of this handle.
        \item Let $f:(\partial_1D^{k}_+\times D^{d-k},\partial_{01}D^k_+\times D^{d-k})\to ( N,\hb N)$ be an embedding disjoint from $\vb N$. Define the $4$-ad $(W',\hb W',\partial_1 W',\vb W')$ given by
        \begin{align*}
            W'&\coloneqq W\cup_{f} D^k_+\times D^{d-k},
            &\hb W'\coloneqq\hb W\cup_f \partial_0D^k_+\times D^{d-k}, \\ 
            \vb W'&\coloneqq\vb W\backslash f(\partial_1D^k_+\times \partial D^{d-k})^\circ\cup D^k_+\times \partial D^{d-k}, &\partial_1 W'\coloneqq\partial_1 W.
        \end{align*}There is a unique triad structure on $N'\coloneqq N\backslash f(\partial_1D^k_+\times D^{d-k})^\circ\cup D^k_+\times \partial D^{d-k}$ making $W'$ into a triad cobordism from $M$ to $N'$. We say $W'$ is obtained by \textit{attaching a left $k$-handle} to $W$ along $f.$ The subspace $(D^k_+\times \{0\},\partial_0D^k_+\times \{0\})\subset (W',\hb W')$ is called the \textit{core} and $\{0\}\times D^{d-k}\subset \hb W'$ the \textit{cocore} of this handle.
    \end{enumerate}A \textit{triad handle decomposition} of a triad cobordism $W:M\leadsto N$ relative to $M$ is a filtration of triples $\cdots \subset (W_i,\hb W_i,\partial_1 W_i)\subset (W_{i+1} ,\hb W_{i+1},\partial_1 W_{i+1})\subset \cdots$ for $i=0,\cdots, k$ for some $k$ where $(W_{i},\hb W_i,\partial_1 W_i,\vb W_i)$ is a triad cobordism from $M$ to some triad $M_i$ satisfying the following properties: \begin{enumerate}
        \item $(W_0,\hb W_0, \partial_{1} W_0, \vb W_0)=(M\times [0,1],\hb M\times [0,1],\vb M\times [0,1], M\times \{0,1\})$
        \item $(W_k,\hb W_k, \partial_1 W_k, \vb W_k)=(W,\hb W,\partial_1 W,\vb W)$
        \item for $i=0,\cdots,k-1$, the cobordism $W_{i+1}$ is diffeomorphic to a cobordism obtained from $W_i$ by attaching an interior, a right or a left handle to $M_i$. This diffeomorphism is assumed to be relative to $M_i.$
    \end{enumerate}Observe that a triad handle decomposition of $W:M\leadsto N$ induces a handle decomposition (in the sense of \cite{WallGeometrical}) of the cobordism $\hb W:\hb  M\leadsto \hb N $ given by $\hb W_i$ for $i=0,\cdots,k.$ More precisely, an interior $k$-handle induces a trivial cobordism on the horizontal boundary, a right $k$-handle induces a $k$-handle and a left $k$-handle induces a $(k-1)$-handle. 
\end{defn}

\begin{figure}
    \centering
    \tikzset{every picture/.style={line width=0.75pt}} 

\begin{tikzpicture}[x=0.75pt,y=0.75pt,yscale=-1,xscale=1]

\draw    (89.88,155.23) .. controls (74.7,160.15) and (68.92,229.52) .. (90.24,224.33) ;
\draw  [dash pattern={on 4.5pt off 4.5pt}]  (89.88,155.23) .. controls (98.91,153.87) and (101.8,225.42) .. (90.24,224.33) ;
\draw    (89.88,155.23) -- (162.15,155.5) ;
\draw    (90.24,224.33) -- (162.51,224.6) ;
\draw  [dash pattern={on 4.5pt off 4.5pt}]  (230.46,198.37) .. controls (231.87,196.78) and (232.21,195.45) .. (231.64,194.36) .. controls (231.07,193.26) and (227.31,191.38) .. (222.1,190.59) ;
\draw    (230.11,181.67) .. controls (234.05,182.54) and (234.05,195.38) .. (230.46,198.37) ;
\draw  [dash pattern={on 4.5pt off 4.5pt}]  (169.28,177.09) .. controls (167.93,176.81) and (165.04,188) .. (170.02,188.82) ;
\draw  [dash pattern={on 4.5pt off 4.5pt}]  (148.61,194.83) .. controls (153.47,197.56) and (151.67,204.94) .. (148.88,204.95) ;
\draw    (96.4,66.47) .. controls (82.64,70.36) and (77.4,125.2) .. (96.73,121.1) ;
\draw  [dash pattern={on 4.5pt off 4.5pt}]  (96.4,66.47) .. controls (104.59,65.39) and (107.21,121.96) .. (96.73,121.1) ;
\draw    (161.91,66.69) .. controls (154.9,68.67) and (150.1,83.88) .. (149.64,97.78) .. controls (149.55,100.53) and (149.63,103.24) .. (149.89,105.78) .. controls (150.91,115.57) and (154.71,122.91) .. (162.24,121.32) ;
\draw    (161.91,66.69) .. controls (164.99,66.28) and (167.28,74) .. (168.38,83.75) .. controls (168.72,86.71) and (168.94,89.86) .. (169.05,93.03) .. controls (169.53,107.22) and (167.59,121.76) .. (162.24,121.32) ;
\draw    (96.4,66.47) -- (161.91,66.69) ;
\draw    (96.73,121.1) -- (162.24,121.32) ;
\draw  [draw opacity=0][fill={rgb, 255:red, 189; green, 52; blue, 27 }  ,fill opacity=1 ] (147.46,102.2) .. controls (147.46,101.54) and (148.28,101) .. (149.28,101) .. controls (150.29,101) and (151.1,101.54) .. (151.1,102.2) .. controls (151.1,102.87) and (150.29,103.4) .. (149.28,103.4) .. controls (148.28,103.4) and (147.46,102.87) .. (147.46,102.2) -- cycle ;
\draw  [draw opacity=0][fill={rgb, 255:red, 189; green, 52; blue, 27 }  ,fill opacity=1 ] (167.12,88.6) .. controls (167.12,87.94) and (167.93,87.4) .. (168.94,87.4) .. controls (169.94,87.4) and (170.76,87.94) .. (170.76,88.6) .. controls (170.76,89.26) and (169.94,89.8) .. (168.94,89.8) .. controls (167.93,89.8) and (167.12,89.26) .. (167.12,88.6) -- cycle ;
\draw    (148.88,204.95) .. controls (149.18,208.16) and (150.94,226.24) .. (162.51,224.6) ;
\draw  [dash pattern={on 4.5pt off 4.5pt}]  (169.78,195.38) .. controls (170.1,201.11) and (169.61,201.33) .. (169.02,205.78) ;
\draw    (169.02,205.78) .. controls (168.43,210.24) and (165.51,224.89) .. (162.51,224.6) ;
\draw    (170.02,188.82) .. controls (170.05,189.74) and (169.46,189.65) .. (169.78,195.38) ;
\draw    (148.88,204.95) .. controls (165.77,206.24) and (221.78,208.21) .. (230.46,198.37) ;
\draw    (222.1,190.59) .. controls (216.88,189.8) and (191.45,189.02) .. (170.02,188.82) ;
\draw    (286.54,155.23) .. controls (271.37,160.15) and (265.58,229.52) .. (286.9,224.33) ;
\draw  [dash pattern={on 4.5pt off 4.5pt}]  (286.54,155.23) .. controls (295.58,153.87) and (298.47,225.42) .. (286.9,224.33) ;
\draw    (286.54,155.23) -- (358.81,155.5) ;
\draw    (286.9,224.33) -- (359.17,224.6) ;
\draw  [dash pattern={on 4.5pt off 4.5pt}]  (427.12,198.37) .. controls (428.53,196.78) and (428.88,195.45) .. (428.31,194.36) .. controls (427.74,193.26) and (423.98,191.38) .. (418.76,190.59) ;
\draw    (426.78,181.67) .. controls (430.72,182.54) and (430.72,195.38) .. (427.12,198.37) ;
\draw    (292.92,68.07) .. controls (279.33,71.96) and (274.15,126.8) .. (293.24,122.7) ;
\draw  [dash pattern={on 4.5pt off 4.5pt}]  (292.92,68.07) .. controls (301,66.99) and (303.59,123.56) .. (293.24,122.7) ;
\draw    (357.62,68.29) .. controls (350.7,70.27) and (345.96,85.48) .. (345.5,99.38) .. controls (345.41,102.13) and (345.49,104.84) .. (345.75,107.38) .. controls (346.76,117.17) and (350.51,124.51) .. (357.95,122.92) ;
\draw    (357.62,68.29) .. controls (360.66,67.88) and (362.92,75.6) .. (364.01,85.35) .. controls (364.34,88.31) and (364.57,91.46) .. (364.67,94.63) .. controls (365.14,108.82) and (363.23,123.36) .. (357.95,122.92) ;
\draw    (292.92,68.07) -- (357.62,68.29) ;
\draw    (293.24,122.7) -- (357.95,122.92) ;
\draw  [draw opacity=0][fill={rgb, 255:red, 189; green, 52; blue, 27 }  ,fill opacity=1 ] (343.35,103.8) .. controls (343.35,103.14) and (344.15,102.6) .. (345.15,102.6) .. controls (346.14,102.6) and (346.95,103.14) .. (346.95,103.8) .. controls (346.95,104.47) and (346.14,105) .. (345.15,105) .. controls (344.15,105) and (343.35,104.47) .. (343.35,103.8) -- cycle ;
\draw  [draw opacity=0][fill={rgb, 255:red, 189; green, 52; blue, 27 }  ,fill opacity=1 ] (362.76,90.2) .. controls (362.76,89.54) and (363.57,89) .. (364.56,89) .. controls (365.55,89) and (366.36,89.54) .. (366.36,90.2) .. controls (366.36,90.86) and (365.55,91.4) .. (364.56,91.4) .. controls (363.57,91.4) and (362.76,90.86) .. (362.76,90.2) -- cycle ;
\draw    (345.55,204.95) .. controls (345.84,208.16) and (347.61,226.24) .. (359.17,224.6) ;
\draw  [dash pattern={on 4.5pt off 4.5pt}]  (366.44,195.38) .. controls (366.76,201.11) and (366.27,201.33) .. (365.69,205.78) ;
\draw    (365.69,205.78) .. controls (365.39,208.03) and (364.5,212.88) .. (363.3,217.09) .. controls (362.13,221.22) and (360.66,224.74) .. (359.17,224.6) ;
\draw  [dash pattern={on 4.5pt off 4.5pt}]  (366.69,188.82) .. controls (366.71,189.74) and (366.13,189.65) .. (366.44,195.38) ;
\draw    (345.55,204.95) .. controls (362.43,206.24) and (418.45,208.21) .. (427.12,198.37) ;
\draw  [dash pattern={on 4.5pt off 4.5pt}]  (418.76,190.59) .. controls (413.55,189.8) and (388.11,189.02) .. (366.69,188.82) ;
\draw    (486.54,152.41) .. controls (471.37,157.33) and (465.58,226.7) .. (486.9,221.51) ;
\draw  [dash pattern={on 4.5pt off 4.5pt}]  (486.54,152.41) .. controls (495.58,151.05) and (498.47,222.6) .. (486.9,221.51) ;
\draw    (486.54,152.41) -- (558.81,152.68) ;
\draw    (486.9,221.51) -- (559.17,221.78) ;
\draw    (489.9,68.05) .. controls (475.66,71.77) and (470.24,124.24) .. (490.24,120.32) ;
\draw  [dash pattern={on 4.5pt off 4.5pt}]  (489.9,68.05) .. controls (498.37,67.02) and (501.08,121.14) .. (490.24,120.32) ;
\draw    (557.7,68.26) .. controls (550.44,70.15) and (545.48,84.7) .. (544.99,98) .. controls (544.9,100.64) and (544.98,103.23) .. (545.26,105.66) .. controls (546.31,115.03) and (550.25,122.05) .. (558.04,120.53) ;
\draw    (557.7,68.26) .. controls (560.88,67.87) and (563.25,75.25) .. (564.39,84.58) .. controls (564.74,87.41) and (564.97,90.43) .. (565.08,93.46) .. controls (565.58,107.03) and (563.58,120.95) .. (558.04,120.53) ;
\draw    (489.9,68.05) -- (557.7,68.26) ;
\draw    (490.24,120.32) -- (558.04,120.53) ;
\draw  [draw opacity=0][fill={rgb, 255:red, 189; green, 52; blue, 27 }  ,fill opacity=1 ] (554.23,107.04) .. controls (554.23,106.4) and (555.08,105.89) .. (556.12,105.89) .. controls (557.16,105.89) and (558,106.4) .. (558,107.04) .. controls (558,107.67) and (557.16,108.19) .. (556.12,108.19) .. controls (555.08,108.19) and (554.23,107.67) .. (554.23,107.04) -- cycle ;
\draw  [draw opacity=0][fill={rgb, 255:red, 189; green, 52; blue, 27 }  ,fill opacity=1 ] (554.17,85.86) .. controls (554.17,85.23) and (555.01,84.71) .. (556.05,84.71) .. controls (557.09,84.71) and (557.94,85.23) .. (557.94,85.86) .. controls (557.94,86.5) and (557.09,87.01) .. (556.05,87.01) .. controls (555.01,87.01) and (554.17,86.5) .. (554.17,85.86) -- cycle ;
\draw    (558.81,152.68) .. controls (562.2,152.17) and (564.45,155.48) .. (565.67,167.81) ;
\draw    (565.92,205.89) .. controls (565.33,210.34) and (562.18,222.07) .. (559.17,221.78) ;
\draw    (566.17,180.72) .. controls (566.19,181.64) and (566.13,186.83) .. (566.44,192.56) ;
\draw [color={rgb, 255:red, 189; green, 52; blue, 27 }  ,draw opacity=1 ]   (345.15,103.8) -- (364.56,90.2) ;
\draw  [draw opacity=0][fill={rgb, 255:red, 189; green, 52; blue, 27 }  ,fill opacity=1 ] (146.2,199.8) .. controls (146.2,198.96) and (147.1,198.28) .. (148.21,198.28) .. controls (149.32,198.28) and (150.22,198.96) .. (150.22,199.8) .. controls (150.22,200.64) and (149.32,201.32) .. (148.21,201.32) .. controls (147.1,201.32) and (146.2,200.64) .. (146.2,199.8) -- cycle ;
\draw  [draw opacity=0][fill={rgb, 255:red, 189; green, 52; blue, 27 }  ,fill opacity=1 ] (166.88,182.06) .. controls (166.88,181.22) and (167.78,180.54) .. (168.89,180.54) .. controls (170,180.54) and (170.89,181.22) .. (170.89,182.06) .. controls (170.89,182.89) and (170,183.57) .. (168.89,183.57) .. controls (167.78,183.57) and (166.88,182.89) .. (166.88,182.06) -- cycle ;
\draw [color={rgb, 255:red, 189; green, 52; blue, 27 }  ,draw opacity=1 ][fill={rgb, 255:red, 189; green, 52; blue, 27 }  ,fill opacity=0.35 ]   (344.88,199.8) .. controls (359.17,200.89) and (401.2,201.31) .. (419.01,195.44) .. controls (422.23,194.38) and (424.66,193.12) .. (425.99,191.61) .. controls (427.84,189.52) and (427.86,187.9) .. (426.38,186.64) .. controls (420.94,182) and (395.09,182.33) .. (365.55,182.06) ;
\draw [color={rgb, 255:red, 189; green, 52; blue, 27 }  ,draw opacity=1 ]   (344.88,199.8) -- (365.55,182.06) ;
\draw  [draw opacity=0][fill={rgb, 255:red, 189; green, 52; blue, 27 }  ,fill opacity=1 ] (554.97,200.47) .. controls (554.97,199.63) and (555.87,198.95) .. (556.97,198.95) .. controls (558.08,198.95) and (558.98,199.63) .. (558.98,200.47) .. controls (558.98,201.31) and (558.08,201.99) .. (556.97,201.99) .. controls (555.87,201.99) and (554.97,201.31) .. (554.97,200.47) -- cycle ;
\draw  [draw opacity=0][fill={rgb, 255:red, 189; green, 52; blue, 27 }  ,fill opacity=1 ] (554.9,174.81) .. controls (554.9,173.97) and (555.8,173.29) .. (556.91,173.29) .. controls (558.01,173.29) and (558.91,173.97) .. (558.91,174.81) .. controls (558.91,175.64) and (558.01,176.32) .. (556.91,176.32) .. controls (555.8,176.32) and (554.9,175.64) .. (554.9,174.81) -- cycle ;
\draw    (558.81,152.68) .. controls (551.08,155.19) and (545.79,174.43) .. (545.27,192.01) .. controls (545.17,195.5) and (545.26,198.92) .. (545.55,202.13) .. controls (546.68,214.52) and (550.87,223.8) .. (559.17,221.78) ;
\draw    (555.28,168.14) .. controls (567.7,167.05) and (588.31,165.79) .. (597.67,169.93) .. controls (607.03,174.07) and (615.42,186.64) .. (603.17,197.21) .. controls (590.92,207.79) and (575.15,204.69) .. (556.92,206.73) ;
\draw    (556.03,181.89) .. controls (588.42,175.22) and (602.67,191.72) .. (556.42,193.62) ;
\draw [color={rgb, 255:red, 74; green, 144; blue, 226 }  ,draw opacity=1 ]   (229.72,186.64) -- (230.46,198.37) ;
\draw [fill={rgb, 255:red, 248; green, 207; blue, 28 }  ,fill opacity=0.25 ]   (345.27,194.83) .. controls (362.16,196.12) and (417.71,196.47) .. (426.38,186.64) .. controls (428.23,184.55) and (428.25,182.93) .. (426.78,181.67) .. controls (421.34,177.03) and (395.49,177.36) .. (365.95,177.09) ;
\draw [color={rgb, 255:red, 74; green, 144; blue, 226 }  ,draw opacity=1 ]   (426.38,186.64) -- (427.12,198.37) ;
\draw [color={rgb, 255:red, 74; green, 144; blue, 226 }  ,draw opacity=1 ][fill={rgb, 255:red, 74; green, 144; blue, 226 }  ,fill opacity=0.38 ]   (229.72,186.64) .. controls (221.93,190.95) and (221.93,198.84) .. (230.46,198.37) ;
\draw [color={rgb, 255:red, 189; green, 52; blue, 27 }  ,draw opacity=1 ]   (148.21,199.8) .. controls (165.09,201.09) and (220.65,201.44) .. (229.32,191.61) .. controls (231.17,189.52) and (231.19,187.9) .. (229.72,186.64) .. controls (224.28,182) and (198.43,182.33) .. (168.89,182.06) ;
\draw    (148.61,194.83) .. controls (165.49,196.12) and (221.05,196.47) .. (229.72,186.64) .. controls (231.57,184.55) and (231.59,182.93) .. (230.11,181.67) .. controls (224.67,177.03) and (198.82,177.36) .. (169.28,177.09) ;
\draw  [color={rgb, 255:red, 74; green, 144; blue, 226 }  ,draw opacity=1 ][fill={rgb, 255:red, 74; green, 144; blue, 226 }  ,fill opacity=0.32 ] (586.37,185.03) .. controls (586.37,183.31) and (591.42,181.92) .. (597.65,181.92) .. controls (603.88,181.92) and (608.93,183.31) .. (608.93,185.03) .. controls (608.93,186.74) and (603.88,188.13) .. (597.65,188.13) .. controls (591.42,188.13) and (586.37,186.74) .. (586.37,185.03) -- cycle ;
\draw [color={rgb, 255:red, 189; green, 52; blue, 27 }  ,draw opacity=1 ]   (558.91,174.81) .. controls (571.34,173.71) and (587.93,173.95) .. (593.68,176.7) .. controls (599.43,179.45) and (599.93,188.12) .. (596.18,191.5) .. controls (592.43,194.89) and (577.18,201.66) .. (558.98,200.47) ;
\draw [fill={rgb, 255:red, 248; green, 207; blue, 28 }  ,fill opacity=0.25 ]   (365.95,177.09) .. controls (365.26,168.92) and (364.59,156.51) .. (358.81,155.5) .. controls (353.03,154.5) and (343.59,181.89) .. (345.27,194.83) ;
\draw    (169.28,177.09) .. controls (168.59,168.92) and (167.93,156.51) .. (162.15,155.5) .. controls (156.36,154.5) and (146.93,181.89) .. (148.61,194.83) ;
\draw  [draw opacity=0][fill={rgb, 255:red, 248; green, 207; blue, 28 }  ,fill opacity=0.25 ] (366.69,188.82) -- (366.14,207.18) -- (363.3,217.09) -- (359.17,224.6) -- (350.97,222.41) -- (347.64,213.95) -- (345.55,204.95) -- cycle ;
\draw  [draw opacity=0][fill={rgb, 255:red, 248; green, 207; blue, 28 }  ,fill opacity=0.25 ] (366.69,188.82) -- (421.14,191.32) -- (428.31,194.36) -- (427.12,198.37) -- (424.14,201.05) -- (408.14,204.01) -- (393.89,205.7) -- (373.64,205.91) -- (345.55,204.95) -- cycle ;
\draw  [draw opacity=0][fill={rgb, 255:red, 248; green, 207; blue, 28 }  ,fill opacity=0.25 ] (549.78,216.7) -- (546.03,205.28) -- (545.53,186.45) -- (550.03,164.45) -- (554.28,155.99) -- (558.81,152.68) -- (563.03,155.57) -- (565.67,167.81) -- (566.44,192.56) -- (587.28,204.43) -- (565.92,205.89) -- (562.78,216.7) -- (559.17,221.78) -- (554.03,221.35) -- cycle ;
\draw  [draw opacity=0][fill={rgb, 255:red, 248; green, 207; blue, 28 }  ,fill opacity=0.25 ] (597.67,169.93) -- (605.28,174.82) -- (609.78,186.03) -- (606.28,194.49) -- (596.78,200.83) -- (587.28,204.43) -- (566.44,192.56) -- (578.53,190.89) -- (586.37,185.03) -- (581.03,181.16) -- (566.06,180.19) -- (565.67,167.81) -- (584.78,167.2) -- cycle ;
\draw  [draw opacity=0][fill={rgb, 255:red, 248; green, 207; blue, 28 }  ,fill opacity=0.25 ] (151.22,202.1) -- (151.97,198.3) -- (148.61,194.83) -- (149.47,184.12) -- (153.11,168.47) -- (157.47,159.17) -- (162.15,155.5) -- (166.47,159.17) -- (169.22,174.18) -- (169.28,177.09) -- (199.47,177.14) -- (220.64,178.63) -- (227.97,181.16) -- (231.72,183.28) -- (224.72,190.47) -- (210.39,189.2) -- (186.72,188.57) -- (170.02,188.82) -- (169.72,193.01) -- (169.62,194.92) -- (148.61,194.83) -- (151.97,198.3) -- (151.22,202.1) -- (148.88,204.95) -- (169.02,205.78) -- (166.22,219.24) -- (162.51,224.6) -- (159.22,224.73) -- (154.97,221.99) -- (151.72,216.06) -- (148.88,204.95) -- cycle ;

\draw (106,28) node [anchor=north west][inner sep=0.75pt]   [align=left] {Right};
\draw (304.8,28.2) node [anchor=north west][inner sep=0.75pt]   [align=left] {Left};
\draw (493.2,28.2) node [anchor=north west][inner sep=0.75pt]   [align=left] {Interior};

\end{tikzpicture}

    \caption{Here are depictions of the three types of handle attachments of \cref{defn of handles} to a trivial triad cobordism for $d=3$. The red portions in the top pictures represent the attaching maps. The red and blue portions in the bottom pictures represent the core and cocores, respectively. The manifold $N'$ (that is, the outgoing boundary of the cobordism $W'$ is depicted in light yellow.}
    \label{pictures of handle}
\end{figure}

We move now to establishing some properties of these objects. These properties were observed in \cite{BNR} for the special case that $\partial_1 W=\emptyset$ and proved using Morse theory. We will deduce our results from theirs by adapting some of the arguments to the more general case when $\partial_1 W\neq \emptyset$. We recall now the relevant notions in loc.cit for the reader's convinience. We assume however the reader's familiarity with this reference for the proofs, as we will only explain the necessary Morse-theoretic modifications.

\begin{rmk}[Morse theory for manifolds with boundary]\label{compare to bnr}
    In \cite{BNR}, the authors study \textit{cobordisms of manifolds with boundary} which are simply triad cobordisms $W$ in the sense of \cref{relative cobordism} where $\partial_1W=\emptyset.$ To avoid confusion, we will use the notation $(\Omega, Y)$ for a cobordism between manifolds with boundary $\Sigma_1$ and $\Sigma_2$ from \cite[Definition 1.1]{BNR}, and $(W,\hb W,\partial_1,\vb W)$ for triad cobordisms. Here are the necessary definitions and facts from \cite{BNR}:
    \begin{enumerate}
        \item A \textit{Morse function} of a cobordism of manifolds with boundary $(\Omega,Y)$ is a smooth function $f:\Omega\to [0,1]$ such that $f(\Sigma_0)=0$ and $f(\Sigma_1)=1$, whose critical points have non degenerate Hessian and are not in $\Sigma_0\cup \Sigma_1$, and $\nabla f$ is everywhere tangent to $Y$. See \cite[Defn. 1.4]{BNR}.
        \item If $(\Omega,Y)$ admits a Morse function with a single critical point $c$, then it is diffeomorphic to a handle attachment in the sense of \cref{defn of handles} relative to $\Sigma_0$ of:
        \begin{enumerate}
            \item Interior type if $c$ is in the interior of $\Omega$,
            \item Left type if $c\in Y$ is \textit{boundary stable} in the sense of \cite[Defn 2.4]{BNR},
            \item Right type if $c\in Y$ is \textit{boundary unstable} in the sense of loc.cit.
        \end{enumerate}This is \cite[Theorem 2.27]{BNR}, by observing that the definitions of handles in \cite[Defs. 2.11 and 2.12]{BNR} are special cases of \cref{defn of handles} for $\partial_1 W=\emptyset$.
    \end{enumerate}To import results from \cite{BNR} to our setting, we see a triad cobordism $(W,\hb W,\partial_1 W,\vb W)$ as the cobordism $(\Omega,Y)=(W,\hb W\cup \partial_1 W)$ between the manifolds with boundary $M$ and $N$. 
\end{rmk}

\begin{prop}[Existence]\label{existence of handle decomposition}
    Let $W:M\leadsto N$ be a triad cobordism, then there exists a triad handle decomposition of $W$ relative to $M.$ In fact, such decomposition can be found where either right or left handles are not present.
\end{prop}
\begin{proof}
    
    By combining \cite[Lemma 2.10]{BNR} with \cite[Thm. 2.27]{BNR}, we see that a Morse function on a cobordism of manifolds with boundary $(\Omega,Y)$ from $\Sigma_0$ to $\Sigma_1$ (in the sense of \cite{BNR}) gives rise to a filtration of the pair $(\Omega,Y)$ where at each step is obtained from the previous by attaching a handle of one of types above (where we see $\Omega$ as a $4$-ad by taking $\partial_1 \Omega=\emptyset$). Moreover, these handles are attached along subsets of the stable and unstable manifolds of each critical point. We now specify to $(\Omega,Y)=(W,\hb W\cup \partial_1 W)$. Since $\partial_1 W$ is diffeomorphic to a product $\vb M\times [0,1]$, we see that there exists a Morse function on $W$ without critical points in $\partial_1 W$, as one can find a classical Morse function on $Y=\hb W\cup \partial_1W$ without critical points in $\partial_1 W$ and then extend it to a Morse function on $W$ by \cite[Lemma 2.1]{BNR}. Moreover, by definition, the gradient flow (for some metric on $W$) of this Morse function is everywhere tangent to $\partial_1 W$. We see that there exists a collar of $\partial_1 W$ which is disjoint from the stable and unstable manifolds of all critical points. Thus, the induced filtration of pairs $(W,\hb W)$ induces the trivial filtration on $\partial_1 W$ and thus corresponds to a triad handle decomposition. We conclude that triad handle decompositions exist. The extra conclusion follows from the proof of \cite[Lemma 2.1]{BNR}, where a Morse function is constructed where the boundary critical points are solely of one type.
\end{proof}

\begin{figure}
    \centering
    \tikzset{every picture/.style={line width=0.75pt}} 

\begin{tikzpicture}[x=0.75pt,y=0.75pt,yscale=-1,xscale=1]

\draw    (252.27,44.74) .. controls (236.73,37.67) and (222.61,74.57) .. (224.46,107.05) .. controls (226.31,139.52) and (233.6,149.48) .. (234.36,150.3) ;
\draw  [dash pattern={on 4.5pt off 4.5pt}]  (252.27,44.74) .. controls (257.35,45.51) and (264.71,67.21) .. (268.24,91.5) .. controls (271.76,115.79) and (270.64,113.95) .. (270.11,124.27) .. controls (269.58,134.6) and (269.51,131.98) .. (268.89,135.49) ;
\draw    (252.27,44.74) .. controls (263.84,46.37) and (301.84,41.87) .. (340.88,31.4) .. controls (379.93,20.92) and (420.89,14.17) .. (440.77,15.03) ;
\draw    (440.77,15.03) .. controls (425.44,19.76) and (414.16,76.15) .. (441.99,75.72) ;
\draw    (440.77,15.03) .. controls (442.08,14.59) and (443.31,14.62) .. (444.46,15.04) .. controls (456.73,19.62) and (459.14,69.81) .. (441.99,75.72) ;
\draw    (440.29,119.61) .. controls (434.01,122.54) and (428.04,136.28) .. (425.39,151.08) ;
\draw    (440.29,119.61) .. controls (445.74,117.77) and (449.68,125.66) .. (451.66,136.88) ;
\draw [color={rgb, 255:red, 74; green, 144; blue, 226 }  ,draw opacity=1 ]   (441.99,75.72) .. controls (435.74,80.02) and (430.63,77.87) .. (407.96,79.73) .. controls (397.52,80.58) and (388.95,90.62) .. (389.5,100.49) .. controls (390.13,112.05) and (403.21,123.39) .. (440.29,119.61) ;
\draw [color={rgb, 255:red, 208; green, 2; blue, 27 }  ,draw opacity=0.67 ]   (224.46,107.05) -- (268.24,91.5) ;
\draw [color={rgb, 255:red, 144; green, 19; blue, 254 }  ,draw opacity=1 ]   (224.46,107.05) .. controls (295.83,114.15) and (383.05,111.77) .. (389.5,100.49) ;
\draw [color={rgb, 255:red, 144; green, 19; blue, 254 }  ,draw opacity=1 ] [dash pattern={on 4.5pt off 4.5pt}]  (268.24,91.5) .. controls (319.18,92.68) and (390.06,92.68) .. (389.5,100.49) ;
\draw  [draw opacity=0][fill={rgb, 255:red, 208; green, 2; blue, 27 }  ,fill opacity=0.33 ] (224.46,107.05) -- (268.24,91.5) -- (340.12,93.18) -- (367.56,95.04) -- (384.81,97.32) -- (389.5,100.49) -- (386.68,102.88) -- (372.43,107.45) -- (351.81,109.74) -- (324.81,110.97) -- (286.18,110.97) -- (258.06,109.67) -- cycle ;
\draw [color={rgb, 255:red, 233; green, 89; blue, 89 }  ,draw opacity=1 ]   (81.46,72.47) -- (141.74,51.81) ;
\draw [color={rgb, 255:red, 155; green, 155; blue, 155 }  ,draw opacity=1 ]   (122.26,67.68) .. controls (146.89,92.3) and (181.28,92.45) .. (210.72,96.63) ;
\draw [shift={(212.51,96.89)}, rotate = 188.48] [color={rgb, 255:red, 155; green, 155; blue, 155 }  ,draw opacity=1 ][line width=0.75]    (10.93,-3.29) .. controls (6.95,-1.4) and (3.31,-0.3) .. (0,0) .. controls (3.31,0.3) and (6.95,1.4) .. (10.93,3.29)   ;
\draw  [draw opacity=0][fill={rgb, 255:red, 144; green, 19; blue, 254 }  ,fill opacity=1 ] (79.46,72.52) .. controls (79.46,71.67) and (80.35,70.98) .. (81.45,70.98) .. controls (82.55,70.98) and (83.44,71.67) .. (83.44,72.52) .. controls (83.44,73.38) and (82.55,74.07) .. (81.45,74.07) .. controls (80.35,74.07) and (79.46,73.38) .. (79.46,72.52) -- cycle ;
\draw  [draw opacity=0][fill={rgb, 255:red, 144; green, 19; blue, 254 }  ,fill opacity=1 ] (139.1,51.89) .. controls (139.1,51.03) and (139.99,50.34) .. (141.09,50.34) .. controls (142.19,50.34) and (143.08,51.03) .. (143.08,51.89) .. controls (143.08,52.74) and (142.19,53.43) .. (141.09,53.43) .. controls (139.99,53.43) and (139.1,52.74) .. (139.1,51.89) -- cycle ;
\draw  [draw opacity=0][fill={rgb, 255:red, 144; green, 19; blue, 254 }  ,fill opacity=1 ] (222.6,106.65) .. controls (222.6,105.8) and (223.49,105.11) .. (224.59,105.11) .. controls (225.69,105.11) and (226.58,105.8) .. (226.58,106.65) .. controls (226.58,107.5) and (225.69,108.2) .. (224.59,108.2) .. controls (223.49,108.2) and (222.6,107.5) .. (222.6,106.65) -- cycle ;
\draw  [draw opacity=0][fill={rgb, 255:red, 144; green, 19; blue, 254 }  ,fill opacity=1 ] (266.34,91.57) .. controls (266.34,90.72) and (267.23,90.03) .. (268.33,90.03) .. controls (269.42,90.03) and (270.31,90.72) .. (270.31,91.57) .. controls (270.31,92.43) and (269.42,93.12) .. (268.33,93.12) .. controls (267.23,93.12) and (266.34,92.43) .. (266.34,91.57) -- cycle ;
\draw  [draw opacity=0][fill={rgb, 255:red, 43; green, 0; blue, 255 }  ,fill opacity=1 ] (440.29,75.7) .. controls (440.29,74.84) and (441.18,74.15) .. (442.28,74.15) .. controls (443.38,74.15) and (444.27,74.84) .. (444.27,75.7) .. controls (444.27,76.55) and (443.38,77.24) .. (442.28,77.24) .. controls (441.18,77.24) and (440.29,76.55) .. (440.29,75.7) -- cycle ;
\draw  [draw opacity=0][fill={rgb, 255:red, 43; green, 0; blue, 255 }  ,fill opacity=1 ] (439.63,119.35) .. controls (439.63,118.5) and (440.52,117.8) .. (441.62,117.8) .. controls (442.71,117.8) and (443.6,118.5) .. (443.6,119.35) .. controls (443.6,120.2) and (442.71,120.9) .. (441.62,120.9) .. controls (440.52,120.9) and (439.63,120.2) .. (439.63,119.35) -- cycle ;
\draw  [draw opacity=0][fill={rgb, 255:red, 43; green, 0; blue, 255 }  ,fill opacity=1 ] (556.92,76.49) .. controls (556.92,75.64) and (557.81,74.95) .. (558.91,74.95) .. controls (560.01,74.95) and (560.9,75.64) .. (560.9,76.49) .. controls (560.9,77.35) and (560.01,78.04) .. (558.91,78.04) .. controls (557.81,78.04) and (556.92,77.35) .. (556.92,76.49) -- cycle ;
\draw  [draw opacity=0][fill={rgb, 255:red, 43; green, 0; blue, 255 }  ,fill opacity=1 ] (556.92,119.35) .. controls (556.92,118.5) and (557.81,117.8) .. (558.91,117.8) .. controls (560.01,117.8) and (560.9,118.5) .. (560.9,119.35) .. controls (560.9,120.2) and (560.01,120.9) .. (558.91,120.9) .. controls (557.81,120.9) and (556.92,120.2) .. (556.92,119.35) -- cycle ;
\draw [color={rgb, 255:red, 155; green, 155; blue, 155 }  ,draw opacity=1 ]   (542.94,96.68) -- (525.73,96.74) -- (468.07,96.94) ;
\draw [shift={(466.07,96.95)}, rotate = 359.8] [color={rgb, 255:red, 155; green, 155; blue, 155 }  ,draw opacity=1 ][line width=0.75]    (10.93,-3.29) .. controls (6.95,-1.4) and (3.31,-0.3) .. (0,0) .. controls (3.31,0.3) and (6.95,1.4) .. (10.93,3.29)   ;
\draw  [dash pattern={on 4.5pt off 4.5pt}]  (234.36,150.3) -- (268.89,135.49) ;
\draw    (425.39,151.08) -- (451.66,136.88) ;
\draw    (234.36,150.3) -- (425.39,151.08) ;
\draw  [dash pattern={on 4.5pt off 4.5pt}]  (268.89,135.49) -- (451.66,136.88) ;
\draw  [draw opacity=0][fill={rgb, 255:red, 233; green, 89; blue, 89 }  ,fill opacity=0.1 ] (451.66,136.88) -- (425.39,151.08) -- (234.36,150.3) -- (268.89,135.49) -- cycle ;
\draw  [draw opacity=0][fill={rgb, 255:red, 43; green, 0; blue, 255 }  ,fill opacity=0.12 ] (268.89,135.49) -- (234.36,150.3) -- (228.23,134.77) -- (225.75,118.08) -- (224.51,100.51) -- (226.25,79.42) -- (230.72,64.48) -- (236.43,52.4) -- (242.15,46.03) -- (247.87,43.62) -- (254.58,45.59) -- (261.53,60.09) -- (266.5,78.98) -- (270.48,112.15) -- cycle ;
\draw  [draw opacity=0][fill={rgb, 255:red, 43; green, 0; blue, 255 }  ,fill opacity=0.11 ] (449.15,69.98) -- (442.28,75.7) -- (435.74,74.15) -- (430.1,70.85) -- (425.13,58.77) -- (425.38,46.91) -- (427.62,35.05) -- (431.59,25.17) -- (435.57,18.8) -- (440.77,15.03) -- (445.26,15.28) -- (449.98,20.99) -- (453.71,35.71) -- (453.71,51.31) -- (452.22,59.65) -- cycle ;
\draw  [draw opacity=0][fill={rgb, 255:red, 43; green, 0; blue, 255 }  ,fill opacity=0.12 ] (451.66,136.88) -- (425.39,151.08) -- (430.85,131.92) -- (437.31,121.59) -- (441.62,119.35) -- (446.01,121.59) -- (449.98,129.72) -- cycle ;

\draw (159.01,65.51) node [anchor=north west][inner sep=0.75pt]    {$\textcolor[rgb]{0.61,0.61,0.61}{f|}\textcolor[rgb]{0.61,0.61,0.61}{_{\partial _{1} D_{+}^{k}}}$};
\draw (50.73,85.43) node [anchor=north west][inner sep=0.75pt]  [font=\footnotesize]  {$\left(\textcolor[rgb]{0.91,0.35,0.35}{\partial }\textcolor[rgb]{0.91,0.35,0.35}{_{1}}\textcolor[rgb]{0.91,0.35,0.35}{D}\textcolor[rgb]{0.91,0.35,0.35}{_{+}^{k}} ,\textcolor[rgb]{0.56,0.07,1}{\partial }\textcolor[rgb]{0.56,0.07,1}{_{01}}\textcolor[rgb]{0.56,0.07,1}{D}\textcolor[rgb]{0.56,0.07,1}{_{+}^{k}}\right)$};
\draw (318.59,111.51) node [anchor=north west][inner sep=0.75pt]  [font=\footnotesize]  {$\left(\textcolor[rgb]{0.82,0.01,0.11}{D}\textcolor[rgb]{0.82,0.01,0.11}{_{+}^{k}} ,\textcolor[rgb]{0.56,0.07,1}{\partial }\textcolor[rgb]{0.56,0.07,1}{_{0}}\textcolor[rgb]{0.56,0.07,1}{D}\textcolor[rgb]{0.56,0.07,1}{_{+}^{k}}\right)$};
\draw (497.29,76.4) node [anchor=north west][inner sep=0.75pt]    {$\textcolor[rgb]{0.61,0.61,0.61}{f'|}\textcolor[rgb]{0.61,0.61,0.61}{_{\partial D^{d-k}}}$};
\draw (543.67,134.06) node [anchor=north west][inner sep=0.75pt]  [font=\footnotesize,color={rgb, 255:red, 43; green, 0; blue, 255 }  ,opacity=1 ]  {$\textcolor[rgb]{0.17,0,1}{\partial D}\textcolor[rgb]{0.17,0,1}{^{d-k}}$};
\draw (383.65,60.78) node [anchor=north west][inner sep=0.75pt]  [font=\footnotesize,color={rgb, 255:red, 43; green, 0; blue, 255 }  ,opacity=1 ]  {$\textcolor[rgb]{0.29,0.56,0.89}{D}\textcolor[rgb]{0.29,0.56,0.89}{^{d-k}}$};
\draw (391.61,137.19) node [anchor=north west][inner sep=0.75pt]  [font=\footnotesize]  {$\textcolor[rgb]{0.91,0.35,0.35}{\partial _{1} W}$};
\draw (238.38,54.65) node [anchor=north west][inner sep=0.75pt]  [font=\footnotesize]  {$\textcolor[rgb]{0.17,0,1}{M}$};
\draw (434.81,31.49) node [anchor=north west][inner sep=0.75pt]  [font=\footnotesize]  {$\textcolor[rgb]{0.17,0,1}{N}$};

\end{tikzpicture}

    \caption{This is a triad cobordism $W:M\leadsto N$ obtained from the trivial triad cobordism $(M\times [0,1],\hb M\times [0,1], \vb M\times [0,1], M\times \{0,1\})$ by attaching a left $k$-handle to $M\times \{1\}$. On the other hand, it is obtained from the trivial cobordism $N\times [0,1]$ by attaching a right $(d-k)$-handle to $N\times \{0\}$ (compare with \cref{dual handle}). We also depict the attaching map $f$, belt sphere $f'$ (see \cref{elementary cobordisms}), the core $(D^k_+,\partial_0D^k_+)$ and cocore $D^{d-k}.$}
    \label{handles}
\end{figure}

\begin{rmk}[Triad handle decompositions from classical]\label{extending handle decompositions}
Let $(M,\hb M)$ be a manifold pair and $W':\hb M\leadsto \hb N$ be a cobordism with trivial boundary (also called $r\partial$-cobordism in \cite{WallGeometrical}) and let $V:\hvb M\leadsto \hvb N$ be the trivial cobordism between the boundaries. Observe that the $4$-ad $W\coloneqq ((W'\cup_{\hb M} M)\times [0,1],W'\times \{0\}\cup_{\hb N} \hb N\times [0,1] , (\vb M\cup V)\times [0,1] ,M\times \{0\}\sqcup (W'\cup_{\partial M} M)\times \{1\})$ is a triad cobordism. Additionally, given any handle decomposition on $W'$ relative to $\hb M$ (in the sense of loc.cit.), one can define a triad handle decomposition on $W$ relative to $M$ by attaching right handles just as prescribed by the handle decomposition on $W'.$ In other words, a filtration on $W'$ relative to $\partial M$ induces one on $(W'\cup_{\hb M}M)\times [0,1]$ relative to $M\times [0,1].$ This filtration is obtained by attaching right handles at each stage along (collars of) the attaching maps of $W'.$ Conversely, given any triad cobordism $W:M\leadsto N$ that admits a triad handle decomposition with only right handles is diffeomorphic to $(\hb W\cup_{\hb M} M)\times [0,1].$ One can see this by using the gradient flow of a Morse function, in the sense of \cite{BNR}, with only boundary critical points to retract $W$ into a collar of $\hb W\cup_{\hb M} M$.
    
\end{rmk}

\begin{defn}\label{elementary cobordisms}
    Let $W:M\leadsto N$ be a $d$-dimensional triad cobordism, we call $W$ \textit{elementary relative to $M$} (resp. $N)$ if it admits a triad handle decomposition relative to $(M,\hb M)$ (resp. $N$) given by a single handle attachment. Such cobordism is called of \textit{interior, left} or \textit{right type} depending on the type of the handle attachment and of index $k$ if the handle attached is a $k$-handle (see \cref{handles}). When dealing with elementary triad cobordisms, we often fix such a triad handle decomposition. In this case, we denote the core and cocore of the handle of $W$ by $\core_W$ and $\cocore_W$ respectively. We denote the intersection $\core_W$ with $M$ by $\attach_W.$ Similarly, we denote the intersection of $\cocore_W$ with $N$ by $\belt_W$, and call it the \textit{belt sphere}.
\end{defn}

 Given the cobordism $W,$ we can define $W^*:N\leadsto M$ given by the same triad seen as a cobordism from $N$ to $M.$ The following lemma follows by unwrapping definitions.

\begin{lemma}\label{dual handle}
Let $k\geq 0$. Let $W:M\leadsto N$ be a $d$-dimensional elementary triad cobordism of left (resp. interior) type of index $k$ relative to $M$. Then $W^*$ is elementary of a right (resp. interior) type and index $(d-k)$ relative to $N$. Moreover, such triad handle decomposition can be arranged so that $\core_{W^*}=\cocore_W$ and $\cocore_{W^*}=\core_W.$
\end{lemma}

\begin{rmk}\label{see submanifolds inside other places}
    In the context above, note that there exists a diffeomorphism $\phi:M\backslash \nu(\attach_W)\cong N\backslash \nu(\belt_W)$, where $\nu(-)$ denotes an open tubular neighborhood. This follows by the definition of handle attachments. Thus, given a submanifold of $M$ which is disjoint from $\attach_W$ we can consider it as a submanifold of $N$ by choosing a small enough tubular neighborhood of $\attach_W$ and taking the image under $\phi$. 
\end{rmk}

Given two triad cobordisms $W:P\leadsto Q$ and $W':Q\leadsto R$, then its \textit{composition} $W'\circ W:P\leadsto R$ is the union $W\cup_Q W'.$ In this case, we call $W$ and $W'$ composable. We recall that, we say that two triad cobordisms $W:M\leadsto N$ and $W':M\leadsto N$ are diffeomorphic if they are diffeomorphic as triads where the diffeomorphism sends $M$ in $W$ to $M$ in $W'$ (and thus $N$ to $N$). Given a triad cobordism $W:P\leadsto Q$, we call $P$ the \textit{ingoing boundary} and $Q$ the \textit{outgoing boundary}.

\begin{prop}[Handle rearrangement]\label{rearragement}
    Let $P\overset{W}{\leadsto} Q\overset{W'}{\leadsto} R$ be two elementary triad cobordisms relative to their ingoing boundaries. If $\belt_W\cap \attach_{W'}=\emptyset$ for some triad handle decomposition of $W$ and $W'$, then there exist two elementary triad cobordisms $\smash{P\overset{M}{\leadsto} Q'\overset{M'}{\leadsto} R}$ such that $M$ (\textit{resp.} $M'$) is of the same type and index as $W'$ (\textit{resp.} $W$) relative to their ingoing boundaries and the composition is diffeomorphic to $W'\circ W$ relative to $P.$ 
\end{prop}
\begin{proof}
    Once again, we deduce this statement by observing that the analogous statement in \cite{BNR} (in this case Thm. $4.1$) can be improved to our setting. By \cite[Prop. 2.35]{BNR}, for any elementary triad cobordism $W$ there exists a Morse function with exactly one critical point. We want to argue that that critical point can be assumed to lie outside $\partial_1 W$. If it lies in the interior of $\partial_1 W$, then the restriction of the Morse function to $\partial_1 W$ is a Morse function. However, since $\partial_1 W$ is a trivial cobordism, it cannot be an elementary one, thus we reach a contradiction. If the critical point lies in the boundary, then using a collar on $\hb W$, we can push the critical point outside of $\partial_1 W$. Thus, $W'\circ W$ has a Morse function with two critical points disjoint from $\partial_1 W$. In \cite[Thm. 4.1]{BNR}, a new Morse function is created where the critical points are swapped. Here the hypothesis are satisfied since $\belt_W\cap \attach_{W'}=\emptyset$ implies that the stable and unstable manifolds of the two critical points are disjoint. Thus, we see that this new Morse function does not have any critical points on $\partial_1(W'\circ W)$ which implies that it gives a triad handle decomposition of $W'\circ W $ given by first attaching a handle of index and type of $W'$ and after one of index and type of $W$. Thus, we obtain the desired result.
\end{proof}

Whenever the conclusion of \cref{rearragement} holds, we say $W$ and $W'$ \textit{can be rearranged}. The following is presented in \cite[Table 1]{BNR} and follows from transversality and dimension counting.

\begin{cor}\label{rearrangement using transversality}
    Let $P\overset{W}{\leadsto} Q\overset{W'}{\leadsto} R$ be two elementary triad cobordisms of index $k$ and $k'$ of any type relative to their ingoing boundaries. Then $W$ and $W'$ can be rearranged if one of the following conditions holds:
    \begin{enumerate}
        \item $W'$ is of interior or left type and $W$ is of any type provided $k\geq k'.$ When $W'$ is of interior type and $W$ is of left type, no condition on $k$ and $k'$ is necessary.
        \item Both cobordisms are of right type and $k\geq k'.$
        \item $W'$ is of right type and $W$ is of left type and $k>k'.$
        \item $W'$ is of right type and $W$ is of interior type.
    \end{enumerate}
\end{cor}
\begin{proof}
    By transversality of embeddings of pairs (see \cite[Lemma 9.2]{genauer}\footnote{Genauer proves that the subspace of smooth maps which are transversal to a fixed smooth map is open and dense in the space of smooth maps of pairs. To conclude that the space of transversal embeddings of pairs is open in the space of smooth maps, it suffices to prove that embeddings of pairs are open in the space of smooth maps. This follows from analogous arguments to \cite[Prop. 2.1.4]{Hirsch1976}.}), we can always isotope the attaching map of the handle of $W'$, resulting in a diffeomorphic cobordism relative to $Q$, so $\belt_W$ and $\attach_{W'}$ intersect transversely. Thus by dimension counting, the hypothesis of \cref{rearragement} is satisfied. The cases where no dimension assumption is needed follow from the fact that either $\belt_W$ is contained in the interior of $Q$ and $\attach_{W'} $ in $\partial Q$, or vice-versa.
\end{proof}

The following propositions constitute the main results of \cite{BNR} and \cite{borodzik2023merging}. Given a cobordism $W:P\leadsto R, $ a factorization of $W$ is a sequence of composable cobordisms $W_i$ such that the composition is diffeomorphic to $W$ relative to $P.$ The following three results follow from \cite[Thm. 3.1]{BNR}, \cite[Thm. 1.2]{borodzik2023merging} and \cite[Thm. 5.1]{BNR}, respectively, by observing that the constructions made in those proofs can be assumed to leave $\partial_1$ untouched (as in \cref{rearragement}) and thus generalize to our context.

\begin{prop}[Splitting]\label{splitting}
    Let $W:P\leadsto R$ be an elementary triad cobordism of interior type of index $0\leq k\leq d-1$ attached to a component of $P$ connected to $\hb P$. Then there exists a factorization $\smash{P\overset{M}{\leadsto} Q\overset{M'}{\leadsto} R}$ of $W$ where $M$ is elementary of left type and $M'$ is elementary of right type and both have index $k.$ Moreover, we have $\belt_W$ and $\attach_{W'}$ intersect transversely at one point.
\end{prop}

\begin{prop}[Merging]\label{merging}
    Let $P\overset{W}{\leadsto} Q\overset{W'}{\leadsto} R$ be two elementary triad cobordisms of index $k$, where $W$ is of left and and $W'$ is of right type relative to their ingoing boundaries. If $\belt_W$ and $\attach_{W'}$ intersect transversely at one point, then $W'\circ W$ is an elementary triad cobordism of interior type and index $k.$
\end{prop}

\begin{prop}[Cancellation]\label{cancellation}
    Let $P\overset{W}{\leadsto} Q\overset{W'}{\leadsto} R$ be two elementary triad cobordisms of the same type and index $k$ and $k+1$ respectively. If $\belt_W$ and $\attach_{W'}$ intersect transversely at one point, then $W'\circ W$ is diffeomorphic to the $(P\times [0,1],\hb P\times [0,1])$ relative to $P.$
\end{prop}

\begin{rmk}\label{cancellation boundary}
    It is important to remark that if $W$ and $W'$ are either of right or of left type, then either $\belt_W$ or $\attach_{W'}$ is completely contained in $\partial Q.$ This implies that the condition in \cref{cancellation} can be checked in the boundary cobordism $\hb W'\circ \hb W.$ In practice, when we wish to cancel a right or left handle, we produce a canceling (classical) handle on the boundary cobordism and then extend it to a triad handle of the same type. 
\end{rmk}

For a triad $(P,\hb P,\vb P)$, we call the cobordism $(P\times [0,1],\hb P\times [0,1],\vb P\times [0,1], P\times \{0\}\sqcup P\times \{1\})$ the \textit{product cobordism}. We say that a cobordism $W:P\leadsto Q$ is \textit{trivial} if it is diffeomorphic to the product cobodism $P\times [0,1]$ relative to $P$, where $P$ is seen in the product cobordism as the inclusion of $P\times \{0\}.$

\begin{prop}[Standard presentation]\label{presentation}
    Let $W:P\leadsto Q$ be a triad cobordism, then $W$ admits a decomposition of the form
    \[W\cong W_0\cup W_0^R\cup W_1^L\cup W_1\cup W_1^R\cup \cdots \cup W_{d-1}^R\cup W_d^L\cup W_d\]where $W_i$ is a composition of elementary interior cobordisms of index $i$ and $W_i^{L/R}$ is a composition of elementary left/right cobordisms of index $i$, relative to their ingoing boundaries. Moreover, for any of the following conditions, there exists such a decomposition satisfying it:
    \begin{enumerate}[label=(\textit{\roman*})]
        \item The cobordisms $W_i^L$ are trivial for all $i\geq 0$.
        \item The cobordisms $W_i^R$ are trivial for all $i\geq 0$.
        \item The cobordisms $W_i$ are trivial for all $0<j<d$.
    \end{enumerate}
\end{prop}
\begin{proof}
    The first result follows from \cref{existence of handle decomposition} and \cref{rearrangement using transversality}. The second result follows from the second claim of \cref{existence of handle decomposition} and from \cref{splitting}. 
\end{proof}

\subsubsection{Strong geometrical connectivity.}\label{strong geo section}

Handle decompositions of cobordisms between closed manifolds are made to be akin to cell structures of CW pairs. Classical algebraic topology tells us that if a map $X\to Y$ of spaces if $k$-connected, then there exists a CW pair $(Z,X)$ made from cells of dimension at least $k+1$ along with a weak equivalence $Z\to Y$ under $X.$ Wall \cite{WallGeometrical} proved a geometric analog of this statement, where a map is replaced by the inclusion of one of the boundary components of a cobordism and cells are replaced by handles. We wish to prove a similar statement for triad cobordisms, triad handle decompositions and  strong connectivity. We start by noticing how different handle attachments behave homotopically. Let $k\geq 0$ be an integer, then the following holds:
\begin{enumerate}[label=(\textsc{\roman*})]
    \item \label{effect of right handles} If $W'$ is a triad cobordism obtained from $W$ by attaching a right $k$-handle, then note that the pair $(W',\hb W')$ is equivalent to $(W\cup_{S^{k-1}}D^k,\hb W'\cup_{S^{k-1}}D^k).$ Clearly, the map $\hb W\to \hb W'$ is $(k-1)$-connected. Moreover, the map $\hb W'\cup_{\hb W}W\to W'$ is an equivalence. This uses the fact that the strict pushout $\hb W'\cup_{\hb W}W$ is the homotopy pushout (since $\hb W\to W$ is a Hurewicz cofibration, see the discussion in the beginning of \cref{strong connectivity section}).  In particular, the map $(W,\hb W)\to (W',\hb W')$ is strongly $(k-1)$-connected. 

    \item \label{effect of left handles} Let $k\geq 1$. If $W'$ is a triad cobordism obtained from $W$ by attaching a left $k$-handle, then note that the pair $(W',\hb W')$ is equivalent to $(W\cup_{\partial_1D^k_+} D^k_+,\hb W'\cup_{\partial_{01}D^k_+}\partial_0D^k_+).$ In particular, the map $W\cup_{\hb W} \hb W'\to W'$ is equivalent to the attachment of a $k$-cell and thus $(k-1)$-connected. Similarly, the map $\hb W\to \hb W'$ is equivalent to the attachment of a $(k-1)$-cell and thus $(k-2)$-connected. Thus the map $(W,\hb W)\to (W',\hb W')$ is strongly $(k-2)$-connected. Additionally, the map $W\to W'$ is an equivalence.
    \item \label{effect of interior handles} If $W'$ is a triad cobordism obtained from $W$ by attaching an interior $k$-handle, then note that the pair $(W',\hb W')$ is equivalent to $(W\cup_{S^{k-1}} D^k,\hb W').$ Since the map on horizontal boundaries is an equivalence, it follows that $(W,\hb W)\to (W',\hb W')$ is strongly $(k-2)$-connected.
\end{enumerate}

We are ready now for the analog of \cite[Thm. 3]{WallGeometrical} for triad cobordisms.

\begin{prop}[Geometrical connectivity]\label{strong geometrical connectivity}
    Let $W:R\leadsto S$ be a triad cobordism of dimension $d$ and $-1\leq k\leq d-5$ be an integer. If the inclusion $(R,\hb R)\hookrightarrow (W,\hb W)$ is strongly $k$-connected, then there exists a triad handle decomposition relative to $(R,\hb R)$ with only right handles of index at least $k+1$, only interior handles of index at least $k+2,$ and no left handles. 
\end{prop}
\begin{proof}
    By \cref{existence of handle decomposition}, there exists a triad handle decomposition of $W$ relative to $R$ consisting of only right and interior handles. By \cref{rearrangement using transversality} and \cref{presentation}, we can assume that any triad cobordism $W:R\leadsto S$ has a triad handle decomposition of the form 
    \[W\cong W_0^R\cup W_1^R\cup \cdots \cup W_{d-1}^R\cup W_0\cup \cdots \cup W_d\]where once again $W_i^R$ is a composite of elementary right cobordisms of index $i$ and $W_i$ a composite of elementary interior cobordisms of index $i$. Denote by $W^R:R\leadsto R'$ and $W^I:R'\leadsto S$ the unions of $W_i^R$ and $W_i$, respectively. Observe that the cobordism $\hb W:\hb R\leadsto \hb S$ is diffeomorphic to $\hb W^R:\hb R\leadsto \hb R'$ relative to $\hb R$, since $\hb W^I$ is a trivial cobordism since $W^I$ is obtained by attaching only interior handles. Moreover, the cobordism $\hb W\cong \hb W^R$ satisfies the hypothesis of \cite[Thm. 3]{WallGeometrical} since it is a classical cobordism between manifolds with trivial boundary (or a $r\partial$-cobordism) and $k\leq (d-1)-4=d-5$. Thus, there exists a handle decomposition of $\hb W$ relative to $\hb R$ whose handles have index at least $k+1$. Given this handle decomposition, we can extend it to a triad handle decomposition of $W^R$ containing only right handles of index at least $k+1$, by \cref{extending handle decompositions}. We proceed now to show that $W^I$ admits a triad handle decomposition given by interior handles of index at least $k+2.$ The cobordism $W^I$ can be viewed as a $r\partial$-cobordism in the sense of \cite{WallGeometrical} (whose boundary cobordism is $\hb W^R\cup \partial_1 W^R$). Hence, a handle decomposition in the sense of loc. cit. induces a triad handle decomposition solely given by interior handles of the same indices as the ones in the handle decomposition. Thus, if we show that the inclusion $R'\to W^I$ is $(k+1)$-connected, applying Wall's result again finishes the claim. By assumption, the map $R\cup_{\hb R} \hb W\to W$ is $(k+1)$-connected. Since $W^R$ is obtained from $R$ by attaching right handles, we conclude from \ref{effect of right handles} above that the maps $R\cup_{\hb R} \hb W^R\to W^R$ and  $R'\to W^R$ are equivalences. Thus it follows that $R'\to W^R\to W$ is a composite of an equivalence with a $(k+1)$-connected map and thus it is $(k+1)$-connected. This finishes the proof.
\end{proof}

\section{Cobordism categories with boundary.}\label{cob cat section}

As mentioned in the introduction, the proof of the main result involves a version of the classical cobordism category (see e.g. \cite{GMTW}) in the context of manifolds with boundary and triad cobordisms between them. In this section, we define the cobordism category $\Cobbt$ of $\Theta$-manifolds with boundary and variations thereof for a map of pairs $\Theta: (B,B^\partial)\to (\BO(d),\BO(d-1))$. The main result of this section (stated below in \cref{surgery on morphisms}) states that the classifying space of $\Cobbt$ is equivalent to that of a certain subcategory of highly connected cobordisms. This result is analogous to \cite[Thm. 3.1]{GRWStableModuli} for cobordism categories of closed manifolds and similar to \cite[Thm A.3]{BP} in the context of manifolds with boundary. Our proof adopts a similar strategy. We assume the reader has prior knowledge of these references, since we will at times explain our proofs by comparison and only indicate the modifications needed in our context. 
\subsection{\texorpdfstring{Cobordism categories of $\Theta$-triads.}{Cobordism categories of Theta-triads.}}
In this subsection, we introduce the cobordism category $\Cobbt$ and state the main result of this section. To do so, we start by defining certain spaces of manifolds with boundary, which will be the object and morphism spaces in our category. For an integer $N\geq 0$, let $\RR^N_+\coloneqq [0,+\infty)\times \RR^{N-1}$ and $\partial \RR^N_+=\{0\}\times \RR^{N-1}.$ We denote the union along the inclusions $\RR^{N}_+\cong\RR^{N}_+\times \{0\}\hookrightarrow \RR^{N+1}_+$ by $\RR^\infty_+$, and $\partial \RR^\infty_+$ for the analogous union of $\partial \RR^N_+.$ Let $n,k\geq 0$ and $U\subset \RR^k\times \RR^{n-k}_+$ be an open subset and observe that $\partial U=U\cap \partial (\RR^k\times \RR^{n-k}_+)$. Let $d\geq 0$ and $\Theta=(\theta,\theta^\partial):(B,B^\partial)\to (\BO(d),\BO(d-1))$ be a map of pairs. We will need to import the following definitions from \cite{GRWStableModuli} and \cite{BP}:
\begin{enumerate}[label=(\alph*)]
    \item The space $\Psi^\partial_\theta(U)$ is the space of pairs $(M,\ell)$ of a $d$-dimensional submanifold with boundary $M$ of $U$ satisfying certain conditions and a bundle map $\ell:TM\to \theta^*\gamma_d$ from \cite[Defn. 2.1]{BP}. From the same reference, for $l\leq d$ the space $\smash{\Psi^\partial_{\theta_l}(U)}$ is the space of pairs $(M,\ell)$ of a $l$-dimensional submanifold with boundary $M$ of $U$ satisfying certain conditions and a bundle map $\ell:TM\oplus \varepsilon^{d-l}\to \theta^*\gamma_d$. 
    \item The space $\Psi_{\theta^\partial}(\partial U)$ is the space of pairs $(M,\ell)$ of a $(d-1)$-dimensional submanifold $M$ of $\partial U$ satisfying certain conditions and a bundle map $\ell:TM\to (\theta^\partial)^*\gamma_{d-1}$ from \cite[Defn. 2.2]{GRWStableModuli}. From the same reference, for $l\leq d-1$ the space $\Psi_{\theta_{l}}(\partial U)$ is the space of pairs of a $l$-dimensional submanifold $M$ of $\partial U$ and a bundle map $\ell:TM\oplus \varepsilon^{d-1-l}\to \theta^*\gamma_{d-1}.$
\end{enumerate}

By definition, there exists a map $\hat{\partial}:\Psi^\partial_{\theta_l}(U)\to \Psi_{\theta_{l-1}}(\partial U)$ taking a pair $(M,\ell)$ to $(\partial M,\ell|_{\partial M}) $ (see also \cite[Lemma B.1]{BP}). We also have a map $\smash{\Psi_{\theta^\partial_l}(\partial U)\to \Psi_{\theta_{l}}(\partial U)}$ by taking $(M,\ell)$ to $(M,\iota\circ \ell)$ where $\iota\circ \ell$ is the composite $TM\oplus \varepsilon^{d-1-l}\oplus \varepsilon^1\to (\theta^\partial)^*\gamma_{d-1}\oplus \varepsilon^1\to \theta^*\gamma_d,$ where $\iota$ is induced by $\Theta.$

\begin{defn}[Spaces of manifolds]\label{spaces of manifolds}
    We define $\Psi_\Theta^\partial(U)$ to be the strict pullback $\smash{\Psi_\theta^\partial(U)\times_{\Psi_{\theta_{d-1}}(\partial U)}\Psi_{\theta^\partial}(\partial U)}$. As a set, $\Psi_\Theta^\partial(U)$ is the set of pairs $(M,\ell)$ where $M$ is a $d$-dimensional neat submanifold with boundary of $U$ and $\ell:(TM,T\partial M)\to (\theta^*\gamma_d,(\theta^{\partial})^*\gamma_{d-1})$ is a map of pairs of vector bundles on $M$ in the sense of \cref{pairs of vector bundles} satisfying the following properties:
    \begin{enumerate}[label=(\textit{\roman*})]
        \item $M$ is closed as a subspace of $U;$
        \item there exists $\epsilon>0$ such that $M\cap( \RR^k\times [0,\epsilon) \times \partial \RR^{n-k}_+)=[0,\epsilon)\times \partial M$ and the restriction of $\ell$ to $[0,\epsilon)\times \partial M\subset M $ is given by
        \[\smash{T([0,\epsilon)\times \partial M)\to T(\partial M)\oplus \varepsilon^1\overset{\ell_{\partial M}\oplus \varepsilon^1}{\to} (\theta^\partial)^*\gamma_{d-1}\oplus \varepsilon^1 \to \theta^*\gamma_d}\]where the leftmost map is the bundle map covering the projection $[0,\epsilon)\times \partial M$ given by the splitting of the source given by the vector $e_1$ in the coordinate $[0,\epsilon).$
    \end{enumerate}We call such a pair $(M,\ell)$ a $\Theta$-manifold and $\ell$ a $\Theta$-structure on $M$. For $l\leq d$, let $\Psi_{\Theta_l}^\partial(U)$ be the strict pullback $\smash{\Psi_{\theta_l}^\partial(U)\times_{\Psi_{\theta_{l-1}}(\partial U)}\Psi_{\theta^\partial_{l-1}}(\partial U)}$. 
\end{defn}

\begin{rmk}\label{the embedded condition gives a collared bundle}
    Let $(M,\ell)\in \Psi^{\partial}_\Theta(U)$, we observe that the first condition of $(ii)$ induces a preferred collar of $(TM,T\partial M)$ (in the sense of \cref{collared pairs}) by taking the inwards pointing vector field given by the basis vector on the $[0,\epsilon)$-coordinate. One can check that, with this collar, the map $\ell$ of collared pairs of vector bundles is collared, by the second part of $(ii)$.
\end{rmk}

\begin{rmk}\label{spaces of manifolds are a pullback}
    By \cite[Lemma B.1]{BP}, the map $\hat{\partial}:\Psi^\partial_{\theta_l}(U)\to \Psi_{\theta_{l-1}}(\partial U)$ is a Serre fibration, so the square
    \[\begin{tikzcd}
        \Psi_{\Theta_l}^\partial(U)\arrow[d]\arrow[r] &\Psi_{\theta_l}^\partial(U)\arrow[d] \\
        \Psi_{\theta^{\partial}_{l-1}}(\partial U)\arrow[r] &\Psi_{\theta_{l-1}}(\partial U).
    \end{tikzcd}\]is a homotopy pullback square and the left vertical map is a Serre fibration.
\end{rmk}

We define now convenient subspaces of the spaces of manifolds defined above.

\begin{defn}
    Define $\psi_\Theta^\partial(n,k)\subset \Psi_\Theta^\partial(\RR^k\times \RR^{n-k}_+)$ consisting of those pairs $(M,\ell)$ such that $M$ is contained in $\RR^k\times [0,1)\times (-1,1)^{n-k-1}.$ Analogously, we define the spaces $\smash{\psi_{\Theta_l}^\partial(n,k)}$ for $l\leq k.$ For $i\in \ZZ$, denote by $\smash{\psi_\Theta^\partial(\infty+i,k)}$ the colimit of these $\psi_\Theta^\partial(n+i,k)$ under the stabilizations over $n$ induced by the inclusion $\RR^k\times [0,1)\times (-1,1)^{n-k-1}\cong \RR^k\times [0,1)\times (-1,1)^{n-k-1}\times \{0\}\subset \RR^k\times [0,1)\times (-1,1)^{n-k}.$
\end{defn}

\begin{nota}
    Let $(W,\ell)\in \smash{\psi_\Theta^\partial(\infty+1,1)}$. Denote by $x_0:\RR\times \RR^\infty_+\to \RR$ the projection onto the first factor. For $K\subset \RR,$ denote $W|_{K}\coloneqq W\cap x_0^{-1}(K)$ and $\hb W|_K\coloneqq \partial W\cap W|_K$ if both subsets are manifolds (possibly with corners). If $K=\{a\}$ is a singleton, we denote $\partial W|_a\coloneqq \hb W|_a$. When $K=[a,b]$ and $W|_{[a,b]}$ is a manifold, then $\partial W|_{[a,b]}=\hb W|_{[a,b]}\cup W|_a\cup W|_b $ so we denote $W|_a\cup W|_b$ by $\vb W|_{[a,b]}.$ We see this triad also as a triad cobordism $W|_{[a,b]}:W|_a\leadsto W|_b$ by setting $\partial_1=\emptyset$ in the sense of \cref{relative cobordism}. 
\end{nota}

\begin{defn}[Cobordism category]\label{cob caties definition original}
    Define a non-unital topological category (in the sense of \cite[Section 2.3]{GRWStableModuli}) $\Cobbt$ having $\smash{\psi^\partial_{\Theta_{d-1}}(\infty,0)}$ as space of objects. The morphism space is the subspace of $\smash{\RR\times 
    \psi_\Theta^\partial(\infty+1,1)}$ consisting of those pairs $(t,(W,\ell))$ where there exists an $\epsilon>0$ such that \[W|_{(-\infty,\epsilon)}= (\RR\times W|_0)|_{(-\infty,\epsilon)}\in \Psi_\Theta^\partial((-\infty,\epsilon)\times \RR^\infty_+) \] and \[W|_{(t-\epsilon,\infty)}=(\RR \times W|_t)|_{(t-\epsilon,\infty)}\in\Psi_\Theta^\partial((t-\epsilon,\infty)\times \RR^\infty_+).\]Here, $\RR\times W|_a$ denotes the manifold $\RR\times W|_a\subset \RR\times \RR^\infty_+$ along with the map $(T(\RR\times W|_a),T(\RR\times \partial W|_a))\to \varepsilon^1\oplus (TW|_a, T\partial W|_a)\to (\theta^*\gamma_d,(\theta^\partial)^*\gamma_{d-1})$. The source and target of $(t,(W,\ell))$ are $(W|_0,\ell|_0)$ and $(W|_t,\ell|_t)$, respectively. Composition is given by "stacking" as defined in \cite[Defn. 2.6]{GRWStableModuli} (see \cref{cob picture}).
\end{defn}

\begin{figure}
    \centering
    \tikzset{every picture/.style={line width=0.75pt}} 

\begin{tikzpicture}[x=0.75pt,y=0.75pt,yscale=-1,xscale=1]

\draw [color={rgb, 255:red, 155; green, 155; blue, 155 }  ,draw opacity=1 ]   (339.28,101.19) -- (302.52,120.35) -- (277.15,133.91) -- (237.48,155.13) -- (167.32,192.81) ;
\draw [shift={(165.56,193.76)}, rotate = 331.76] [color={rgb, 255:red, 155; green, 155; blue, 155 }  ,draw opacity=1 ][line width=0.75]    (10.93,-3.29) .. controls (6.95,-1.4) and (3.31,-0.3) .. (0,0) .. controls (3.31,0.3) and (6.95,1.4) .. (10.93,3.29)   ;
\draw [shift={(341.06,100.26)}, rotate = 152.47] [color={rgb, 255:red, 155; green, 155; blue, 155 }  ,draw opacity=1 ][line width=0.75]    (10.93,-3.29) .. controls (6.95,-1.4) and (3.31,-0.3) .. (0,0) .. controls (3.31,0.3) and (6.95,1.4) .. (10.93,3.29)   ;
\draw [color={rgb, 255:red, 155; green, 155; blue, 155 }  ,draw opacity=1 ]   (277.55,20.23) -- (277.15,133.91) ;
\draw [shift={(277.56,18.23)}, rotate = 90.2] [color={rgb, 255:red, 155; green, 155; blue, 155 }  ,draw opacity=1 ][line width=0.75]    (10.93,-3.29) .. controls (6.95,-1.4) and (3.31,-0.3) .. (0,0) .. controls (3.31,0.3) and (6.95,1.4) .. (10.93,3.29)   ;
\draw [color={rgb, 255:red, 155; green, 155; blue, 155 }  ,draw opacity=1 ]   (131.56,126.72) -- (589.56,148.68) ;
\draw [shift={(591.56,148.77)}, rotate = 182.74] [color={rgb, 255:red, 155; green, 155; blue, 155 }  ,draw opacity=1 ][line width=0.75]    (10.93,-3.29) .. controls (6.95,-1.4) and (3.31,-0.3) .. (0,0) .. controls (3.31,0.3) and (6.95,1.4) .. (10.93,3.29)   ;
\draw [color={rgb, 255:red, 155; green, 155; blue, 155 }  ,draw opacity=1 ]   (183.56,54.4) -- (182.56,184.5) ;
\draw [color={rgb, 255:red, 155; green, 155; blue, 155 }  ,draw opacity=1 ]   (326.56,10.74) -- (326.06,107.76) ;
\draw [color={rgb, 255:red, 155; green, 155; blue, 155 }  ,draw opacity=1 ]   (74.06,179.65) -- (476.56,199.49) ;
\draw [color={rgb, 255:red, 155; green, 155; blue, 155 }  ,draw opacity=1 ]   (177.06,100.7) -- (579.56,120.55) ;
\draw [color={rgb, 255:red, 74; green, 144; blue, 226 }  ,draw opacity=1 ][line width=0.75]    (237.48,155.13) .. controls (232.96,104.65) and (253.85,70.92) .. (273.15,68.16) .. controls (287.65,66.08) and (301.25,81.47) .. (302.52,120.35) ;
\draw    (490.56,116.14) -- (441.65,142.29) -- (347.06,192.88) ;
\draw [color={rgb, 255:red, 245; green, 166; blue, 35 }  ,draw opacity=1 ][line width=1.5]  [dash pattern={on 5.63pt off 4.5pt}]  (302.52,120.35) -- (466.1,129.22) ;
\draw [color={rgb, 255:red, 74; green, 144; blue, 226 }  ,draw opacity=1 ][line width=0.75]    (450.56,137.31) .. controls (449.85,105.67) and (454.49,92.31) .. (459.01,91.63) .. controls (463.97,90.89) and (468.79,105.42) .. (466.27,127.85) .. controls (466.21,128.3) and (466.16,128.76) .. (466.1,129.22) ;
\draw [color={rgb, 255:red, 74; green, 144; blue, 226 }  ,draw opacity=1 ]   (369.56,180.09) .. controls (370.22,154.18) and (375.13,147.03) .. (381.02,146.78) .. controls (386.92,146.52) and (389.4,148.87) .. (389.82,170.41) ;
\draw [color={rgb, 255:red, 245; green, 166; blue, 35 }  ,draw opacity=1 ][line width=1.5]    (237.48,155.13) .. controls (265.56,156.71) and (331.56,181.41) .. (369.56,180.09) ;
\draw    (273.15,68.16) .. controls (315.82,62.38) and (415.06,109.97) .. (459.01,91.63) ;
\draw  [dash pattern={on 4.5pt off 4.5pt}]  (389.82,170.41) .. controls (378.22,168.61) and (351.82,155.56) .. (373.02,146.74) ;
\draw    (356.91,125.68) .. controls (352.91,129.65) and (350.51,146.5) .. (381.02,146.78) ;
\draw [color={rgb, 255:red, 155; green, 155; blue, 155 }  ,draw opacity=1 ]   (348.06,62.78) -- (347.06,192.88) ;
\draw [color={rgb, 255:red, 155; green, 155; blue, 155 }  ,draw opacity=1 ]   (491.06,19.12) -- (490.56,116.14) ;
\draw    (373.02,146.74) .. controls (383.15,143.26) and (418.08,137.03) .. (450.56,137.31) ;
\draw  [draw opacity=0][fill={rgb, 255:red, 233; green, 89; blue, 89 }  ,fill opacity=0.1 ] (151.81,65.3) -- (270.81,67.73) -- (292.14,68.41) -- (310.89,71.72) -- (329.39,76.13) -- (361.89,84.29) -- (407.39,93.77) -- (426.89,95.98) -- (448.89,94.87) -- (459.01,91.63) -- (548.01,96.48) -- (543.64,102.15) -- (540.14,115.6) -- (539.56,142.16) -- (537.39,141.62) -- (508.14,140.3) -- (498.39,139.86) -- (469.42,138.97) -- (450.56,137.31) -- (426.62,138.27) -- (408.5,139.86) -- (391,142.72) -- (377,146.03) -- (381.02,146.78) -- (506.01,154.44) -- (500,159.04) -- (494.5,170.95) -- (493.06,187.58) -- (359.14,179.77) -- (337.7,177.32) -- (312.7,171.36) -- (279.7,162.99) -- (251.45,156.37) -- (237.48,155.13) -- (110.56,148.77) -- (111.81,117.78) -- (114.81,91.1) -- (131.81,71.92) -- cycle ;
\draw  [draw opacity=0][fill={rgb, 255:red, 233; green, 89; blue, 89 }  ,fill opacity=0.1 ] (512.25,159.92) -- (514.06,175.24) -- (389.82,170.41) -- (375.25,165.88) -- (364.25,157.06) -- (363,151.32) -- (362.88,150.86) -- (360.75,142.72) -- (367.5,145.59) -- (381.02,146.78) -- (506.01,154.44) -- (510,156.84) -- cycle ;
\draw  [draw opacity=0][fill={rgb, 255:red, 233; green, 89; blue, 89 }  ,fill opacity=0.08 ] (273.15,68.16) -- (288,68.63) -- (309.5,70.84) -- (339.75,78.11) -- (374.75,87.6) -- (409.75,94.43) -- (426.89,95.98) -- (441.25,95.98) -- (459.01,91.63) -- (548.01,96.48) -- (553,100.83) -- (558,117.36) -- (559.56,132.9) -- (302.52,120.35) -- (157.06,113.05) -- (153.06,65.86) -- cycle ;

\draw (586.4,154.34) node [anchor=north west][inner sep=0.75pt]  [font=\footnotesize]  {$\mathbb{\textcolor[rgb]{0.61,0.61,0.61}{R}}$};
\draw (258.8,1.4) node [anchor=north west][inner sep=0.75pt]  [font=\footnotesize]  {$\textcolor[rgb]{0.61,0.61,0.61}{[ 0,\infty )}$};
\draw (167.56,195.8) node [anchor=north west][inner sep=0.75pt]  [font=\footnotesize,color={rgb, 255:red, 155; green, 155; blue, 155 }  ,opacity=1 ]  {$\partial \mathbb{R}_{+}^{\infty }$};
\draw (179.2,41.79) node [anchor=north west][inner sep=0.75pt]  [font=\scriptsize]  {$\textcolor[rgb]{0.29,0.56,0.89}{0}$};
\draw (345.6,48.5) node [anchor=north west][inner sep=0.75pt]  [font=\scriptsize]  {$\textcolor[rgb]{0.29,0.56,0.89}{t}$};
\draw (248.4,50.73) node [anchor=north west][inner sep=0.75pt]  [font=\footnotesize]  {$\textcolor[rgb]{0.29,0.56,0.89}{W|_{0}}$};
\draw (446.8,74.32) node [anchor=north west][inner sep=0.75pt]  [font=\footnotesize]  {$\textcolor[rgb]{0.29,0.56,0.89}{W|}\textcolor[rgb]{0.29,0.56,0.89}{_{t}}$};
\draw (363.8,62.42) node [anchor=north west][inner sep=0.75pt]  [font=\footnotesize]  {$\textcolor[rgb]{0.82,0.01,0.11}{W|_{[ 0,t]}}$};
\draw (252,168.51) node [anchor=north west][inner sep=0.75pt]  [font=\footnotesize]  {$\textcolor[rgb]{0.96,0.65,0.14}{\partial ^{h} W|_{[ 0,t]}}$};

\end{tikzpicture}

    \caption{This is an example of the underlying submanifold of a morphism $(t,(W,\ell))$ in $\Cobbt$.}
    \label{cob picture}
\end{figure}

It is convenient to define a variant of the cobordism category above where the objects and morphisms contain a fixed submanifold. Let $L$ be a compact $(d-1)$-dimensional submanifold of $[0,1)\times (-\frac{1}{2},0]\times (-1,1)^{\infty-2}$ (possibly with corners) such that $\smash{\partial L= L\cap \partial([0,1)\times (-\frac{1}{2},0]\times (-1,1)^{\infty-2})}$. Let $\smash{\hb L\coloneqq L\cap \{0\}\times (-\frac{1}{2},0]\times (-1,1)^{\infty-2}} $ and $\smash{\vb L\coloneqq L\cap [0,1)\times \{0\}\times (-1,1)^{\infty-2}}$. We require that near $\{0\}\times (-\frac{1}{2},0]\times \RR^{\infty-2}$ it agrees with $[0,1)\times\hb L$ and that near $[0,1)\times \{0\}\times \RR^{\infty-2}$ it agrees with $(-1,0]\times \vb L.$ Additionally, let $(\ell_L,\ell_{\hb L}):(TL,T\hb L)\oplus \varepsilon^1\to (\theta^*\gamma_{d},(\theta^\partial)^*\gamma_{d-1}) $ be a map of pairs of vector bundles such that the restriction of $\ell_L$ to $T[0,\epsilon)\times \hb L$ agrees with the composition \[\smash{T[0,\epsilon)\times \hb L\to TL\oplus \varepsilon^1\overset{\ell_{\hb L}}{\to} (\theta^\partial)^*\gamma_{d-1}\oplus \varepsilon^1\to \theta^*\gamma_d},\]and the restriction to $T(-\epsilon,0]\times \vb L$ is a product structure, that is, the family $t\mapsto \ell_L|_{\{t\}\times \vb L}$ is a constant homotopy of maps of pairs $(TL,T\partial L)\to (\theta^*\gamma_d,(\theta^\partial)^*\gamma_{d-1})$ for $t<\epsilon$ for some $\epsilon.$  

\begin{defn}\label{cob with L}
    Define the non-unital topological category $\Cobbtl$ whose space of objects is the space $\smash{\psi_{\Theta_{d-1},L}^\partial(\infty,0)} \subset \smash{\psi_{\Theta_{d-1}}^\partial(\infty,0)}$ consisting of those $(M,\ell)$ such that $M\cap [0,1)\times (-\infty,0]\times \RR^{\infty-2}=L$ as $\Theta$-manifolds, that is, $\ell$ agrees with $\ell_L$ in this subspace. Its space of morphisms is the subspace $\smash{\psi_{\Theta_{d-1},L}^\partial(\infty+1,1) \subset \psi_{\Theta_{d-1}}^\partial(\infty+1,1)}$ of those $(t,(W,\ell))$ such that $W\cap (\RR\times[0,1)\times(-\infty,0]\times \RR^{\infty-2})=\RR\times L$ as $\Theta$-manifolds (see \cref{cob picture l}).  
\end{defn}

\begin{figure}
    \centering
    \tikzset{every picture/.style={line width=0.75pt}} 

\begin{tikzpicture}[x=0.75pt,y=0.75pt,yscale=-1,xscale=1]

\draw [color={rgb, 255:red, 155; green, 155; blue, 155 }  ,draw opacity=1 ]   (339.28,101.19) -- (302.52,120.35) -- (277.15,133.91) -- (237.48,155.13) -- (167.32,192.81) ;
\draw [shift={(165.56,193.76)}, rotate = 331.76] [color={rgb, 255:red, 155; green, 155; blue, 155 }  ,draw opacity=1 ][line width=0.75]    (10.93,-3.29) .. controls (6.95,-1.4) and (3.31,-0.3) .. (0,0) .. controls (3.31,0.3) and (6.95,1.4) .. (10.93,3.29)   ;
\draw [shift={(341.06,100.26)}, rotate = 152.47] [color={rgb, 255:red, 155; green, 155; blue, 155 }  ,draw opacity=1 ][line width=0.75]    (10.93,-3.29) .. controls (6.95,-1.4) and (3.31,-0.3) .. (0,0) .. controls (3.31,0.3) and (6.95,1.4) .. (10.93,3.29)   ;
\draw [color={rgb, 255:red, 155; green, 155; blue, 155 }  ,draw opacity=1 ]   (277.55,20.23) -- (277.15,133.91) ;
\draw [shift={(277.56,18.23)}, rotate = 90.2] [color={rgb, 255:red, 155; green, 155; blue, 155 }  ,draw opacity=1 ][line width=0.75]    (10.93,-3.29) .. controls (6.95,-1.4) and (3.31,-0.3) .. (0,0) .. controls (3.31,0.3) and (6.95,1.4) .. (10.93,3.29)   ;
\draw [color={rgb, 255:red, 155; green, 155; blue, 155 }  ,draw opacity=1 ]   (131.56,126.72) -- (589.56,148.68) ;
\draw [shift={(591.56,148.77)}, rotate = 182.74] [color={rgb, 255:red, 155; green, 155; blue, 155 }  ,draw opacity=1 ][line width=0.75]    (10.93,-3.29) .. controls (6.95,-1.4) and (3.31,-0.3) .. (0,0) .. controls (3.31,0.3) and (6.95,1.4) .. (10.93,3.29)   ;
\draw [color={rgb, 255:red, 155; green, 155; blue, 155 }  ,draw opacity=1 ]   (183.56,54.4) -- (182.56,184.5) ;
\draw [color={rgb, 255:red, 155; green, 155; blue, 155 }  ,draw opacity=1 ]   (326.56,10.74) -- (326.06,107.76) ;
\draw [color={rgb, 255:red, 155; green, 155; blue, 155 }  ,draw opacity=1 ]   (74.06,179.65) -- (476.56,199.49) ;
\draw [color={rgb, 255:red, 155; green, 155; blue, 155 }  ,draw opacity=1 ]   (177.06,100.7) -- (579.56,120.55) ;
\draw [color={rgb, 255:red, 74; green, 144; blue, 226 }  ,draw opacity=1 ][line width=0.75]    (237.48,155.13) .. controls (232.96,104.65) and (253.85,70.92) .. (273.15,68.16) ;
\draw    (490.56,116.14) -- (441.65,142.29) -- (347.06,192.88) ;
\draw [color={rgb, 255:red, 114; green, 188; blue, 42 }  ,draw opacity=1 ][line width=1.5]  [dash pattern={on 5.63pt off 4.5pt}]  (302.52,120.35) -- (466.1,129.22) ;
\draw [color={rgb, 255:red, 74; green, 144; blue, 226 }  ,draw opacity=1 ][line width=0.75]    (436.18,144.54) .. controls (435.47,112.91) and (437.38,78.54) .. (442.58,78.94) ;
\draw [color={rgb, 255:red, 74; green, 144; blue, 226 }  ,draw opacity=1 ]   (369.56,180.09) .. controls (370.22,154.18) and (375.13,147.03) .. (381.02,146.78) .. controls (386.92,146.52) and (389.4,148.87) .. (389.82,170.41) ;
\draw [color={rgb, 255:red, 245; green, 166; blue, 35 }  ,draw opacity=1 ][line width=1.5]    (237.48,155.13) .. controls (265.56,156.71) and (331.56,181.41) .. (369.56,180.09) ;
\draw  [dash pattern={on 4.5pt off 4.5pt}]  (389.82,170.41) .. controls (378.22,168.61) and (351.82,155.56) .. (373.02,146.74) ;
\draw    (356.91,125.68) .. controls (352.91,129.65) and (350.51,146.5) .. (381.02,146.78) ;
\draw [color={rgb, 255:red, 155; green, 155; blue, 155 }  ,draw opacity=1 ]   (348.06,62.78) -- (347.06,192.88) ;
\draw [color={rgb, 255:red, 155; green, 155; blue, 155 }  ,draw opacity=1 ]   (491.06,19.12) -- (490.56,116.14) ;
\draw    (373.02,146.74) .. controls (383.15,143.26) and (403.7,144.26) .. (436.18,144.54) ;
\draw  [draw opacity=0][fill={rgb, 255:red, 233; green, 89; blue, 89 }  ,fill opacity=0.1 ] (153.61,65.8) -- (272.61,68.22) -- (293.94,68.91) -- (312.69,72.22) -- (444.38,79.44) -- (544.78,87.44) -- (543.98,98.64) -- (541.94,116.1) -- (541.36,142.66) -- (539.19,142.12) -- (509.94,140.79) -- (500.19,140.35) -- (471.22,139.47) -- (452.36,137.8) -- (428.42,138.77) -- (410.3,140.35) -- (392.8,143.22) -- (378.8,146.53) -- (382.82,147.27) -- (507.81,154.93) -- (501.8,159.54) -- (496.3,171.44) -- (494.86,188.08) -- (360.94,180.26) -- (339.49,177.81) -- (314.49,171.86) -- (281.49,163.48) -- (253.24,156.87) -- (239.28,155.62) -- (112.36,149.27) -- (113.61,118.28) -- (116.61,91.6) -- (133.61,72.41) -- cycle ;
\draw  [draw opacity=0][fill={rgb, 255:red, 233; green, 89; blue, 89 }  ,fill opacity=0.1 ] (512.25,159.92) -- (514.06,175.24) -- (389.82,170.41) -- (375.25,165.88) -- (364.25,157.06) -- (363,151.32) -- (362.88,150.86) -- (360.75,142.72) -- (367.5,145.59) -- (381.02,146.78) -- (506.01,154.44) -- (510,156.84) -- cycle ;
\draw  [draw opacity=0][fill={rgb, 255:red, 142; green, 223; blue, 58 }  ,fill opacity=0.14 ] (273.15,68.16) -- (288,68.63) -- (309.5,70.84) -- (442.58,78.94) -- (546.58,86.14) -- (551.53,92.16) -- (556.18,102.94) -- (558,117.36) -- (559.56,132.9) -- (302.52,120.35) -- (157.06,113.05) -- (153.06,65.86) -- cycle ;
\draw [color={rgb, 255:red, 142; green, 223; blue, 58 }  ,draw opacity=1 ][line width=0.75]    (273.15,68.16) .. controls (287.65,66.08) and (301.25,81.47) .. (302.52,120.35) ;
\draw [color={rgb, 255:red, 29; green, 223; blue, 11 }  ,draw opacity=1 ]   (273.15,68.16) -- (442.58,78.94) ;
\draw [color={rgb, 255:red, 142; green, 223; blue, 58 }  ,draw opacity=1 ][line width=0.75]    (442.58,78.94) .. controls (457.08,76.87) and (464.83,90.34) .. (466.1,129.22) ;
\draw [color={rgb, 255:red, 142; green, 223; blue, 58 }  ,draw opacity=1 ][line width=0.75]    (53.95,59.06) .. controls (68.45,56.98) and (82.05,72.37) .. (83.32,111.25) ;
\draw  [draw opacity=0][fill={rgb, 255:red, 29; green, 223; blue, 11 }  ,fill opacity=1 ] (51.97,59.06) .. controls (51.97,58.2) and (52.86,57.51) .. (53.95,57.51) .. controls (55.05,57.51) and (55.94,58.2) .. (55.94,59.06) .. controls (55.94,59.91) and (55.05,60.6) .. (53.95,60.6) .. controls (52.86,60.6) and (51.97,59.91) .. (51.97,59.06) -- cycle ;
\draw  [draw opacity=0][fill={rgb, 255:red, 102; green, 151; blue, 53 }  ,fill opacity=1 ] (81.33,111.25) .. controls (81.33,110.39) and (82.23,109.7) .. (83.32,109.7) .. controls (84.42,109.7) and (85.31,110.39) .. (85.31,111.25) .. controls (85.31,112.1) and (84.42,112.79) .. (83.32,112.79) .. controls (82.23,112.79) and (81.33,112.1) .. (81.33,111.25) -- cycle ;

\draw (586.4,154.34) node [anchor=north west][inner sep=0.75pt]  [font=\footnotesize]  {$\mathbb{\textcolor[rgb]{0.61,0.61,0.61}{R}}$};
\draw (258.8,1.4) node [anchor=north west][inner sep=0.75pt]  [font=\footnotesize]  {$\textcolor[rgb]{0.61,0.61,0.61}{[}\textcolor[rgb]{0.61,0.61,0.61}{0,\infty }\textcolor[rgb]{0.61,0.61,0.61}{)}$};
\draw (167.56,195.8) node [anchor=north west][inner sep=0.75pt]  [font=\footnotesize,color={rgb, 255:red, 155; green, 155; blue, 155 }  ,opacity=1 ]  {$\partial \mathbb{R}_{+}^{\infty }$};
\draw (179.2,41.79) node [anchor=north west][inner sep=0.75pt]  [font=\scriptsize]  {$\textcolor[rgb]{0.29,0.56,0.89}{0}$};
\draw (345.6,48.5) node [anchor=north west][inner sep=0.75pt]  [font=\scriptsize]  {$\textcolor[rgb]{0.29,0.56,0.89}{t}$};
\draw (248.4,50.73) node [anchor=north west][inner sep=0.75pt]  [font=\footnotesize]  {$\textcolor[rgb]{0.29,0.56,0.89}{W|}\textcolor[rgb]{0.29,0.56,0.89}{_{0}}$};
\draw (434.4,58.92) node [anchor=north west][inner sep=0.75pt]  [font=\footnotesize]  {$\textcolor[rgb]{0.29,0.56,0.89}{W|}\textcolor[rgb]{0.29,0.56,0.89}{_{t}}$};
\draw (357.5,45.02) node [anchor=north west][inner sep=0.75pt]  [font=\footnotesize]  {$\textcolor[rgb]{0.82,0.01,0.11}{W|}\textcolor[rgb]{0.82,0.01,0.11}{_{[ 0,t]}}$};
\draw (252,168.51) node [anchor=north west][inner sep=0.75pt]  [font=\footnotesize]  {$\textcolor[rgb]{0.96,0.65,0.14}{\partial }\textcolor[rgb]{0.96,0.65,0.14}{^{h}}\textcolor[rgb]{0.96,0.65,0.14}{W|}\textcolor[rgb]{0.96,0.65,0.14}{_{[ 0,t]}}$};
\draw (287.6,51.53) node [anchor=north west][inner sep=0.75pt]  [font=\footnotesize]  {$\textcolor[rgb]{0.56,0.87,0.23}{L}$};
\draw (565.33,99) node [anchor=north west][inner sep=0.75pt]  [font=\footnotesize]  {$\textcolor[rgb]{0.56,0.87,0.23}{L\times \mathbb{R}}$};
\draw (564.8,122.33) node [anchor=north west][inner sep=0.75pt]  [font=\footnotesize]  {$\textcolor[rgb]{0.48,0.68,0.26}{\partial ^{h} L\times \mathbb{R}}$};
\draw (81.6,67.38) node [anchor=north west][inner sep=0.75pt]  [font=\footnotesize]  {$\textcolor[rgb]{0.56,0.87,0.23}{L}$};
\draw (43.6,38.13) node [anchor=north west][inner sep=0.75pt]  [font=\footnotesize]  {$\textcolor[rgb]{0.11,0.87,0.04}{\partial ^{v} L}$};
\draw (61.1,114.63) node [anchor=north west][inner sep=0.75pt]  [font=\footnotesize]  {$\textcolor[rgb]{0.44,0.68,0.18}{\partial ^{h} L}$};

\end{tikzpicture}

    \caption{This is an example of the underlying submanifold of a morphism $(t,W)$ in $\Cobbtl$.}
    \label{cob picture l}
\end{figure}

Recall from \cite[Section 2.3]{GRWStableModuli} that for a non-unital topological category $\C$, we define its classifying space $\B\C$ as the geometric realization $||N_\bullet \C||$ of its nerve $N_\bullet \C$, which is a semi-simplicial space. The following is proved exactly as \cite[Cor. 2.17]{GRWStableModuli}.

\begin{prop}\label{have or not L in BCOb}
    For any $L$ as above, the map induced by inclusion \(\B\Cobbtl\to \B\Cobbt\) is a weak equivalence.
\end{prop}

We will be interested in certain subcategories of these categories. A \textit{subcategory} $\C$ of $\Cobbtl$ is a pair of a collection of objects $\ob(\C)\subset \ob(\Cobbtl)$ and a collection of morphisms $\mor(\C)\subset \mor(\Cobbtl)$ such that:
\begin{enumerate}[label=(\textit{\alph*})]
    \item $\ob(\C)$ is a union of path components of $\ob(\Cobbtl)$;
    \item the images of source and target maps $s_0,s_1:\mor(\C)\subset \mor(\Cobbtl)\to \ob(\Cobbtl)$ lie in the $\ob(\C);$ 
    \item for every $a,b\in \C$, the morphism space $\C(a,b)$ defined by the strict pullback $ \mor(\C)\times_{\ob(\C)^{\times 2}}\{(a,b)\}$, given by the source and target maps, is a union of path components of $\Cobbtl(a,b);$
    \item for every $(P,\ell_P)\in \C$ and $t>0$, the morphism $(t,[0,t]\times P)$ with the \textit{cylindrical $\Theta$-structure} $\ell_P\oplus \varepsilon^1$ is in $\C.$ 
\end{enumerate}
Unwrapping the definition, $\C$ is a non-unital topological category and the inclusion $\C\to \Cobbt$ is a functor of non-unital topological categories. The following definition is analogous to \cite[Def. 2.9]{GRWStableModuli}.
    
\begin{defn}\label{cobordism categories definition}
    Let $k\geq -1 $ be an integer. Let $\Cobbtlbk$ be the subcategory of $\Cobbtl$ consisting of the same objects and those morphisms $(t,(W,\ell))$ where the inclusion of $\partial W|_t\to \hb W|_{[0,t]}$ is $k$-connected. When $k=\smash{\lfloor \frac{d-3}{2}\rfloor}$, we denote $\Cobbtlbk$ by $\Cobbtlb$. Let now $k\geq -2$. Let $\Cobbtlk$ be the subcategory of $\Cobbtlb$ with the same objects and those morphisms $(t,(W,\ell))$ where the inclusion $(W|_t,\partial W|_t)\to (W|_{[0,t]},\hb W|_{[0,t]})$ is strongly $k$-connected in the sense of \cref{strong conn defn}. 
\end{defn}

The main result of this section (and the only one used in later sections) is the following result, whose proof ocupies the remainder of this section. Recall the definition of a triad handle decomposition from \cref{defn of handles}. Once and for all, we see $(L,\hb L,\emptyset, \vb L)$ as a triad cobordism from $(\vb L,\hvb L,\emptyset)$ to $\emptyset.$

\begin{teo}\label{surgery on morphisms}
    Let $d\geq 5$ and $k\geq -2$ be integers such that $2k\leq d-3$. Assume $(L,\hb L)$ admits a triad handle decomposition only using handles of index at most $d-k-3$ and of any type relative to $\emptyset$. Then the map 
    \[\B\Cobbtlk\to \B\Cobbt\]is a weak equivalence.
\end{teo}

\subsubsection{Flexible models.} To prove \cref{surgery on morphisms}, it is convenient to model the homotopy types of the classifying spaces of these subcategories of $\Cobbtl$ in a more tractable way, similar to \cite[Defs. 2.13 and 2.18]{GRWStableModuli}. We start by introducing two models for these classifying spaces. 

\begin{defn}\label{flexible model def}
    Let $\D_\Theta^\partial\subset \RR\times \RR_{>0}\times \psi_{\Theta}^\partial(\infty+1,1)$ be the subspace of tuples $(t,\epsilon,(W,\ell)) $ such that $[t-\epsilon,t+\epsilon]$ consists of \textit{regular values} for $x_0:W\to \RR$, that is, points $s\in \RR$ such that the maps $x_0:W\to \RR$ and $x_0|_{\partial W}$ have $s$ has a regular value. We define a partial order on $\Dt$ by $(t,\epsilon,(W,\ell))<(t',\epsilon',(W',\ell'))$ if and only if $(W,\ell)=(W',\ell')$ and $t+\epsilon<t'-\epsilon'.$ Denote by $\Dl$ the full subposet of $\Dt$ consisting of those $(t,\epsilon,(W,\ell))$ such that $W\cap (\RR\times [0,+\infty)\times (-\infty,0]\times \RR^{\infty-2})=\RR\times L$ as $\Theta$-manifolds. If $\C\subset \Cobbtl$ is a subcategory, we denote by $\D_{\Theta,L}^{\partial,\C}\subset \Dl$ the smallest subposet containing the objects $(t,\epsilon, (W,\ell))$ where $W|_t$ is an object of $\C$ and of morphisms $(t,\epsilon,(W,\ell))<(t',\epsilon',(W',\ell'))$ where $W'|_{[t,t']}=W|_{[t,t']}:W|_t\leadsto W|_{t'}$ is a morphism of $\C.$
\end{defn}

\begin{prop}\label{flexible model}
    Let $\C\subset \Cobbtl$ be a subcategory of $\Cobbtl,$ then there is a weak equivalence \(\B\C\simeq \B\D_{\Theta,L}^{\partial,\C}.\)
\end{prop}
\begin{proof}
    The proof follows analogously as the proof of \cite[Prop. 2.14]{GRWStableModuli} by replacing the appropriate terms and using \cite[Lemma 10.1]{genauer} instead of Ehresmann's fibration lemma.
\end{proof}


We introduce now a "more flexible model" in the spirit of \cite[Defn. 2.18]{GRWStableModuli}. Roughly speaking, these models have better "mapping in"-properties than the previous model.

\begin{defn}\label{more flexible model def}
    Let $p\geq 0$ be an integer and define $\X^{\partial}_{\Theta,L}$ to be the semisimplicial space where the space of $p$-simplices consists of tuples $a\in \RR^{p+1}$, $\epsilon\in (\RR_{>0})^{p+1}$ and $(W,\ell) \in\Psi_{\Theta}^\partial((a_0-\epsilon_0,a_p+\epsilon_p)\times \RR^{\infty}_+)$ satisfying:
    \begin{enumerate}[label=(\roman*)]
        \item For all $i=0,\cdots, p-1$, we have $a_i+\epsilon_i<a_{i+1}-\epsilon_{i+1};$
        \item $W$ lies in $(a_0-\epsilon_0,a_p+\epsilon_p)\times [0,1)\times (-1,1)^{\infty-1};$
        \item $W$ agrees with $(a_0-\epsilon_0,a_p+\epsilon_p)\times L$ as $\Theta$-manifold pairs on the subspace $U\coloneqq (a_0-\epsilon_0,a_p+\epsilon_p)\times [0,+\infty)\times (-\infty,0]\times \RR^{\infty-2},$ that is, they agree on $\Psi^\partial_{\Theta}(U).$
    \end{enumerate}Given $\C\subset \Cobbtl$ a subcategory, let $\X^{\partial,\C}_{\Theta,L}$ be the subspace of $\X^{\partial}_{\Theta,L}$ consisting of those tuples where for each regular value $t$ of $x_0$, $W|_t\in \C$ and for each pair of regular values $t_0<t_1$ contained in the union over $i$ of $(a_i-\epsilon_i,a_i+\epsilon_i),$ the cobordism $W|_{[t_0,t_1]}$ is a morphism in $\C.$ 
\end{defn}

We can map the nerve of the topological poset $\D_{\Theta,L}^{\partial,\C}$ to $\X^{\partial,\C}_{\Theta,L}$ taking a $p$-simplex $(a_0,\cdots,a_p;\epsilon_0,\cdots \epsilon_p,(W,\ell))$ to the same tuple where we consider the restriction of $(W,\ell)$ to $(a_0-\epsilon_0,a_p+\epsilon_p)\times \RR^{\infty}_+.$ The following result is proved using exactly the same method as in \cite[Prop. 2.20]{GRWStableModuli}.

\begin{prop}\label{more flexible model}
    Let $\C\subset \Cobbtl$ be a subcategory of $\Cobbtl,$ then the map \(\smash{N_\bullet\D_{\Theta,L}^{\partial,\C}\to \X^{\partial,\C}_{\Theta,L}}\) is a weak equivalence after geometric realization.
\end{prop}

\subsection{Surgery on morphisms.}

In this section, we prove \cref{surgery on morphisms}. To do so, we consider the following filtrations from \cref{cobordism categories definition}
\[ \Cobbtlk\subset \Cobbtlkk\subset \cdots \subset \Cob^{\partial,0}_{\Theta,L}\subset \Cob^{\partial,-1}_{\Theta,L}=\Cobbtlb \]and 
\[ \Cobbtlbk\subset \Cobbtlbkk\subset \cdots \subset \Cob^{\partial,(-1)}_{\Theta,L}\subset \Cob^{\partial,(-2)}_{\Theta,L}=\Cobbtl. \]The main inputs to the proof of \cref{surgery on morphisms} are the two results stated below, which give conditions to the inclusions of the subcategories in the filtrations above to induce an equivalence on classifying spaces. 

\begin{prop}[Surgery on the interior]\label{surgery on interior morphisms}
    Assume $d\geq 5$. Suppose that the following conditions are satisfied:
    \begin{enumerate}[label=(\roman*)]
        \item $2k\leq d-3;$
        \item $(L,\hb L)$ admits a triad handle decomposition where interior and left handles have index at most $d-k-3$ relative to $\emptyset$.
    \end{enumerate}Then the map induced by inclusion \(\B\Cobbtlk\to \B\Cobbtlkk\) is a weak equivalence.
\end{prop}


\begin{prop}[Surgery on the boundary]\label{surgery on boundary morphisms}
    Assume $d\geq 3$. Suppose that the following conditions are satisfied:
    \begin{enumerate}[label=(\roman*)]
        \item $2k\leq d-3;$
        \item $(L,\hb L)$ admits a triad handle decomposition where right handles have index at most $d-k-3$ and left handles have index at most $d-k-2$ relative to $\emptyset$.
    \end{enumerate}Then the map induced by inclusion \(\B\Cobbtlbk\to \B\Cobbtlbkk\) is a weak equivalence.
\end{prop}

These results are analogous to \cite[Thm. 3.1]{GRWStableModuli} and \cite[Thms. A.2 and B.4]{BP}. Assuming these results, we can now prove our main statement of this section.

\begin{proof}[Proof of \cref{surgery on morphisms}]
    This follows by applying \cref{surgery on boundary morphisms} to conclude that the inclusions in the second filtration induce equivalences on classifying spaces for $k\leq \lfloor\frac{d-3}{2}\rfloor$ and thus the inclusion $\Cobbtlb\to \Cobbtl$ induces an equivalence of classifying spaces. Applying \cref{surgery on interior morphisms}, we deduce that the inclusions in the first filtration are equivalences after taking classifying spaces if $2k\leq d-3$. Thus, we obtain that the maps $\B\Cobbtlk\to \B\Cobbtlb\to \B\Cobbtl$ are equivalences. The result now follows by post-composing this with the equivalence of \cref{have or not L in BCOb}.
\end{proof}

\subsubsection{The approach due to Galatius and Randal-Williams.}
The remainder of this section is devoted to proving \cref{surgery on interior morphisms,surgery on boundary morphisms}. We start by explaining the overall strategy of the proofs. We employ a strategy originally due to Galatius and Randal-Williams in \cite{GRWStableModuli}. Let us consider a general context. Let $\A\subset \BB$ be subcategories of $\Cobbtl.$ We wish to prove that the inclusion $\A\hookrightarrow \BB$ induces an equivalence on classifying spaces. To do so, we will construct an augmented bi-semisimplicial space\footnote{Here, we mean that this bi-semisimplicial space is augmented in the second coordinate. In other words, it is a semi-simplicial object in the category of augmented semi-simplicial spaces.} $\smash{(\Da)_{\bullet,\bullet}}$ where $\smash{(\Da)_{\bullet,-1}=N_\bullet\Db}$ along with maps $\smash{\iota:N_\bullet\Da\to (\Da)_{\bullet,0}}$ of semi-simplicial spaces fitting in a commutative square of solid arrows

\begin{equation}\label{oscar sorens strategy}
    \begin{tikzcd}
        \B\Da \arrow[d, "\iota"]\arrow[rr, "\ref{more flexible model}", "\simeq"'] & &\left|\Xa\right| \arrow[d, "\text{inc}"]\\
        \left|\left(\Da\right)_{\bullet,\bullet}\right| \arrow[r,"p" ]\arrow[rru, dashed,"\K^1"] &\B\Db \arrow[r,"\ref{more flexible model}", "\simeq"'] &\left|\Xb\right|
    \end{tikzcd}
\end{equation}where the map $p$ is induced by the augmentation, satisfying the following properties:
\begin{enumerate}
    \item\label{1 grw} (\textit{Contractability of surgery data}) The map $p$ is a weak equivalence.
    \item\label{2 grw} (\textit{Parameterized surgery move}) There exists a dashed lift $\K^1$ as above making the upper triangle strictly commute and a homotopy $\K^t$ making the lower triangle commute up to homotopy.
\end{enumerate}By \ref{1 grw} and commutativity of the lower triangle, the map $\K^1$ induces an injection on all homotopy groups. By the commutativity of the upper triangle, the map $\K^1$ induces also a surjection on all homotopy groups. Thus $\K^1$ is a weak equivalence and thus so are the vertical maps. By \cref{flexible model} and \cref{more flexible model}, the right vertical map is equivalent to the map $\B\A\to \B\BB$ and hence the latter is a weak equivalence. The remaining sections are devoted to constructing this bisimplical space and this diagram for $\A$ and $\BB$ the steps of the filtrations above.


\subsubsection{Proof of \cref{surgery on interior morphisms}.}


Recall the definition of the submanifolds $V\coloneqq(-2,0)\times \RR^k\times \RR^{d-k} $ and $\bar{V}\coloneqq[-2,0]\times \RR^k\times \RR^{d-k}$ of $\RR^{k+1}\times \RR^{d-k}$ from \cite[Section A.1]{BP}. We start with some notation.

\begin{nota}
    We fix once and for all an infinite set $\Omega$. Let $(\Lambda,\delta,e)$ be a triple where $\Lambda\subset \Omega$ is a finite and $\delta:\Lambda\to [p]^\vee$ is a function and $e:\Lambda\times \bar{V}\hookrightarrow \RR\times [0,1)\times (0,1)\times (-1,1)^{\infty-2}$. Here $[p]^\vee$ means the set of monotone non-decreasing maps from $[p]\coloneqq \{0,\cdots,p\}$ to $[1]$ (see \cite[149]{BP}). We say that such triple is a \textit{(interior) surgery datum} for a $p$-simplex $(a,\epsilon,(W,\ell))$ in the nerve of $\smash{\D_{\Theta,L}^{\partial,\C}}$ if the triple satisfies conditions $(i)$ to $(v)$ in \cite[Def. A.3.]{BP}. We call the subspaces $D_i=e(\delta^{-1}(i)\times \partial_-D^{k+2}\times \{0\})$ the \textit{cores} of (the trace of) the surgery datum $(\Lambda,\delta,e)$ for $i=0,\cdots,p,$ where $\partial_-D^{k+2}\coloneqq \{(x_1,\cdots, x_{k+2})\in \RR^{k+2}|~\sum_{i=1}^{k+2}x_i^2=1,~ x_1\geq 0\}$. Notice that since the third coordinate is restricted to $(0,1)$, the embedding $e$ is automatically disjoint from $L.$
\end{nota}

We define now the bi-semisimplicial spaces promised above in the context of the first filtration. The following definition is essentially \cite[Def. A.3]{BP}.

\begin{defn}\label{defn of space of surgery}
    Given a $p$-simplex $x=(a,\epsilon,(W,\ell))$ of $\Dkk$. Let $\Z_0(x)$ to be set of surgery data $(\Lambda,\delta,e)$ for $x$ such that the map 
    \[\left(W|_{a_i}\cup D_i|_{[a_{i-1},a_i]},\partial W|_{a_i}\right)\to \left(W_{[a_{i-1},a_i]},\hb W|_{[a_{i-1},a_i]}\right)\]is strongly $k$-connected for each $i=1,\cdots, p.$ Define $\Z_p(x)$ in the same way as \cite[Def.A.5]{BP}. Similarly to \cite{BP} define the bi-semisimplicial space given by $\smash{(\Dk)_{p,q}}$ to be the space of pairs $(x,y)$ where $x$ is a $p$-simplex of $\smash{(\Dkk)_\bullet\coloneqq N_\bullet \Dkk}$ and $y\in \Z_q(x)$ topologized as a subspace of
    \[\left(\Dkk\right)_p\times \left( \coprod_{\Lambda\subset \Omega} C^\infty(\Lambda\times \bar{V},\RR\times \RR^\infty_+) \right)^{(p+2)(q+1)}. \]
\end{defn}

    We can define the vertical left map $\iota$ for the square \eqref{oscar sorens strategy} in this setting as the inclusion of the empty surgery data, that is, $\Lambda=\emptyset$. We prove now the first property of the square \eqref{oscar sorens strategy} in this context. The following result should be compared to \cite[Thm. A.6]{BP} and \cite[Thm. 3.4]{GRWStableModuli}.

\begin{prop}\label{contractability of surgery data}
    Assume the hypothesis of \cref{surgery on interior morphisms}, then the augmentation map \[\left(\Dk\right)_{\bullet,\bullet}\to \left(\Dkk\right)_\bullet\] is a weak equivalence after geometric realization.
\end{prop}
\begin{proof}
    This proof is very similar to the proof of \cite[Theorem 3.4]{GRWStableModuli}. However, there are steps where extra care is needed. We assume the reader's familiarity with this reference. We do not give a complete and detailed proof, but instead we sketch the strategy, which follows by analogy to loc.cit., and focus on the points where it differs. For every $p\geq 0$, the semi-simplicial space $\smash{(\Dk)_{p,\bullet}}$ is augmented over $\smash{(\Dk)_p}.$ We will show that the augmentation induces a weak equivalence on geometric realizations for every $p$. This implies the claim by taking geometric realization in the $p$-coordinate. To prove that the augmentation map induces a weak equivalence, we use the criterion of \cite[Thm. 6.2]{GRWStableModuli}. For that, one must check that this augmented semi-simplicial space is an \textit{augmented topological flag complex} in the sense of \cite[Def. 6.1]{GRWStableModuli}. We sketch the argument here for completeness. Just as in the case of this reference (see the remark in p.334 of \cite{GRWI}), our augmented semi-simplicial space does not satisfy this criteria of \cite[Thm. 6.2]{GRWStableModuli}. However, one can find an equivalent model which weakens some of the conditions in the definition of this space and satisfies this criteria, analogous to $\tilde{D}_{\theta,L}^\kappa(\RR^N)$ in \cite[Defn. 6.8]{GRWStableModuli}, by only assuming that the embeddings in a surgery datum are immersions which restrict to embeddings on the cores. We proceed implicitly with this weaker model. To apply \cite[Thm. 6.2]{GRWStableModuli}, we must check that, for every $p,$ this augmented topological flag complex satisfies conditions $(i)-(iii)$ of loc.cit. To verify condition $(i),$ one proceeds exactly as \cite[Prop. 6.10]{GRWStableModuli}, using that the embeddings in our surgery data are contained in the interior. 

    Condition $(ii)$ is equivalent to the claim that $\Z_0(x)$ is non-empty for every $x\in \Dkk.$ We will proceed as in \cite[Prop. 6.13]{GRWStableModuli} to check this condition. We explain the part of the argument which is different in our case in detail now. Let $x=(a,\epsilon,(W,\ell))$ be in $\Dkk.$ It follows from the assumptions, the pair $(W|_{[a_{i-1},a_i]},W|_{a_{i}})$ is $(k-1)$-connected for every $i$ and the pair $(\hb W|_{[a_{i-1},a_i]},\partial W|_{a_{i+1}})$ is $b\coloneqq\lfloor \frac{d-3}{2}\rfloor$-connected. Thus, we see that $(W|_{[a_{i-1},a_i]},W|_{a_{i}})$ is already $k$-connected, since $k\leq b$ and $(W|_{a_i},\partial W|_{a_i})\to (W|_{[a_{i-1},a_i]},\hb W|_{[a_{i-1},a_i]})$ is strongly $(k-1)$-connected. By Whitehead's theorem, $W|_{[a_{i-1},a_i]}$ is equivalent to a CW complex $X$ obtained from $W|_{a_{i}}$ by attaching cells of dimension at least $k+1.$ These cells can be assumed to be embedded in $W|_{[a_{i-1},a_i]}$ since $2(k+1)< d,$ so we get an embedding
    \[\hat{e}_{i,0}:\Lambda_{i,0}\times (D^{k+1},\partial D^{k+1})\to (W|_{[a_{i-1}+\epsilon_{i-1},a_i+\epsilon_i]},W|_{a_i+\epsilon_i})\]for finite set $\Lambda_{i,0}\subset \Omega$ such that \(W|_{a_i+\epsilon_i}\cup \text{im}(\hat{e}_{i,0})\to W|_{[a_{i-1}+\epsilon_{i-1},a_i+\epsilon_i]}\)is $(k+1)$-connected for every $i$. This implies that the map \[{(W|_{a_i+\epsilon_i}}\cup {\text{im}(\hat{e}_{i,0})}, \partial W|_{a_i+\epsilon_i})\to(W|_{[a_{i-1}+\epsilon_{i-1},a_i+\epsilon_i]},\hb W|_{[a_{i-1}+\epsilon_{i-1},a_i+\epsilon_i]})\] is strongly $(k+1)$-connected. These embeddings can be assumed to be disjoint from $\hb W|_{[a_{i-1}+\epsilon_{i-1},a_i+\epsilon_i]},$ since they can be defined in the interior of $W|_{a_{i}+\epsilon_{i}}.$
    These embeddings can be isotoped to be disjoint from $L\times \RR$ by the following argument. By taking a triad handle decomposition of $L$ (where $L$ is seen as a triad cobordism from $(\vb L,\hvb L,\emptyset)$ to $\emptyset$ by taking $\partial_0 L=\hb L$ and $\partial_1L=\emptyset$) as in the hypothesis, it suffices to make these embeddings disjoint from the products of the cores by $\RR.$ The cores of the right handles are automatically disjoint since they are collars of submanifolds of $\hb W|_{[a_{i-1}+\epsilon_{i-1},a_i+\epsilon_i]}.$ The left and interior cores have dimension at most $d-k-2$ so $\hat{e}_{i,0}$ is generically disjoint since $(d-k-2)+(k+1)<d$, by transversality. Once again by the fact that $2(k+1)<d,$ $\hat{e}_{i,0}$ and $\hat{e}_{j,0}$ can be assumed to be disjoint for all $i\neq j.$ One proceeds now as in \cite[Prop. 6.13]{GRWStableModuli} to extend these embeddings to a surgery datum. Thus, we conclude that $\Z_0(x)$ is non-empty and so condition $(ii)$ is verified. 
    
    The verification of $(iii)$ goes exactly as \cite[Prop. 6.12]{GRWStableModuli} since any tuple $(k+1)$-dimensional submanifolds of a $d$-dimensional manifold can be made pairwise disjoint by transversality. This finishes the proof by using \cite[Theorem 6.2]{GRWI}.  
\end{proof}

We move now to establish the second property of \eqref{oscar sorens strategy} in this context. Before that, we record an application of the Blakers-Massey theorem needed for this proof.

\begin{lemma}\label{blackers massey lemma}
    Let $(M,K)$ be a $1$-connected compact manifold pair of dimension $d\geq 3$. Let $\phi:(P,\partial P)\hookrightarrow (M,K)$ be an embedding of pairs, where $P$ is a simply connected manifold with boundary of dimension $k\geq 0$ with no closed components and path-connected boundary. Let $(M',K')$ be the manifold pair obtained from $(M,K)$ by removing an open tubular neighborhood of $\phi.$ If $d\geq k+3,$ then the induced map \(\pi_i(M',K')\to \pi_i(M,K)\) is an isomorphism for $i\leq d-k-1$ and a surjection for $i=d-k$ for every basepoint. 
\end{lemma}
\begin{proof}
    Start by observing that the claim is equivalent to proving that the square
    \[\begin{tikzcd}
        K'\arrow[d]\arrow[r] & M'\arrow[d] \\
        K \arrow[r] & M
    \end{tikzcd}\]is homotopy $(d-k-1)$-cartesian in the sense of \cite[Defn. 3.3.1]{Munson_Volić_2015}. We establish this by proving that this square is homotopy $(d-k)$-cocartesian and applying the Blakers-Massey theorem. In other words, we prove now that the $M'\cup_{K'} K\to M$ is $(d-k)$-connected. This map is a pushout of the map \(S(\nu_{P,M})\cup_{S(\nu_{\partial P,K})}D(\nu_{\partial P,K})\hookrightarrow D(\nu_{P,M})\) where $\nu_{P,M}$ and $\nu_{\partial P,K}$ are the normal bundles  of $P$ in $M$ and $\partial P$ in $K$ respectively (which are oriented since $P$ is simply-connected), and $D(-)$ and $S(-)$ are the disc and sphere bundles of a vector bundle (with a metric). The right-hand side is a $d$-dimensional manifold with boundary homotopy equivalent to $P$ and the left-hand side is its boundary. By Poincaré duality, we have that $\H_i(D(\nu_{P,M}),\partial D(\nu_{P,M}))\cong\H^{d-i}(P)\cong \H_{k-d+i}(P,\partial P)$ which vanishes for $k-d+i\leq 0$, since $P$ has no closed components and hence $\partial P\hookrightarrow P$ is $0$-connected. Hence $\H_i(D(\nu_{P,M}),\partial D(\nu_{P,M}))$ vanishes for $i\leq d-k.$ Since $d\geq k+3$ and $P$ is simply connected, we deduce that $\partial D(\nu_{P,M})$ is simply connected. Thus, the map $\partial D(\nu_{P,M})\hookrightarrow D(\nu_{P,M})$ is $(d-k)$-connected. Hence, the square above is homotopy $(d-k)$-cocartesian in the sense of \cite[Defn. 3.7.1]{Munson_Volić_2015}. Additionally, the map $K'\to K$ is $(d-1-k)$-connected by transversality, since $\partial P$ is $(k-1)$-dimensional. We prove now that $K'\to M'$ is $1$-connected. Since $d\geq k+2,$ we have $\pi_0(K')=\pi_0(K)$ and $\pi_0(M')\to \pi_0(M)$ so $\pi_0(M')=\pi_0(K').$ Moreover, since $d\geq k+3$, the same applies for fundamental groups and thus $\pi_1(K')$ surjects onto $\pi_1(M').$ We conclude by applying the Blakers-Massey theorem, as stated in \cite[Thm. 4.2.3]{Munson_Volić_2015}, to deduce that the square above is homotopy $(d-k-1)$-cartesian.
\end{proof}

The following result is similar to \cite[Lemma A.9]{BP} and uses \cite[Prop. 3.6]{GRWStableModuli} (see \cref{why this is different} below for more details on this similarity).

\begin{prop}\label{standard interior}
    Assume the hypothesis of \cref{surgery on interior morphisms}, then there exists a homotopy

    \[\K:I\times \left(\Dk\right)_{p,0}\to \left(\Xkk\right)_p\]such that the image of $(t,(W,\ell))$ lies in $\left(\Xk\right)_p$ if $t=1$ or if $(W,\ell)\in \left(\Dk\right)_p$.
\end{prop}
\begin{proof}
    We use the map $\K$ constructed in \cite[(A.1)]{BP}, that uses the standard family of \cite[Prop. 3.6]{GRWStableModuli}, and verify the claim in the statement above. Once again, we assume familiarity of the reader with both constructions. Denote the image of $(t,(a,\epsilon,(W,\ell),e))$ under $\K$ by $(a,\epsilon,(W_t,\ell))$ (since this homotopy is constant on the parameters $a$ and $\epsilon).$ The manifold $W_t$ is obtained from $W_0$ by removing the interior of $e$ and attaching the manifold $P_t\in \Psi_{\theta}(V)$ from \cite[Prop 3.6]{GRWStableModuli}. Since, by \cite[Rmk. A.8]{BP} the family is constant on its horizontal boundary, we immediately have a $\theta^{\partial}$-structure in $\hb W_t$ from $\hb W$ by taking the constant family of $\theta^\partial$ structures $\ell|_{\hb W_t}\coloneqq \ell|_{\hb W}$. Thus, we can see $W_t\in \Psi_{\Theta,L}^\partial(\RR^\infty).$ It remains to check that the image of $\K$ lies in $\smash{(\Xk)_p}$ if $t=1$ or $(W,\ell)\in \smash{(\Dk)_p}$. Let us start by proving the case $(a,\epsilon,(W,\ell))\in \smash{(\Dk)}_p.$ Let $t\in [0,1],$ we have to show that, for two regular values $a<b$ of $x_0,$ the map
    \[(W_t|_{b},\partial W_t|_{b})\to (W_t|_{[a,b]},\hb W_t|_{[a,b]})\]is strongly $k$-connected, provided $(W|_b,\partial W|_b)\to (W|_{[a,b]},\hb W|_{[a,b]})$ is strongly $k$-connected. Since the family $\hb W_t$ is constant in $t$, it suffices to prove that 
    \[W_t|_b\cup_{\partial W_t|_b}\hb W_t|_{[a,b]}\to W_t|_{[a,b]}\]is $(k+1)$-connected. Denote the manifold $W|_{[a,b]}\backslash \text{im}(e)^\circ$ by $X$. Observe that the pair $(W_t|_{[a,b]},\hb W_t|_{[a,b]})$ is the union of $(X,\hb W_t|_{[a,b]})$ with $P_t|_{[a_\lambda,b_\lambda]}$ for all $\lambda\in \Lambda_{i,0}$ for some $i.$ We start by proving that \[Y\coloneqq X|_b\cup_{\partial X|_b}\hb W_t|_{[a,b]}\to X\]is $(k+1)$-connected. The pair $(X,Y)$ is obtained from $(W|_{[a,b]}, W|_b\cup_{\partial W|_b}\hb W|_{[a,b]})$ by cutting out a neighborhood of embedded copies of $(D^{k+1},\partial D^{k+1}).$  By applying \cref{blackers massey lemma} for each copy iteratively and $2k\leq d-3$, we conclude that $\pi_{i}(X,Y)$ vanishes for $i\leq k+1\leq d-(k+1)-1,$ since it is isomorphic to $\pi_i(W|_{[a,b]}, W|_b\cup_{\partial W|_b}\hb W|_{[a,b]}),$ which vanishes by assumption. Since connectivity of maps is preserved under pushouts, the map \(W_t|_b\cup_{\partial W_t|_b} \hb W_t|_{[a,b]}=W_t|_b\cup_{X|_b} Y\to W_t|_b\cup_{X|_b} X\) is also $(k+1)$-connected, since $\partial X|_b=\hb W_t|_b$. By construction, the map $X\cup_Y W_t|_b\to W_t|_{[a,b]}$ is a pushout of a $(k+1)$-connected map, by Property $(iv)$ of \cite[Lemma A.7]{BP} (see square $(\text{A}.4)$ in \cite{BP}). We conclude that, the map $W_t|_b\cup \hb W_t|_{[a,b]}\to W_t|_{[a,b]}$ is $(k+1)$-connected. This establishes the case $(a,\epsilon,(W,\ell))\in N_p\Dk.$ 
    
    We now concern ourselves with the case $t=1$ and $(a,\epsilon,(W,\ell))\in \smash{(\Dkk)_p}.$ In \cite[Lemma 3.7]{GRWStableModuli}, the proof is divided in three steps. We focus on Step $2$ of loc.cit., since it is the non-trivial step. More precisely, we have to show that, for two regular values $a\in (a_{i-1}-\epsilon_{i-1},a_{i-1}+\epsilon_{i-1})$ and $b\in(a_i-\epsilon_i,a_i+\epsilon_i)$ of $x_0,$ the map
    \[W_1|_{b}\cup_{\partial W_1|_{b}} \hb W_1|_{[a,b]}\to W_1|_{[a,b]}\]is $(k+1)$-connected. It suffices to prove the analogous statement for $\widetilde{W}_1$ obtained from $W$ by only doing the surgeries for $\Lambda_{i,0}.$ Isotope the submanifold $D_i$ to a submanifold $\widetilde{D_i}$ so it is disjoint from the image of $e_{i}$, using the fact that it has trivial normal bundle. Once again by \cref{blackers massey lemma}, since $k+1\leq d-(k+1)-1$, we deduce that the map
    \[X|_b\cup \widetilde{D}_i|_{[a,b]}\cup_{\partial W|_{b}}\hb W|_{[a,b]}\to X,\]is $(k+1)$-connected, since it is obtained from $W|_b\cup \widetilde{D}_i|_{[a,b]}\cup_{\partial W|_{b}}\hb W|_{[a,b]}\to W|_{[a,b]}$, which is $(k+1)$-connected by assumption, by removing a neighborhood of an $(k+1)$-disc. Once again by Property $(iv)$ of \cite[Lemma A.7]{BP}, we deduce that 
    \[\widetilde{W}_1|_{b}\cup_{\partial \widetilde{W}_1|_{b}} \hb \widetilde{W}_1|_{[a,b]}\cup \widetilde{D}_i|_{[a,b]}\to \widetilde{W}_1|_{[a,b]}\]is $(k+1)$-connected. However, by Property $(v)$ of \cite[Lemma A.7]{BP}, $\widetilde{D}_i$ retracts to $\widetilde{W}_1|_b$ so the map $\widetilde{W}_1|_b\hookrightarrow \widetilde{W}_1|_b\cup \widetilde{D}_i|_{[a,b]}$ is an equivalence. This implies
    \[\widetilde{W}_1|_{b}\cup_{\partial \widetilde{W}_1|_{b}} \hb \widetilde{W}_1|_{[a,b]}\to \widetilde{W}_1|_{[a,b]}\]is $(k+1)$-connected. This finishes the proof.
\end{proof}

\begin{rmk}\label{why this is different}
    The difference between this statement and \cite[Lemma A.9]{BP} is two-fold. Firstly, we work with a more general notion of tangential structure and therefore must produce a family of manifold pairs with such structure. Secondly, our definition of $\Dk$ and $\Xk$ differs from the analogous one in that paper in two ways. First, we do not assume that $\partial W|_t$ is $(n-2)$-connected for every regular value $t$. Second, their connectivity condition of the morphisms is dual to the notion of strong connectivity in the sense that the square in \cref{cocartesian definition} is assumed to be homotopy $k$-cartesian instead of homotopy $(k+1)$-cocartesian (see also \cref{1 conn: old same as new} for a relation between the two).
\end{rmk}

\begin{proof}[Proof of \cref{surgery on interior morphisms}]
    Consider the square \eqref{oscar sorens strategy} for $\A= \Cobbtlk$ and $\BB=\Cobbtlkk$, where the bottom left space is defined in \cref{defn of space of surgery} and the left vertical map as the inclusion of the empty surgery datum. By \cref{contractability of surgery data}, the map $p$ is a weak equivalence and by \cref{standard interior} and an analogous extension as in \cite[301]{GRWStableModuli} to a homotopy between geometric realizations gives a lift $\K$. By the discussion below \eqref{oscar sorens strategy}, we conclude that the map $\A\hookrightarrow B$ induces an equivalence on classifying spaces.
\end{proof}

\subsubsection{Proof of \cref{surgery on boundary morphisms}.}


Recall now the definition of $V\coloneqq (-2,0)\times \RR^k\times \RR^{d-k+1}_+$ and $\bar{V}\coloneqq [-2,0]\times \RR^k\times \RR^{d-k+1}_+$ from \cite[Section B.2]{BP}.

\begin{nota}
    We fix again an infinite set $\Omega$. Let $(\Lambda,\delta,e)$ be a triple where $\Lambda\subset \Omega$ is a finite and $\delta:\Lambda\to [p]^\vee$ (see \cite[149]{BP}) is a function and $e:\Lambda\times \bar{V}\hookrightarrow \RR\times [0,1)\times (0,1)\times (-1,1)^{\infty-1}$. We say that such triple is a \textit{boundary surgery datum} for a $p$-simplex $(a,\epsilon,(W,\ell))$ in the nerve of $\smash{\D_{\Theta,L}^{\partial,\C}}$ if the triple satisfies conditions $(i)$ to $(vi)$ in \cite[Def. B.5]{BP}. We call the subspace $D_i=e(\delta^{-1}(i)\times \partial_-D^{k+1}\times \{0\})$ the \textit{cores} of the trace of the surgery datum $(\Lambda,\delta,e)$ for $i=0,\cdots,p.$ The space of such data is topologized similarly to \cite[Def B.5]{BP}. Notice that since the third coordinate is restricted to $(0,1)$, the embedding $e$ is automatically disjoint from $L.$
\end{nota}

\begin{defn}
    Given a $p$-simplex $x=(a,\epsilon,(W,\ell))$ of $\Dbkk$. Let $\Z_0(x)$ to be space of surgery data $(\Lambda,\delta,e)$ for $x$ such that the map 
    \[\partial W|_{a_i}\cup D_i\to \hb W|_{[a_{i-1},a_i]}\]is $k$-connected for each $i=1,\cdots, p.$ Define $\Z_p(x)$ in the same way as \cite[Def.B.5]{BP}. Similarly to \cite{BP} define the bi-semisimplicial space given by $\smash{(\Dbk)_{p,q}}$ to be the space of pairs $(x,y)$ where $x$ is a $p$-simplex of $\smash{ \Dbkk}$ and $y\in \Z_q(x).$
\end{defn}


The following result should be compared to \cite[Lemma B.6]{BP} and \cite[Thm. 3.4]{GRWStableModuli} and is proved in the same way as \cref{contractability of surgery data}.

\begin{prop}\label{contractability of boundary surgery data}
    Assume the hypothesis of \cref{surgery on boundary morphisms}, then the augmentation map \[\left(\Dbk\right)_{\bullet,\bullet} \to \left(\Dbkk\right)_\bullet\] is a weak equivalence after geometric realization.
\end{prop}

The following result is essentially \cite[Prop. B7]{BP}.

\begin{prop}\label{standard boundary}
    Assume the hypothesis of \cref{surgery on boundary morphisms}, there exists a homotopy

    \[\K:I\times \left(\Dbk\right)_{p,0}\to \left(\Xbkk\right)_p\]such that the image of $(t,(W,\ell))$ lies in $\left(\Xbk\right)_p$ if $t=1$ or $(W,\ell)\in \left(\Dbk\right)_p$.
\end{prop}
\begin{proof}
    Such family was constructed in \cite[Prop. B.7]{BP} as a $\theta$-manifold pair. We can follow exactly the same proof by replacing the map $\smash{\Psi^{\partial}_{\theta}(V)\to \Psi^{\partial}_{\theta_{d-1}}(V)}$ by $\smash{\Psi^{\partial}_{\Theta}(V)\to \Psi^{\partial}_{\theta^\partial}(V)}$. The former is a Serre fibration by \cite[Lemma B.1]{BP}. By definition of $\Psi^{\partial}_{\Theta}(V)$ as a pullback, so is the latter. This property is the only necessary ingredient for the proof of \cite[Prop. B.7]{BP}. The rest of the argument follows from \cite[Lemma 3.7]{GRWStableModuli} applied to $\hb W.$
\end{proof}

\begin{proof}[Proof of \cref{surgery on boundary morphisms}]
    This follows analogously as the proof of \cref{surgery on interior morphisms} by replacing \cref{contractability of surgery data} by \cref{contractability of boundary surgery data} and \cref{standard interior} by \cref{standard boundary} and considering \eqref{oscar sorens strategy} for $\A= \Cobbtlbk$ and $\BB=\Cobbtlbkk$.
\end{proof}


\section{Stable stability.}\label{sta sta section}

One of the key steps in the proof of \cref{main no tang} and a generalization of it allowing tangential structures is establishing some form of "stable stability" in the context of cobordism categories of manifolds with boundary. The first goal of this section is to make this wish into a precise claim (see \cref{stable stability} below). The second goal of this section is to reduce this statement to two closure properties of the subcategory of "stably stable cobordisms" (see \cref{main closure property} and \cref{main closure property high} below). The proof of these statements is deferred to the next section. From now on, we specialize our study to odd-dimensional manifolds and cobordism categories whose morphisms are odd-dimensional cobordisms. Our definitions and proofs are close in spirit to \cite{GRWII}, in particular the proof of Theorem $2.15$ in loc.cit. Nevertheless, we follow a slightly different strategy, but we refer to the analogous concepts and proofs in \cite{GRWII} as they appear.

\subsection{The statement.} The goal of this subsection is to state stable stability in the form of \cref{stable stability} below. The first important definition towards this goal is that of a $\Theta$-end for a map of pairs $\Theta:(B,B^\partial)\to (\BO(d),\BO(d-1))$. To define this, consider the following definitions. Let $(M,K)$ be a $d$-dimensional manifold pair (see \cref{section of pairs of manifolds}), recall that $(TM,TK)$ is a collared vector bundle pair, where the collar is given by an inwards-pointing vector field along $K$. A \textit{$\Theta$-structure} on $(M,K)$ is a collared bundle map $\ell:(TM,TK)\to (\theta^*\gamma_d,(\theta^{\partial})^*\gamma_{d-1}).$ The \textit{space of $\Theta$-structures} on $(M,K)$ is $\Bunc(TM,\Theta^*\gamma_d)$ as in \cref{collared pairs}. A \textit{framing} on a pair $(M,K)$ is a $\text{fr}_d$-structure for \(\text{fr}_d:(\EO(d),\EO(d-1))\to (\BO(d),\BO(d-1)).\) The pair $(\RR^{d}_+,\partial\RR^{d}_+)$ admits the canonical framing coming from the basis $\{e_1,\cdots, e_{d}\}.$ 

\begin{ass}\label{the assumptions for stable stability}
    For the entirety of this subsection, we fix the following choices:
    \begin{enumerate}[label=\underline{$\mathrm{\Roman*}$}]
        \item An odd integer $d=2n+1\geq 7$.
        \item A map of pairs $\Theta:(B,B^\partial)\to (\BO(2n+1),\BO(2n))$ such that $B^\partial$ is path-connected.
        \item \label{the basepoint framing map} A collared bundle map $\tau:\text{fr}^*_d(\gamma_{2n+1},\gamma_{2n})\to \Theta^*(\gamma_{2n+1},\gamma_{2n}).$ This induces a preferred $\Theta$-structure for every framed $(2n+1)$-manifold pair. In particular, we have a preferred choice of $\Theta$-structure on $(\RR^{2n+1}_+,\partial \RR^{2n+1}_+),$ called the \textit{basepoint $\Theta$-structure.} Moreover, up to homotopy $\tau$ is determined by the choice of basepoint $\Theta$-structure on $(\RR^{2n+1}_+,\partial \RR^{2n+1}_+).$
        \item \label{the L} A compact $(d-1)$-dimensional submanifold $L$ of $[0,1)\times (-\frac{1}{2},0]\times (-1,1)^{\infty-2}$ and $\Theta$-structure $\ell_L$ on $L$ satisfying the assumptions of \cref{cob with L}. We assume that $\hvb L\neq \emptyset.$
    \end{enumerate} 
\end{ass}

We now introduce the notion of a \textit{standard framing}. The following definitions are necessary for the definition of $\Theta$-ends and are analogous to \cite[Defn.2.9]{GRWII}. For every integer $k\geq 0,$ we consider the embedding
\[\rho:(D^k,\partial D^k) \hookrightarrow (\RR^{k+1}_+,\partial \RR^{k+1}_+)\]taking $x=(x_1,\cdots, x_k)\in D^k\subset \RR^k$ to $(\sqrt{1-|x|^2},x_1,\cdots,x_k).$ Define $(V_1,W_1)$ to be the submanifold pair $(S^n\times \rho(D^{n+1}),S^n\times \rho(\partial D^{n+1}))\subset (\RR^{n+1}\times \RR^{n+2}_+,\RR^{n+1}\times \partial\RR^{n+2}_+)$ induced by $\rho$ and the inclusion $S^n\subset \RR^{n+1}.$ Recall the triad $(D^{k}_+,\partial_0 D^k_+,\partial_1D^k_+)$ from \cref{triad homotopy}. 

\begin{defn}\label{right standard}
    The \textit{right standard framing} $\xi^r$ on $(S^n\times D^{n+1}_+,S^n\times \partial_0D^{n+1}_+)$ is the framing induced by the codimension $0$ embedding of pairs
    \begin{align*}
        (S^n\times D^{n+1}_+,S^n\times \partial_0D^{n+1}_+) & \hookrightarrow (\RR^{n+1}\times \RR^n_+,\RR^{n+1}\times \partial \RR^n_+)\cong (\RR^{2n+1}_+,\partial \RR^{2n+1}_+)\\
        (x;y_0,\cdots,y_{n}) & \mapsto (2e^{-\frac{y_1}{2}} x;y_0,y_2,\cdots,y_{n})
    \end{align*} We say that a $\Theta$-structure on $S^n\times D^{n+1}_+$ is \textit{standard} if it is homotopic to $\ell_r\coloneqq \tau\circ\xi^r,$ that is, in the same path component in $\Bunc(T(S^n\times D^{n+1}_+),\Theta^*\gamma).$
\end{defn}

\begin{defn}\label{left standard}
    The \textit{left standard framing} $\xi^l$ on $(D^{n+1}\times D^{n},\partial D^{n+1}\times D^{n})$ is the framing induced by the codimension $0$ embedding of pairs
    \begin{align*}
        (D^{n+1}\times D^{n},\partial D^{n+1}\times D^{n}) & \hookrightarrow (\RR^{n+2}_+\times \RR^{n-1},\partial \RR^{n+2}_+\times \RR^{n-1})=(\RR^{2n+1}_+,\partial \RR^{2n+1}_+)\\
        (x;y_1,\cdots,y_{n}) & \mapsto (2e^{-\frac{y_1}{2}} \rho(x);y_2,y_3,\cdots,y_{n})
    \end{align*} We say that a $\Theta$-structure on $D^{n+1}\times D^{n}$ is \textit{standard} if it is homotopic to $\ell_l\coloneqq \tau\circ\xi^l,$ that is, in the same path component in $\Bunc(T(D^{n+1}\times D^{n}),\Theta^*\gamma).$
\end{defn}

\begin{defn}[Standard $\Theta$-structure]\label{standard structure on V1}
    Fix a closed $2n$-disc $D$ inside the interior of the product of the lower hemispheres $D^n_-\times D^n_-\subset S^n\times S^n$ and denote the complement $W_1\backslash \text{int}(D)$ by $W_{1,1}.$ A $\Theta$-structure $\ell$ on the pair $(V_1,W_{1,1})$ is \textit{standard} if both structures $\bar{e}^*\ell$ and $\bar{f}^*\ell$ are standard, where $\bar{e}$ and $\bar{f}$ are the embeddings defined by
    \begin{align*}
        \bar{e}:(S^n\times D^{n+1}_+,S^n\times \partial_0D^{n+1}_+) & \hookrightarrow (S^{n}\times \rho(D^{n+1}),S^n\times \rho(\partial D^{n+1}))\subset \RR^{n+1}\times \RR^{n+2}_+
    \end{align*}induced by the inclusion $D^{k+1}_+\subset D^{k+1}$ and $\rho,$ and 
    \begin{align*}
        \bar{f}:(D^{n+1}\times D^{n},S^n\times D^{n}) & \hookrightarrow (S^{n}\times \rho(D^{n+1}),S^n\times \rho(\partial D^{n+1}))\subset \RR^{n+1}\times \RR^{n+2}_+\\
        (x;y) & \mapsto \left(y,\sqrt{1-|y|^2};\rho(x)\right)
    \end{align*}
\end{defn}


\begin{lemma}\label{the space of standard is connected}
    The space of standard $\Theta$-structures on $(V_1,W_{1,1})$ extending the basepoint $\Theta$-structure (see \ref{the basepoint framing map}) on a disc in $\partial W_{1,1}$ is non-empty and connected.
\end{lemma}
\begin{proof}
    This is an adaptation of \cite[Lemma 7.6/7]{GRWI}. Denote the images of $\bar{e}$ and $\bar{f}$ by $E$ and $F$, respectively. Observe that the intersection $E\cap F$ is a contractible pair, \textit{i.e.} equivalent to $(*,*)$ as a pair. Two framings on a contractible manifold pair are homotopic if and only if they induce the same orientation. Since both $\bar{e}$ and $\bar{f}$ are orientation preserving embeddings, it follows that the restrictions of a right standard framing on $E$ and a left standard framing on $F$ to $E\cap F$ are homotopic. Thus, we can homotope them to agree on this subspace pair. This defines a standard framing on $E\cup F,$ which is isotopy equivalent to $V_1$ as pairs. This proves the first claim.

    Let us denote the space of standard $\Theta$-structures on $(V_1,W_{1,1})$ by $X$. Denote by $X_0$ (resp. $X_1$) the path component of $\Bunc(T(S^n\times D^{n+1}_+),\Theta^*\gamma)$ (resp. $\Bunc(T(D^{n+1}\times D^{n}),\Theta^*\gamma)$) of standard $\Theta$-structures. Denote also by $Z$ the path component of $\Bunc(T(E\cap F),\Theta^*\gamma)$ given by the basepoint $\Theta$-structure (here we are identifying $E\cap F$ with $(\RR^{2n+1}_+,\partial \RR^{2n+1}_+)$). One can show that the restriction maps $X_i\to Z$ are Serre fibrations (see \cref{bundle map fibration}), and thus $X$ is equivalent to the homotopy pullback $X_0\times_Z X_1$. Denote by $X^*$ the space of standard $\Theta$-structures on $(V_1,W_{1,1})$ extending the basepoint $\Theta$-structure on a disc $D$ in $\partial W_{1,1}$. We want to show that $X^*$ is path-connected. Denote also by $X_0^*$ and $X_1^*$ the analogous spaces of $\Theta$-structures extending the basepoint $\Theta$-structure on $D$. Observe that the previous pullback description of $X$ implies that $X^*$ is equivalent to the product $X_0^*\times X_1^*$, since $Z$ is equivalent to the space of $\Theta$-structures on $D$. Thus, to prove that $X^*$ is path-connected, it suffices to show that $X_i^*$ is path-connected for $i=0,1.$ For $i=0$, this now follows by the fact that $X_0$ is path-connected and the fact that we have a fiber sequence \(X^*_0\to X_0\to Z\) that admits a section $s:Z\to X_0$ by pre-composing a $\Theta$-structure on $E\cap F$ with the embedding from \cref{right standard}. For $i=1$, this follows verbatim.
\end{proof}

As in \cite[Defn. 2.5]{GRWII}, it will be convenient to consider a version of the cobordism category $\Cobbtl$ defined in \cref{cob with L} where the submanifold $L$ from \ref{the L} is removed in the following way.

\begin{defn}\label{cobordism categories definition L cut}
    Let $L$ be as in \cref{the assumptions for stable stability}. Define the non-unital topological category $\CobbtL$ to have objects $(M^\circ,\ell|_{M^\circ})$ for $(M,\ell)\in \Cobbtl$, where $(M^\circ,\hb M^\circ) \coloneqq(M\backslash \int(L),\partial M\backslash \int(\hb L)) $. Here $\int(L)\coloneqq L\backslash \vb L $ and $\int(\hb L)=\hb L\backslash \hvb L.$ The morphisms are triples $(t,(W^\circ,\ell|_{W^\circ}))$ for $(t,(W,\ell))$ a morphism in $\Cobbtl$, where $(W^\circ,\hb W^\circ) \coloneqq (W\backslash \int(L)\times [0,t],\hb W\backslash \int(\hb L)\times [0,t])$. The topology on the space of objects and morphisms is the unique one such that the assignment $(M\mapsto M^\circ)$ from $\Cobbtl$ to $\CobbtL$ is an isomorphism of topological categories. Similarly, define $\smash{\CobbtLo}$ to be the subcategory with the same objects and those morphisms $W:M\leadsto M'$ such that $(M',\hb M')\to (W,\hb W)$ is strongly $(n-1)$-connected.
\end{defn}

\begin{rmk}
    Observe that given $W:M\leadsto N$ in $\Cobbtl$, then the $4$-ad $(W^\circ,\hb W^\circ, \vb L\times [0,1], M^\circ\sqcup N^\circ)$ is a triad cobordism in the sense of \cref{relative cobordism}.
\end{rmk}

\begin{rmk}
Recall that $\Cobbtls$ is a subcategory of $\Cobbtl$ of those morphisms $W:M\leadsto N$ such that $(N,\partial N)\to (W,\hb W)$ is strongly $(n-1)$-connected. It turns out that $\smash{\Cobbtls}$ is not, in general, isomorphic to $\smash{\CobbtLs}$. It will be clear later that $\smash{\CobbtLs}$ is more convenient for our purposes. However, under some conditions on $L$ (which will be satisfied in our case), we will see that $\smash{\CobbtLs}$ is isomorphic to $\smash{\Cobbtls}$ (see \cref{cut or not cut is the same} below). \end{rmk}

We are now ready for the definition of a $\Theta$-end. From now on, we directly work with the underlying manifolds and $\Theta$-structures where $L$ is cut out. Therefore, we drop the notation of $M^\circ$ unless necessary, that is, when we talk about an object $(M,\ell)\in \CobbtL$, then $(M,\hb M)=(N^\circ,\hb N^\circ)$ for some $N\in \Cobbtl$.

\begin{defn}\label{theta end in l}
    Let $L$ be as in \cref{the assumptions for stable stability}, a \textit{$\Theta$-end} $K$ in $\CobbtL$ is a sequence of composable morphisms
    \[\{K|_{[i,i+1]}:K|_i\leadsto K|_{i+1}\}_{i\geq 0}\]in $\CobbtL$ such that:
    \begin{enumerate}[label=(\roman*)]
        \item the inclusions $(K|_i,\hb K|_i)\hookrightarrow (K|_{[i,i+1]},\hb K|_{[i,i+1]}) \hookleftarrow(K|_{i+1},\hb K|_{i+1}) $ are strongly $(n-1)$-connected, for every $i\geq 0$;
        \item for every $i\geq 0,$ there exists an embedding $\omega:(V_1,W_{1,1})\hookrightarrow (K|_{[i,i+1]},\hb K|_{[i,i+1]})$ in a path component which intersects $\vb L$. We require that $\im(\omega)$ is disjoint from $\vb L$ and that $\omega^*\ell_{K_{[i,i+1]}}$ is a standard $\Theta$-structure. 
    \end{enumerate}
\end{defn}

For any $\Theta$-end $K$ in $\CobbtL,$ define
\[\F(P,K|_i)\subset \CobbtL(P,K|_i)\]to be the subspace of those $\Theta$-cobordisms $(s,W)$ such that $\ell_W:(W,\hb W)\to (B,B^\partial)$ is strongly $n$-connected. The following is similar to \cite[Lemma 2.14]{GRWII}.

\begin{lemma}\label{subfunctor}
    For every $i\geq 0$, the subspace $\F(P,K|_i)\subset \CobbtL(P,K|_i)$ determines a subfunctor \[\F(-,K|_i):\left(\CobbtLs\right)^{\text{op}}\to \Top\]of $\CobbtL(-,K|_i).$ Moreover, post-composition with $K|_{[i,i+1]}$ defines a natural transformation
    \[\F(-,K|_i)\Rightarrow \F(-,K|_{i+1}).\]
\end{lemma}
\begin{proof}
    The first claim is equivalent to the following: given $M:Q\leadsto P$ in $\CobbtLs$ and $W:P\leadsto K|_i$ such that $\ell_W$ is strongly $n$-connected, then $\ell_{W\circ M}$ is strongly $n$-connected. To prove this, it suffices to observe that $(W,\hb W)\hookrightarrow (M\cup_P W,\hb M\cup_{\hb P} \hb W)$ is strongly $(n-1)$-connected, by \cref{closure properties strongly}. This follows since $(P,\hb P)\hookrightarrow (M,\hb M)$ is strongly $(n-1)$-connected, using \cref{strong conn and pushouts}. The second claim follows exactly by the same argument using that $(K|_i,\hb K|_i)\hookrightarrow (K|_{[i,i+1]}, \hb K|_{[i,i+1]})$ is strongly $(n-1)$-connected.
\end{proof}

Denote by $\F(P,K|_\infty) $ the homotopy colimit of $ \F(P,K|_i)$ under the stabilization maps $K|_{[i,i+1]}\circ (-)$. By \cref{subfunctor}, this defines a presheaf of spaces on $\CobbtLs,$ by taking a functorial model of the homotopy colimit, e.g. the telescope. We recall from the introduction that a map of spaces $f:X\to Y$ is called an \textit{abelian homology equivalence} if $f_*:\H_k(X;f^*\L)\to \H_k(Y;\L)$ is an isomorphism for all \textit{abelian} local systems $\L: \Pi_1(Y)\to \Ab$, that is, those $\L$ such that for every $x\in X$ the action of the commutator subgroup of $\pi_1(X,x)$ on $\L(x)$ is trivial. We are now ready to state the main result of this section. Recall the \cref{the assumptions for stable stability}.

\begin{teo}[Stable stability]\label{stable stability}
    For any $\Theta$-end $K$, the functor $\F(-,K|_\infty)$ takes morphisms to abelian homology equivalences, provided $(B,B^\partial)$ is $1$-connected.
\end{teo}

The path to the proof of this statement is inspired by the proof of \cite[Thm. 2.15]{GRWII} and will roughly go as follows. First, we will reduce this statement to a claim about the behaviour of the functor $\F(-,K|_\infty)$ under certain simpler types of morphisms. Second, for each of these types, we further reduce this statement to a certain \textit{closure property}, which is then proved in the next section. \cref{stable stability} is the only result we use in later sections. The next subsection establishes the first step in the proof of this statement. 

\subsection{\texorpdfstring{Elementary simplifications of \cref{stable stability}.}{Elementary simplifications.}}\label{simplification}

From now on we fix a $\Theta$-end $K$. The goal of this subsection is to reduce the proof of \cref{stable stability} for $K$ to the verification that certain "elementary" morphisms in $\smash{\CobbtLs}$ are sent to abelian homology equivalences by $\smash{\F(-,K|_\infty)}$. To do so, consider the following definition. Let $\W$ be the collection of morphisms $M:P\leadsto Q$ of $\CobbtLs$ such that $\smash{\F(Q,K|_\infty)\to \F(P,K|_\infty)}$ is an abelian homology equivalence. Notice that \cref{stable stability} is equivalent to $\W=\smash{\CobbtLs}.$ We start by showing that $\W$ is a subcategory which contains all isomorphisms in the sense of the definition below. Analogously to the definition below \cref{have or not L in BCOb}, a subcategory of $\CobbtL$ is a collection of morphisms closed under composition, which is the union of entire path components and contains $[0,t]\times P$ with the cylindrical $\Theta$-structure $\ell_P\oplus \varepsilon^1$. 
    \begin{defn}\label{isomorphisms}
    Let $\C$ be a subcategory of $\CobbtL$ and $W:M\leadsto N$ be a morphism in $\C.$ We say that $W$ is an \textit{isomorphism} if there exists another morphism $W':N\leadsto M$ in $\C$ such that $W\circ W'$ and $W'\circ W$ are in the same path components of $[0,s]\times N$ and $ [0,s]\times M$ (both with the cylindrical $\Theta$-structure) in $\C(N,N)$ and in $\C(M,M)$, respectively, for some $s>0$.
    \end{defn}

\begin{exam}\label{cylinders are isomorphisms}
    Let $t>0$ and $P\in \CobbtL.$ For any $\Theta$-structure $\ell$ on $ [0,t]\times P$, the morphism $W\coloneqq(t, [0,t]\times P,\ell)$ is an isomorphism: let $W^*$ be the dual $\Theta$-cobordism given by precomposing $\ell$ with map on tangent bundles induced by the reflection around the axis $\{\frac{t}{2}\}\times \RR^{\infty}_+;$ observe that $W$ and $W^*$ are composable morphisms, whose compositions have $\Theta$-structure homotopic to $[0,2t]\times \ell_P$ relative to both ends. Hence, $W$ and $W^*$ are mutual inverses and hence are isomorphisms. In general, any morphism $(t,W)$ whose underlying manifold is diffeomorphic to $[0,s]\times P$ for some $s>0$ is an isomorphism.
\end{exam}

\begin{lemma}\label{W is a subcategory and isos}
    The collection $\W$ is a subcategory of $\CobbtLs$ which contains all isomorphisms. 
\end{lemma}
\begin{proof}
    Start by observing that $\W$ is closed under composition since abelian homology equivalences are closed under composition. To show that $\W$ contains entire path components of the morphism spaces $\CobbtL(P,Q)$, observe that a path in $\CobbtLs(P,Q)$ from $(t,(W,\ell))$ to $(t',(W',\ell'))$ induces a homotopy between the induced pre-composition maps by these morphisms, since the composition law is continuous. Thus, $\W$ contains entire path components since abelian homology equivalences are homotopy invariant. To show that $[0,s]\times P$ with the cylindrical $\Theta$-structure is in $\W$, it suffices to observe that the induced pre-composition map is isotopic to the identity by scaling the $s$-parameter and thus a homotopy equivalence.
    
    We are left to show that all isomorphisms are in $\W$. Let $W:M\leadsto N$ be an isomorphism in $\CobbtLs$ and $W':N\leadsto M$ be another morphism in $\CobbtLs$ such that $W\circ W'$ and $W'\circ W$ are in the same path component as the corresponding cylinders $[0,s]\times N$ or $[0,s]\times M.$ This implies that the induced maps of $W$ and $W'$ are homotopy inverses to each other. Thus, $(-)\circ W$ is a homotopy equivalence and hence an abelian homology equivalence.
\end{proof}

\begin{lemma}\label{two out of three}
    The subcategory $\W\subset \CobbtLs$ has the $2$-out-of-$3$ property.
\end{lemma}
\begin{proof}
    This follows from the $2$-out-of-$3$ property for abelian homology equivalences, which we prove now. Let $g\circ f:X\to Y\to Z$ be two maps of spaces, where at least two out of the three maps $f,g$ or $g\circ f$ are abelian homology equivalences. Thus, in particular, all maps are homology equivalences since homology equivalences satisfy the $2$-out-of-$3$ property. In particular $\H_1(X)\cong \H_1(Y)\cong \H_1(Z).$ Suppose $g$ and $g\circ f$ are abelian homology equivalences and let $\L$ be an abelian system on $Y.$ Then, since $\H_1(Y)\cong \H_1(Z),$ there exists a local system $\L'$ on $Z$ such that $\L=g^*\L'.$ This implies that for $f$ induces an isomorphism for all abelian systems on $Y.$ The remaining cases are simpler since one always considers abelian local systems on $Z$ and systems which are pulled back from abelian systems are themselves abelian.
\end{proof}

Recall the definition of elementary triad cobordisms, their type and index from \cref{elementary cobordisms}. Often we will use results in \cref{section of pairs of manifolds} to deduce consequences for morphisms in $\CobbtL$. We will hinge on the following remark.

\begin{rmk}[Factorizations = Triad handle decompositions]\label{factorizations and handle structures}
    Let $W:M\leadsto N$ be a morphism in $\CobbtL$. Given a triad handle decomposition of the underlying triad cobordism of $W$ in the sense of \cref{defn of handles}, there exists a path in $\CobbtL(M,N)$ from $W$ to a composition
    \[\smash{M\overset{W_1}{\leadsto} M_1\overset{W_2}{\leadsto} M_2\leadsto \cdots \leadsto M_k=N}\]whose underlying triad cobordisms $W_i$ are elementary of the types and indices present in the triad handle decomposition. This is seen by the following argument: The triad handle decomposition induces a decomposition $W=W_1\cup_{M_1}W_2\cup_{M_2}\cup \cdots\cup_{M_{k-1}}W_k$ where $W_i$ is a triad elementary cobordism from $M_i$ to $M_{i+1}$. We can embed the underlying triad $W$ such that $M_i$ lies in the subspace $\{i\}\times \RR^\infty_+$, the cobordisms $W_i$ lie in $[i,i+1]\times \RR^\infty_+$ and satisfy the appropriate collaring conditions to be seen as morphisms in $\CobbtL(M_i,M_{i+1})$ for all $i\geq 0$. By the Whitney's embedding theorem, one can see that this embedding of $W$ is isotopic to the original inclusion of $W$ into $\RR\times \RR^\infty_+$. Moreover, such an isotopy can be promoted to a path in the morphism space $\CobbtL(M,N),$ hence justifying the original claim.
\end{rmk}

We move now to the first simplification of \cref{stable stability}.

    \begin{lemma}\label{no need top}
    If $\W$ contains every morphism $M:P\leadsto Q$ in $\CobbtLs$ whose underlying triad cobordism is elementary relative to $Q$ of the following type:
    \begin{enumerate}[label=(\roman*)]
        \item left type and index $n+1\leq k< 2n+1;$
        \item right type and index $n\leq k<2n;$
        \item interior type and index $2n,$
    \end{enumerate}then $\W=\CobbtLs.$
\end{lemma}
\begin{proof}
    By combining \cref{strong geometrical connectivity,splitting} and \cref{factorizations and handle structures}, every morphism $M$ in $\CobbtLo$ admits a factorization into elementary triad cobordisms of index at least $n$ for right type, $n+1$ for left type and $2n+1$ for interior type relative to $Q$, up to a path in the morphism space $\CobbtLs(P,Q)$. By \cref{two out of three}, $M$ is in $\W$ if all elementary pieces are. By hypothesis, all the elementary pieces are in $\W$ except the ones of left or interior type of index $2n+1$ and right type of index $2n.$ We prove now that our assumptions imply that these cases lie in $\W$ too.

    Let $M:P\leadsto Q$ be a morphism of $\CobbtLs$ whose underlying cobordism is elementary of interior type and index $2n+1$ relative to $Q.$ Let $\phi:S^{2n}\hookrightarrow Q$ be the attaching map of the interior $(2n+1)$-handle of such a triad handle decomposition. By invariance of domain, this embedding is a diffeomorphism onto one component of $Q$. In particular, this component is disjoint from the boundary and therefore disjoint from $\vb L.$ Choose an embedding $\varphi:S^{0}\times D^{2n}\hookrightarrow Q$ sending one disc to the image of $\phi$ and the other to the interior of any other component (which exists since $\hvb L\neq \emptyset$ and by the discussion above). Let $U$ be the result of attaching an interior handle to $Q$ along this attaching map. This cobordism admits a $\Theta$-structure extending $\ell_Q$ since $B$ is path-connected, and thus we can see it as a morphism $U:Q\leadsto R$ in $\CobbtL$. Moreover, the inclusion $(R,\hb R)\to (U,\hb U)$ is strongly $(2n-2)$-connected and thus it lies in $\CobbtLs$. The composition $U\circ M$ consists of a $2n+1$ and a $2n$-handle of interior type relative to $R$ such that the belt (recall \cref{elementary cobordisms}) of the latter intersects the attaching map of the former exactly at one point, by construction. Thus, by \cref{cancellation}, $U\circ M$ is diffeomorphic to $R\times [0,1]$ as a cobordism relative to $R,$ which is in $\W$ by \cref{W is a subcategory and isos,cylinders are isomorphisms}. By hypothesis, $U$ is in $\W$ and thus so is $M$, by \cref{two out of three}.

    Let $M:P\leadsto Q$ be a morphism in $\smash{\CobbtLs}$ whose underlying cobordism is elementary of right type of index $2n$ relative to $Q.$ By invariance of domain, the attaching map $\phi:S^{2n-1}\hookrightarrow \hb Q$ is a diffeomorphism onto one path component of $Q$ disjoint from $\vb L.$ Choose similarly an embedding $\varphi:S^0\times D^{2n-1}\hookrightarrow \hb Q$ sending one disc to the image of $\phi$ and the other to any other component of $\hb Q$ and let $U:Q\leadsto R$ be the result of adding a right handle along this embedding. Since $B^\partial$ is path-connected, $U$ admits a $\Theta$-structure and is strongly $(n-1)$-connected relative to its target, so it lies in $\CobbtLs$. The composition $U\circ M$ consists of a right $2n$-handle and a left $2n$-handle such that the belt of the latter intersects the attaching sphere of the former in exactly one point. By \cref{merging}, this is isomorphic to an elementary interior cobordism of index $2n$ relative to $Q.$ This is in $\W$ by hypothesis and since $U$ is by assumption, so is $M.$

    Let $M:P\leadsto Q$ be a morphism in $\smash{\CobbtLs}$ whose underlying cobordism is elementary of left type of index $2n+1$ relative to $Q.$ Consider the attaching map $\phi: (D^{2n},S^{2n-1})\hookrightarrow (Q,\hb Q).$ Once again by invariance of domain, this embedding is a diffeomorphism onto one component of $Q$. In particular, this component does not contain the boundary and therefore is disjoint from $\vb L.$ The map $\phi|_{S^{2n-1}}$ is closed and open and thus a diffeomorphism onto a component of $\hb Q.$ Let $\varphi:S^0\times D^{2n-1}\hookrightarrow \hb Q$ be an embedding sending one disc to the image of $\phi$ and the other one to another path component (once again exists by the discussion above, since $\hvb L$ is non-empty). Let $U:Q\leadsto R$ be the result of a $1$-right handle attachment at $\varphi.$ By the same reason as before, $U$ admits a $\Theta$-structure and is strongly $(n-1)$-connected relative to its target, so it lies in $\CobbtLs$. By \cref{cancellation}, the composite $U\circ M$ is a cancelling pair of a left $2n$ and left $(2n+1)$-handle and thus is in $\W$. Since $U\in \W$ by hypothesis, we see that $M\in \W.$ This finishes the proof.
\end{proof}

We introduce the following definition and notation, which will be convenient for the proof of \cref{stable stability}.

\begin{nota}[Basepoint component]\label{basepoint component def}
    Notice that $\hvb L$ has finitely many connected components (by compactness of $L$) and is non-empty. By definition, each $K|_{[i,i+1]}$ contains an embedded copy of $V_1$ in some path component intersecting $\vb L.$ There exists a path component $\partial_0\vb L$ of $\hvb L$ such that for infinitely many $i,$ $\hb K|_{[i,i+1]}$ contains an embedded copy of $\hb V_1=W_{1,1}$ in a path component of $\hb K|_{[i,i+1]}$ intersecting $\partial_0\vb L.$ Choose a path component $\vb_0L$ of $\vb L$ such that $\partial(\vb_0L)$ contains $\partial_0\vb L$ and call it the \textit{basepoint component}. Observe that $(\vb_0 L,\partial_0\vb L)$ is a pair in the sense of \cref{section of pairs of manifolds}. It follows that for infinitely many $i,$ $K|_{[i,i+1]}$ contains an embedded copy of $V_1$ in a path component intersecting $\vb_0 L.$
\end{nota}

\begin{ass}
    Fix a basepoint component $(\vb_0L,\partial_0\vb L)$ as above. For the purposes of \cref{stable stability}, we may assume that for every $i\geq 0,$ the cobordism $K|_{[i,i+1]}$ contains a $V_1$ in a path component intersecting $\vb_0L$ by the following argument: By composing the appropriate number of $K|_{[i,i+1]}$ and rescaling the first coordinate, we produce a $\Theta$-end $K'$ with the same underlying manifold (that is, the union of all $K'_{[i,i+1]})$ and isotopy equivalent submanifold of $\RR\times \RR^{\infty}_+,$ that has this property. This assumption will not affect the conclusion of \cref{stable stability}, since it takes the colimit over $i$. Such colimit for $K$ agrees with the one for $K'.$ 
\end{ass}

Given a morphism $M:P\leadsto Q$ such that the underlying cobordism is elementary of any type and of index $k\geq n$ relative to $Q,$ we say that it is \textit{attached to the basepoint component} if the attaching map of the $k$-handle lies in a component of $Q$ intersecting the basepoint component $\vb_0 L$ and the horizontal boundary of the attaching map of $k$-handle lies in a component of $\hb Q$ intersecting the component $\partial_0\vb L\subset \partial(\vb_0L)$ for some handle decomposition of $M$ relative to $Q$.

\begin{lemma}\label{basepoint component statement}
    The conclusion of \cref{no need top} still holds if additionally $M$ is assumed to be attached to the basepoint component. 
\end{lemma}
\begin{proof}
    It suffices to prove that an elementary triad cobordism $M:P\leadsto Q$ of type and index as in the statement but not attached in the basepoint component lies in $\W.$ We proceed by checking this for each type.
    Let $M:P\leadsto Q$ be elementary of right type and index $n\leq k<2n$ relative to $Q$ and let $\phi:\partial D^{k}\times D^{2n-k}\hookrightarrow \hb Q$ be its attaching map. Choose an embedding $\varphi:S^0\times D^{2n-1}\hookrightarrow \hb Q$ sending one disc to the component of the image of $\phi$ but disjoint to it (possible because $k<2n$) and the other disc to the basepoint component $\partial_0\vb L.$ Define $T:Q\leadsto S$ to be an elementary triad cobordism relative to $Q$ given by attaching a right handle along $\varphi$ with any $\Theta$-structure extending $\ell_Q$ (this is again possible since $B^\partial$ is path connected, as in \cref{no need top}). Since $\hb T$ is the trace of the surgery of $\hb Q$ at $\varphi,$ the belt of this handle lies in the basepoint component in $S.$ Thus $T$ is an elementary left cobordism of index $2n$ relative to $S$ attached to the basepoint component in $S$ and thus in $\W$. By \cref{rearragement,factorizations and handle structures}, there exists a factorization $M'\circ T':P\leadsto R\leadsto S$ of $T\circ M$, where $T'$ is right elementary with index $1$ relative to $P$ and $M'$ is right elementary of index $k.$ The attaching map of $T'$ in $\hb P$ is one disc in the basepoint component and another disc, thus be the same argument above, $T'$ is left elementary of index $2n$ attached to the basepoint component at $R.$ Thus, $T'\in \W.$ However, the component of the image of the attaching map of $M'$ is the basepoint component, since $T'$ performed the $1$-surgery between the component of $\phi$ and the basepoint component. Thus, $M\in \W.$ 

    When $M:P\leadsto Q$ is elementary of left type and index $n+1\leq k<2n+1$ relative to $Q$, we proceed similarly by defining $\varphi:S^0\times D^{2n-1}\hookrightarrow \hb Q$ to send one disc to the basepoint component and another to the component of the image of the attaching map of the $(k-1)$-handle in $\hb M.$ Define $T:Q\leadsto S$ to be the elementary triad cobordism given by attaching a right handle along $\varphi$ along with any $\Theta$-structure extending $\ell_Q$ (this is again possible since $B^\partial$ is path connected, as above). We obtain a factorization $T\circ M\cong M'\circ T':P\leadsto S$ where $T':P\leadsto R$ is left elementary of index $2n$ relative to $R$ and $M'$ left elementary of index $k$ relative to $S.$ Proceeding exactly as in the last paragraph, we see that the morphisms $T$ and $T'$ are in $\W$ since they are both attached to the basepoint component.  By the same argument as above, $M'$ is attached to the basepoint component and thus lies in $\W.$ This implies that $M\in \W.$

    When $M:P\leadsto Q$ is of interior type and index $2n$, we proceed similarly by defining $\varphi:S^0\times D^{2n}\hookrightarrow  Q$ to send one disc to (interior of) the basepoint component and another to the component of (but disjoint from) the image of the attaching map of the $2n$-handle in $\hb M.$ We define $T:Q\leadsto S $ to be the elementary interior cobordism of index $1$ relative to $Q$ along with any $\Theta$-structure extending $\ell_Q$ (this is now possible since $B$ is path connected, as in \cref{no need top}). Since $T$ is interior elementary of index $2n$ relative to $S$ and attached to the basepoint component, it lies in $\W.$ Proceeding similarly to before, we deduce that $M\in \W.$ This finishes the proof.  
\end{proof}

\subsection{The closure property for middle handles.}\label{middle section}

We move now to the second goal of this section: reducing the proof of \cref{stable stability} to two \textit{closure properties} of the subcategory $\W$, as mentioned in the introduction. This will be done in the next two subsections, where we will state these properties and prove that they imply \cref{stable stability}. By \cref{basepoint component statement}, it suffices to prove that certain elementary triad cobordisms are in $\W$ to show that $\W=\CobbtLs$. Roughly speaking, we will check this condition by induction on the index of the elementary triad cobordisms. In this subsection, we concern ourselves with the base case of this induction, which is proved using the first closure property (see \cref{main closure property} below).

Roughly speaking, the main strategy can be seen as a generalization of the strategy of \cref{basepoint component statement} without assuming $T$ and $T'$ are in $\W$ but having $W'\in \W$ to deduce that $W\in \W.$ The next two subsubsections introduce models for "rearrangements" of morphisms and "traces of surgeries". These notions are also used in \cref{higher section} so we phrase them in higher generality.

\subsubsection{Models for rearrangements.} We start by modelling rearrangements (see \cref{rearragement} for the terminology and also \cref{rearrangement name} below for the precise relation) of morphisms in $\CobbtL$. This is inspired by \cite[Section 4.1]{GRWII}.


\begin{defn}\label{support}
    Let $(t,W):M\leadsto N$ be a morphism in $\CobbtL.$ The support $\supp(W)$ of $W$ is the smallest closed subset $A$ of $\RR^\infty_+$ such that
    \[W\cap ([0,t]\times \RR^\infty_+\backslash A)=[0,t]\times M\backslash A\]as $\Theta$-manifolds (recall \cref{spaces of manifolds} for the terminology).
\end{defn}

\begin{figure}
\centering

\tikzset{every picture/.style={line width=0.75pt}} 

\begin{tikzpicture}[x=0.75pt,y=0.75pt,yscale=-1,xscale=1]

\draw    (165.47,85.94) .. controls (157.83,81.22) and (153.4,256.56) .. (168.69,243.98) ;
\draw  [dash pattern={on 4.5pt off 4.5pt}]  (165.47,85.94) .. controls (173.12,81.22) and (175.94,235.33) .. (168.69,243.98) ;
\draw    (165.47,85.94) .. controls (214.57,81.22) and (291.83,82.79) .. (317.58,84.37) ;
\draw    (168.69,243.98) .. controls (218.19,246.34) and (291.42,243.98) .. (317.18,241.62) ;
\draw    (317.58,84.37) .. controls (313.44,82.97) and (311.05,135.81) .. (319.33,132.1) ;
\draw    (317.58,84.37) .. controls (321.72,82.97) and (323.25,129.55) .. (319.33,132.1) ;
\draw    (316.76,205.78) .. controls (314.03,204.74) and (311.72,244.82) .. (317.18,242.07) ;
\draw    (316.76,205.78) .. controls (319.49,204.74) and (319.76,240.18) .. (317.18,242.07) ;
\draw    (317.92,158.88) .. controls (314.86,158.17) and (311.86,196.08) .. (317.98,194.2) ;
\draw    (317.92,158.88) .. controls (320.98,158.17) and (320.88,192.9) .. (317.98,194.2) ;
\draw    (243.14,82.8) .. controls (239,81.41) and (236.61,134.25) .. (244.89,130.54) ;
\draw  [dash pattern={on 4.5pt off 4.5pt}]  (243.14,82.8) .. controls (247.27,81.41) and (248.81,127.98) .. (244.89,130.54) ;
\draw    (244.89,130.54) .. controls (265.27,130.73) and (299.07,131.25) .. (319.33,132.1) ;
\draw    (243.13,157.32) .. controls (263.51,157.51) and (297.31,158.03) .. (317.57,158.88) ;
\draw    (272.51,194.35) .. controls (292.89,194.54) and (297.72,193.35) .. (317.98,194.2) ;
\draw    (272.91,205.75) .. controls (293.29,205.94) and (296.5,204.92) .. (316.76,205.78) ;
\draw    (242.65,157.36) .. controls (235.1,154.77) and (230.72,250.85) .. (245.83,243.95) ;
\draw  [dash pattern={on 4.5pt off 4.5pt}]  (242.65,157.36) .. controls (250.2,154.77) and (252.99,239.21) .. (245.83,243.95) ;
\draw    (244.89,130.54) .. controls (239.51,130.21) and (229.05,152.09) .. (242.65,157.36) ;
\draw    (272.51,194.35) .. controls (260.44,194.29) and (258.43,205.23) .. (272.91,205.75) ;
\draw [color={rgb, 255:red, 208; green, 2; blue, 27 }  ,draw opacity=1 ]   (160.24,117.53) -- (313.96,117.53) ;
\draw [color={rgb, 255:red, 208; green, 2; blue, 27 }  ,draw opacity=1 ] [dash pattern={on 4.5pt off 4.5pt}]  (170.3,108.15) -- (321.6,108.67) ;
\draw [color={rgb, 255:red, 208; green, 2; blue, 27 }  ,draw opacity=1 ] [dash pattern={on 4.5pt off 4.5pt}]  (160.24,117.53) -- (170.3,108.15) ;
\draw [color={rgb, 255:red, 208; green, 2; blue, 27 }  ,draw opacity=1 ] [dash pattern={on 4.5pt off 4.5pt}]  (239.92,117.53) -- (247.16,108.15) ;
\draw [color={rgb, 255:red, 208; green, 2; blue, 27 }  ,draw opacity=1 ]   (158.23,176.4) -- (314.76,176.4) ;
\draw [color={rgb, 255:red, 208; green, 2; blue, 27 }  ,draw opacity=1 ] [dash pattern={on 4.5pt off 4.5pt}]  (172.72,165.98) -- (319.99,167.02) ;
\draw [color={rgb, 255:red, 208; green, 2; blue, 27 }  ,draw opacity=1 ] [dash pattern={on 4.5pt off 4.5pt}]  (158.23,176.4) -- (172.72,165.98) ;
\draw [color={rgb, 255:red, 208; green, 2; blue, 27 }  ,draw opacity=1 ] [dash pattern={on 4.5pt off 4.5pt}]  (237.5,176.4) -- (246.76,166.5) ;
\draw [color={rgb, 255:red, 74; green, 144; blue, 226 }  ,draw opacity=1 ]   (159.04,190.99) -- (315.16,190.47) ;
\draw [color={rgb, 255:red, 74; green, 144; blue, 226 }  ,draw opacity=1 ] [dash pattern={on 4.5pt off 4.5pt}]  (172.72,180.57) -- (319.99,181.09) ;
\draw [color={rgb, 255:red, 74; green, 144; blue, 226 }  ,draw opacity=1 ] [dash pattern={on 4.5pt off 4.5pt}]  (159.04,190.99) -- (172.72,180.57) ;
\draw [color={rgb, 255:red, 74; green, 144; blue, 226 }  ,draw opacity=1 ] [dash pattern={on 4.5pt off 4.5pt}]  (236.3,190.99) -- (248.77,180.57) ;
\draw [color={rgb, 255:red, 74; green, 144; blue, 226 }  ,draw opacity=1 ]   (160.65,223.29) -- (313.96,223.29) ;
\draw [color={rgb, 255:red, 74; green, 144; blue, 226 }  ,draw opacity=1 ] [dash pattern={on 4.5pt off 4.5pt}]  (172.32,212.35) -- (318.38,212.35) ;
\draw [color={rgb, 255:red, 208; green, 2; blue, 27 }  ,draw opacity=1 ]   (321.6,108.67) -- (313.96,117.53) ;
\draw [color={rgb, 255:red, 74; green, 144; blue, 226 }  ,draw opacity=1 ] [dash pattern={on 4.5pt off 4.5pt}]  (160.65,223.29) -- (172.32,212.35) ;
\draw [color={rgb, 255:red, 74; green, 144; blue, 226 }  ,draw opacity=1 ] [dash pattern={on 4.5pt off 4.5pt}]  (249.98,212.35) -- (237.3,223.29) ;
\draw [color={rgb, 255:red, 208; green, 2; blue, 27 }  ,draw opacity=1 ]   (319.99,167.02) -- (314.76,176.4) ;
\draw [color={rgb, 255:red, 74; green, 144; blue, 226 }  ,draw opacity=1 ]   (319.99,181.09) -- (315.16,190.47) ;
\draw [color={rgb, 255:red, 74; green, 144; blue, 226 }  ,draw opacity=1 ]   (318.38,212.35) -- (313.96,223.29) ;
\draw  [color={rgb, 255:red, 208; green, 2; blue, 27 }  ,draw opacity=0.67 ][fill={rgb, 255:red, 208; green, 2; blue, 27 }  ,fill opacity=0.28 ] (170.3,108.15) -- (171.31,122.03) -- (172.72,165.98) -- (165.47,171.19) -- (158.23,176.4) -- (158.57,163.98) -- (158.97,148.35) -- (159.17,135.59) -- (159.81,121.09) -- (160.24,117.53) -- cycle ;
\draw  [color={rgb, 255:red, 74; green, 144; blue, 226 }  ,draw opacity=0.66 ][fill={rgb, 255:red, 74; green, 144; blue, 226 }  ,fill opacity=0.29 ] (159.04,190.99) -- (172.72,180.57) -- (172.32,212.35) -- (160.65,223.29) -- (159.37,206.01) -- cycle ;
\draw    (425.1,85.7) .. controls (417.77,80.99) and (413.53,256.32) .. (428.19,243.74) ;
\draw  [dash pattern={on 4.5pt off 4.5pt}]  (425.1,85.7) .. controls (432.43,80.99) and (435.13,235.09) .. (428.19,243.74) ;
\draw    (425.1,85.7) .. controls (472.17,80.99) and (546.25,82.56) .. (570.94,84.13) ;
\draw    (428.19,243.74) .. controls (475.64,246.1) and (545.86,243.74) .. (570.56,241.38) ;
\draw    (570.94,84.13) .. controls (566.97,82.74) and (564.69,135.58) .. (572.62,131.86) ;
\draw    (570.94,84.13) .. controls (574.91,82.74) and (576.38,129.31) .. (572.62,131.86) ;
\draw    (570.16,205.54) .. controls (567.54,204.51) and (565.32,244.58) .. (570.56,241.83) ;
\draw    (570.16,205.54) .. controls (572.78,204.51) and (573.04,239.94) .. (570.56,241.83) ;
\draw    (571.27,158.64) .. controls (568.33,157.93) and (565.46,195.85) .. (571.33,193.96) ;
\draw    (571.27,158.64) .. controls (574.2,157.93) and (574.11,192.66) .. (571.33,193.96) ;
\draw    (499.56,82.57) .. controls (495.6,81.17) and (492.97,197.02) .. (500.9,193.3) ;
\draw  [dash pattern={on 4.5pt off 4.5pt}]  (499.56,82.57) .. controls (503.53,81.17) and (504.66,190.75) .. (500.9,193.3) ;
\draw    (467.53,194.11) .. controls (487.07,194.3) and (551.9,193.11) .. (571.33,193.96) ;
\draw    (467.92,205.51) .. controls (487.46,205.71) and (550.73,204.69) .. (570.16,205.54) ;
\draw    (572.62,131.86) .. controls (567.47,131.53) and (557.44,153.41) .. (570.47,158.68) ;
\draw    (467.53,194.11) .. controls (455.96,194.05) and (454.03,204.99) .. (467.92,205.51) ;
\draw [color={rgb, 255:red, 208; green, 2; blue, 27 }  ,draw opacity=1 ]   (420.09,117.29) -- (567.47,117.29) ;
\draw [color={rgb, 255:red, 208; green, 2; blue, 27 }  ,draw opacity=1 ] [dash pattern={on 4.5pt off 4.5pt}]  (429.73,107.91) -- (574.8,108.43) ;
\draw [color={rgb, 255:red, 208; green, 2; blue, 27 }  ,draw opacity=1 ] [dash pattern={on 4.5pt off 4.5pt}]  (420.09,117.29) -- (429.73,107.91) ;
\draw [color={rgb, 255:red, 208; green, 2; blue, 27 }  ,draw opacity=1 ] [dash pattern={on 4.5pt off 4.5pt}]  (496.48,117.29) -- (503.42,107.91) ;
\draw [color={rgb, 255:red, 208; green, 2; blue, 27 }  ,draw opacity=1 ]   (418.16,176.16) -- (568.24,176.16) ;
\draw [color={rgb, 255:red, 208; green, 2; blue, 27 }  ,draw opacity=1 ] [dash pattern={on 4.5pt off 4.5pt}]  (432.05,165.74) -- (573.26,166.79) ;
\draw [color={rgb, 255:red, 208; green, 2; blue, 27 }  ,draw opacity=1 ] [dash pattern={on 4.5pt off 4.5pt}]  (418.16,176.16) -- (432.05,165.74) ;
\draw [color={rgb, 255:red, 74; green, 144; blue, 226 }  ,draw opacity=1 ]   (418.93,190.75) -- (568.63,190.23) ;
\draw [color={rgb, 255:red, 74; green, 144; blue, 226 }  ,draw opacity=1 ] [dash pattern={on 4.5pt off 4.5pt}]  (432.05,180.33) -- (573.26,180.85) ;
\draw [color={rgb, 255:red, 74; green, 144; blue, 226 }  ,draw opacity=1 ] [dash pattern={on 4.5pt off 4.5pt}]  (418.93,190.75) -- (432.05,180.33) ;
\draw [color={rgb, 255:red, 74; green, 144; blue, 226 }  ,draw opacity=1 ]   (420.47,223.06) -- (567.47,223.06) ;
\draw [color={rgb, 255:red, 74; green, 144; blue, 226 }  ,draw opacity=1 ] [dash pattern={on 4.5pt off 4.5pt}]  (431.66,212.11) -- (571.71,212.11) ;
\draw [color={rgb, 255:red, 208; green, 2; blue, 27 }  ,draw opacity=1 ]   (574.8,108.43) -- (567.47,117.29) ;
\draw [color={rgb, 255:red, 74; green, 144; blue, 226 }  ,draw opacity=1 ] [dash pattern={on 4.5pt off 4.5pt}]  (420.47,223.06) -- (431.66,212.11) ;
\draw [color={rgb, 255:red, 208; green, 2; blue, 27 }  ,draw opacity=1 ]   (573.26,166.79) -- (568.24,176.16) ;
\draw [color={rgb, 255:red, 74; green, 144; blue, 226 }  ,draw opacity=1 ]   (573.26,180.85) -- (568.63,190.23) ;
\draw [color={rgb, 255:red, 74; green, 144; blue, 226 }  ,draw opacity=1 ]   (571.71,212.11) -- (567.47,223.06) ;
\draw  [color={rgb, 255:red, 208; green, 2; blue, 27 }  ,draw opacity=0.67 ][fill={rgb, 255:red, 208; green, 2; blue, 27 }  ,fill opacity=0.28 ] (429.73,107.91) -- (430.7,121.79) -- (432.05,165.74) -- (425.1,170.95) -- (418.16,176.16) -- (418.48,163.75) -- (418.86,148.12) -- (419.06,135.35) -- (419.67,120.86) -- (420.09,117.29) -- cycle ;
\draw  [color={rgb, 255:red, 74; green, 144; blue, 226 }  ,draw opacity=0.66 ][fill={rgb, 255:red, 74; green, 144; blue, 226 }  ,fill opacity=0.29 ] (418.93,190.75) -- (432.05,180.33) -- (431.66,212.11) -- (420.47,223.06) -- (419.25,205.78) -- cycle ;
\draw    (500.83,206.04) .. controls (498.21,205.01) and (494.06,246.62) .. (499.29,243.87) ;
\draw  [dash pattern={on 4.5pt off 4.5pt}]  (500.83,206.04) .. controls (503.44,205.01) and (501.77,241.98) .. (499.29,243.87) ;
\draw [color={rgb, 255:red, 208; green, 2; blue, 27 }  ,draw opacity=1 ] [dash pattern={on 4.5pt off 4.5pt}]  (502.65,166.27) -- (497.07,176.11) ;

\draw (192.38,87.76) node [anchor=north west][inner sep=0.75pt]    {$W$};
\draw (279.63,88.24) node [anchor=north west][inner sep=0.75pt]    {$W'$};
\draw (440.14,88.52) node [anchor=north west][inner sep=0.75pt]    {$\L_W(W')$};
\draw (512.51,88.01) node [anchor=north west][inner sep=0.75pt]    {$\R_{W'}(W)$};

\end{tikzpicture}

\caption{In this picture, the red and blue shaded areas represent the supports of $W$ and $W'$ respectively. Notice that the support of $\L_W(W')$ is the blue shaded area and the support of $\R_{W'}(W)$ is the red shaded area.}
\label{rearrangement picture}

\end{figure}
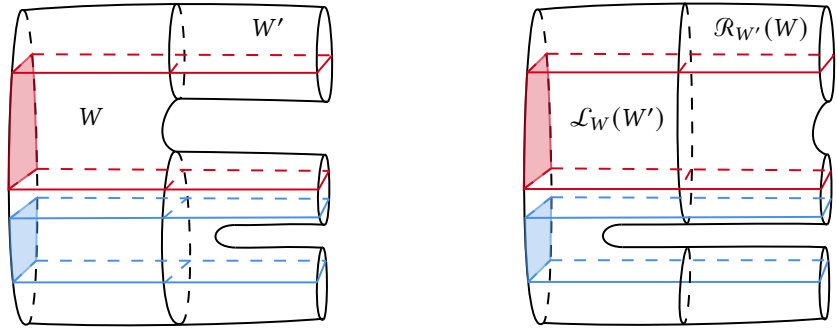

\begin{defn}\label{rearragement model}
    Let $(t,W):M\leadsto N$ and $(t',W'):N\leadsto R$ be morphisms in $\CobbtL$ such that $\supp(W)\cap \supp(W')=\emptyset.$ We define the \textit{rearrangement} (or \textit{interchange of support}) to be the morphisms
    \[\R_{W'}(W)=(W\backslash [0,t]\times \supp(W'))\cup ([0,t]\times R\backslash \supp(W))\]
    and 
    \[\L_{W}(W')=(W'\backslash [0,t']\times \supp(W))\cup ([0,t']\times M\backslash \supp(W'))\]See \cref{rearrangement picture} for an example.
\end{defn}

\begin{rmk}\label{path of rearrangement}
    Notice that $\R_{W'}(W)$ is a morphism in $\CobbtL$ from $N'\coloneqq(M\backslash \supp(W'))\cup (R\backslash\supp(W))$ to $(N\backslash \supp(W'))\cup (R\backslash\supp(W))=R,$ since $N\backslash \supp(W')=R\backslash\supp(W').$ Similarly, $\L_W(W')$ is a morphism from $M$ to $N'.$ Moreover, there exists a path \(\tau(W,W'):[0,1]\to \CobbtL(M,R)\) from $W'\circ W$ to $\R_{W'}(W)\circ \L_W(W').$ See \cite[148]{GRWII} for more details. In particular, both compositions are diffeomorphic relative to their ends.
\end{rmk}

\begin{rmk}\label{rearrangement name}
    The term interchange of support is used in \cite{GRWII}. We will also use the term \textit{rearrangement} since when $W$ and $W'$ are elementary, then $\R_{W'}(W)$ and $\L_W(W')$ are a model for the rearrangement of $W$ and $W'$ (in the sense of \cref{rearragement}). That is, $\R_{W'}(W)$ and $\L_W(W')$ are elementary of type and index of $W$ and $W'$, respectively, and $W'\circ W\cong \R_{W'}(W)\circ \L_W(W')$. More precisely, the underlying cobordism of $\R_{W'}(W)$ is obtained from $N'$ by attaching a handle at the attaching map of $W$ in $ M\backslash \supp(W')\subset N'$ and the underlying cobordism of $\L_W(W')$ is obtained from $M$ by attaching a handle at the attaching map of $W'$ in $ N\backslash\supp(W)\subset M.$ 
\end{rmk}

\subsubsection{Models for traces of surgeries.} In \cref{basepoint component statement}, we saw that it suffices to prove that certain elementary triad cobordisms are in $\W$. We now provide models for elementary triad cobordisms (recall \cref{elementary cobordisms}), which come with choices of a core, cocore, and attaching and belt spheres. Moreover, the support of these models is easy to describe. We prove below in \cref{every elementary is (reverse) trace} that any elementary triad cobordism can be realized as such model up to composition with an isomorphism in the sense of \cref{isomorphisms}.


\begin{cons}[Interior traces]\label{trace of interior surgery}
    Let $k\geq 0$ be an integer. Let $A^I\coloneqq\{0\}\times \partial D^k\times \RR^{2n+1-k}\subset [0,1]\times \RR^k\times \RR^{2n+1-k} $ and let $T^I$ be a submanifold of $[0,1]\times \RR^k\times \RR^{2n+1-k}$ seen as a triad with $\hb T^I=\emptyset$ satisfying the following properties:
    \begin{enumerate}[label=(\roman*)]
        \item $T^I$ agrees with $[0,1]\times A^I$ outside 
        \([0,1]\times D^k\times D^{2n+1-k},\)
        \item $T^I\cap ([0,\epsilon)\times \RR^{2n+1})=[0,\epsilon)\times A^I, $ for some $0<\epsilon<1,$
        \item $T^I\cap ((1-\epsilon,1] \times \RR^{2n+1}_+)=(1-\epsilon,1]\times P $ for a manifold $P$ diffeomorphic to $\RR^k\times \partial D_+^{2n+1-k}$, for some $0<\epsilon<1.$
    \end{enumerate}such that the triad $T^I$ is obtained from $A^I=T^I\cap \{0\}\times \RR^k\times \RR^{2n+1-k}$ by attaching an interior $k$-handle along the inclusion
    \[\partial D^k\times D^{2n+1-k}\hookrightarrow \partial D^k\times \RR^{2n+1-k}.\]Such submanifold exists: take the graph of the function $f:D^k\times D^{2n+1-k}\to [-1,1]$ given by $f(x,y)=\sum_{i=1}^kx_i^2-\sum_{j=1}^{2n+1-k}y_j^2,$ scale and translate it to have image $[0,1].$ In order to satisfy the properties above, one must add a $\epsilon$-collar to it. The core and cocore of $\smash{T^I}$ are the core and cocore of the interior $k$-handle in the sense of \cref{defn of handles}.
    
\end{cons}

\begin{cons}[Right traces]\label{trace of right surgery}
    Let $k\geq 0$ be an integer. Let $(A^R,\partial A^R)\coloneqq\{0\}\times \partial D^k\times (\RR^{2n+1-k}_+, \partial\RR^{2n+1-k}_+)$ and let $(T^R,\hb T^R)$ be a submanifold pair of $[0,1]\times \RR^k\times (\RR^{2n+1-k}_+,\partial \RR^{2n+1-k}_+)$ satisfying the following properties:
    \begin{enumerate}[label=(\roman*)]
        \item $(T^R,\hb T^R)$ agrees with $[0,1]\times A^R$ outside 
        \([0,1]\times D^k\times (D^{2n+1-k}_+,\partial_0D^{2n+1-k}_+),\)
        \item $(T^R,\hb T^R)\cap ([0,\epsilon)\times \RR^k\times (\RR^{2n+1-k}_+,\partial \RR^{2n+1-k}_+)) =[0,\epsilon)\times (A^R,\partial A^R), $ for some $0<\epsilon<1,$
        \item $(T^R,\hb T^R)\cap ((1-\epsilon,1]\times \RR^k\times (\RR^{2n+1-k}_+,\partial \RR^{2n+1-k}_+) )=(1-\epsilon,1]\times (P,\hb P) $ for a manifold pair $(P,\hb P)$ diffeomorphic to $(\RR^k\times \partial_1D_+^{2n+1-k},\RR^k\times \partial_{01}D_+^{2n+1-k})$, for some $0<\epsilon<1.$
    \end{enumerate}such that $(T^R,\hb T^R)$ is obtained from $(A^R,\partial A^R)$ by attaching a right $k$-handle along the inclusion
    \[(\partial D^k\times D^{2n+1-k}_+,\partial D^k\times \partial_0 D^{2n+1-k}_+)\hookrightarrow (\partial D^k\times \RR^{2n+1-k}_+,\partial D^k\times \partial\RR^{2n+1-k}_+).\]Such submanifold triad exists: we can use the already constructed $T^I$ and restrict it to the subspace $[0,1]\times \RR^k\times \RR^{2n+1-k}_+\subset [0,1]\times \RR^k\times \RR^{2n+1-k}.$ The core and cocore of $T^R$ are the core and cocore of the right $k$-handle in the sense of \cref{defn of handles}.

\end{cons}

\begin{cons}[Left traces]\label{trace of left surgery}
    Let $k\geq 0.$ Let $(A^L,\partial A^L)\coloneqq \{0\}\times(\partial_1 D^k_+,\partial_{01} D_+^k)\times \RR^{2n+1-k}$ and let $(T^L\hb T^L)$ be a submanifold pair of $[0,1]\times (\RR^k_+,\partial \RR^k_+)\times \RR^{2n+1-k}$ satisfying the following properties:
    \begin{enumerate}[label=(\roman*)]
        \item $(T^L,\hb T^L)$ agrees with $[0,1]\times (A^L,\partial A^L)$ outside 
        \([0,1]\times (D^k_+,\partial_0 D^k_+)\times D^{2n+1-k},\)
        \item $(T^L,\hb T^L)\cap ([0,\epsilon)\times (\RR^k_+,\partial \RR^k_+)\times \RR^{2n+1-k})=[0,\epsilon)\times (A^L,\partial A^L), $ for some $0<\epsilon<1,$
        \item $(T^L,\hb T^L)\cap ((1-\epsilon,1]\times (\RR^k_+,\partial \RR^k_+)\times \RR^{2n+1-k})=(1-\epsilon,1]\times (P',\partial P') $ for a manifold pair $P'$ diffeomorphic to $(\RR^k_+,\partial \RR^k_+)\times \partial D^{2n+1-k}$, for some $0<\epsilon<1.$
    \end{enumerate}such that the $(T^L,\hb T^L)$ is obtained from $(A^L,\partial A^L)$ by attaching a left $k$-handle along the inclusion
    \[(\partial_1 D^k_+\times D^{2n+1-k},\partial_{01} D_+^k\times D^{2n+1-k})\hookrightarrow (\partial_1 D^k_+\times \RR^{2n+1-k},\partial_{01} D_+^k\times \RR^{2n+1-k}).\]Such submanifold triad exists: we can use the already constructed $T^I$ and restrict it to the subspace $[0,1]\times \RR^k_+\times \RR^{2n+1-k}\subset [0,1]\times \RR^k\times \RR^{2n+1-k}.$ The core and cocore of $T^L$ are the core and cocore of the left $k$-handle in the sense of \cref{defn of handles}.
    
\end{cons}

\begin{defn}[Trace of a surgery]\label{traces}
    Let $N\in \CobbtL$ and $\sigma:\RR^k\times (\RR^{2n+1-k}_+,\partial \RR^{2n+1-k}_+) \hookrightarrow (\RR^{\infty}_+,\partial \RR^{\infty}_+)$ be an embedding of pairs such that $\sigma^{-1}(N)=A^R$ as manifold pairs. The \textit{trace} $\tr(\sigma):N\leadsto N_\sigma$ of $\sigma$ is the elementary triad cobordism of right type and index $k$ relative to $N$ given by
    \[([0,1]\times N\backslash \sigma(A^R))\cup ((\id_{[0,1]}\times \sigma)(T^R)).\]For an embedding $\sigma:\RR^k\times \RR^{2n+1-k}\hookrightarrow \RR^\infty_+\backslash \partial \RR^{\infty}_+$ such that $\sigma^{-1}(N\backslash \hb N)=A^I,$ denote by $\tr(\sigma):N\leadsto N_\sigma$ the elementary interior cobordism of index $k$ relative to $N$ contructed analogously using $T^I.$ For an embedding of pairs $\sigma:(\RR^k_+,\partial \RR^k_+)\times \RR^{2n+1-k}\hookrightarrow (\RR^\infty_+,\partial \RR^\infty_+)$ such that $\sigma^{-1}(N)=A^L$, denote by $\tr(\sigma):N\leadsto N_\sigma$ the elementary left cobordism of index $k$ relative to $N$ constructed analogously using $T^L.$ Let $r:[0,1]\times \RR^\infty_+\to [0,1]\times \RR^{\infty}_+$ given by $r(t,x)=(1-t,x).$ We define the \textit{reverse trace} $\rtr(\sigma):N_\sigma\leadsto N$ of $\sigma$ to be the image of $\tr(\sigma)$ under the diffeomorphism $r.$ In other words, $\rtr(\sigma)$ is given by \[([0,1]\times N\backslash \sigma(A^R))\cup ((\id_{[0,1]}\times \sigma)(r(T^R))).\] Denote by $\core_\sigma$ and $\cocore_\sigma$ to be the image of core and cocores of $T^R, T^I$ or $T^L$ under $\sigma.$ We call $\sigma$ an \textit{attaching map} of interior, right or left type on $N$.
\end{defn}

\begin{rmk}[How to construct $\Theta$-structures on $\tr(\sigma)$]\label{theta of traces}\label{support of traces}
    In general, the cobordism $\tr(\sigma)$ admits no preferred $\Theta$-structure making it a morphism of $\Cobbtl.$ Let us provide a general recipe to endow this cobordism with such a structure. We will focus on the case when $\sigma$ is of interior type. Let $T_0^I$ be the submanifold of $T^I$ given by $A^I\cup [0,1]\times \partial D^k\times \RR^{2n+1-k}\backslash D^{2n+1-k}.$ By pulling back the structure $\ell_N$ along $\sigma|_{\sigma^{-1}(N)}$, we can endow $T_0^I$ with a $\Theta$-structure. Thus, to endow $\tr(\sigma)$ with a $\Theta$-structure it suffices to extend the latter structure from $T_0^I$ to $T^I.$ In particular, given such an extension $\ell$, the morphism $(\tr(\sigma),\ell)$ has support at most $\sigma(D^k\times D^{2n+1-k}_+).$ One proceeds analogously for the right and left types.
    \end{rmk}

When considering $\Theta$-structures on $\tr(\sigma)$, we will always implicitely assume that they are constructed as above in \cref{support of traces}. In particular, a morphism with underlying cobordism $\tr(\sigma)$ will always have support at most $\sigma(D^k\times D^{2n+1-k}_+),$ $\sigma(D^k\times D^{2n+1-k})$ or $\sigma(D^k_+\times D^{2n+1-k})$ depending if $\sigma$ is of right, interior or left type.

\begin{lemma}\label{every elementary is (reverse) trace}
    Let $W:M\leadsto N$ be a morphism in $\CobbtL$ whose underlying cobordism is elementary (in the sense of \cref{elementary cobordisms}). Then there exists an embedding $\sigma$ of the same type (see  \cref{traces}) of $W$ relative to $M$ and an extension $\ell$ of  the $\Theta$-structure $\ell_M$ to \[\tr(\sigma):M\leadsto M_\sigma\in \CobbtL\] such that there exists an isomorphism $U:M_\sigma\leadsto N$ in $\smash{\CobbtL}$ such that $U\circ \tr(\sigma)$ and $W$ are in the same component of $\smash{\CobbtL(M,N)}.$ Dually, there exists $\sigma'$ and an extension of $\ell_N$ such that $\rtr(\sigma'):N_{\sigma'}\leadsto N$ after pre-composed with an isomorphism is in the same component of $W.$
\end{lemma}
\begin{proof}
    The proof is completely analogous to \cite[Lemma 4.5]{GRWII}.
\end{proof}

\subsubsection{The closure property.}

We move now to stating the first closure property, which will be essential for the base case of the induction alluded to in the beginning of this subsection. 

\begin{defn}\label{the surgery construction}
    Let $W:M\leadsto N$ be a morphism in $\CobbtL$ and let $\sigma$ be an attaching map on $N$ in the sense of \cref{traces} disjoint from the support of $W$ and $\ell$ be a $\Theta$-structure on $\tr(\sigma).$ Define the \textit{surgery on $W$ along $\sigma$}, denoted by \[\chi(W,\sigma,\ell):\chi(M,\sigma)\leadsto \chi(N,\sigma),\]to be the morphism $\R_{\tr(\sigma)}(W),$ which is well defined since $\supp(\tr(\sigma))$ is disjoint from the support of $W.$
\end{defn}
 
\begin{rmk}\label{the surgery is in fact surgery}
    Let $W:M\leadsto N$ be a morphism in $\CobbtL$ whose underlying triad cobordism is elementary relative to $N$, and let $\sigma$ be an attaching map in $N$ disjoint from the attaching map of the handle of $W.$ Then the underlying cobordism of $\chi(W,\sigma)$ is an elementary triad cobordism relative to $\chi(N,\sigma)$ between the results of the surgery $\chi(M,\sigma)$ and $\chi(N,\sigma)$ on $M$ and $N$ along $\sigma.$ This cobordism is built by attaching a handle of the type and index of $W$ to the attaching map of $W.$ This is well defined because $\sigma$ is disjoint from the attaching map of the handle of $W.$ In particular, the underlying cobordism $\chi(W,\sigma,\ell)$ does not depend on the choice of $\ell$ (see also \cref{rearrangement name}).
\end{rmk}

\begin{nota}[Translation of $\sigma$]
    Given an right attaching map $\sigma:(\partial \RR^k\times \RR^{2n+1-k}_+,\partial \RR^k\times \partial\RR^{2n+1-k}_+)\hookrightarrow (\RR^\infty_+,\partial \RR^\infty_+)$ in $N,$ we denote by $\sigma_p=\sigma(-,3 p\cdot e_1+(-))$ for $p\geq 1,$ where $e_1$ is the first basis vector of $\partial \RR^{2n+1-k}\cong \RR^{2n-k}.$ Observe that this is an attaching map on $N$. Similarly, for an interior or left attaching map $\sigma$, one defines $\sigma_p$ by restricting in the second factor to the disc translated by the first basis vector of $\RR^{2n+1-k}.$ 
\end{nota}

    One sees that $\sigma_1$ is disjoint from the support of $\tr(\sigma).$ Therefore, it is disjoint from the support of $\chi(W,\sigma,\ell).$ To define a $\Theta$-structure on $\tr(\sigma_p)$ we proceed similarly to \cref{support of traces}. Assume $\sigma$ is of interior type, let $T_{0,t}^I\coloneqq A^I\cup [0,1]\times \partial D^k\times \RR^{2n+1-k}\backslash (D^{2n+1-k}+t\cdot e_1)$ and pull back the structure $\ell_N$ to obtain a $\Theta$-structure $\ell_t$ on it. For every $t>0,$ translation along the $e_1$-coordinate, we have a diffeomorphism $\smash{T_0^I\cong T_{0,t}^I},$ which induces a structure $\ell_t$ on $\smash{T_0^I}.$ Assume we are given a $(0,\infty)$-parameter family\footnote{Here for a submanifold $T\in \Psi^\partial(U)$ (see \cref{spaces of manifolds}), a \textit{$(0,\infty)$-parameter family of $\Theta$-structures} is simply a path $(0,\infty)\to \Psi^\partial_{\Theta}(U)$ lifting the constant path at $T$ along the projection $\Psi^\partial_{\Theta}(U)\to \Psi^\partial(U)$ which forgets the $\Theta$-structure.} $\smash{\hat{\ell}_t}$ of $\Theta$-structures on $T$ such that $\smash{\hat{\ell}_t}$ extends $\ell_t$ in $\smash{T_0^I}.$ This induces a structure on each $\smash{\tr(\sigma_p)}$ by taking $\smash{\hat{\ell}_{3p}}$ making it into a morphism in $\smash{\CobbtL}.$ Proceed analogously for right and left types.

\begin{defn}\label{the multiple surgery}
    Let $W:M\leadsto N$ be a morphism in $\CobbtL$ and let $\sigma$ an attaching map disjoint from the support of $W$ and $\ell_t$ be a $(0,\infty)$-family of $\Theta$-structures as above. Define the \textit{$p$-th (translated) iterated surgery on $W$ along $\sigma$} to be the morphism given by the inductive formula
    \[\chi^{p}(W,\sigma,\ell)\coloneqq\chi(\chi^{p-1}(W,\sigma,\ell),\sigma_{p},\ell_{3p})\]and the initial value $\chi^0(W,\sigma,\ell)=W.$
    
\end{defn}

\begin{rmk}
    If $W$ lies in $\CobbtLs$, then so does $\chi(W,\sigma,\ell).$ This follows by using \cref{strong geometrical connectivity} and \cref{the surgery is in fact surgery}: for any handle decomposition of $W$, there exists a handle decomposition of $\chi(W,\sigma,\ell)$ with same the number of handles in each index as one of $W.$
\end{rmk}

Recall that we fix the choices of \cref{the assumptions for stable stability}.

\begin{teo}\label{main closure property}
    Assume $(B,B^\partial)$ is $1$-connected. Let $W:M\leadsto N$ be a morphism in $\CobbtLs$ and $\sigma$ an attaching map of index $n$ for right type or $n+1$ for left type disjoint from the support of $W$ and $\ell_t$ be a $(0,\infty)$-family of $\Theta$-structures as above. If $\chi^p(W,\sigma,\ell)\in \W$ for every $p\geq 1,$ then so is $W.$
\end{teo}

As mentioned before, we will defer the proof of this statement to the next section.

\subsubsection{Deducing stable stability from \cref{main closure property}.}\label{middle section ded}

In this subsubsection, we deduce the following result.  

\begin{prop}\label{middle dim stability}
    Let $W:M\leadsto N$ be a morphism whose underlying cobordism is elementary relative to $N$ of right type and index $n$ or elementary of left type and index $n+1$ attached to the basepoint component. Then $W\in \W.$ 
\end{prop}

This will be the base case in an induction proving \cref{stable stability}. We proceed in a similar way to \cite[Sections 3 and 4]{GRWII}. Roughly speaking, we first start by proving this claim directly in the case where the handles are attached \textit{trivially}. Secondly, we use \cref{main closure property} to reduce to this case, by constructing an attaching map $\sigma$ such that the handles in $\chi^p(W,\sigma,\ell)$ are trivially attached. 

\begin{figure}
    \centering
    \tikzset{every picture/.style={line width=0.75pt}} 

\begin{tikzpicture}[x=0.75pt,y=0.75pt,yscale=-1,xscale=1]

\draw    (197.56,99.73) .. controls (185.68,92.11) and (183.8,199.92) .. (197.21,201.38) ;
\draw    [dash pattern={on 4.5pt off 4.5pt}] (197.56,99.73) .. controls (201.46,95.37) and (208.75,198.43) .. (197.21,201.38) ;
\draw    (197.56,99.73) .. controls (226.84,97.4) and (299.22,99.07) .. (334.24,97.69) ;
\draw    (197.21,201.38) .. controls (225.49,201.69) and (296.42,204.66) .. (331.44,203.29) ;
\draw    (334.24,97.69) .. controls (329.96,95.01) and (323.16,139.38) .. (334.23,135.98) ;
\draw    (331.38,171.16) .. controls (327.1,168.47) and (320.37,206.69) .. (331.44,203.29) ;
\draw  [dash pattern={on 4.5pt off 4.5pt}]  (334.24,97.69) .. controls (340.97,96.89) and (339.03,135.15) .. (334.23,135.98) ;
\draw  [dash pattern={on 4.5pt off 4.5pt}]  (331.38,171.16) .. controls (337.86,169.93) and (337.41,204.91) .. (331.44,203.29) ;
\draw    (334.23,135.98) .. controls (311.7,135.32) and (294.31,172.76) .. (331.38,171.16) ;
\draw    (429.49,201.34) .. controls (439.52,200.14) and (444.16,101.26) .. (430.77,99.69) ;
\draw    (429.49,201.34) .. controls (425.55,205.67) and (419.2,102.56) .. (430.77,99.69) ;
\draw    (429.49,201.34) .. controls (400.19,203.44) and (366.47,202.18) .. (331.44,203.29) ;
\draw    (430.77,99.69) .. controls (402.49,99.16) and (369.28,96.58) .. (334.24,97.69) ;
\draw    (331.38,171.16) .. controls (353.9,171.99) and (371.3,134.67) .. (334.23,135.98) ;
\draw [color={rgb, 255:red, 208; green, 2; blue, 27 }  ,draw opacity=1 ]   (326.25,183.97) .. controls (244.62,183.31) and (257.21,107.73) .. (329.13,112.35) ;
\draw [color={rgb, 255:red, 74; green, 144; blue, 226 }  ,draw opacity=1 ]   (328.41,122.58) .. controls (371.2,114) and (405.36,181.66) .. (327.33,179.35) ;
\draw [color={rgb, 255:red, 208; green, 2; blue, 27 }  ,draw opacity=1 ] [dash pattern={on 4.5pt off 4.5pt}]  (329.13,112.35) .. controls (364.01,111.36) and (391.34,111.69) .. (426.22,111.69) ;
\draw [color={rgb, 255:red, 208; green, 2; blue, 27 }  ,draw opacity=1 ] [dash pattern={on 4.5pt off 4.5pt}]  (326.25,183.97) .. controls (361.13,182.98) and (390.62,183.31) .. (425.5,183.31) ;
\draw [color={rgb, 255:red, 65; green, 117; blue, 5 }  ,draw opacity=1 ][line width=1.5]    (426.22,111.69) .. controls (423.7,124.89) and (423.7,164.16) .. (425.5,183.31) ;

\draw (206.44,108.76) node [anchor=north west][inner sep=0.75pt]    {$W\cong \rtr(\sigma)$};
\draw (384.34,123.94) node [anchor=north west][inner sep=0.75pt]    {$\tr(\sigma) $};
\draw (186.16,74.43) node [anchor=north west][inner sep=0.75pt]    {$M$};
\draw (326.82,73.68) node [anchor=north west][inner sep=0.75pt]    {$N$};
\draw (418.61,72.12) node [anchor=north west][inner sep=0.75pt]    {$N_{\sigma }$};

\end{tikzpicture}

    \caption{In this picture, we see a morphism $W:M\leadsto N$ whose underlying cobordism is elementary. We represent its core relative to $N$ as the red arc. Observe that if perform surgery on a translation of the attaching map of $W$, the attaching map becomes trivial. In the picture, we see that the attaching map of $W$ seen in $N_\sigma$ bounds a disc, here depicted in green. However by \cref{the surgery construction}, $\chi(W,\sigma,\ell)$ can be obtained by taking $N_\sigma$ and attaching a handle precisely along this attaching map. Thus, it is a trivial handle attachment.}
    \label{figure reaarang middle}
\end{figure}
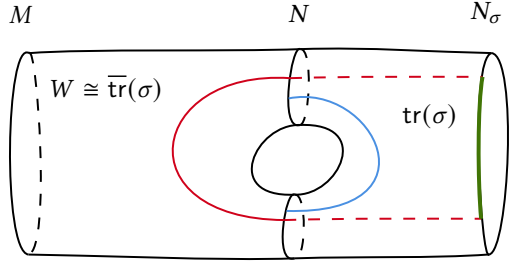

Let us give a few more words about this construction. We will construct $\sigma$ to be an attaching map in $N$ of the same type and index $W$ where $\sigma|_{\sigma^{-1}(N)}$ is "parallel" to an attaching map of the handle of $W$. By \cref{the surgery is in fact surgery}, $\chi(W,\sigma,\ell)$ is elementary of the type of $W$ with unique handle attached to $\chi(N,\sigma).$ However, since surgery along its attaching map is performed, this handle is attached trivially. The same holds for $\chi^p(W,\sigma,\ell).$ Based on this idea, we shall use \cref{main closure property} to deduce stable stability for $W.$ (See \cref{figure reaarang middle}). 

\paragraph{\textit{Stability for $V_1$.}}

We start by proving that $\W$ contains a certain restricted class of morphisms. Let us first define the following notion.

\begin{defn}[Triad connected sum]\label{triad connect sum}
    Let $(W,\hb W,\vb W)$ and $(W',\hb W',\vb W')$ be two $d$-dimensional manifold triads such that $\vb W$ and $\vb W' $ are non-empty. Let $e:(D^{d-1}_+,\partial_0D^{d-1}_+)\hookrightarrow (\vb W,\hvb W)$ and $e':(D^{d-1}_+,\partial_0D^{d-1}_+)\hookrightarrow (\vb W',\hvb W')$ be two embeddings. Define the \textit{triad connected sum} $W\natural W'$ as the triad obtained by attaching a right $1$-handle to $W\sqcup W'$ along $e\sqcup e'.$ The result on horizontal and vertical boundaries are boundary connected sums. 
\end{defn}

\begin{rmk}If $\hvb W$ and $\hvb W$ are connected, then the diffeomorphism type of $W\natural W'$ only depends on the local orientations induced by the embeddings $e$ and $e'$. If either $W$ or $W'$ admit a orientation-reversing diffeomorphism, then the diffeomorphism type of the triad $W\natural W'$ is independent of $e$ and $e'$. For example, this is the case for the triad $V_1$ (see the definition in the introduction or above \cref{right standard}). Unless otherwise specified, if $W$ is a cobordism from $P$ to $Q$ then we assume that the connected sum is done in $\partial Q$ and the result is a cobordism from $P$ to $Q\natural \vb W.$\end{rmk}

Let $P\in \CobbtL,$ let \[H_P:P'\leadsto P\]be any morphism in $\CobbtL$
whose underlying triad cobordism is given by taking the triad connected sum (see \cref{triad connect sum}) of $P\times [0,1]$ with $V_1$ on the basepoint component (recall \cref{basepoint component def}) in $P\times \{0\} $ away from $\vb L$ and whose $\Theta$-structure is such that the restriction to $V_1$ is standard in the sense of \cref{standard structure on V1}. More explicitly, the underlying triad cobordism of $H_P$ is $(P\times [0,1]\natural V_1, \hb P\times [0,1]\natural W_{1,1}, \vb L\times [0,1], P\times \{0\}\natural \vb V_1\sqcup P\times \{1\})$, so in particular $P'\cong P$. This $\Theta$-structure exists since we are in the basepoint component and by \cref{the space of standard is connected}. This is in $\CobbtLs$ since the pair $(V_1,\hb V_1)$ is $n$-connected. Similarly, let \(\prescript{}{P}{H}:P\leadsto P'\) be any morphism whose underlying triad cobordism is given by the triad connected sum of $P\times [0,1]$ with $V_1$ along $P\times\{1\}$ in the basepoint component with standard $\Theta$-structure. 

\begin{prop}\label{stability for V_1}
    For any $P\in \CobbtL,$ then $H_P$ and $\prescript{}{P}{H}$ are in $\W.$  
\end{prop}
\begin{proof}
    This is proven analogously to \cite[Thm. 3.1]{GRWII} (see also \cite[Lemma 7.15]{GRWStableModuli} for the same argument) by using \cref{the space of standard is connected} in place of \cite[Lemma 2.12]{GRWII}.
\end{proof}

\paragraph{\textit{Reflections of cobordisms.}} We will now see that the morphism $\prescript{}{P}{H}$ can be expressed as a composition $W\circ r(W)$ of an morphism $W$ whose underlying triad cobordism is elementary and its "reflection" $r(W)$, where the unique handle of $W$ is attached trivially. We shall apply \cref{main closure property} to a general morphism $M$ whose underlying cobordism is elementary and its reflection to reduce stable stability to the trivially attached case, i.e. \cref{stability for V_1}. We start by making the notion of the reflection of a $\Theta$-cobordism precise.

\begin{cons}[Reflection]\label{reflection}
    Let $W:M\leadsto N$ be $\rtr(\sigma)$ for some $\sigma:\RR^n\times \RR^{n+1}_+\hookrightarrow \RR^\infty_+$ (which we call of \textit{right type}) or $\sigma:\RR^{n+1}_+\times \RR^n\hookrightarrow \RR^{\infty}_+$ (which we call of \textit{left type}) as in \cref{traces}. Let $r:[0,1]\times \RR^\infty_+\to [0,1]\times \RR^\infty_+$ given by $r(t,x)=(1-t,x).$ Let $r(W):N\leadsto M$ whose underlying submanifold of $[0,1]\times \RR^\infty_+$ is $r(W)$. The manifold $r(W)$ is left (resp. right) elementary of index $n+1$ (resp. $n$) relative to $M$ if $\sigma$ is of right (resp. left)  type. More specifically: if $\sigma $ is of right type, the manifold $W$ contains the single handle
\[\cocore_\sigma=(\id_{[0,1]}\circ \sigma)(\cocore_{T^R}):(D^{n+1}_+\times \RR^{n},\partial_0 D^{n+1}_+\times \RR^{n})\to W\]relative to $M$. If $\sigma $ is of left type, the manifold $M$ contains the single handle
\[\cocore_\sigma=(\id_{[0,1]}\circ \sigma)(\cocore_{T^L}):(D^n\times \RR^{n+1}_+,D^n\times \partial \RR^{n+1}_+)\to W\]relative to $M$. Thus, it suffices to determine a $\Theta$-structure on $r(\cocore_\sigma)$ which extends $\ell_M$ in $r(\belt_\sigma)\subset M,$ where $\belt_\sigma=\cocore_\sigma\cap N.$ In the right case, we define a $\Theta$-structure on $r(\cocore_\sigma)$ by insisting that the union
\[(D^{n+1}\times D^{n},\partial D^{n+1}\times D^{n})\cong r(\cocore_\sigma)(D^{n+1}_+\times D^n)\cup_M \cocore_\sigma (D^{n+1}_+\times D^n)\subset r(W)\cup_M W \]has a standard $\Theta$-structure (see below the proof of existence of such structure). In the left case, we define a $\Theta$-structure on $r(\cocore_\sigma)$ by insisting that the union
\[(S^n\times D^{n+1}_+,S^n\times \partial_0D^{n+1}_+)\cong r(\cocore_\sigma)(D^{n}\times D^{n+1}_+)\cup_M \cocore_\sigma (D^{n}\times D^{n+1}_+)\subset r(W)\cup_M W\]has a standard $\Theta$-structure (see below the proof of existence of such structure).
\end{cons}
\begin{proof}[Proof of existence of such structure]
    We prove that there exists a $\Theta$-structure on $r(\cocore_\sigma)$ such that the union above has standard $\Theta$-structure. We will focus on the first case, as the second is completely analogous. Pick a standard $\Theta$-structure $\ell^{\std}$ on $(D_R,\hb D_R)\coloneqq (D^{n+1}\times D^n, \partial D^{n+1}\times D^n)$ and let $\ell\coloneqq  (\cocore_\sigma|_{D_R})^*\ell_M$. We start by observing that $\ell^{\std}$ and $\ell$ are homotopic as bundle maps. We will show that $\ell$ is homotopic to the basepoint $\Theta$-structure fixed in \cref{the assumptions for stable stability}, as $\ell^{\std}$ is by definition. This will use that $\sigma$ is in the basepoint component. We prove in the case that $(\theta^*\gamma_{2n+1},(\theta^\partial)^*\gamma_{2n})$ is orientable, in which case all $\Theta$-manifolds are orientable. If this assumption is not satisfied, there exists only one $\Theta$-structure on $D_R$ as $B^\partial$ is path-connected so the claim follows trivially. By assumption, there exists a path in $\hb W$ from a point in the image of $\cocore_\sigma$ and a point in the basepoint component $\vb_0 L.$ Thus, for any embedding $e$ of $D_R$ inducing the same orientation as $\cocore_\sigma$ into the basepoint component, we have that $e^*\ell_{\vb L}$ is homotopic to $\ell$. Moreover, by the definition of the basepoint component (see \cref{basepoint component def}), we see that there exists a path from it to an embedded copy of $V_1$ with standard $\Theta$-structure. One key feature of the standard $\Theta$-structure is that when restricted to an orientation preserving embedding of the disc it is the basepoint $\Theta$-structure. We assume that the embedding $\sigma$ induces the same orientation as the one induced by the embedding of this disc after composing with the embedding of the copy of $V_1$. (If this is not the case, we can change $\sigma$ by an orientation reversing diffeomorphism, since it induces the same underlying manifold after taking $\rtr.$) We conclude that $e^*\ell_{\vb L} $ is homotopic to the basepoint $\Theta$-structure, which implies that $\ell$ is homotopic to the basepoint $\Theta$-structure. By the (bundle) homotopy extension property of the pair $(S^n\times D^n,D_R)$, we can extend the homotopy from $\ell$ to the basepoint $\Theta$-structure from $D_R$ to a homotopy from a $\Theta$-structure extending $\ell$ to $S^n\times D^n$ to a standard $\Theta$-structure.  
\end{proof}

\begin{rmk}
    As in \cref{support of traces}, the support of $\tr(\sigma)$ is $\sigma(D^n\times D^{n+1}_+)$ or $\sigma(D^{n+1}_+\times D^{n})$ if $\sigma$ is of right or left type. It is clear to see $\supp(M)=\supp(r(M)).$
\end{rmk}

\paragraph{\textit{Identifying $\chi^p(W,\sigma,\ell)$.}}


Let $W:M\leadsto N$ be $\rtr(\sigma)$ as in \cref{reflection}, where $\sigma$ lies in the basepoint component. It is easy to see that $\supp(W)$ is disjoint from $\supp(\tr(\sigma)).$ Moreover, it is also clear that $\supp(W)$ is disjoint from $\supp(\tr(\sigma_p))$ (recall \cref{the multiple surgery}) for all $p\geq 1$. The trace $\tr(\sigma)$ is elementary of the type of $\sigma$ relative to $N.$ We fix a $(0,\infty)$-family of $\Theta$-structures $\ell_t$ as in \ref{the multiple surgery} by insisting, as before, that the union of the core of $W$ relative to $N$ with the core of $\tr(\sigma)$ relative to $N$ is standard.

\begin{lemma}\label{after surgery it is H}
    Let $W =\rtr (\sigma)$ be as above. Then $\chi^p( W\circ r(W),\sigma,\ell)$ is the composition of $H_{\chi^p(N,\sigma,\ell)}$ with an isomorphism, for every $p\geq 1$. In particular, it is in $\W.$
\end{lemma}
\begin{proof}

    Let $X\coloneqq\chi^p(W\circ r(W),\sigma,\ell),$ we show that there exists an embedding of pairs $h:(V_1,W_{1,1})\hookrightarrow (X,\hb X)$ with standard induced $\Theta$-structure and a path $\gamma$ from a point in $h(W_{1,1})$ to a point in $\hb \chi^p(N,\sigma,\ell)$ in $\hb X$ disjoint from $h$ in its interior such that the complement of $h$ and a tubular neighborhood of $\gamma$ (i.e. an embedding $\widetilde{\gamma}:([0,1]\times D^{2n}_+,[0,1]\times \partial_0 D^{2n}_+)\to (X,\hb X)$ such that $\widetilde{\gamma}|_{[0,1]\times\{0\}}=\gamma$) is diffeomorphic to an interval cobordism (and thus, an isomorphism by \cref{isomorphisms}). We also note that the complement of $h$ and a tubular neighborhood of $\gamma$ is diffeomorphic to the union of the complement of $h$ with a disc $(D^{2n+1}_+,\partial_1D^{2n+1}_+)$ along $h(\vb V_1,\hvb V_1)$. 

    Assume $\sigma$ is of right type. The union of $\cocore_\sigma$ with $r(\cocore_\sigma)$ induces an embedding
    \[f':(D^{n+1}\times D^{n},\partial D^{n+1}\times D^n)\hookrightarrow r(W)\cup_P W\]with standard $\Theta$-structure, by definition. Since the image of $f'$ is inside the support of $W,$ it is disjoint from the support of $\tr(\sigma_p).$ The rearrangement (recall \cref{rearragement model}) defines an embedding $f$ into $X.$ (Moreover this embedding is isotopic to $f'$ in $X\circ \L_{W\circ r(W)}(\tr(\sigma_p))=\tr(\sigma_p)\circ W\circ r(W)$, by \cref{path of rearrangement}). The core of $W=\rtr(\sigma)$ is disjoint from the support of $\tr(\sigma_p)$ and intersects the cocore of $W$ transversely at one point (by definition of handle attachment). In other words, it intersects the image of $f'$ transversely at one point and thus, its image after the rearrangement will have the same intersection behavior with $f.$ However, the embedding of the attaching map of $W$ can be seen also in $\chi^p(N,\sigma,\ell)$ giving an embedding
    \[g:(S^{n-1}\times D^{n+1}_+,S^{n-1}\times \partial_0D^{n+1}_+)\hookrightarrow (\chi^p(N,\sigma,\ell),\hb \chi^p(N,\sigma,\ell))\]since it is disjoint from the support of $\tr(\sigma_p).$ However, this embedding extends to an embedding
    \[e':(D^n\times D^{n+1}_+,D^n\times \partial_0D^{n+1}_+)\hookrightarrow (\chi^p(N,\sigma,\ell),\hb \chi^p(N,\sigma,\ell))\]since the target is the result of right surgery at a translate of $g.$ Moreover, the induced $\Theta$-structure on $e$ after glued with the one of the core of $W$ (after rearrangement) is standard, by definition. This produces an embedding 
    \[e:(S^n\times D^{n+1}_+,S^n\times \partial_0D^{n+1}_+)\hookrightarrow (X,\hb X)\]with standard $\Theta$-structure. By the discussion above, the cores of $e$ and $f$ intersect transversely at one point in the horizontal boundary. Therefore, the union of these embeddings is diffeomorphic to $(V_1,W_{1,1})$ with standard $\Theta$-structure. Thus, we define $h$ as the gluing of these embeddings. We can choose a path $\gamma$ in $\hb X$ connecting $\chi^p(N,\sigma,\ell)$ and this embedding and disjoint from such image, since $n\geq 2.$
    

    We finish the proof by showing that the union of the complement of $h$ with a disc $(D^{2n+1}_+,\partial_1D^{2n+1}_+)$ along $h(\vb V_1,\hvb V_1)$ is diffeomorphic to an interval cobordism\footnote{We thank Oscar Randal-Williams for explaining the argument underlying the analagous claim in the proof of \cite[Lemma 4.8]{GRWII}, which inspired the present proof.}. We claim that this complement is diffeomorphic to
    \[(Y,\hb Y)\coloneqq (X\backslash f'(D^{n+1}\times \int(D^{n}))\cup D^{n+2}_+\times \partial D^n, \hb X\backslash f'(\partial D^{n+1}\times D^n)\cup \partial_1D^{n+2}_+\times \partial D^n)\]which is the result of doing surgery along $f'$ of left type and index $n+1$, that is, the outgoing boundary of an elementary triad cobordism of left type and index $n+2$ relative to $(X,\hb X)$. This can be seen from the fact that $(V_1,W_{1,1})$ is obtained as the outgoing boundary of an elementary triad cobordism relative to $(D^{2n+1}_+,\partial_1D^{2n+1}_+)$ of left type and index $n+2$ along the embedding $\bar{f}$ from \cref{standard structure on V1}. To show now that $(Y,\hb Y)$ is the trivial triad cobordism, we show that the inverse surgery move applied to the trivial cobordism recovers $(X,\hb X).$ We use the fact that $(X,\hb X)$ is the result of composing an elementary triad cobordism $W$ of right type and index $n$ relative to $\chi^p(N,\sigma,\ell)$ with its reflection $r(W)$, which can be seen by unravelling definitions, and that $f'$ is the inclusion of the union of the cocore of $W$ with its reflection in $r(W)$, together with the following claim. Let $d,k\geq 0$ be integers.

    \begin{claim}\label{claim that GRW didn't prove}
        Let $W:R\leadsto S$ be an elementary triad cobordism of dimension $d$ of right type and index $k$ relative to $S$. Let $\phi:(\partial D^k\times D^{d-k}_+,\partial D^k\times \partial_1D^k_+)\hookrightarrow (S,\hb S)$ be an attaching map of the unique handle of $W$. Then $W\circ r(W)$ is the result of doing surgery of right type and index $k-1$ to $S\times [-2,2]$ along the embedding $\tilde{\phi}$ 
        \[\begin{tikzcd}
            (\partial D^k\times D^{d-k+1}_+, \partial D^k\times \partial_1D^{d-k+1}_+) \arrow[r, "\cong"] & (\partial D^k\times D^{d-k}_+\times [-1,1], \partial D^k\times \partial_1D^{d-k}_+\times [-1,1]) \arrow[d, "{\phi\times[-1,1] }"] \\
            & ~~~~~~~~~~~~~~~~ S\times [-1,1]\subseteq S\times [-2,2] 
        \end{tikzcd}\]
    \end{claim}
    \begin{proof}[Proof of \cref{claim that GRW didn't prove}]
        By definition, the result of doing surgery of right type and index $k-1$ to $S\times [-2,2]$ along the embedding $\tilde{\phi}$ is the outgoing boundary of the elementary triad cobordism relative to $S\times [-2,2]$ of right type and index $k$ given by attaching a handle along $\tilde{\phi}$. Unravelling definitions (see \cref{defn of handles}), we see that this result is the pair
        \[(S\times [-2,2]\backslash \tilde{\phi}(\partial D^k\times D^{d-k+1}_+) \cup D^k\times \partial_0 D^{d-k+1}_+,\hb S\times [-2,2]\backslash \tilde{\phi}(\partial D^k\times D^{d-k+1}_+) \cup D^k\times \partial_{01} D^{d-k+1}_+).\]We will see this pair as a union of two pairs $M_0$ and $M_1$, where $M_0\cong W$, $M_1\cong r(W)$ and $M_0\cap M_1\cong R$, where the last identification is compatible with the previous two. We take $M_0$ as 
        \[(S\times [-2,0]\backslash \tilde{\phi}(\partial D^k\times D^{d-k+1}_{+-}) \cup D^k\times \partial_0 D^{d-k+1}_{+-},\hb S\times [-2,0]\backslash \tilde{\phi}(\partial D^k\times D^{d-k+1}_{+-}) \cup D^k\times \partial_{01} D^{d-k+1}_{+-})\]where $D^{d-k+1}_{+-}$ is the subspace of $D^{d-k+1}_+\cong D^{d-k}_+\times [-1,1]$ of the points where the last coordinate is non-positive. One can see that this is diffeomorphic to $(S\times [-2,0] \cup_\phi D^k\times \partial_0 D^{d-k+1}_{+-},\hb S\times [-2,0] \cup_\phi D^k\times \partial_{01} D^{d-k+1}_{+-}),$ hence diffeomorphic to $W$. We take $M_1$ as \[(S\times [0,2]\backslash \tilde{\phi}(\partial D^k\times D^{d-k+1}_{++}) \cup D^k\times \partial_0 D^{d-k+1}_{++},\hb S\times [0,2]\backslash \tilde{\phi}(\partial D^k\times D^{d-k+1}_{++}) \cup D^k\times \partial_{01} D^{d-k+1}_{++})\]where $D^{d-k+1}_{++}$ is the subspace of $D^{d-k+1}_+\cong D^{d-k}_+\times [-1,1]$ of the points where the last coordinate is non-negative. The same argument as before identifies $M_1$ with $r(W)$. The union $M_0\cup M_1$ is thus diffeomorphic to $W\circ r(W)$, which completes the proof.
    \end{proof}This finishes the case when $\sigma$ is of right type. The case when $\sigma$ is of left type follows analogously by dualizing the arguments above.
\end{proof}

We are now ready to prove \cref{middle dim stability}.

\begin{proof}[Proof of \cref{middle dim stability}.]
    By \cref{every elementary is (reverse) trace}, it suffices to assume that $W$ is $\rtr(\sigma)$ for some $\sigma,$ where $\sigma$ lies in the basepoint component. Let $\ell$ be the constant path on the chosen $\Theta$-structure for $\tr(\sigma).$ By combining \cref{main closure property} and \cref{after surgery it is H}, we see that $W\circ r(W)\in \W.$ Since the reflection of an elementary of right type and index $n$ is an elementary of left type and index $n+1$ and vice-versa, we have that $r(W)\circ r(r(W))\in \W.$ Consider now the equality $W\circ (r(W)\circ r(r(W)))=(W\circ r(W))\circ r(r(W)).$ One concludes that all morphisms induce isomorphisms on $\H_1(-,\ZZ).$ On the other hand, pre-composition with $r(W)$ induces surjections on homology with abelian coefficients from its target, since $W\circ r(W)\in \W$, and injections on homology with abelian coefficients on the target of $r(r(W))$, since $r(W)\circ r(r(W))\in \W.$ Since the abelian coefficients in all three spaces are the same (as the maps above induce isomorphisms on $\H_1(-,\ZZ)$), it follows that $r(W)\in \W.$ Thus, by \cref{two out of three}, we see $W\in \W.$ 
\end{proof}

\subsection{The closure property for higher handles.}\label{higher section}

We move now to stating the main closure property \cref{main closure property high} for higher dimensional handles, which will be the main input for the induction step to prove \cref{stable stability}. This strategy is very similar to the previous subsection but instead of parallel surgery, we apply "meridian" surgery. To make this precise, we consider the following definition. Let $\sigma$ be an attaching map to $N$ of one of the following types:
\begin{enumerate}[label=(\textit{\alph*})]
    \item\label{bound 1}  $\sigma:\RR^k\times \RR^{2n+1-k}_+\hookrightarrow \RR^{\infty}_+$ for $n<k<2n;$ \hspace*{\fill}(Right)
    \item $\sigma:\RR^{k}_+\times \RR^{2n+1-k}\hookrightarrow \RR^\infty_+$ for $n+1<k<2n+1;$ \hspace*{\fill}(Left)
    \item\label{bound 3} $\sigma:\RR^k\times \RR^{2n+1-k}\hookrightarrow \RR^\infty_+\backslash \partial \RR^{\infty}_+$ for $n+1<k<2n+1.$ \hspace*{\fill}(Interior)
\end{enumerate}Let $\iota:D^{k-1}\to \partial D^k$ for the embedding onto the lower hemisphere which is inverse to the stereographic projection $(y_1,\cdots, y_{k})\mapsto \frac{1}{1-y_k}(y_1,\cdots, y_{k-1}).$ Consider $D^k$ as a subspace of $D^{k-1}\times D^1$ with coordinates $(y,z)=(y_1,\cdots, y_{k-1},z)$. Consider the family of embedding
\[\mu_t:\partial D^{2n+1-k}\times D^k\hookrightarrow \partial D^k\times \RR^{2n+1-k}\]given by $\mu_t(x;y,z)=(\iota(y),t(1+\frac{1}{3}z)x)$ for $t\in (2,\infty).$ From this family, we will only need the following properties:
\begin{enumerate}[label=(\textit{\roman*})]
    \item The image $\mu_t(\partial D^{2n+1-k}\times \{0\})$ is $\{-e_k\}\times t\cdot \partial D^{2n+1-k}$, the \textit{meridian sphere} of $\sigma$ to $e_k\in \partial D^k$ of radius $t,$ and the image of $\partial D^{2n+1-k}\times D^k$ lies in $\partial D^k\times \RR^{2n+1-k}\backslash D^{2n+1-k}.$
    \item If $t'> 2t,$ then the images of $\mu_t$ and $\mu_{t'}$ are disjoint.
    \item For all $t$, we have $\mu_t^{-1}(\partial_1D^k_+\times \RR^{2n+1-k})=\partial D^{2n+1-k}\times D^k_+$ and $\mu_t^{-1}(\partial D^k\times \RR^{2n+1-k}_+)=\partial_1D^{2n+1-k}_+\times D^k.$
\end{enumerate}

\begin{defn}[The family of meridian attaching maps]\label{definition of phi}
For the right type, let $\chi:\partial D^k\times \RR^{2n+1-k}_+\to \RR^k\times \RR^{2n+1-k}_+$ be the inclusion induced by $\partial D^k\hookrightarrow \RR^k.$ For the left type, let $\chi:\partial_1 D^k_+\times \RR^{2n+1-k}\to \RR^k_+\times \RR^{2n+1-k}$ be the inclusion induced by $\partial_1 D^k_+\hookrightarrow \RR^k_+.$ For interior type, let $\chi:\partial D^k\times \RR^{2n+1-k}\to \RR^k\times \RR^{2n+1-k}$ be the inclusion induced by $\partial D^k\hookrightarrow \RR^k.$ For each of the types above, denote $\phi_t^\sigma=\sigma\circ \chi\circ \mu_t$ using property $(iii)$ of $\mu$ above. We can extend $\phi_t^\sigma$ to an attaching map disjoint from support of $\rtr(\sigma)$ of index $2n+1-k$ of left, right and interior type for $\sigma$ of right, left and interior type and index $k$, respectively such that the images of the extensions of $\phi_t^\sigma$ and $\phi_{t'}^\sigma$ are disjoint if $t'>2t$ by $(ii).$ For simplicity, we denote such family of extensions also by $\phi_t^\sigma.$ 
\end{defn}
\begin{rmk}[$\Theta$-structures for meridian surgery]\label{theta for expanded phi}
    To define $\Theta$-structures on $\tr(\phi_t^\sigma),$ we can consider the extension $\hat{\mu}_t:D^{2n+1-k}\times D^k\hookrightarrow [0,2]\times \partial D^k\times \RR^{2n+1-k}$ of $\mu_t$ given by $\hat{\mu}_t(x,y,z)=(\frac{3}{2}(1-|x|^2)(1+\frac{1}{3}z),\iota(y),t(1+\frac{1}{3}z)x).$ Let $\ell$ be a $\Theta$-structure on $[0,2]\times \partial D^k\times \RR^{2n+1-k}$ if $\sigma $ is of interior type, on $[0,2]\times \partial D^k\times \RR^{2n+1-k}_+$ if $\sigma$ is of right type and on $[0,2]\times \partial_1 D^k_+\times \RR^{2n+1-k}$ if $\sigma$ is of left type. We can extend $\mu_t^*\ell$ to get a $\Theta$-structure on $\tr(\phi_t^\sigma).$ Moreover, the space of such extensions is contractible.
\end{rmk}

\begin{defn}\label{the multiple expanded surgery}
    Let $W:M\leadsto N$ be a morphism in $\CobbtL$ and let $\sigma$ an attaching map whose support contains the support of $W$ and $\phi_t^\sigma$ be as in \cref{definition of phi}. Let $\ell$ be a $\Theta$-structure as in \ref{theta for expanded phi}. Define the \textit{$p$-th iterated (meridian) surgery on $W$ along $\phi$} $\chi^p(W,\phi^\sigma,\ell)$ to be given by the inductive formula
    \[\chi^{p}(W,\phi^\sigma,\ell)\coloneqq\chi(\chi^{p-1}(W,\phi^\sigma,\ell),\phi_{3p}^\sigma,\mu_{3p}^*\ell)\]and the initial value $\chi^0(W,\phi^\sigma,\ell)=W.$ This is well defined since $\phi_{3p}^\sigma$ is disjoint from the support of $\sigma$.
    
\end{defn}

We state now the second closure property we use. This will play the same role as \cref{main closure property} but now for the induction step in the proof of \cref{stable stability}. We will finish by assembling all the pieces and giving the proof of the later. Recall that we conisder the choices fixed in \cref{the assumptions for stable stability}.

\begin{teo}\label{main closure property high}
    Assume $(B,B^\partial)$ is $0$-connected. Let $W:M\leadsto N$ be a morphism in $\CobbtLs$ and $\sigma$ an attaching map of index $k$ satisfying the bounds of \ref{bound 1}--\ref{bound 3} above for each type whose support contains the support of $W$ and $\ell$ be a $\Theta$-structure as in \ref{theta for expanded phi}. If $\chi^p(W,\phi^\sigma,\ell)\in \W$ for every $p\geq 1,$ then so is $W.$
\end{teo}


\subsubsection{Deducing stable stability from \cref{main closure property high}.}\label{higher section ded}

\begin{figure}
    \centering
    \tikzset{every picture/.style={line width=0.75pt}} 

\begin{tikzpicture}[x=0.75pt,y=0.75pt,yscale=-1,xscale=1]

\draw    (150.96,120.55) .. controls (141.3,114.47) and (139.78,200.64) .. (150.68,201.81) ;
\draw  [dash pattern={on 4.5pt off 4.5pt}]  (150.96,120.55) .. controls (154.13,117.07) and (160.06,199.45) .. (150.68,201.81) ;
\draw    (150.96,120.55) .. controls (174.76,118.7) and (233.59,120.03) .. (262.05,118.93) ;
\draw    (150.68,201.81) .. controls (173.66,202.06) and (231.31,204.44) .. (259.77,203.33) ;
\draw    (262.05,118.93) .. controls (258.57,116.78) and (253.04,152.25) .. (262.04,149.53) ;
\draw [color={rgb, 255:red, 65; green, 117; blue, 5 }  ,draw opacity=1 ]   (259.72,177.65) .. controls (256.24,175.5) and (250.77,206.05) .. (259.77,203.33) ;
\draw    (262.05,118.93) .. controls (267.52,118.28) and (265.94,148.87) .. (262.04,149.53) ;
\draw [color={rgb, 255:red, 65; green, 117; blue, 5 }  ,draw opacity=1 ]   (259.72,177.65) .. controls (264.99,176.67) and (264.63,204.64) .. (259.77,203.33) ;
\draw    (262.04,149.53) .. controls (243.73,149.01) and (229.6,178.93) .. (259.72,177.65) ;
\draw    (150.68,201.81) .. controls (126.86,203.49) and (72.46,200.05) .. (43.99,200.94) ;
\draw    (150.96,120.55) .. controls (127.98,120.13) and (74.74,115.64) .. (46.27,116.53) ;
\draw    (43.94,175.25) .. controls (62.24,175.92) and (75.98,145.68) .. (45.85,146.74) ;
\draw    (46.27,116.53) .. controls (42.79,114.38) and (36.85,149.45) .. (45.85,146.74) ;
\draw  [dash pattern={on 4.5pt off 4.5pt}]  (45.86,116.13) .. controls (51.32,115.48) and (49.75,146.07) .. (45.85,146.74) ;
\draw    (43.12,175.25) .. controls (39.64,173.11) and (34.99,203.65) .. (43.99,200.94) ;
\draw  [dash pattern={on 4.5pt off 4.5pt}]  (43.12,175.25) .. controls (48.4,174.28) and (48.03,202.24) .. (43.17,200.94) ;
\draw [color={rgb, 255:red, 208; green, 2; blue, 27 }  ,draw opacity=1 ]   (62.28,162.93) .. controls (55.78,173.32) and (111.04,175.72) .. (143.6,173.53) ;
\draw [color={rgb, 255:red, 208; green, 2; blue, 27 }  ,draw opacity=1 ] [dash pattern={on 4.5pt off 4.5pt}]  (62.28,162.93) .. controls (67.97,152.54) and (149.65,152.54) .. (154.93,152.94) ;
\draw [color={rgb, 255:red, 74; green, 144; blue, 226 }  ,draw opacity=1 ]   (143.6,160.74) .. controls (185,154.93) and (239.05,161.73) .. (241.9,166.12) ;
\draw [color={rgb, 255:red, 74; green, 144; blue, 226 }  ,draw opacity=1 ] [dash pattern={on 4.5pt off 4.5pt}]  (241.9,166.12) .. controls (240.68,180.91) and (202.48,202.9) .. (187.85,202.5) ;
\draw [color={rgb, 255:red, 74; green, 144; blue, 226 }  ,draw opacity=1 ]   (144,183.92) .. controls (166.72,188.91) and (171.59,194.9) .. (187.85,202.5) ;
\draw [color={rgb, 255:red, 208; green, 2; blue, 27 }  ,draw opacity=1 ] [dash pattern={on 4.5pt off 4.5pt}]  (143.6,173.53) .. controls (176.06,169.72) and (232.55,186.51) .. (255.71,186.91) ;
\draw [color={rgb, 255:red, 208; green, 2; blue, 27 }  ,draw opacity=1 ] [dash pattern={on 4.5pt off 4.5pt}]  (154.93,152.94) .. controls (187.85,147.74) and (237.43,151.74) .. (265.06,138.55) ;
\draw    (444.77,123.3) .. controls (435.02,117.16) and (433.49,204.12) .. (444.48,205.3) ;
\draw  [dash pattern={on 4.5pt off 4.5pt}]  (444.77,123.3) .. controls (447.96,119.79) and (453.94,202.91) .. (444.48,205.3) ;
\draw    (444.77,123.3) .. controls (468.77,121.43) and (528.1,122.78) .. (556.8,121.66) ;
\draw    (444.48,205.3) .. controls (467.66,205.55) and (525.8,207.94) .. (554.5,206.83) ;
\draw    (556.8,121.66) .. controls (553.29,119.5) and (547.71,155.29) .. (556.79,152.55) ;
\draw [color={rgb, 255:red, 65; green, 117; blue, 5 }  ,draw opacity=1 ]   (554.45,180.92) .. controls (550.94,178.75) and (545.43,209.57) .. (554.5,206.83) ;
\draw    [dash pattern={on 4.5pt off 4.5pt}](556.8,121.66) .. controls (562.31,121.01) and (560.73,151.88) .. (556.79,152.55) ;
\draw [color={rgb, 255:red, 65; green, 117; blue, 5 }  ,draw opacity=1 ]   (554.45,180.92) .. controls (559.77,179.93) and (559.4,208.15) .. (554.5,206.83) ;
\draw    (556.79,152.55) .. controls (538.32,152.01) and (524.07,182.21) .. (554.45,180.92) ;
\draw    (444.48,205.3) .. controls (420.46,206.99) and (365.6,203.52) .. (336.88,204.41) ;
\draw    (444.77,123.3) .. controls (421.59,122.88) and (367.89,118.35) .. (339.18,119.24) ;
\draw    (336.83,178.5) .. controls (355.29,179.17) and (369.15,148.66) .. (338.76,149.72) ;
\draw    (339.18,119.24) .. controls (335.67,117.08) and (329.68,152.47) .. (338.76,149.72) ;
\draw  [dash pattern={on 4.5pt off 4.5pt}]  (338.77,118.84) .. controls (344.28,118.19) and (342.7,149.05) .. (338.76,149.72) ;
\draw    (336.01,178.5) .. controls (332.5,176.33) and (327.81,207.16) .. (336.88,204.41) ;
\draw  [dash pattern={on 4.5pt off 4.5pt}]  (336.01,178.5) .. controls (341.33,177.51) and (340.96,205.73) .. (336.06,204.41) ;
\draw [color={rgb, 255:red, 208; green, 2; blue, 27 }  ,draw opacity=1 ]   (549.95,192.12) .. controls (474.45,188.62) and (504.45,134.62) .. (552.5,136.5) ;
\draw [color={rgb, 255:red, 74; green, 144; blue, 226 }  ,draw opacity=1 ]   (334,144) .. controls (373.45,134.62) and (401.95,193.62) .. (331.45,192.12) ;
\draw [color={rgb, 255:red, 74; green, 144; blue, 226 }  ,draw opacity=1 ] [dash pattern={on 4.5pt off 4.5pt}]  (41,136.5) .. controls (81.45,134.12) and (80.95,168.12) .. (143.6,160.74) ;
\draw [color={rgb, 255:red, 74; green, 144; blue, 226 }  ,draw opacity=1 ] [dash pattern={on 4.5pt off 4.5pt}]  (39,189.5) .. controls (73.45,180.62) and (124.45,184.62) .. (144,183.92) ;
\draw    (554.45,180.92) .. controls (574.15,178.21) and (584.04,175.81) .. (590.75,181.22) .. controls (597.45,186.62) and (597.45,188.62) .. (594.45,199.62) .. controls (591.45,210.62) and (571.96,207.29) .. (554.5,206.83) ;
\draw    (556.8,121.66) .. controls (574.95,121.12) and (606.45,122.62) .. (617.45,120.62) ;
\draw    (617.45,120.62) .. controls (613.94,118.46) and (608.36,154.25) .. (617.43,151.51) ;
\draw    (617.45,120.62) .. controls (622.96,119.97) and (621.37,150.84) .. (617.43,151.51) ;
\draw    (556.79,152.55) .. controls (577.45,150.62) and (610.95,156.62) .. (617.43,151.51) ;

\draw (56.09,122.15) node [anchor=north west][inner sep=0.75pt]    {$W\cong \rtr(\sigma)$};

\draw (180.14,122.15) node [anchor=north west][inner sep=0.75pt]    {$\tr(\phi ^{\sigma })$};

\draw (349.26,126.15) node [anchor=north west][inner sep=0.75pt]    {$\R_{W}(\tr(\phi^\sigma))$};

\draw (450.79,126.15) node [anchor=north west][inner sep=0.75pt]    {$\chi \left( W,\phi ^{\sigma },l\right)$};

\draw (570.5,126.15) node [anchor=north west][inner sep=0.75pt]    {$V_{1}$};

\end{tikzpicture}

    \caption{In this picture, we see an morphism $W:M\leadsto N$ whose underlying triad cobordism is elementary of index $k$ relative to $N$ composed with the trace of a surgery along the "meridian sphere". Observe that $W$ does not admit a left inverse, that is, there does not exist any bordism $W':N\leadsto M$ where the composite $W'\circ W$ is an isomorphism. However, as the picture illustrates, $\chi(W,\phi^\sigma,\ell)$ does. This right inverse $V_1$ is obtained by attaching a handle $2n-k+2$ relative to $\chi(N,\phi^\sigma)$ along the green circle. This morphism is elementary of index $k-1$ relative to its target and thus, by induction, is in $\W.$}
    \label{figure higher}
\end{figure}
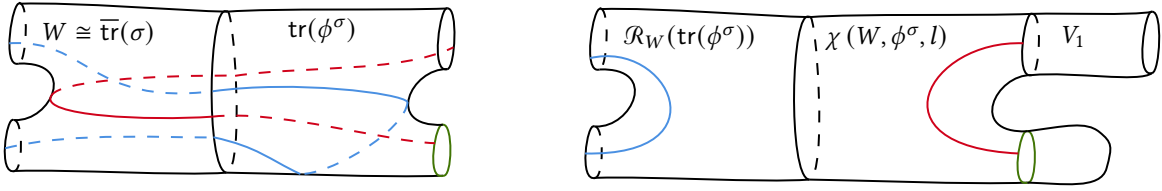

Let $W$ be $\rtr(\sigma)$ for an attaching map $\sigma$ as in \cref{main closure property high}. We have produced maps $\phi_p^\sigma$ disjoint from the support of $W$, whose attaching map are isotopic to the meridian of $\sigma$ at $x.$ Geometrically, doing "meridian surgery" using $\phi^\sigma$ produces a left-invertible morphism, whose inverse is already in $\W$ by induction. (See \cref{figure higher} for a depiction of this fact in the classical case without boundary\footnote{The picture is slightly misleading as the red arc on the right seems to be homotopic to the cocore of $\tr(\phi^\sigma)$ in the left picture, which is false. Unfortunately, this is a disadvantage of low dimensions where not all rearrangements can be realized.}.) The next proposition establishes exactly this.

\begin{lemma}\label{higher are left cancellable after rearranging}
     There exist a $\Theta$-structure $\ell$ on $\tr(\phi^\sigma)$ as in \cref{theta for expanded phi} and a morphism $V'_p:\chi^p(N,\phi^\sigma,\ell)\leadsto S_{\!p}$ in $\CobbtLs$ for every $p\geq 1$, such that the underlying cobordism is elementary of the same type of $\sigma$ and strictly lower index relative to $S_{\!p}$, and that $V'_p\circ \chi^p(W,\phi^\sigma,\ell)$ is an isomorphism. 
\end{lemma}
\begin{proof}
    We focus on the case where $\sigma$ is of right type as the remaining cases are analogous. We start by choosing the $\Theta$-structure $\ell$ on $\tr(\phi^\sigma)$. To do so, start by observing that there exists an embedding of underlying pairs $e:\tr(\phi^\sigma)\hookrightarrow N\times [0,1]$ relative to $N$ by the following argument: It suffices to find an embedding of the core of the unique handle of $\tr(\phi^\sigma)$ extending its attaching map. More precisely, we must find an embedding $(D^{2n+1-k}_+\times D^k,\partial_0D^{2n+1-k}_+\times D^k)\hookrightarrow (N\times [0,1],\hb N\times [0,1])$ extending the map $\phi^\sigma:(\partial_1D^{2n+1-k}_+\times D^k,\partial_{01}D^{2n+1-k}_+\times D^k)\hookrightarrow (N\times \{0\},\hb N\times \{0\}).$ This is possible since the attaching map $\phi^\sigma|_{(\partial_1D^{2n+1-k}_+,\partial_{01}D^{2n+1-k})}$ is trivial, that is, it extends to $(D^{2n+1-k}_+,\partial_0D^{2n+1-k}_+)$ using $\sigma$. Observe that the pullback $e^*(\ell_N\oplus \varepsilon^1)$ of the cylindrical structure may not be cylindrical away from the support of $\tr(\phi^\sigma)$, however we can modify the cylindrical structure on $N\times [0,1]$ to another structure $\ell'$ such that $e^*\ell'$ is cylindrical away from the support of $\tr(\phi^\sigma)$. Let $\ell$ be $e^*\ell'$ for such a choice of $\Theta$-structure on $N\times [0,1].$
    
    We will now construct the morphism $V'_p$. We will restrict ourselves to the case $p=1$ as the other cases are given by iterating the construction for this case. We define the underlying cobordism of $V_1'$ to be complement of the embedding $e$ and we choose the $\Theta$-structure $\ell_V$ on $V_1'$ given by pulling back $\ell'$ along the inclusion $V_1'\hookrightarrow N\times [0,1].$  As the embedding $e$ can be chosen to be disjoint from $\vb L$, we can see $(V_1',\ell_V)$ as a morphism of $\CobbtL$ with source $\chi^1(N,\phi^\sigma,\ell)$. We denote its target by $S_1$. Since the underlying cobordism $V'_1$ is the complement of $\tr(\phi^\sigma)$ in $N\times [0,1]$, we make the following observations:
    \begin{enumerate}
        \item The cobordism $V'_1$ is elementary of left type and index $2n-k+2$ relative to $\chi^1(N,\phi^\sigma,\ell)$. Thus, it is elementary of right type and index $k-1$ relative to $S_1$ and thus lies in the subcategory $\CobbtLs$ and satisfies the first condition of the claim.
        \item There exists a handle decomposition of $V_1'$ relative to $\chi^1(N,\phi^\sigma,\ell)$ whose attaching map intersects the belt sphere of $\tr(\phi^\sigma)$ exactly at one point, as $V_1'$ cancels $\tr(\phi^\sigma)$ in the sense of \cref{cancellation}. 
    \end{enumerate}However, the belt sphere of $\phi^\sigma$ is, by definition, isotopic to the attaching map of the unique handle (see \cref{the surgery is in fact surgery}) of $\chi^1(W,\phi^\sigma,\ell)$. Thus, the handle decomposition of $V_1'$ can be arranged so that its attaching map relative to $\chi^1(N,\phi^\sigma,\ell)$ intersects the attaching map of $\chi^1(W,\phi^\sigma,\ell)$ relative to $\chi^1(N,\phi^\sigma,\ell)$ exactly at one point. Thus $V'_1\circ \chi^p(W,\phi^\sigma,\ell)$ is diffeomorphic to a product cobordism by \cref{cancellation} and hence an isomorphism by \cref{cylinders are isomorphisms}. We conclude that $V'_1$ satisfies our desired properties, hence we have finished the proof.
    \end{proof}

We can now state and prove "stable stability" for higher elementary triad cobordisms.

\begin{prop}\label{stability for higher}
    Let $W:M\leadsto N$ be a morphism whose underlying cobordism is elementary of index $n<k<2n$ for right type, $n+1<k<2n+1$ for interior type and $n+1<k<2n+1$ for left type attached to the basepoint component, then $W\in \W.$ 
\end{prop}
\begin{proof}
    Once again, we can assume without loss of generality that $W$ is $\rtr(\sigma)$ as above. We proceed by induction in the following way for each type. We start with right type. Assume by induction that the result is true for $k-1.$ The base case is \cref{middle dim stability}. Define $\phi_p^\sigma$ as above and note that $V_p$ in \cref{higher are left cancellable after rearranging} is in $\W$ by induction for every $p\geq 1.$ Thus, we have $\chi^p(W,\phi^\sigma,\ell)\in \W.$ Apply now \cref{main closure property} to deduce the result. Thus, we conclude that the result holds for elementary triad cobordism of right type. The same argument applies for the left type, where the base case is index $n+1$ and is established by \cref{middle dim stability}. For interior type, the base case (index $n+1$) is established by using \cref{splitting} to factor this morphism as the composite of a pair of morphisms whose underlying triad cobordisms are elementary of left and right type of index $n+1$. Since we have proved this claim for right and left elementary triad cobordisms of index $n+1$, we conclude that both factors are in $\W$. Thus, we establish the base case for interior handles. The induction step goes exactly as the previous cases. This finishes the proof. 
\end{proof}

We finish by putting everything together and prove \cref{stable stability}.

\begin{proof}[Proof of \cref{stable stability}.]
    This follows by combining \cref{middle dim stability} and \cref{stability for higher} with \cref{basepoint component statement} to deduce that $\W=\CobbtLo.$
\end{proof}



\section{Proof of the closure properties.}\label{proof of closure section}

The goal of this section is to prove \cref{main closure property,main closure property high}. We start by sketching the overall strategy, which is inspired by the proof of \cite[Thm. 2.15]{GRWII}. Let $W:M\leadsto N$, $\sigma$ and $\ell_t$ be as in the hypothesis of either claim. We shall construct \textit{augmented semi-simplicial spaces} $X(N)_\bullet\to \F(N,K|_\infty)$ and $X(M)_\bullet\to \F(M,K|_\infty)$ satisfying the following properties:
\begin{enumerate}[label=(\textit{\roman*})]
    \item\label{step 1} The map $(-)\circ W: \F(N,K|_\infty)\to \F(M,K|_\infty) $ lifts to a map $(-)\circ W: X(N)_\bullet\to X(M)_\bullet$ of semi-simplicial spaces (see \cref{post functoriality,post functoriality high}).
    \item\label{step 3} The induced map $||X(N)_\bullet||\to \F(N,K|_\infty)$ is an equivalence and similarly for the version with $M$ instead of $N$ (see \cref{contractability middle right,contractability of higher,final contractability of middle left}).
    \item\label{step 2} The map $X(N)_{p-1}\to X(M)_{p-1}$ is an abelian homology equivalence, for any $p\geq 1$ (see \cref{de are cuts,de are cuts high}).
\end{enumerate}Since abelian homology equivalences are preserved under homotopy colimits, then \ref{step 2} implies that $||X(N)_\bullet||\to ||X(M)_\bullet||$ is an abelian homology equivalence. By \ref{step 3}, we conclude that $W\in \W.$ The construction of such augmented semi-simplicial spaces takes the form of \cref{right resolution,left resolution} for the middle-dimensional case, and \cref{right resolution high,left resolution high,interior resolution} for the higher-dimensional case. We fix again the data present in the \cref{the assumptions for stable stability} and a $\Theta$-end $K$.

\subsection{The middle-dimensional case.}\label{proof closure middle section}
The goal of this section is to prove \cref{main closure property}. Although the hypotheses of \cref{main closure property} require $(B,B^\partial)$ to be $1$-connected, we do not start by imposing such a condition. The only step that requires this extra assumption is \cref{contractibility easy}. We start by defining the augmented semi-simplicial space mentioned in the sketch above.


\begin{defn}[Right resolution] \label{right resolution}
    Let $n\geq 1$, $P,Q\in \CobbtL$ and $(s,W)\in \F(P,Q).$ Let \[\phi:(\partial D^n\times (0,\infty)\times \RR^{n}_+,\partial D^n\times (0,\infty)\times \partial\RR^{n}_+)\hookrightarrow (P,\hb P) \]be an embedding and $\ell:(0,\infty)\to \Thetaspace(T(D^n\times D^{n+1}_+),\Theta^*\gamma_{2n+1})$ a map, such that
    \(\smash{\ell_t|_{\partial D^n\times D^{n+1}_+}=(\phi^*\ell_P)|_{\partial D^n\times (t\cdot e_1+D^{n+1}_+)}}\), where $\Thetaspace(T(D^n\times D^{n+1}_+),\Theta^*\gamma_{2n+1})$ denotes the space of collared bundle maps $(T(D^n\times D^{n+1}_+),T(D^n\times \partial_0D^{n+1}_+))\to (\theta^*\gamma_{2n+1},(\theta^\partial)^*\gamma_{2n})$. Define $\de(W,\phi,\ell)_0$ to be the space of triples $(t,c,\L)$ where $t\in (0,\infty),$ 
    \(c:(D^n\times D^{n+1}_+, D^n\times \partial_0D^{n+1}_+)\hookrightarrow (W,\hb W)\) an embedding of pairs (recall \cref{section of pairs of manifolds}) and $\L:[0,1]\to \Thetaspacev(T(D^n\times D^{n+1}_+),\Theta^*\gamma_{2n+1})$ satisfying the following properties:
    \begin{enumerate}[label=(\textit{\roman*})]
        \item \textit{(Collared translation of $\phi$ by $t$ near $P$)} For some $\delta>0,$ we have
        $c(x,v)=\phi(\frac{x}{|x|},v+t\cdot e_1)+(1-|x|)\cdot e_0$, where $e_0$ is the first coordinate vector of $\RR\times \RR^{\infty}_+$ for all $x$ such that $1-|x|<\delta.$ 
        \item The image $(C,\hb C)=c(D^n\times D^{n+1}_+,D^n\times \partial_0D^{n+1}_+)$ is disjoint from $([0,s]\times \vb L)\cup (\{s\}\times Q)$ and $c^{-1}(P)=\partial D^n\times D^{n+1}_+.$
        \item $\L(0)=c^*\ell_W$ and $\L(1)=\ell_t.$
        \item The map $\ell_W|_{W\backslash C}:(W\backslash C,\hb W\backslash \hb C)\to (B,B^\partial)$ is strongly $n$-connected.
    \end{enumerate}This space is topologized as a subspace of 
    \[ \RR\times \Emb(D^n\times D^{n+1}_+,[0,\infty)\times  \RR^\infty_+)\times \Thetaspacev(D^n\times D^{n+1}_+,\Theta^*\gamma)^{[0,1]}\]where the middle space denotes the space of embeddings of pairs $e:(D^n\times D^{n+1}_+, D^n\times \partial_0D^{n+1}_+)\hookrightarrow ([0,\infty)\times \RR^\infty_+,[0,\infty)\times \partial \RR^\infty_+).$ Define $\de(W,\phi,\ell)_p$ to be the subspace of $\de(W,\phi,\ell)_0^{\times (p+1)}$ consisting of those tuples $(t_i;c_i;\L_i)_{i=0,\cdots,p}$ such that:
    \begin{enumerate}[label=(\textit{\roman*})]
        \item\label{right t} $t_0<t_1<\cdots<t_p.$
        \item\label{right pairwise disjoint} The embeddings $c_i$ are pairwise disjoint.
        \item\label{right strongly n} The map $\ell_W|_{W\backslash C}:(W\backslash C,\hb W\backslash \hb C)\to (B,B^\partial)$ is strongly $n$-connected, for $C$ the union of the images of all $c_i.$
    \end{enumerate}
    
\end{defn}

\begin{lemma}\label{right destab are semi simplicial}
    The subspaces $\smash{\de(W,\phi,\ell)_p\subset \de(W,\phi,\ell)_0^{\times (p+1)}}$ are preserved under forgetting factors. Therefore, it induces the structure of a semi-simplicial space to $\de(W,\phi,\ell)_\bullet.$ 
\end{lemma}
\begin{proof}
    We want to show that if we restrict the domain of the projection map $\smash{\text{proj}_i:\de(W,\phi,\ell)_0^{\times(p+1)}} \to \allowbreak  \smash{\de(W,\phi,\ell)_0^{\times p}}$ to the $\de(W,\phi,\ell)_p,$ the image lies in $\de(W,\phi,\ell)_{p-1}.$ It is obvious to see that the conditions \ref{right t} and \ref{right pairwise disjoint} of \cref{right resolution} are closed under forgetting some factors. It remains to prove that, if $M=W\backslash\cup_{\!j}C_{\!j}\to B$ is strongly $n$-connected, then so is $M'=W\backslash \cup_{\!j\neq i} C_{\!j}\to B$. We prove that the inclusion $M\to M'$ is strongly $(n-1)$-connected and apply \cref{closure properties strongly}. The cobordism $M$ is obtained from $M'$ by removing the core of a right $n$-handle. Therefore, $M'$ is obtained from $M$ by attaching its dual handle, which is a left $(n+1)$-handle. Thus by \ref{effect of left handles} in \cref{strong geo section}, the map $M\to M'$ is $(n-1)$-connected. 
\end{proof}

\begin{defn}[Left resolution]\label{left resolution}
    Let $n\geq 1$, $P,Q\in \CobbtL$ and $(s,W)\in \F(P,Q).$ Let \[\phi:(\partial_1 D^{n+1}_+\times (0,\infty)\times \RR^{n-1},\partial_{01} D^{n+1}_+\times (0,\infty)\times \RR^{n-1})\hookrightarrow (P,\hb P) \]be an embedding and $\ell:(0,\infty)\to \Thetaspace(T(D^{n+1}_+\times D^{n}),\Theta^*\gamma_{2n+1})$ a map, such that
    \(\smash{\ell_t|_{\partial_1 D^{n+1}_+\times D^{n}}=\phi^*\ell_P|_{\partial_{1} D^{n+1}_+\times (t\cdot e_1+D^{n})}}\), where $\Thetaspace(T(D^{n+1}_+\times D^{n}),\Theta^*\gamma_{2n+1})$ denotes the space of collared bundle maps $(T(D^{n+1}\times D^{n}),T(D^{n+1}_+\times \partial_0D^{n}))\to (\theta^*\gamma_{2n+1},(\theta^\partial)^*\gamma_{2n})$. Define $\de(W,\phi,\ell)_0$ to be the space of triples $(t,c,\L)$ where $t\in (0,\infty),$ 
    \(c:(D^{n+1}_+\times D^{n}, \partial_0D^{n+1}_+\times D^{n})\hookrightarrow (W,\hb W)\) an embedding of pairs and $\L:[0,1]\to \Thetaspacev(T(D^{n+1}_+\times D^{n}),\Theta^*\gamma_{2n+1})$ satisfying the following properties:
    \begin{enumerate}[label=(\textit{\roman*})]
        \item \textit{(Collared translation of $\phi$ by $t$ near $P$)} For some $\delta>0,$ we have
        \(c(x,v)=\phi(\frac{x}{|x|},v+t\cdot e_1)+(1-|x|)\cdot e_0\) where $e_0$ is the first coordinate vector of $\RR\times \RR^{\infty}_+$ for all $x$ such that $1-|x|<\delta.$ 
        \item The image $(C,\hb C)=c(D^{n+1}_+\times D^{n},\partial_0D^{n+1}_+\times D^{n})$ is disjoint from $([0,s]\times \vb L)\cup (\{s\}\times Q)$ and $c^{-1}(P)=\partial_1 D^{n+1}_+\times D^{n}.$
        \item $\L(0)=c^*\ell_W$ and $\L(1)=\ell_t.$
        \item The map $\ell_W|_{W\backslash C}:(W\backslash C,\hb W\backslash \hb C)\to (B,B^\partial)$ is strongly $n$-connected.
    \end{enumerate}This space is topologized in the same way as \cref{right resolution}. Define $\de(W,\phi,\ell)_p$ to be the subspace of $\smash{\de(W,\phi,\ell)_0^{\times (p+1)}}$ consisting of those tuples $(t_i;c_i;\L_i)_{i=0,\cdots,p}$ such that:
    \begin{enumerate}[label=(\textit{\roman*})]
        \item\label{left t} $t_0<t_1<\cdots<t_p;$
        \item\label{left pairwise disjoint} the embeddings $c_i$ are pairwise disjoint;
        \item\label{left strongly n} the map $\ell_W|_{W\backslash C}:(W\backslash C,\hb W\backslash\hb C)\to (B,B^\partial)$ is strongly $n$-connected, for $C$ the union of the images of all $c_i.$
    \end{enumerate}
    
\end{defn}

\begin{lemma}\label{left destab are semi simplicial}
    The subspaces $\smash{\de(W,\phi,\ell)_p\subset \de(W,\phi,\ell)_0^{\times (p+1)}}$ are preserved under forgetting factors. Therefore, it induces the structure of a semi-simplicial space to $\de(W,\phi,\ell)_\bullet.$
\end{lemma}
\begin{proof}
    As before, it suffices to prove that, if $M=W\backslash\cup_jC_j\to B$ is strongly $n$-connected, then so is $M'=W\backslash \cup_{j\neq i} C_j\to B$ also is. We prove that the inclusion $M\to M'$ is strongly $(n-1)$-connected. The cobordism $M$ is obtained from $M'$ by removing the core of a left $(n+1)$-handle. Therefore, $M'$ is obtained by attaching its dual handle, which is a right $n$-handle. Thus by \ref{effect of right handles} in \cref{strong geo section}, the map $M\to M'$ is strongly $(n-1)$-connected. The claim follows now by \cref{closure properties strongly}.
\end{proof}

We now introduce the augmented semi-simplicial spaces as promissed in the introduction to this section.

\begin{defn}
    Given $P,Q\in \CobbtL$ and $\phi$ and $\ell$ as in one of the above definitions. Define $\de_{\phi,\ell}(P,Q)_p$ to be the space of pairs $(W,x)$ where $W\in \F(P,Q)$ and $x\in \de(W,\phi,\ell)_p.$ For every $p\geq 0,$ this space maps to $\F(P,Q)$ by forgetting $x$.
\end{defn}

\begin{lemma}\label{quasi fibration}
    The homotopy fiber of the map $||\de_{\phi,\ell}(P,Q)_\bullet||\to \F(P,Q)$ over $W\in \F(P,Q)$ is equivalent to $||\de(W,\phi,\ell)_\bullet||$.
\end{lemma}
\begin{proof}
    By \cite[Lemma 2.14]{semisimplicial}\footnote{To use this result, one is forced to prove that all spaces forming the semi-simplicial space are compactly generated. However, the statement is invariant under changing $X_p$ by its Kellyfication $kX_p$ since the map $X\to kX$ is a weak equivalence (and the identity on underlying sets) for every space $X$.}, it suffices to show that $\de_{\phi,\ell}(P,Q)_p\to \F(P,Q)$ is a quasi-fibration for each $p\geq 0.$ This follows analogously to \cite[Lemma 7.17]{GRWI}. 
\end{proof}

\subsubsection{Functoriality of the resolutions.}
We now establish "functoriality" properties of these resolutions with respect to pre- and post-composition. In this section, we take $\de_{\phi,\ell}(P,Q)$ to be of any of the types above. Start by observing that post-composition with a cobordism $W':Q\leadsto S$ and a triple $(t,e,\L)\in \de(W,\phi,\ell)$ induces a triple $(t,\iota(c),\L)$ in $\de(W\cup_Q W',\phi,\ell)$ since $c$ is disjoint from $Q$, where $\iota:W\hookrightarrow W\cup_Q W'$ is the inclusion.

\begin{lemma}\label{post functoriality}
    Let $W':Q\leadsto S$ be a morphism in $\CobbtL$ which is strongly $(n-1)$-connected relative to $Q$, then post-composition defines a map of augmented semi-simplicial spaces
    \[W'\circ(-):\de_{\phi,\ell}(P,Q)_\bullet\to \de_{\phi,\ell}(P,S)_\bullet\]by taking a tuple $(W,(t,c,\L))$ to the tuple $(W\cup_Q W',(t,\iota(c),\L)).$ 
\end{lemma}
\begin{proof}
    It suffices to prove that, given $(W,(t,c,\L))\in \de_{\phi,\ell}(P,Q)_p$ then the restriction of the $\Theta$-structure to $(W\cup W')\backslash C$ is strongly $n$-connected. This manifold is equal to $(W\backslash C)\cup W'$, since $C$ is disjoint from $Q=W\cap W',$ so it suffices to prove that $W\backslash C\to (W\backslash C)\cup W'$ is strongly $(n-1)$-connected, by \cref{closure properties strongly}. This follows by \cref{strong conn and pushouts} since $(Q,\hb Q)\hookrightarrow (W',\hb W')$ is strongly $(n-1)$-connected.
\end{proof}

As a consequence of this lemma, the maps $K_{[i,i+1]}\circ(-):\F(N,K|_i)\to \F(N,K|_{i+1})$ lift to maps of augmented semi-simplicial spaces, and also for $M$. Thus, we can define the augmented semi-simplicial spaces $X(M)$ and $X(N)$ from the introduction by levelwise taking the homotopy colimit along these lifts. We denote these spaces by $\de_{\phi,\ell}(M,K|_\infty)$ amd $\de_{\phi,\ell}(N,K|_\infty)$, respectively. We will return to this closer to the end of the proof of \cref{main closure property}. We move now to constructing a similar lift of the pre-composition map.

\begin{cons}\label{reparameterization drama}
    Fix a diffeomorphism $\epsilon_s:(D^n\times D^{n+1}_+,D^n\times \partial_0 D^{n+1}_+)\cong ((s\cdot e_1+D^n)\times D^{n+1}_+,(s\cdot e_1+D^n)\times \partial_0 D^{n+1}_+)\cup [0,s]\times (\partial D^n\times D^{n+1}_+,\partial D^n\times \partial_0 D^{n+1}_+)\eqqcolon B_s\cup A_s=D_s$. Given $(t,c,\L)\in \de(W,\phi,\ell)$ for $\phi$ of right type and $(s,W'):P'\leadsto P$ such that $\supp(W)\cap \phi=\emptyset,$ define $c_s:D_s\hookrightarrow W'\cup W$ by taking $B_s$ to $c$ and $A_s$ to $[0,s]\times c|_{\partial D^n\times D^{n+1}_+}.$ The triple $(t,\epsilon^{-1}_s(c_s),\epsilon^{-1}_s(\L_s))$, where $\L_s$ is given by extending the path of $\Theta$-structure $\L$ along the collar $A_t,$ satisfies properties $(i)-(iii)$ of \cref{right resolution,left resolution}. We prove in the next lemma that such definition satisfies $(iv)$. When $\phi$ is of left type, one proceeds in a similar way by adding now a collar $[0,s]\times (\partial_1D^{n+1}_+\times D^{n},\partial_{01}D^{n+1}_+\times D^{n})$ to $(D^{n+1}_+\times D^{n},\partial_{0}D^{n+1}_+\times D^{n}).$ (For more details, see \cite[Defn. 4.13]{GRWII}.)
\end{cons}

\begin{lemma}\label{pre functoriality}
    Let $(s,W'):P'\leadsto P$ be a morphism in $\CobbtLs$ such that $\supp(W)\cap \phi=\emptyset$, then pre-composition as defined above 
    \[ (-)\circ W':\de_{\phi,\ell}(P,Q)_\bullet\to \de_{\phi,\ell}(P',Q)_\bullet\]is a map of augmented semi-simplicial spaces. Moreover, for $W'':Q\leadsto S,$ then the pre and post-composition maps commute strictly in the natural way. 
\end{lemma}
\begin{proof}
    As mentioned before, it suffices to prove that the construction above satisfies (iv). Let $(W,(t,c,\L))$ be a $p$-simplex of $\de_{\phi,\ell}(P,Q).$ We wish to prove that the $\Theta$-structure of $(W'\cup W)\backslash([-s,0]\times \vb C\cup C)$ is strongly $n$-connected and apply \cref{closure properties strongly}. It suffices to prove that the inclusion 
    \[W\backslash C\to (W'\cup W)\backslash([-s,0]\times \vb C\cup C)\]is strongly $(n-1)$-connected. By \cref{strong conn and pushouts}, this follows if $P\backslash \hb C\to W'\backslash C$ is strongly $(n-1)$-connected. However, this follows by \cref{strong conn and pushouts}, since $P\to W'$ is strongly $(n-1)$-connected and $(P\backslash \vb C,\hb P\backslash \hvb C)\to (P,\hb P)$ induces an isomorphism of fundamental groupoids, since $n\geq 3$. 
\end{proof}

This establishes \ref{step 1} from the introduction, namely that, we have a map of semi-simplicial spaces $(-)\circ W: \de_{\phi,\ell}(M,K|_\infty)_\bullet\to \de_{\phi,\ell}(N,K|_\infty)_\bullet$ lifting the map $(-)\circ W: \F(N,K|_\infty)\to \F(M,K|_\infty)$ .

\subsubsection{Contractability of left core complexes of middle dimension.}


The goal of this subsubsection is to establish \ref{step 3} from the introduction for the left type. More precisely, we prove the following result.
\begin{teo}\label{final contractability of middle left}
    Let $P\in \smash{\CobbtL}$ and $K$ a $\Theta$-end. Let $\phi$ and $\ell:(0,\infty)\to \Thetaspace(T(D^{n+1}_+\times D^{n}),\Theta^*\gamma_{2n+1})$ be as in \cref{left resolution}. Then the map induced by augmentation
    \(\smash{||\de_{\phi,\ell}(P,K|_\infty)_\bullet||\to \F(P,K|_\infty)}\) is an equivalence, provided $(B,B^\partial)$ is $1$-connected.
\end{teo}
By \cref{quasi fibration}, to prove this statement it suffices to prove that the homotopy colimit of the realizations of $\de(W\cup K_{[i,i+j]},\phi,\ell)_\bullet$ is contractible, for any $W\in \Fo(P,K|_i)$. To do so, we define the following variations of the semi-simplicial spaces $\de(-).$ Fix $W,\phi$ and $\ell$ as in \cref{left resolution}.

\begin{defn}\label{left resolution immersion}
    Define $\deb(W,\phi,\ell)_0$ to be the space of triples $(t,c,\L)$ where $t\in (0,\infty),$ 
    \[c:(D^{n+1}_+\times D^{n}, \partial_0D^{n+1}_+\times D^{n})\to (W,\hb W)\]an immersion of pairs and $\L:[0,1]\to \Thetaspacev(D^{n+1}_+\times D^{n})$ satisfying the following properties:
    \begin{enumerate}[label=(\roman*)]
        \item $c$ is an embedding when restricted to $(D^{n+1}_+\times \{0\},\partial_0D^{n+1}_+\times \{0\})$ and for some $\delta>0,$ we have
        \(c(x,v)=\phi(\frac{x}{|x|},v+t\cdot e_1)+(1-|x|)\cdot e_0\) where $e_0$ is the first coordinate vector of $\RR\times \RR^{\infty}_+$ for all $x$ such that $1-|x|<\delta.$ 
        \item The image $C=c(D^{n+1}_+\times \{0\},\partial_0D^{n+1}_+\times\{0\})$ is disjoint from $([0,s]\times \vb L)\cup (\{s\}\times Q)$ and $c^{-1}(P)=\partial_1 D^{n+1}_+\times D^{n}.$
        \item $\L(0)=c^*\ell_W$ and $\L(1)=\ell_t.$
        \item The map $\ell_W|_{W\backslash C}:(W\backslash C,\hb W\backslash \hb C)\to (B,B^\partial)$ is strongly $n$-connected.
    \end{enumerate}Define $\deb(W,\phi,\ell)_p$ to be the subspace of $\deb(W,\phi,\ell)_0^{\times (p+1)}$ constisting of those tuples $(t_i;c_i;\L_i)_{i=0,\cdots,p}$ such that:
    \begin{enumerate}[label=(\roman*)]
        \item $t_0<t_1<\cdots<t_p.$
        \item The embeddings $c_i|_{D^{n+1}\times \{0\}}$ are pairwise disjoint.
        \item The map $\ell_W|_{W\backslash C}:(W\backslash C,\hb W\backslash \hb C)\to (B,B^\partial)$ is strongly $n$-connected, for $C$ the union of the images of all $c_i|_{D^{n+1}\times \{0\}}.$
    \end{enumerate}This defines a semi-simplicial space $\deb(W,\phi,\ell)_\bullet.$ Denote by $\debd(W,\phi,\ell)_\bullet$ the semi-simplicial space where $\debd(W,\phi,\ell)_p$ is the set $\deb(W,\phi,\ell)_p$ with the discrete topology. The identity induces a map of semi-simplicial spaces
    \[\debd(W,\phi,\ell)_\bullet\to \deb(W,\phi,\ell)_\bullet.\]
    
\end{defn}

\begin{defn}\label{left resolution embedded core no (v)}
    Define $\dec(W,\phi,\ell)_0$ to be the set of triples $(t,c,\L)$ where $t\in (0,\infty),$ 
    \[c:(D^{n+1}_+\times D^{n}, \partial_0D^{n+1}_+\times D^{n})\to (W,\hb W)\]an immersion of pairs and $\L:[0,1]\to \smash{\Thetaspacev}(D^{n+1}_+\times D^{n})$ satisfying properties $(i)-(iii)$ of \cref{left resolution immersion}. Define $\smash{\dec(W,\phi,\ell)_p}$ to be the subset of $\smash{\dec(W,\phi,\ell)_0^{\times (p+1)}}$ constisting of those tuples $(t_i;c_i;\L_i)_{i=0,\cdots,p}$ such that:
    \begin{enumerate}[label=(\textit{\roman*})]
        \item $t_0<t_1<\cdots<t_p.$
        \item the embeddings $c_i|_{D^{n+1}\times \{0\}}$ are pairwise disjoint.
    \end{enumerate}This defines a semi-simplicial space $\dec(W,\phi,\ell)_\bullet,$ where all spaces have the discrete topology. The identity induces a map 
    \[\debd(W,\phi,\ell)_\bullet\to \dec(W,\phi,\ell)_\bullet.\]
\end{defn}

\begin{defn}\label{left resolution immersion no (v)}
    Define $\deh(W,\phi,\ell)_0$ to be the set of triples $(t,c,\L)$ where $t\in (0,\infty),$ 
    \[c:(D^{n+1}_+\times D^{n}, \partial_0D^{n+1}_+\times D^{n})\to (W,\hb W)\]an immersion of pairs and $\L:[0,1]\to \Thetaspacev(D^{n+1}_+\times D^{n})$ satisfying:
    \begin{enumerate}[label=(\textit{\roman*'})]
        \item $c$ is a self-transverse immersion with no triple points when restricted to $(D^{n+1}_+\times \{0\},\partial_0D^{n+1}_+\times \{0\})$ and for some $\delta>0,$ we have
        \(c(x,v)=\phi(\frac{x}{|x|},v+t\cdot e_1)+(1-|x|)\cdot e_0\)where $e_0$ is the first coordinate vector of $\RR\times \RR^{\infty}_+$ for all $x$ such that $1-|x|<\delta.$  
    \end{enumerate}and $(ii)-(iii)$ of \cref{left resolution immersion}. Define $\deh(W,\phi,\ell)_p$ to be the subset of $\deh(W,\phi,\ell)_0^{\times (p+1)}$ constisting of those tuples $(t_i;c_i;\L_i)_{i=0,\cdots,p}$ such that:
    \begin{enumerate}[label=(\textit{\roman*})]
        \item $t_0<t_1<\cdots<t_p.$
        \item the embeddings $c_i|_{D^{n+1}_+\times \{0\}}$ are pairwise transverse and without triple intersections.
    \end{enumerate}This defines a semi-simplicial space $\deh(W,\phi,\ell)_\bullet,$ where all spaces have the discrete topology. The identity induces a map 
    \[\dec(W,\phi,\ell)_\bullet\to \deh(W,\phi,\ell)_\bullet.\]
    
\end{defn}

Putting all these variations together, we have the following string of maps
\[\de(W,\phi,\ell)_\bullet\to\deb(W,\phi,\ell)_\bullet\leftarrow\debd(W,\phi,\ell)_\bullet\to \dec(W,\phi,\ell)_\bullet\to \deh(W,\phi,\ell)_\bullet.\]The strategy to prove that the leftmost semi-simplicial space has contractible realization after taking the homotopy colimit of the maps $K_{[i,i+1]}\circ (-)$ will go as follows. First, we prove that the rightmost semi-simplicial set has contractible realization. Secondly, we prove that all arrows above induce a weak equivalence after realization and taking homotopy colimits. The first step does not require taking this homotopy colimit, however this is necessary for the second step. We start with the first step. This is the first point in our proof where the assumption that $(B,B^\partial)$ is $1$-connected is necessary.

\begin{lemma}\label{contractibility easy}
    The space $||\deh(W,\phi,\ell)_\bullet||$ is contractible, provided $(B,B^\partial)$ is $1$-connected.
\end{lemma}
\begin{proof}
    We start by showing that $\deh(W,\phi,\ell)_0$ is non-empty. Since $(B,B^\partial)$ and thus also $(W,\hb W)$ (by $n\geq 2$) are $1$-connected, \cref{strong implies lifting} assures the existence of a lift
    \[\begin{tikzcd}
	{(\partial_1 D^{n+1}_+\times D^{n},\partial_{01} D^{n+1}_+\times D^{n})} & {(P,\hb P)} & {(W, \hb W)} \\
	{(D^{n+1}_+\times D^{n},\partial_{0} D^{n+1}_+\times D^{n})} && {(B,B^\partial)}
	\arrow["\phi", from=1-1, to=1-2]
	\arrow[hook, from=1-1, to=2-1]
	\arrow[hook, from=1-2, to=1-3]
	\arrow["{\ell_W}", from=1-3, to=2-3]
	\arrow["\hat{c}", dashed, from=2-1, to=1-3]
	\arrow["{\ell_t}", from=2-1, to=2-3]
\end{tikzcd}\]for any $t\in (0,\infty)$, which can be assumed to be covered by a bundle map between the tangent bundles, making the upper triangle commute and the lower one commute up to homotopy of bundle maps, since $\ell_t$ is covered by a bundle map. Thus by Smale-Hirsch, we can homotope it to be an immersion of pairs. We can now define a triple $(t,\hat{c},\L)$ where $\L$ is a bundle homotopy witnessing that the bottom triangle is homotopy commutative through bundle maps, such that it defines a point in $\deh(W,\phi,\ell)_0.$ We proceed by proving contractability assuming non-emptyness. For $k\geq 0$, let $f:\partial I^k\to ||\deh(W,\phi,\ell)_\bullet||$ be a map, which we can assume to be simplicial with respect to a triangulation of $\partial I^k$ by the simplicial approximation theorem. Thus, for each vertex $v_i\in \partial I^k$ of this triangulation, we have an element $f(v_i)=(t_i,c_i,\L_i)\in \deh(W,\phi,\ell)_0.$ Fix such a $v_i$. Let $(t,c',\L)$ be a triple such that $t\neq t_j$ for any $j$, $c'$ is the result of an isotopy of $c_i$ which extends the translation of $\phi_{t_i}$ to $\phi_{t}$ and $\L$ be the path extending the path $\ell$ on $P$. By transversality, we can find a regular homotopy of pairs from $\smash{c'|_{D^{n+1}_+\times \{0\}}}$ relative to its vertical boundary (and hence, constant on $\phi_t$) such the result $c$ is transverse to all $c_j$ and without triple intersections. Therefore, $f$ is in the star of $f(v_j)$ and so $f$ extends to the cone of $\partial I^k.$
\end{proof}

 We now move to comparing the different variations of the resolutions defined above. Let $W:P\leadsto K|_i$ for some $i$ and denote by $\de(W_j,\phi,\ell)_p\coloneqq\de(W_i\cup K|_{[i,i+j]},\phi,\ell)_p$ for $j\geq 0.$ Denote by $\de(W_\infty,\phi,\ell)_p$ the (strict) colimit of $\de(W_j,\phi,\ell)_p$ along post composition map $K|_{[i,i+j]}\circ (-)$ (see \cref{post functoriality}). One can check that this is equivalent to the homotopy colimit, since the maps are inclusions. Since the geometric realization commutes with homotopy colimits, we have an equivalence
\[\hocolim_j||\de(W_j,\phi,\ell)_\bullet||\to ||\de(W_\infty,\phi,\ell)_\bullet||.\]We use the analogous notation for the variations of $\de(W,\phi,\ell)$ defined previously.

\begin{prop}\label{from immersed to embedded cores}
    The map 
    \[\dec(W_\infty,\phi,\ell)_\bullet\to \deh(W_\infty,\phi,\ell)_\bullet\]induces an equivalence on geometric realizations.
\end{prop}
\begin{proof}
    
     Denote the source of this map by $X_\bullet$ and the target by $Y_\bullet.$ Let $k\geq 0$ be an integer and $F:(I^k,\partial I^k)\to (||Y_\bullet||,||X_\bullet||)$ be a map, which we can assume to be simplicial with respect to some triangulation of $I^k.$ We want to show that $F$ is homotopic to a map whose image is contained in $||X_\bullet||.$ We argue now that it suffices to find a  map of pairs $F'$ homotopic to $F$ where $F'(\sigma)\in X_p$ for every simplex $\sigma$ of dimension $p$ for $p=\{0,1\}$: This follows by two properties. First, the map $Y_p\to \smash{Y_0^{\times (p+1)}}$ given by the differentials is injective. Second, for $p\geq 1$, $v=(v_0,v_1,\cdots, v_p)\in Y_p$ belongs to $X_p$ if and only if $(v_i,v_j)\in X_1$ for every $0\leq i,j\leq p.$ We start by finding a homotopic map where the images of all vertices $v$ of $I^k$ are in $X_0$. In other words, given $F(v)=(t,c,\L),$ it suffices to find another vertex $w_v=(t+\epsilon,c',\L')\in Y_0$ for arbitrarily small enough $\epsilon>0$ where $c'$ is embedded, once restricted to its core, and transverse to $c$ (and to all immersions given by the images of the other vertices, which $c$ was already transverse) without triple intersections: by applying this technique to all interior vertices $v$ and rechoosing $F$ on those vertices to be $w_v$, we produce a new map, since $(F(v),w_{v'})\in Y_1$ whenever $(F(v),F(v'))\in Y_1$. This map is homotopic to $F$ since for each $v,$ $(F(v),w_v)\in Y_1.$ 
     \begin{figure}
         \centering
        \tikzset{every picture/.style={line width=0.75pt}} 

\begin{tikzpicture}[x=0.75pt,y=0.75pt,yscale=-1,xscale=1]

\draw  [color={rgb, 255:red, 155; green, 155; blue, 155 }  ,draw opacity=0.69 ][fill={rgb, 255:red, 155; green, 155; blue, 155 }  ,fill opacity=0.11 ] (237.36,119.22) -- (261.56,95.02) -- (427.67,95.02) -- (427.67,151.5) -- (403.46,175.7) -- (237.36,175.7) -- cycle ; \draw  [color={rgb, 255:red, 155; green, 155; blue, 155 }  ,draw opacity=0.69 ] (427.67,95.02) -- (403.46,119.22) -- (237.36,119.22) ; \draw  [color={rgb, 255:red, 155; green, 155; blue, 155 }  ,draw opacity=0.69 ] (403.46,119.22) -- (403.46,175.7) ;
\draw [color={rgb, 255:red, 155; green, 155; blue, 155 }  ,draw opacity=1 ] [dash pattern={on 4.5pt off 4.5pt}]  (263.5,151.23) -- (427.67,151.5) ;
\draw [color={rgb, 255:red, 155; green, 155; blue, 155 }  ,draw opacity=1 ] [dash pattern={on 4.5pt off 4.5pt}]  (237.36,175.7) -- (263.5,151.23) ;
\draw [color={rgb, 255:red, 155; green, 155; blue, 155 }  ,draw opacity=1 ] [dash pattern={on 4.5pt off 4.5pt}]  (263.5,151.23) -- (263.09,94.76) ;
\draw  [color={rgb, 255:red, 74; green, 144; blue, 226 }  ,draw opacity=0.7 ][fill={rgb, 255:red, 74; green, 144; blue, 226 }  ,fill opacity=0.14 ] (247.34,109.83) -- (253.64,103.52) -- (418.94,103.52) -- (418.94,159.31) -- (412.64,165.61) -- (247.34,165.61) -- cycle ; \draw  [color={rgb, 255:red, 74; green, 144; blue, 226 }  ,draw opacity=0.7 ] (418.94,103.52) -- (412.64,109.83) -- (247.34,109.83) ; \draw  [color={rgb, 255:red, 74; green, 144; blue, 226 }  ,draw opacity=0.7 ] (412.64,109.83) -- (412.64,165.61) ;
\draw [color={rgb, 255:red, 74; green, 144; blue, 226 }  ,draw opacity=0.72 ] [dash pattern={on 4.5pt off 4.5pt}]  (254,103.52) -- (254.13,159.8) ;
\draw [color={rgb, 255:red, 74; green, 144; blue, 226 }  ,draw opacity=0.54 ] [dash pattern={on 4.5pt off 4.5pt}]  (254.13,159.8) -- (418.94,159.31) ;
\draw  [color={rgb, 255:red, 208; green, 2; blue, 27 }  ,draw opacity=0.32 ][fill={rgb, 255:red, 208; green, 2; blue, 27 }  ,fill opacity=0.07 ] (237.17,140.09) -- (262.5,114.76) -- (428.17,114.76) -- (428.17,127.53) -- (402.83,152.87) -- (237.17,152.87) -- cycle ; \draw  [color={rgb, 255:red, 208; green, 2; blue, 27 }  ,draw opacity=0.32 ] (428.17,114.76) -- (402.83,140.09) -- (237.17,140.09) ; \draw  [color={rgb, 255:red, 208; green, 2; blue, 27 }  ,draw opacity=0.32 ] (402.83,140.09) -- (402.83,152.87) ;
\draw [color={rgb, 255:red, 208; green, 2; blue, 27 }  ,draw opacity=0.31 ] [dash pattern={on 4.5pt off 4.5pt}]  (237.17,152.87) -- (262.68,128.72) ;
\draw [color={rgb, 255:red, 208; green, 2; blue, 27 }  ,draw opacity=0.28 ] [dash pattern={on 4.5pt off 4.5pt}]  (262.68,128.72) -- (427.44,128.72) ;

\end{tikzpicture}

         \caption{In this picture, we represent the region $U_0\times D^n$ in blue and $U_1\times D^n$ in red in a chart $x$.}
         \label{chart picture}
     \end{figure}
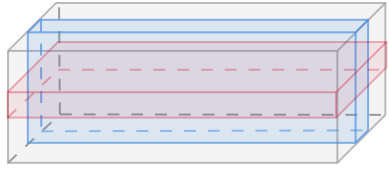
    \\\\
    The embedding $c'$ will be constructed in a few steps via local modifications around the self-intersection points. Firstly, we fix a "chart" around the self-intersection points and all modifications will essentially be supported in this chart. Secondly, we define a local modification to $c$ eliminating intersection points, provided $c$ satisfies a certain hypothesis (D). Thirdly, we define a local modification that takes any embedding into one satisfying (D). This last step will the only place where the infinite supply of embedded $V_1$'s is necessary. This strategy is inspired by the constructions in the proof of \cite[Prop. 5.5]{GRWII}.
    \\\\
    \textit{Removing intersection circles.} Before implementing the aforementioned construction, we make a preliminary modification to $c$ removing certain pathological self-intersections, which uses that $n\geq 3$ and $(B,B^\partial)$ is $1$-connected. Let $v$ be an interior vertex of $I^k$ and $(t,c,\L)$ its image under $F$. Let $j\geq 0$ be an integer such that the image of $c$ lies in $W_j.$ Denote the self-intersection locus of $\smash{c|_{D^{n+1}_+\times \{0\}}}$ by $(K,\partial K)$, that is the $1$-dimensional submanifold pair of $(W_j,\hb W_j)$ of those points where the pre-image of $c$ is more than one point. This is a disjoint union of manifolds diffeomorphic to $([0,1],\{0,1\})$ and $(S^1,\emptyset).$ Fix a path component $K'$ and observe that $c:c^{-1}(K')\to K'$ is a $2$-fold cover. Hence, if $K'$ is a compact interval, then $c^{-1}(K')\cong K'\times \{0,1\} $. However, if $K'$ is a circle, then $c|_{c^{-1}(K')}$ might be the non-trivial $2$-fold cover of $K'$. Since our construction will only apply to the case where $c|_{c^{-1}(K')}$ is trivial covering, we first find an embedding $c'$ which is isotopic to $c$ relative to $c^{-1}(P)$ such that its self-intersection locus is a disjoint union of intervals. This is possible as an application of \cite[Lemma C.6]{BP}, by inspecting that the proof only uses that $n\geq 3$, that $(D^{n+1}_+,\partial_0D^{n+1}_+)$ is $0$-connected, and that $(W,\hb W)$ is $1$-connected (which follows from $(B,B^\partial)$ being $1$-connected and $n\geq 2$). Thus, we assume from now on that our original embedding $c$ self-intersects in a disjoint union of compact intervals.  
    \\\\
     \textit{Parameterizing self-intersection loci.} We proceed now with the aforementioned move. Recall that we denote the self-intersection locus of $c$ by $K$. Fix a path component $K'$ of $K$. We present a method to rechoose $c$ on a neighborhood of the pre-image of $K'$ so its self-intersection locus is $K\backslash K'.$ By hypothesis, $c^{-1}(K')$ has two connected components diffeomorphic to $K'$, which we denote by $K_0$ and $K_1$. Observe that the normal bundle of $K_i$ in $D\coloneqq D^{n+1}_+\times \{0\}$ is trivial: The tangent bundle of $K_i$ is trivial and its sum with the normal bundle is the pullback of $(TD^{n+1}_+,T\partial_0D^{n+1}_+)$, which is trivial. Thus, the normal bundle is an $n$-dimensional collared vector bundle pair which is stably trivial, thus it is trivial. Moreover, the normal bundle of $K'$ in $W_j$ splits as the sum of the normal bundles of $K_0$ and $K_1$ in $D.$ Thus, we can find a (framed) tubular neighborhood $x:(K'\times D^n\times D^n,\partial K'\times D^n\times D^n)\hookrightarrow (W_j,\hb W_j)$ disjoint from $\vb W_j$ such that, up to scaling the $\{0\}\times D^n$-direction, there exists an immersion $c_0$ which agrees with $c$ outside of a neighborhood of $\im(x)$ such that:
    \begin{enumerate}[label=(\roman*)]
        \item there exists a neighborhood $U_0$ of $K_0$ in $D^{n+1}_+\times \{0\}$ such that $\im(x)\cap c_0(U_0)=x(K'\times D^n\times \{0\}).$ Additionally, $\im(x)\cap c_0(U_0\times D^n)=x(K\times D^n\times \frac{1}{2} D^n).$
        \item there exists a neighborhood $U_1$ of $K_1$ in $D^{n+1}_+\times \{0\}$ such that $\im(x)\cap c_0(U_1)=x(K'\times \{0\}\times D^n).$ Additionally, $\im(x)\cap c_0(U_1\times D^n)=x(K\times \frac{1}{2} D^n\times D^n).$
    \end{enumerate}By replacing $c$ by $c_0$, we can assume without loss of generality that $c=c_0,$ since this change is supported away from $\partial_1D^{n+1}_+\times D^{n}.$ (See \cref{chart picture} for the restriction of the immersion to the chart $x$ for the case $K'$ is diffeomorphic to $[0,1].$)
        \\\\
    \textit{The hypothesis} (D). Let us introduce the following hypothesis on the immersion $c$: there exists an embedding $f:(K'\times S^n\times D^n,\partial K'\times S^n\times D^n)\hookrightarrow (W,\hb W)$ whose core (that is, $\im(f|_{K'\times S^n\times \{0\}})$) intersects $c$ exactly in $f(K'\times \{x\}\times \{0\})$ for some $x\in S^n,$ and an embedding $\gamma: (K'\times [0,1],\partial K'\times [0,1])\hookrightarrow (D^{n+1}_+,\partial_0D^{n+1}_+) $ disjoint from $c^{-1}(K')\backslash K_0$ such that $\im(\gamma(-,0))=K_0$ and $c(\gamma(-,1))=f(K'\times \{x\}\times \{0\}).$ Assume additionally that the framings of the normal bundle of $c(\gamma(-,0))$ and $c(\gamma(-,1))$ induced by the identifications above are in the same path component. We say $c$ satisfies (D) if this hypothesis holds. We shall describe a construction which can be applied to $c$ satisfying (D) and later modify $c$ to satisfy (D). Assume $c$ satisfies (D). Choose a chart $x_f$ as before for the intersection of $c$ and $f$ satisfying the analogous properties where $K_0$ lies in the source of $c$ and $K_1$ on the source of $f$. We can find a framed tubular neighborhood \[\hat{\gamma}:(K'\times [0,1]\times D^{n}\times D^{n-1},\partial K'\times [0,1]\times D^{n}\times D^{n-1})\hookrightarrow (W_j,\hb W_j) \] of $c(\im(\gamma))$ satisfying the following properties:
    \begin{enumerate}[label=(\roman*)]
        \item when restricted to $K'\times \{0\}\times D^n\times D^{n-1}$, $\hat{\gamma}$ agrees with $x$ restricted to $K'\times D^n\times D^{n-1}_+\subset K'\times D^n\times \partial D^{n},$ where $D^{n-1}_+$ is the upper hemisphere of $\partial D^n,$ under the clear identification $K'\times D^{n}\times D^{n-1}\cong K'\times D^{n}\times D^{n-1}_+.$
        \item when restricted to $K'\times \{1\}\times D^n\times D^{n-1}$, $\hat{\gamma}$ agrees with $x_f$ restricted to $K'\times D^n\times D^{n-1}_-\subset K'\times D^n\times \partial D^{n},$ where $D^{n-1}_-$ is the lower hemisphere of $\partial D^n.$ 
    \end{enumerate}Denote by $\hat{\gamma}_0$ the restriction of $\hat{\gamma}$ to $K'\times [0,1]\times \partial(\frac{1}{2}D^n)\times D^{n-1}\cong [0,1]\times \partial(\frac{1}{2}D^n)\times D^n$. Let $c':(D^{n+1}_+\times D^{n}, \partial_0D^{n+1}_+\times D^{n})\hookrightarrow (W_j,\hb W_j)$ be the following immersion: Outside the $U_1\times D^n,$ it agrees with $c.$ Inside $\frac{1}{2}U_1\times D^n$ it is given by gluing $\hat{\gamma}_0$ with the image of $f$ except the pre-image of $x.$ This construction requires isotopies of $x(K'\times \frac{1}{2}D^n\times \partial D^n)$ to agree with $\hat{\gamma}(K'\times \{0\}\times \partial(\frac{1}{2} D^n)\times D^{n-1})$ and similarly for $x_f.$ This produces an immersion $c'$ which agrees with $c$ outside $U_1\times D^n$, whose core has intersection locus $K\backslash K'.$ We remark that in this step no additional copies of $V_1$ were used (see \cref{dual move picture} for a picture of this move).
    \begin{figure}
        \centering
        
\tikzset{every picture/.style={line width=0.75pt}} 

\begin{tikzpicture}[x=0.75pt,y=0.75pt,yscale=-1,xscale=1]

\draw  [color={rgb, 255:red, 155; green, 155; blue, 155 }  ,draw opacity=0.81 ][fill={rgb, 255:red, 155; green, 155; blue, 155 }  ,fill opacity=0.11 ] (82.79,73.59) -- (104.99,56.54) -- (159.77,90.8) -- (159.83,139.39) -- (137.64,156.44) -- (82.85,122.18) -- cycle ; \draw  [color={rgb, 255:red, 155; green, 155; blue, 155 }  ,draw opacity=0.81 ] (159.77,90.8) -- (137.58,107.85) -- (82.79,73.59) ; \draw  [color={rgb, 255:red, 155; green, 155; blue, 155 }  ,draw opacity=0.81 ] (137.58,107.85) -- (137.64,156.44) ;
\draw  [color={rgb, 255:red, 184; green, 233; blue, 134 }  ,draw opacity=1 ][fill={rgb, 255:red, 184; green, 233; blue, 134 }  ,fill opacity=0.32 ] (289.1,132.45) -- (229.21,133.23) .. controls (212.68,133.45) and (199.09,112.97) .. (198.86,87.49) .. controls (198.63,62.02) and (211.86,41.19) .. (228.39,40.97) -- (288.27,40.19) .. controls (304.81,39.97) and (318.4,60.45) .. (318.63,85.93) .. controls (318.85,111.4) and (305.63,132.23) .. (289.1,132.45) .. controls (272.56,132.67) and (258.97,112.19) .. (258.74,86.71) .. controls (258.52,61.23) and (271.74,40.4) .. (288.27,40.19) ;
\draw [color={rgb, 255:red, 208; green, 2; blue, 27 }  ,draw opacity=1 ]   (44.9,143) .. controls (63.61,125.98) and (64.9,116.9) .. (82.68,96.02) ;
\draw [color={rgb, 255:red, 208; green, 2; blue, 27 }  ,draw opacity=1 ]   (110.38,180.66) .. controls (118.91,167.61) and (123.4,153.06) .. (138.35,137.68) ;
\draw [color={rgb, 255:red, 208; green, 2; blue, 27 }  ,draw opacity=1 ][line width=0.75]    (138.35,137.68) -- (159.89,118.16) ;
\draw [color={rgb, 255:red, 208; green, 2; blue, 27 }  ,draw opacity=1 ] [dash pattern={on 4.5pt off 4.5pt}]  (82.68,96.02) -- (102.33,78.54) ;
\draw [color={rgb, 255:red, 208; green, 2; blue, 27 }  ,draw opacity=1 ]   (159.89,118.16) .. controls (178.6,101.14) and (240.24,97.38) .. (258.96,80.35) ;
\draw [color={rgb, 255:red, 155; green, 155; blue, 155 }  ,draw opacity=0.81 ] [dash pattern={on 4.5pt off 4.5pt}]  (102.9,55.23) -- (102.52,104.19) ;
\draw [color={rgb, 255:red, 155; green, 155; blue, 155 }  ,draw opacity=0.81 ] [dash pattern={on 4.5pt off 4.5pt}]  (82.85,122.18) -- (102.52,104.19) ;
\draw [color={rgb, 255:red, 208; green, 2; blue, 27 }  ,draw opacity=1 ] [dash pattern={on 4.5pt off 4.5pt}]  (102.33,78.54) .. controls (107.76,75.82) and (122.16,72.84) .. (130.4,72.41) ;
\draw [color={rgb, 255:red, 208; green, 2; blue, 27 }  ,draw opacity=1 ]   (130.4,72.41) .. controls (149.86,71.28) and (174.19,85.35) .. (199.07,80.13) ;
\draw [color={rgb, 255:red, 0; green, 0; blue, 0 }  ,draw opacity=1 ]   (199.07,80.13) -- (258.96,80.35) ;
\draw [color={rgb, 255:red, 208; green, 2; blue, 27 }  ,draw opacity=1 ] [dash pattern={on 4.5pt off 4.5pt}]  (199.07,80.13) .. controls (226.96,75.59) and (250.16,54.25) .. (266.25,55.61) ;
\draw [color={rgb, 255:red, 208; green, 2; blue, 27 }  ,draw opacity=1 ]   (258.96,80.35) .. controls (278.04,67.19) and (278.6,56.07) .. (266.25,55.61) ;
\draw [color={rgb, 255:red, 184; green, 233; blue, 134 }  ,draw opacity=1 ] [dash pattern={on 4.5pt off 4.5pt}]  (228.39,40.97) .. controls (241.84,39.75) and (257.64,58.9) .. (258.58,83.64) .. controls (259.51,108.38) and (248.57,132.4) .. (229.21,133.23) ;
\draw [color={rgb, 255:red, 0; green, 0; blue, 0 }  ,draw opacity=1 ]   (92.5,87.28) -- (150.79,126.21) ;
\draw [color={rgb, 255:red, 74; green, 144; blue, 226 }  ,draw opacity=1 ]   (92.51,13.63) -- (92.22,65.22) ;
\draw [color={rgb, 255:red, 74; green, 144; blue, 226 }  ,draw opacity=1 ] [dash pattern={on 4.5pt off 4.5pt}]  (92.5,87.28) -- (92.22,65.22) ;
\draw [color={rgb, 255:red, 74; green, 144; blue, 226 }  ,draw opacity=1 ] [dash pattern={on 4.5pt off 4.5pt}]  (92.5,87.28) -- (92.04,186.43) ;
\draw [color={rgb, 255:red, 74; green, 144; blue, 226 }  ,draw opacity=1 ]   (150.29,13.92) -- (150.61,97.85) ;
\draw [color={rgb, 255:red, 155; green, 155; blue, 155 }  ,draw opacity=0.81 ] [dash pattern={on 4.5pt off 4.5pt}]  (102.52,104.19) -- (159.83,139.39) ;
\draw  [draw opacity=0][fill={rgb, 255:red, 74; green, 144; blue, 226 }  ,fill opacity=0.19 ] (150.29,13.92) -- (150.98,195.05) -- (92.04,186.43) -- (92.51,13.63) -- cycle ;
\draw  [draw opacity=0][fill={rgb, 255:red, 208; green, 2; blue, 27 }  ,fill opacity=0.1 ] (51.41,136.72) -- (60.2,126.53) -- (75.84,104.51) -- (82.68,96.02) -- (102.33,78.54) -- (115.14,74.32) -- (130.4,72.41) -- (142.27,72.96) -- (153.69,76.14) -- (167.91,78.86) -- (181.94,81.13) -- (190.55,81.58) -- (201.59,79.77) -- (221.99,73.41) -- (233.91,67.69) -- (239.02,65.24) -- (250.43,59.11) -- (260.91,56.16) -- (267.27,55.71) -- (272.32,56.84) -- (274.2,60.02) -- (274.01,64.11) -- (271.39,68.87) -- (265.78,75) -- (254.92,82.95) -- (245.19,88.39) -- (229.85,93.61) -- (199.91,100.42) -- (184.94,104.96) -- (168.09,111.32) -- (154.43,123.12) -- (140.77,136.06) -- (131.6,145.37) -- (125.8,153.31) -- (110.38,180.66) -- (105.68,190.06) -- (45.57,190.85) -- (44.9,143) -- cycle ;
\draw  [draw opacity=0][fill={rgb, 255:red, 144; green, 19; blue, 254 }  ,fill opacity=0.18 ] (91.3,87.41) -- (99.45,80.05) -- (113.93,74.45) -- (130.4,72.41) -- (141.06,73.09) -- (152.48,76.27) -- (166.7,78.99) -- (180.74,81.26) -- (189.34,81.72) -- (197.87,80.26) -- (257.75,80.49) -- (252.04,84.46) -- (242.31,89.91) -- (226.97,95.13) -- (197.03,101.94) -- (182.06,106.48) -- (165.21,112.83) -- (162.96,114.78) -- (149.59,126.34) -- cycle ;
\draw  [color={rgb, 255:red, 155; green, 155; blue, 155 }  ,draw opacity=0.81 ][fill={rgb, 255:red, 155; green, 155; blue, 155 }  ,fill opacity=0.11 ] (381.9,72.58) -- (404.65,54.9) -- (458.99,90.72) -- (459.05,142.05) -- (436.3,159.74) -- (381.96,123.92) -- cycle ; \draw  [color={rgb, 255:red, 155; green, 155; blue, 155 }  ,draw opacity=0.81 ] (458.99,90.72) -- (436.25,108.4) -- (381.9,72.58) ; \draw  [color={rgb, 255:red, 155; green, 155; blue, 155 }  ,draw opacity=0.81 ] (436.25,108.4) -- (436.3,159.74) ;
\draw [color={rgb, 255:red, 208; green, 2; blue, 27 }  ,draw opacity=1 ]   (343.96,145.88) .. controls (362.7,127.9) and (363.99,118.31) .. (381.79,96.26) ;
\draw [color={rgb, 255:red, 208; green, 2; blue, 27 }  ,draw opacity=1 ]   (409.52,185.65) .. controls (418.89,170.55) and (423.39,153.55) .. (442.12,135.57) ;
\draw [color={rgb, 255:red, 208; green, 2; blue, 27 }  ,draw opacity=1 ]   (442.12,135.57) -- (461.8,117.12) ;
\draw [color={rgb, 255:red, 208; green, 2; blue, 27 }  ,draw opacity=1 ] [dash pattern={on 4.5pt off 4.5pt}]  (381.79,96.26) -- (401.47,77.8) ;
\draw [color={rgb, 255:red, 208; green, 2; blue, 27 }  ,draw opacity=1 ]   (461.8,117.12) .. controls (480.54,99.14) and (539.56,97.7) .. (558.29,79.72) ;
\draw [color={rgb, 255:red, 155; green, 155; blue, 155 }  ,draw opacity=0.81 ] [dash pattern={on 4.5pt off 4.5pt}]  (402.04,53.18) -- (401.65,104.89) ;
\draw [color={rgb, 255:red, 155; green, 155; blue, 155 }  ,draw opacity=0.81 ] [dash pattern={on 4.5pt off 4.5pt}]  (381.96,123.89) -- (401.65,104.89) ;
\draw [color={rgb, 255:red, 208; green, 2; blue, 27 }  ,draw opacity=1 ] [dash pattern={on 4.5pt off 4.5pt}]  (401.47,77.8) .. controls (406.9,74.92) and (421.32,71.78) .. (429.57,71.33) ;
\draw [color={rgb, 255:red, 208; green, 2; blue, 27 }  ,draw opacity=1 ] [dash pattern={on 4.5pt off 4.5pt}]  (429.57,71.33) .. controls (449.06,70.13) and (473.42,84.99) .. (498.34,79.48) ;
\draw [color={rgb, 255:red, 208; green, 2; blue, 27 }  ,draw opacity=1 ] [dash pattern={on 4.5pt off 4.5pt}]  (498.34,79.48) .. controls (526.25,74.68) and (549.49,52.15) .. (565.6,53.59) ;
\draw [color={rgb, 255:red, 208; green, 2; blue, 27 }  ,draw opacity=1 ]   (558.29,79.72) .. controls (577.41,65.81) and (577.97,54.07) .. (565.6,53.59) ;
\draw [color={rgb, 255:red, 74; green, 144; blue, 226 }  ,draw opacity=1 ] [dash pattern={on 4.5pt off 4.5pt}]  (527.69,38.13) .. controls (541.15,36.84) and (556.98,57.05) .. (557.92,83.19) .. controls (558.85,109.32) and (547.89,134.69) .. (528.51,135.56) ;
\draw [color={rgb, 255:red, 74; green, 144; blue, 226 }  ,draw opacity=1 ]   (391.53,12.19) -- (391.35,63.73) ;
\draw [color={rgb, 255:red, 74; green, 144; blue, 226 }  ,draw opacity=1 ] [dash pattern={on 4.5pt off 4.5pt}]  (391.16,72.64) -- (391.35,63.73) ;
\draw [color={rgb, 255:red, 74; green, 144; blue, 226 }  ,draw opacity=1 ] [dash pattern={on 4.5pt off 4.5pt}]  (391.55,100.89) -- (391.16,191.74) ;
\draw [color={rgb, 255:red, 74; green, 144; blue, 226 }  ,draw opacity=1 ]   (449.64,12.21) -- (449.99,113.15) ;
\draw [color={rgb, 255:red, 155; green, 155; blue, 155 }  ,draw opacity=0.81 ] [dash pattern={on 4.5pt off 4.5pt}]  (401.65,104.89) -- (459.04,142.07) ;
\draw [color={rgb, 255:red, 74; green, 144; blue, 226 }  ,draw opacity=1 ]   (450.56,141.44) -- (450.18,200.85) ;
\draw [color={rgb, 255:red, 74; green, 144; blue, 226 }  ,draw opacity=1 ]   (449.99,113.15) .. controls (468.73,95.17) and (541.99,84.39) .. (560.73,66.41) ;
\draw [color={rgb, 255:red, 74; green, 144; blue, 226 }  ,draw opacity=1 ]   (450.56,141.44) .. controls (467.98,115.07) and (539.37,108.84) .. (558.11,90.86) ;
\draw [color={rgb, 255:red, 74; green, 144; blue, 226 }  ,draw opacity=1 ] [dash pattern={on 4.5pt off 4.5pt}]  (419.08,64.71) .. controls (428.5,64.57) and (439.09,65.43) .. (449.58,66.34) ;
\draw [color={rgb, 255:red, 74; green, 144; blue, 226 }  ,draw opacity=1 ] [dash pattern={on 4.5pt off 4.5pt}]  (391.16,72.64) .. controls (401.28,65.21) and (416.27,65.45) .. (419.08,64.71) ;
\draw [color={rgb, 255:red, 74; green, 144; blue, 226 }  ,draw opacity=1 ]   (449.58,66.34) .. controls (468.17,66.65) and (491.96,75.04) .. (500.77,66.41) ;
\draw [color={rgb, 255:red, 74; green, 144; blue, 226 }  ,draw opacity=1 ] [dash pattern={on 4.5pt off 4.5pt}]  (391.55,100.89) .. controls (407.46,80.41) and (429.11,79.81) .. (449.59,82.79) .. controls (470.07,85.77) and (486.51,94.54) .. (497.96,90.38) ;
\draw [color={rgb, 255:red, 74; green, 144; blue, 226 }  ,draw opacity=1 ]   (560.73,66.41) .. controls (566.78,45.73) and (574.29,38.63) .. (587.65,37.29) .. controls (601.01,35.96) and (618.63,57.06) .. (618.04,85.6) .. controls (617.46,114.14) and (603.68,134.94) .. (588.47,134.73) .. controls (573.27,134.53) and (559.23,114.83) .. (558.11,90.86) ;
\draw [color={rgb, 255:red, 74; green, 144; blue, 226 }  ,draw opacity=1 ]   (502.4,113.42) .. controls (508.3,126.89) and (517.72,135.81) .. (528.51,135.56) ;
\draw [color={rgb, 255:red, 74; green, 144; blue, 226 }  ,draw opacity=1 ]   (500.77,66.41) .. controls (510.33,36.2) and (525.5,38.84) .. (527.69,38.13) ;
\draw [color={rgb, 255:red, 74; green, 144; blue, 226 }  ,draw opacity=1 ]   (527.69,38.13) -- (587.65,37.29) ;
\draw [color={rgb, 255:red, 74; green, 144; blue, 226 }  ,draw opacity=1 ]   (528.51,135.56) -- (588.47,134.73) ;
\draw  [draw opacity=0][fill={rgb, 255:red, 208; green, 2; blue, 27 }  ,fill opacity=0.1 ] (350.48,139.24) -- (359.29,128.48) -- (374.94,105.23) -- (381.79,96.26) -- (401.47,77.8) -- (414.29,73.35) -- (429.57,71.33) -- (441.46,71.91) -- (452.89,75.26) -- (467.13,78.14) -- (481.18,80.54) -- (489.8,81.02) -- (500.85,79.1) -- (521.28,72.39) -- (533.22,66.34) -- (538.33,63.76) -- (549.76,57.28) -- (560.25,54.17) -- (566.62,53.69) -- (571.68,54.89) -- (573.55,58.24) -- (573.37,62.56) -- (570.74,67.59) -- (565.12,74.06) -- (554.26,82.45) -- (544.51,88.21) -- (529.15,93.72) -- (499.17,100.91) -- (484.18,105.71) -- (467.32,112.42) -- (453.64,124.89) -- (439.96,138.55) -- (430.78,148.38) -- (424.97,156.77) -- (409.52,185.65) -- (404.82,195.58) -- (343.96,197.1) -- (343.96,145.88) -- cycle ;
\draw [color={rgb, 255:red, 74; green, 144; blue, 226 }  ,draw opacity=1 ] [dash pattern={on 4.5pt off 4.5pt}]  (496.51,90.67) .. controls (497.44,99) and (499.49,106.76) .. (502.4,113.42) ;
\draw  [draw opacity=0][fill={rgb, 255:red, 74; green, 144; blue, 226 }  ,fill opacity=0.19 ] (391.53,12.19) -- (391.16,72.64) -- (449.99,113.15) -- (449.64,12.21) -- cycle ;
\draw  [draw opacity=0][fill={rgb, 255:red, 74; green, 144; blue, 226 }  ,fill opacity=0.19 ] (391.55,100.89) -- (391.16,191.74) -- (450.18,200.85) -- (450.56,141.44) -- cycle ;
\draw  [draw opacity=0][fill={rgb, 255:red, 74; green, 144; blue, 226 }  ,fill opacity=0.19 ] (391.16,72.64) -- (402.87,67.08) -- (414.73,64.88) -- (428.47,64.68) -- (439.4,65.48) -- (453.3,66.48) -- (473.76,69.07) -- (490.93,69.87) -- (497.33,68.47) -- (500.77,66.41) -- (504.08,57.69) -- (508.45,49.5) -- (513.92,43.3) -- (521.1,39.11) -- (528.75,38.11) -- (541.09,37.71) -- (556.08,37.91) -- (569.78,37.69) -- (583.05,37.29) -- (581.02,38.29) -- (574.47,41.28) -- (569.78,46.48) -- (565.6,53.59) -- (560.73,66.41) -- (555.11,71.07) -- (544.8,75.87) -- (528.87,81.86) -- (500.61,91.05) -- (485.62,96.04) -- (469.34,101.73) -- (449.99,113.15) -- cycle ;
\draw  [draw opacity=0][fill={rgb, 255:red, 74; green, 144; blue, 226 }  ,fill opacity=0.19 ] (391.55,100.89) -- (398.84,92.85) -- (410.24,85.85) -- (423.98,82.06) -- (438.35,81.26) -- (451.77,83.26) -- (463.48,85.66) -- (475.2,88.85) -- (488.94,91.65) -- (496.51,90.67) -- (558.11,90.86) -- (552.49,95.44) -- (543.27,99.64) -- (526.57,105.43) -- (496.43,115.62) -- (473.79,123.81) -- (457.86,133.2) -- (450.56,141.44) -- cycle ;
\draw  [draw opacity=0][fill={rgb, 255:red, 74; green, 144; blue, 226 }  ,fill opacity=0.19 ] (558.11,90.86) -- (558.89,99.84) -- (559.98,105.23) -- (563.88,114.42) -- (565.85,118.05) -- (567.79,121.61) -- (572.63,127.41) -- (580.75,132.8) -- (588.47,134.73) -- (528.51,135.56) -- (521.88,134.6) -- (514.7,130.6) -- (508.3,123.61) -- (500.96,110.83) -- (497.99,100.84) -- (496.51,90.67) -- cycle ;
\draw  [draw opacity=0][fill={rgb, 255:red, 74; green, 144; blue, 226 }  ,fill opacity=0.19 ] (537.65,133.2) -- (544.83,127.61) -- (550.61,120.02) -- (554.2,112.03) -- (557.33,98.44) -- (557.92,83.19) -- (556.54,70.47) -- (552.17,57.69) -- (547.33,48.9) -- (542.65,44.7) -- (536.87,39.91) -- (530.78,38.51) -- (587.65,37.29) -- (593.08,37.91) -- (601.05,42.5) -- (607.13,49.1) -- (612.76,58.89) -- (616.66,70.87) -- (617.91,83.46) -- (616.97,97.24) -- (614.01,109.63) -- (609.79,119.02) -- (605.11,125.61) -- (599.48,131.2) -- (593.86,133.6) -- (588.47,134.73) -- (528.51,135.56) -- cycle ;
\draw [color={rgb, 255:red, 74; green, 144; blue, 226 }  ,draw opacity=1 ][line width=0.75]    (150.61,97.85) -- (150.79,146.1) ;
\draw [color={rgb, 255:red, 74; green, 144; blue, 226 }  ,draw opacity=1 ]   (150.79,146.1) -- (150.98,195.05) ;

\end{tikzpicture}

        \caption{In this picture, we represent the move explained assuming the hypothesis (D) for the case $K'\cong [0,1].$ On the left, we see the two strands of the immersion $c$ around the chart $x$ (colored in grey) around the intersection submanifold $c(K')$ (colored in black) colored in blue and red as in \cref{chart picture}. In this case, $U_0$ (resp. $U_1$) is the intersection of the red (resp. blue) strand with the grey box. To prevent visual cluttering, we represent only the core of $c$. In green, we see the "geometric dual" embedding $f$ and its intersection with $c$ (colored in black). In light purple, we see the embedded path $\gamma$ of embeddings of $K'$ between the black strips. On the right, we see the effect of the move. We see that the restriction to $\hb W$ is simply an embedded (double) connected sum with the spheres $f(\{0\}\times S^n\times \{0\})$ and $f(\{1\}\times S^n\times \{0\})$. This move can be thought of a $K'$-parameterized connected sum between $c$ and $f$ along the intersection $c(K').$} 
        \label{dual move picture}
    \end{figure}
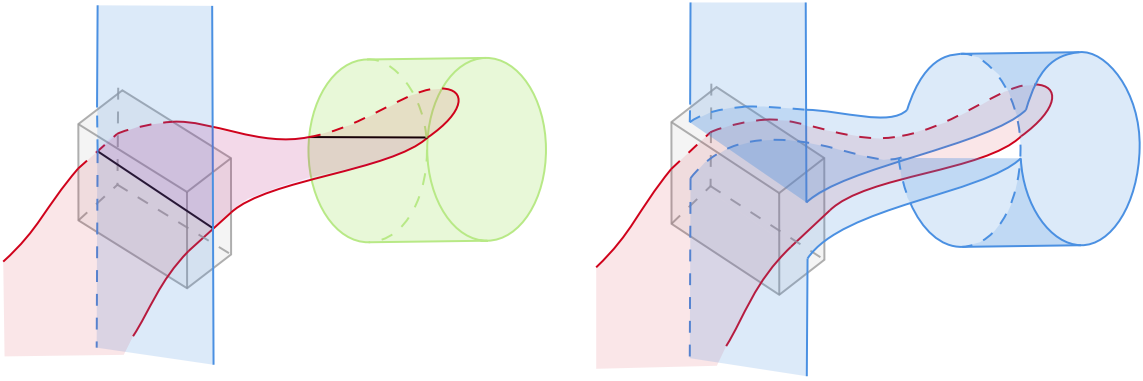
    \\\\
    \textit{Modifying $c$ to satisfy} (D). We shall now modify $c$ locally so that the result satisfies (D). In this step, we use the available copies of $V_1$ and the hypothesis $ n\geq 3$. Recall the embedding $\bar{e}:(S^n\times D^{n+1}_+,S^n_+\times \partial_0 D^{n+1}_+)\hookrightarrow (V_1,W_{1,1})$ from \cref{standard structure on V1}. Let $e:(S^n\times K'\times D^n,S^n\times \partial K'\times D^n)\hookrightarrow (V_1,W_{1,1})$ be the restriction of $\bar{e}$ along the inclusion $K'\times D^n\subset D^{n+1}$ given by $[0,1]\times D^n\hookrightarrow D^{n+1}$ (given by the first coordinate). The intersection of $\bar{f}(D^{n+1}\times \{0\})$ from loc.cit. with $e(S^n\times K'\times \{0\})$ is exactly $e(\{N\}\times K'\times \{0\}), $ where $N=(0,\cdots,0,1).$ By assumption, the cobordism $W_{j+1}$ contains an embedded copy of $V_1$ disjoint from the immersion $c.$ We proceed by defining an embedded boundary connected sum of the core of the immersion $c$ with the embedding $\bar{f}(D^{n+1}\times \{0\})$. To do so, we pick an embedded path $\alpha$ in $\hb W_{j+1}$ from a point $p$ in the image of $c$ disjoint from the self-intersection locus $c^{-1}(K)$ and $P$ to a point in the image of $\bar{f}$ disjoint from $e(\{N\}\times K'\times \{0\}),$ which only intersects the image of the immersion $c$ or the image of $\bar{f}$ at the endpoints. This is possible by transversality, using $n\geq 2$. The result of this embedded boundary connected sum is an immersion $c'\coloneqq c\natural_\alpha \bar{f} :(D^{n+1}\times \{0\},\partial_0 D^{n+1}_+\times \{0\})\hookrightarrow (W_{j+1},\hb W_{j+1})$ satisfying the following properties by construction: 
    \begin{enumerate}[label=(\roman*)]
        \item It has trivial normal bundle: this follows from the fact that an embedded connected sum of two embeddings with trivial normal bundle has trivial normal bundle.
        \item It differs from $c$ in a neighborhood of $c^{-1}(p)$.
        \item Its self-intersection locus is the same as the one of $c$.
        \item Its intersection with the embedding $\bar{e}$ is exactly $e(\{N\}\times K'\times \{0\}).$  
    \end{enumerate}Thus, it suffices to find an embedding $\gamma: (K'\times [0,1],\partial K'\times [0,1])\hookrightarrow (D^{n+1}_+,\partial_0D^{n+1}_+) $ satisfying the conditions of the hypothesis (D)  with respect to the immersion $c'$ and the embedding $e$. However, since $(D^{n+1}_+,\partial_0D^{n+1}_+)$ is $1$-connected (as $n\geq 0$), it follows that both embeddings $(K',\hb K')\to (D^{n+1}_+,\partial_0D^{n+1}_+)$ induced by the self-intersection locus of $c$ and $e(\{N\}\times K'\times \{0\})$ are both homotopic as maps of pairs to the standard inclusion $[0,1]\subset D^{n+1}$. Since $n\geq 3$, these embeddings are isotopic (by \textit{homotopy implies isotopy} using that $n+1\geq 2\dim(K')+2$) and moreover the embedding $(K'\sqcup K',\hb K'\sqcup \hb K')\hookrightarrow (D^{n+1}_+,\partial_0D^{n+1}_+) $, given by both embeddings of $K'$ above, is isotopic to the disjoint union of the standard embedding $[0,1]\subset D^{n+1}_+$ and a translate. By isotopy extension, one finds an ambient isotopy of $(D^{n+1}_+,\partial_0D^{n+1}_+)$ relative to $\partial_1D^{n+1}_+$ extending the isotopy of embeddings of $K'\sqcup K'$ from above. Now, this ambient isotopy allows us to pullback an embedding $\gamma_0:([0,1]\times [0,1],(\partial [0,1])\times [0,1])\hookrightarrow (D^{n+1}_+,\partial_0D^{n+1}_+)$ which restricts to the standard embedding on $[0,1]\times \{0\}$ and its translate on $[0,1]\times \{1\}$, hence obtaining the desired embedding $\gamma$. By transversality, one can assume that this embedding is disjoint from the remaining components of the preimage of the self-intersection locus, since $n\geq 3$. We conclude that $c'$ satisfies hypothesis (D).
    \\\\
    \textit{Rechoosing $F$ on vertices up to homotopy.} We shall define a map $F_v:(I^k,\partial I^k)\to (||Y_\bullet||,||X_\bullet||)$ homotopic to $F$. Choose a compactly supported vector field on $W_{j}$ extending the translation in $\phi(\partial_1D^{n+1}_+\times (0,\infty)\times \RR^{n-1})\subset P$ by the unit vector of $(0,\infty).$ By taking its flow, we produce a family $\psi_\epsilon:W_j\to W_j$ of diffeomorphism with the property that $\psi_\epsilon(\phi(x,t,v))=\phi(x,t+\epsilon,v) $ for $\epsilon>0.$ Let $c_\epsilon$ be the immersion $\psi_\epsilon(c)$ and notice that $(t+\epsilon, c_\epsilon, \psi_\epsilon^*\L)\in Y_0.$ Consider the immersion $c'$ obtained from $c_\epsilon$ by first modifying it so it satisfies (D) for every path component of the self-intersection loci and further modifying it to remove all of those components, by applying the construction above iteratively. In other words, $c'$ is an embedding to $W_{j'}$ for some $j'\geq j$ which agrees with $c_\epsilon$ outside a neighborhood of $K.$ Moreover, the modifications above can be chosen to be disjoint from $c$ and thus $c'$ is disjoint from $c.$ To produce a path $\L'$ from $(c')^*\ell_{W_{j'}}$ to $\ell_{t+\epsilon},$ we shall construct a path from $c_\epsilon^*\ell_{W_{j'}}$ to $(c')^*\ell_{W_{j'}}.$ When modifying $c_\epsilon$ to satisfy (D), we replaced a part of this immersion by a part of the embedding $\smash{\bar{f}}.$ However, $\smash{\bar{f}^*\ell_{W_{j'}}}$ is standard in the sense of \cref{left standard}, so extends to $(\RR^{n-1}\times \RR^{n+2}_+,\RR^{n-1}\times \partial \RR^{n+2}_+)$. Thus the homotopy class of the $\Theta$-structure does not change after this modification. Similarly, the second modification is done by gluing $c$ with the embedding $f.$ However, in this particular case, $f^*\ell_{W_j}$ extends over a contractible space since it is a subspace of the embedding $\bar{e}$ in $(V_1,W_{1,1})$ which again has standard structure. We conclude that $c_\epsilon^*\ell_{W_{j'}}$ and $(c')^*\ell_{W_{j'}}$ are homotopic relative to $\partial_1D^{n+1}_+\times D^n.$ We let $\L'$ be the composite of the path $\L_\epsilon $ with such a homotopy. Let $F_v(w)=F(w)$ for all vertices $w\neq v$ and $F_v(v)=(t+\epsilon, c',\L')\in X_0.$ This assignment extends to $1$-simplices uniquely since $c'$ stays transverse to any other immersions, to which $c$ was transverse. Moreover, it extends to higher simplices uniquely as they are determined by the $0-$ and $1$-simplices. Moreover, $F_v$ is homotopic to $F$ since $(F_v(v),F(v))\in Y_1.$ We let $F'$ be the result of iterating this modification to all vertices of the triangulation of $I^k.$ This produces a homotopic map whose vertices map to $X_\bullet.$ \\\\ 
    \textit{Rechoosing $F'$ to lie in $X_\bullet.$} By construction, the images of all the vertices along $F'$ lie in $X_0.$ In other words, all immersions $c$ given by $F'(v)=(t,c,\L)$ are actually embeddings when restricted to its cores. To produce a map $F''$ where all the $1$-simplices lie in $X_1$, it suffices to modify these embeddings $c$ to be pairwise disjoint. To do so, we simply apply the same technique above for removing self-intersections of the immersion
    \[\underset{(t,c,\L)\in F'((I^k)_0)}{\coprod} c: \coprod (D^{n+1}_+\times D^n,\partial_0 D^{n+1}_+\times D^n)\to (W_{\infty},\hb W_\infty).\]The self-intersection locus of this immersion is precisely all the pairwise intersections of the embeddings $c.$ Thus, modifying such immersion to an embedding is making all $1$-simplices lie in $X_1.$ Proceeding as above, the result will be a homotopic map $F''$ sending $0$- and $1$-simplices to $X_\bullet$ and thus, sending any simplex of any dimension to $X_\bullet.$ This can be applied to any such map $F$ and thus, concludes the proof.\end{proof}

\begin{prop}\label{from no (v) to (v)}
    The map 
    \[\debd(W_\infty,\phi,\ell)_\bullet\to \dec(W_\infty,\phi,\ell)_\bullet\]induces an equivalence on geometric realizations.
\end{prop}

\begin{lemma}\label{only need to make boundary n con}
    Let $x=\{(t_i,c_i,\L_i)\}_{i=0,\cdots,p}$ be a $p$-simplex of $\dec(W,\phi,\ell)$. Then $x\in \debd(W,\phi,\ell)_p$ if and only if the map $\hb W\backslash \hb C\to B^{\partial}$ is $n$-connected.
\end{lemma}
\begin{proof}
    The "only if" direction is immediate from the definition. Given a $p$-simplex $x=\{(t_i,c_i,\L_i)\}_{i=0,\cdots,p}$ in $\dec(W,\phi,\ell)_\bullet$, the cobordism $W$ is obtained from $W\backslash C$ by attaching a right handle of index $n$ at the belt of the each left handle $c_i.$ For the "if" direction, it suffices to show $W\backslash C\cup_{\hb W\backslash\hb C}B^\partial\to B$ is $(n+1)$-connected. We can replace the source by the equivalent $W\backslash C\cup_{\hb W\backslash\hb C} \hb W \cup_{\hb W} B^\partial.$ Since $W\backslash C\to W$ is equivalent to the attachment of a right handle, we have $W\backslash C\cup_{\hb W\backslash\hb C} \hb W \simeq W.$ Thus, the map above is equivalent to the map $W\cup_{\hb W}B^\partial\to B$ which is $(n+1)$-connected by hypothesis.\end{proof}

\begin{proof}[Proof of \cref{from no (v) to (v)}]
    Denote the source and target of the map in the statement by $Z_\bullet$ and $X_\bullet$ respectively. Let $F:(I^k,\partial I^k)\to (||X_\bullet||,||Z_\bullet||)$ be a map, which we can assume to be simplicial with respect to some triangulation of $I^k,$ and $\sigma$ be a simplex of $I^k.$ We shall modify $F$ up to homotopy so that the cores of $F(\sigma)$ satisfy the condition of \cref{only need to make boundary n con}. Let $C_\sigma$ be the union of the images of $c_i$ for $F(\sigma)=\{(t_i,c_i,\L_i)\}.$ Let $K_\sigma(F)$ be the kernel of the map $\pi_{n-1}(\hb W\backslash\hb C_\sigma)\to \pi_{n-1}(B^\partial).$ Let $I_\sigma(F)$ be the image of $\pi_{n}(\hb W\backslash\hb C_\sigma)\to \pi_{n}(B^\partial).$ Since the map $\hb W\backslash \hb C_\sigma\to \hb W$ is $(n-1)$-connected, it suffices to replace $F$ so that $K_\sigma(F)$ vanishes and $I_\sigma(F)=\pi_n(B^\partial).$ If $n\geq 3,$ then $\pi_1(\hb W\backslash \hb C_\sigma)\cong \pi_1(\hb W)\cong \pi_1(B^\partial).$ In this case, we apply the Hurewicz theorem to deduce that $\pi_n(\hb W,\hb W\backslash \hb C_\sigma)$ is generated as a $\pi_1(B^\partial)$-module by the meridian spheres $c_i|_{\{0\}\times \partial D^n}.$ This module surjects onto $K_\sigma(F)$, which proves that the latter is a finitely generated $\ZZ[\pi_1(B^\partial)]$-module. If $n=2$, then the Seifert-Van Kampen theorem shows that $K_\sigma(F)$ is normally generated by the meridian circles of the removed handles (as pointed out in \cite[166]{GRWII}). Moreover since $n\geq 2$, if $K_\sigma (F)$ vanishes, then $\pi_1(\hb W\backslash \hb C_\sigma)\cong \pi_1(B^\partial)$ and the following exact sequence of $\pi_1(B^\partial)$-modules
    \[\cdots\to \pi_n(\hb W\backslash \hb C_\sigma)\to \pi_n(B^\partial)\to \pi_n(B^\partial, \hb W\backslash \hb C_\sigma)\to 0 \]shows that $\pi_n(B^\partial)$ is generated by $I_\sigma(F)$ and finitely many elements lifted from the rightmost module, since the latter is finitely generated.
    \\\\
    \textit{Killing $K_\sigma(F).$} We describe a modification on the embeddings $c_i$ and describe its main properties. This is analogous to \cite[166]{GRWII} and similar to the modification of $c$ to satisfy (D) in the proof of \cref{from immersed to embedded cores}. We assume $n\geq 2.$ Let $v$ be a vertex of $\sigma$ and $F(v)=(t,c,\L)$ and let $j\geq 0$ be such that $\im(c)\subset W_j.$ We start by considering the perturbation $c_\epsilon$ for some $\epsilon>0$ as defined in \cref{from immersed to embedded cores}. Recall the embedding $\bar{f}:(D^{n+1}\times D^n,\partial D^{n+1}\times D^n)\hookrightarrow (V_1,W_{1,1}).$ Denote also by $\bar{f}$ the composition of this embedding with the inclusion of a fixed copy of $V_1$ in $K|_{[j,j+1]}$ seen now in $W_{j+1}.$ Choose an embedded path $\gamma$ in $\hb W_{j+1}$ from $c_\epsilon(x,0)$ for some point $x\in \partial_0D^{n+1}$ and $\bar{f}(x,0)$ disjoint from $\hvb W_j$ and $C_\sigma$, and from $\im (c_\epsilon)$ and $\im(\bar{f})$ except at times $t=\{0,1\}.$ Using a thickening of this path, we can construct an embedded connected sum of the cores of $c_\epsilon$ and $\bar{f}$ and further thicken it to an embedding $c':(D^{n+1}_+\times D^n,\partial_0 D^{n+1}\times D^n)\hookrightarrow (W_{j+1},\hb W_{j+1})$ which agrees with $c_\epsilon$ outside of a neighborhood of $x.$ As in \cref{from immersed to embedded cores}, since the $\Theta$-structure $\bar{f}^*\ell_{W_{j+1}}$ is standard, the structures $(c')^*\ell_{W_{j+1}}$ and $c_\epsilon^*\ell_{W_{j+1}}$ are homotopic. By choosing such homotopy and composing it with the path $\L_\epsilon$, we construct a path $\L'$ from $(c')^*\ell_{W_{j+1}}$ to $\ell_{t+\epsilon}.$ Denote $F_v:(I^k,\partial I^k)\to (||X_\bullet||,||Z_\bullet||)$ the map definedd on vertices by $F_v(w)=F(w)$ for $w\neq v$ and $F_v(v)=(t+\epsilon,c',\L')$. For $\epsilon$ small enough, we see that this map is well-defined and homotopic to $F$ relative to $\partial I^k.$ (See \cref{adding dual picture 1} for a picture portraying this move.) The following claim establishes properties of this construction.
    \begin{figure}
        \centering
        \tikzset{every picture/.style={line width=0.75pt}} 

\begin{tikzpicture}[x=0.75pt,y=0.75pt,yscale=-1,xscale=1]

\draw [color={rgb, 255:red, 208; green, 2; blue, 27 }  ,draw opacity=1 ]   (170.09,189.4) .. controls (168.89,172.16) and (179.32,116.32) .. (209.15,114.73) ;
\draw [color={rgb, 255:red, 208; green, 2; blue, 27 }  ,draw opacity=1 ]   (209.15,114.73) .. controls (231.08,109.62) and (253.01,189.07) .. (256.25,189.39) ;
\draw [color={rgb, 255:red, 126; green, 211; blue, 33 }  ,draw opacity=1 ]   (209.15,114.73) .. controls (208.75,105.98) and (209.51,91.62) .. (220.47,91.84) .. controls (231.44,92.07) and (240.79,115.68) .. (235.04,134.51) ;
\draw [color={rgb, 255:red, 126; green, 211; blue, 33 }  ,draw opacity=1 ]   (209.15,114.73) .. controls (209.51,129.72) and (215.27,150.46) .. (231.08,139.61) ;
\draw  [draw opacity=0][fill={rgb, 255:red, 208; green, 2; blue, 27 }  ,fill opacity=0.14 ] (256.25,189.39) -- (170.09,189.4) -- (170.03,180.99) -- (173.5,160.83) -- (178.54,144.56) -- (185.19,131.32) -- (192.14,122.54) -- (199.51,117.44) -- (205.8,115.2) -- (212.27,114.57) -- (217.12,116.32) -- (223.23,121.75) -- (230.78,131.96) -- (238.69,148.55) -- cycle ;
\draw  [fill={rgb, 255:red, 0; green, 0; blue, 0 }  ,fill opacity=1 ] (209.81,109.92) .. controls (209.81,109.56) and (209.54,109.27) .. (209.21,109.27) .. controls (208.88,109.27) and (208.61,109.56) .. (208.61,109.92) .. controls (208.61,110.29) and (208.88,110.58) .. (209.21,110.58) .. controls (209.54,110.58) and (209.81,110.29) .. (209.81,109.92) -- cycle ;
\draw  [fill={rgb, 255:red, 0; green, 0; blue, 0 }  ,fill opacity=1 ] (210.17,120.14) .. controls (210.17,119.77) and (209.9,119.48) .. (209.57,119.48) .. controls (209.24,119.48) and (208.97,119.77) .. (208.97,120.14) .. controls (208.97,120.5) and (209.24,120.79) .. (209.57,120.79) .. controls (209.9,120.79) and (210.17,120.5) .. (210.17,120.14) -- cycle ;
\draw [color={rgb, 255:red, 208; green, 2; blue, 27 }  ,draw opacity=1 ]   (310.29,189.88) .. controls (309.09,172.64) and (319.52,116.8) .. (349.35,115.2) ;
\draw [color={rgb, 255:red, 208; green, 2; blue, 27 }  ,draw opacity=1 ]   (349.35,115.2) .. controls (359.54,112.83) and (369.72,128.7) .. (378.03,146.4) ;
\draw  [draw opacity=0][fill={rgb, 255:red, 208; green, 2; blue, 27 }  ,fill opacity=0.14 ] (396.45,189.87) -- (310.29,189.88) -- (310.23,181.47) -- (313.7,161.31) -- (318.74,145.04) -- (325.39,131.8) -- (332.34,123.02) -- (339.71,117.92) -- (346,115.68) -- (352.47,115.04) -- (357.32,116.8) -- (363.43,122.22) -- (370.98,132.43) -- (378.03,146.4) -- (392.19,142.86) -- (409.99,138.22) -- (425.67,132.06) -- (425.62,123.54) -- (427.96,112.86) -- (431.02,106.63) -- (435.87,102.96) -- (438.6,102.52) -- (443.06,103.28) -- (446.83,104.72) -- (451.15,108.23) -- (454.74,113.49) -- (457.44,120.51) -- (459.06,127.05) -- (457.98,135.35) -- (454.92,142.21) -- (451.33,145.4) -- (445.76,147.95) -- (439.38,149.32) -- (434.43,147.79) -- (430.66,144.76) -- (427.08,138.76) -- (420.23,141.73) -- (402.8,146.84) -- (388.78,152.26) -- (381.84,154.91) -- cycle ;
\draw [color={rgb, 255:red, 126; green, 211; blue, 33 }  ,draw opacity=1 ]   (438.6,102.52) .. controls (436.95,88.77) and (466.34,89.25) .. (456.63,111.75) ;
\draw  [fill={rgb, 255:red, 0; green, 0; blue, 0 }  ,fill opacity=1 ] (439.79,97.72) .. controls (439.79,97.36) and (439.53,97.06) .. (439.2,97.06) .. controls (438.87,97.06) and (438.6,97.36) .. (438.6,97.72) .. controls (438.6,98.08) and (438.87,98.38) .. (439.2,98.38) .. controls (439.53,98.38) and (439.79,98.08) .. (439.79,97.72) -- cycle ;
\draw [color={rgb, 255:red, 126; green, 211; blue, 33 }  ,draw opacity=1 ]   (438.6,102.52) .. controls (438.84,111.03) and (446.39,118.69) .. (453.93,114.86) ;
\draw  [fill={rgb, 255:red, 0; green, 0; blue, 0 }  ,fill opacity=1 ] (439.79,106.34) .. controls (439.79,105.97) and (439.53,105.68) .. (439.2,105.68) .. controls (438.87,105.68) and (438.6,105.97) .. (438.6,106.34) .. controls (438.6,106.7) and (438.87,106.99) .. (439.2,106.99) .. controls (439.53,106.99) and (439.79,106.7) .. (439.79,106.34) -- cycle ;
\draw [color={rgb, 255:red, 208; green, 2; blue, 27 }  ,draw opacity=1 ]   (438.6,102.52) .. controls (429.24,102.32) and (424.17,118.64) .. (425.67,132.06) ;
\draw [color={rgb, 255:red, 208; green, 2; blue, 27 }  ,draw opacity=1 ]   (438.6,102.52) .. controls (458.25,100.74) and (472,146.69) .. (439.38,149.32) ;
\draw [color={rgb, 255:red, 208; green, 2; blue, 27 }  ,draw opacity=1 ]   (378.03,146.4) .. controls (400.01,139.75) and (408.37,140.46) .. (425.67,132.06) ;
\draw [color={rgb, 255:red, 208; green, 2; blue, 27 }  ,draw opacity=1 ]   (381.84,154.91) .. controls (403.25,145.49) and (417.54,143.58) .. (427.08,138.76) ;
\draw [color={rgb, 255:red, 208; green, 2; blue, 27 }  ,draw opacity=1 ]   (381.84,154.91) .. controls (385.72,164.16) and (391.38,180.91) .. (396.45,189.87) ;
\draw [color={rgb, 255:red, 208; green, 2; blue, 27 }  ,draw opacity=1 ]   (427.08,138.76) .. controls (429.13,144.05) and (433.71,148.6) .. (439.38,149.32) ;

\end{tikzpicture}

        \caption{In this picture, we see the effect of an embedded $D^{n+1}_+$ in $W$, represented in red, having a "geometrically dual sphere", represented in green, on the meridian class. The meridian class, depicted in black in this picture, lives in $\pi_{n-1}(\hb W\backslash \partial_0 D^{n+1}_+)$. Notice that when a geometrically dual sphere exists, in particular it follows that the meridian class vanishes. Observe that taking an embedded connected sum with a disc $D^{n+1}$ admitting a dual sphere creates a dual sphere for the result. The move described is a special case of this fact, where the source of this additional disc comes from the existing $V_1$. Additionally, this dual sphere allows us to modify the initial disc to be disjoint from a fixed embedded $n$-sphere in $\hb W$. This is a $0$-dimensional version of hypothesis (D) in \cref{from immersed to embedded cores}.}
        \label{adding dual picture 1}
    \end{figure}
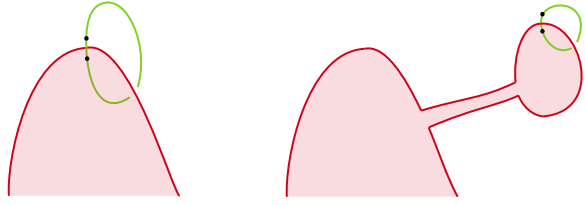
    
    \begin{claim}\label{claim 5.6}
     Given $k,F$ and $v$ as above, the map $F_v$ is homotopic to $F$ relative to $\partial I^k$ satisfies the following properties: 
     \begin{enumerate}[label=(\roman*)]
         \item for every simplex $\sigma\in I^k,$ there exists a surjection $K_\sigma(F)\to K_\sigma(F_v),$ sending the class of the meridian of each simplex of $F$ to the meridian of each simplex of $F_v$ except $v$ and whose kernel contains the class of the meridian of $v.$ 
         \item for every simplex $\sigma\in I^k,$ we have $I_\sigma(F)\subseteq I_\sigma(F_v).$
     \end{enumerate}
\end{claim}\begin{proof}
    Apply \cite[Claim 5.6]{GRWStableModuli} to the horizontal boundary of $c.$ The only difference on our construction that we take a connected sum with the sphere $\bar{f}$ instead of $\bar{e}$, but the proof follows verbatim.
\end{proof}
Apply the construction above and \cref{claim 5.6} for every interior vertex $v\in I^k.$ This produces a homotopic map $F'$ where the meridian of all vertices are trivial in $K_\sigma(F').$ Since this kernel is generated by such meridians, we have $K_\sigma(F')=0.$ By induction on the dimension of the simplices, we can obtain $F'$ such that $K_\sigma(F')=0$ for every simplex $\sigma$ of $I^k.$
    \\\\
    \textit{Adding elements to $I_\sigma(F').$} Since $K_\sigma(F)=0,$ it follows that $\pi_n(B^\partial)$ is generated by $I_\sigma(F')$ along with finitely many other elements. Choose finitely many elements 
    \[\{g_\alpha:S^n\to \hb W_\infty\}_\alpha\]such that $[l_{\hb W_\infty}\circ g_\alpha]\in \pi_n(B^\partial)$ generate $\pi_n(B^\partial)$ along with $I_\sigma(F').$ If $g_\alpha$ is disjoint from $C_\sigma$, then this class lifts to $\pi_n(\hb W_\infty\backslash C_\sigma)$ and thus, lies in $I_\sigma(F').$ We shall modify $F'$ and $g_\alpha$ so this is the case. This modification is heavily inspired by the proof of \cite[Prop. 5.5]{GRWII} and in fact it restricts to it when restricted to the horizontal boundary. We may suppose that all the $g_\alpha$ are smooth and transverse to the embeddings given by the images of the vertices of $\sigma$ along $F'.$ Let $v$ be a vertex of $\sigma$ and $F'(v)=(t,c,\L)$ and let $(t+\epsilon,c_\epsilon,\L_\epsilon)$ be an $\epsilon$-perturbation as before. By compactness, the intersection of $c_\epsilon(D^{n+1}_+\times \{0\})$ with $g_\alpha$ is a finite set of points. Let $p$ be such an intersection point, then we can find a chart $x:(D^{n+1}_+\times D^n,\partial_0D^{n+1}_+\times D^n)\hookrightarrow (W_\infty, \hb W_\infty)$ disjoint from $\vb W_\infty$ where $x(0,0)=p$ such that:
\begin{enumerate}[label=(\roman*)]
    \item $x(D^{n+1}_+\times D^n)\cap c_\epsilon(D^{n+1}_+\times \{0\})=x(D^{n+1}_+\times \{0\})$ and $x(D^{n+1}_+\times D^n)\cap g_\alpha=x(\{0\}\times D^n);$
    \item additionally, $x(D^{n+1}_+\times D^n)\cap c_\epsilon(D^{n+1}_+\times D^n)=x(D^{n+1}_+\times \frac{1}{2}D^n).$
\end{enumerate}

\begin{figure}
    \centering
    \tikzset{every picture/.style={line width=0.75pt}} 

\begin{tikzpicture}[x=0.75pt,y=0.75pt,yscale=-1,xscale=1]

\draw    (82,98.39) .. controls (106.86,74.74) and (153.66,74.74) .. (171.39,98.39) ;
\draw  [dash pattern={on 4.5pt off 4.5pt}]  (82,188.7) .. controls (106.86,165.05) and (153.66,165.05) .. (171.39,188.7) ;
\draw [color={rgb, 255:red, 208; green, 2; blue, 27 }  ,draw opacity=1 ]   (118.52,98.03) -- (118.56,188.66) ;
\draw [color={rgb, 255:red, 208; green, 2; blue, 27 }  ,draw opacity=1 ]   (118.52,98.03) -- (128.13,80.83) ;
\draw [color={rgb, 255:red, 208; green, 2; blue, 27 }  ,draw opacity=1 ] [dash pattern={on 4.5pt off 4.5pt}]  (128.13,80.83) -- (127.78,171.11) ;
\draw [color={rgb, 255:red, 208; green, 2; blue, 27 }  ,draw opacity=1 ] [dash pattern={on 4.5pt off 4.5pt}]  (118.56,188.66) -- (127.78,171.11) ;
\draw [color={rgb, 255:red, 74; green, 144; blue, 226 }  ,draw opacity=1 ]   (82,141.73) -- (170.68,141.73) ;
\draw    (82,98.39) -- (82,188.7) ;
\draw    (82,188.7) -- (171.39,188.7) ;
\draw    (82,98.39) -- (171.39,98.39) ;
\draw    (171.39,98.39) -- (171.39,188.7) ;
\draw    (479.83,110.2) .. controls (475.96,122.77) and (540.26,125.28) .. (535.05,110.2) ;
\draw    (486.97,117.32) .. controls (493.92,105.39) and (524.62,111.67) .. (526.36,118.58) ;
\draw  [dash pattern={on 4.5pt off 4.5pt}]  (416.3,85.91) .. controls (429.05,75.23) and (437.74,155.02) .. (420.36,158.79) ;
\draw    (272.6,97.9) .. controls (297.45,74.25) and (344.26,74.25) .. (361.98,97.9) ;
\draw  [dash pattern={on 4.5pt off 4.5pt}]  (272.6,188.21) .. controls (297.45,164.57) and (344.26,164.57) .. (361.98,188.21) ;
\draw [color={rgb, 255:red, 208; green, 2; blue, 27 }  ,draw opacity=1 ]   (309.11,97.54) -- (318.73,80.34) ;
\draw [color={rgb, 255:red, 208; green, 2; blue, 27 }  ,draw opacity=1 ] [dash pattern={on 4.5pt off 4.5pt}]  (318.73,80.34) -- (318.37,170.62) ;
\draw [color={rgb, 255:red, 208; green, 2; blue, 27 }  ,draw opacity=1 ] [dash pattern={on 4.5pt off 4.5pt}]  (309.15,188.17) -- (318.37,170.62) ;
\draw    (272.6,97.9) -- (272.6,188.21) ;
\draw    (272.6,188.21) -- (361.98,188.21) ;
\draw    (272.6,97.9) -- (361.98,97.9) ;
\draw    (361.98,97.9) -- (361.98,114.59) ;
\draw    (361.98,114.59) .. controls (379.85,109.72) and (398.57,113.62) .. (410.45,113.17) ;
\draw    (361.98,136.04) .. controls (373.9,130.19) and (404.53,137.99) .. (410.75,134.47) ;
\draw    (361.98,136.04) -- (361.98,188.21) ;
\draw    (416.3,85.91) .. controls (411.76,98.99) and (411.76,103.38) .. (410.45,113.17) ;
\draw    (410.75,134.47) .. controls (412.67,144.6) and (413.81,151.1) .. (420.36,158.79) ;
\draw    (420.36,158.79) .. controls (438.48,167.36) and (482.44,171.91) .. (529.23,169.73) ;
\draw    (529.23,169.73) .. controls (558.74,167.36) and (591.35,166.71) .. (603.97,153.76) ;
\draw    (589.49,89.68) .. controls (619.53,99.1) and (621.51,140.7) .. (603.97,153.76) ;
\draw    (416.3,85.91) .. controls (454.46,68.55) and (557.41,80.25) .. (589.49,89.68) ;
\draw [color={rgb, 255:red, 74; green, 144; blue, 226 }  ,draw opacity=1 ]   (338.17,141.24) .. controls (373.06,122.5) and (437.73,135.5) .. (450.49,139.73) ;
\draw [color={rgb, 255:red, 74; green, 144; blue, 226 }  ,draw opacity=1 ]   (450.49,139.73) .. controls (487.64,147.2) and (571.03,135.83) .. (571.6,120.55) ;
\draw [color={rgb, 255:red, 74; green, 144; blue, 226 }  ,draw opacity=1 ]   (448.79,103.65) .. controls (482.82,74.4) and (567.62,94.55) .. (571.6,120.55) ;
\draw [color={rgb, 255:red, 74; green, 144; blue, 226 }  ,draw opacity=1 ]   (291.66,141.25) .. controls (291.66,124.13) and (325.45,122.66) .. (374.48,120.55) .. controls (423.51,118.44) and (428.08,120.23) .. (448.79,103.65) ;
\draw [color={rgb, 255:red, 208; green, 2; blue, 27 }  ,draw opacity=1 ]   (409.74,128.03) .. controls (447.94,126.73) and (479.42,125.75) .. (486.97,117.32) ;
\draw [color={rgb, 255:red, 208; green, 2; blue, 27 }  ,draw opacity=1 ] [dash pattern={on 4.5pt off 4.5pt}]  (434.04,80.57) .. controls (460.98,95.2) and (464.95,108.2) .. (486.97,117.32) ;
\draw [color={rgb, 255:red, 208; green, 2; blue, 27 }  ,draw opacity=1 ]   (309.14,154.03) .. controls (313.37,138.69) and (338.74,131.66) .. (365.97,129) .. controls (393.2,126.35) and (399.23,127.32) .. (409.74,128.03) ;
\draw [color={rgb, 255:red, 74; green, 144; blue, 226 }  ,draw opacity=1 ]   (272.59,141.25) -- (291.66,141.25) ;
\draw [color={rgb, 255:red, 74; green, 144; blue, 226 }  ,draw opacity=1 ]   (338.17,141.24) -- (361.27,141.24) ;
\draw [color={rgb, 255:red, 208; green, 2; blue, 27 }  ,draw opacity=1 ]   (396.88,116) .. controls (417.31,116.33) and (429.22,112.1) .. (434.04,80.57) ;
\draw [color={rgb, 255:red, 208; green, 2; blue, 27 }  ,draw opacity=1 ]   (309.12,117.95) .. controls (336.19,115.35) and (381.29,117.63) .. (396.88,116) ;
\draw [color={rgb, 255:red, 208; green, 2; blue, 27 }  ,draw opacity=1 ]   (309.11,97.54) -- (309.12,117.95) ;
\draw [color={rgb, 255:red, 208; green, 2; blue, 27 }  ,draw opacity=1 ]   (309.14,154.03) -- (309.15,188.17) ;
\draw [color={rgb, 255:red, 208; green, 2; blue, 27 }  ,draw opacity=1 ] [dash pattern={on 0.84pt off 2.51pt}]  (309.12,117.95) .. controls (317.19,131.6) and (312.65,143.95) .. (309.14,154.03) ;
\draw [color={rgb, 255:red, 208; green, 2; blue, 27 }  ,draw opacity=1 ] [dash pattern={on 0.84pt off 2.51pt}]  (410.02,115.68) .. controls (413.05,116.98) and (414.47,123.8) .. (409.74,128.03) ;
\draw  [draw opacity=0][fill={rgb, 255:red, 208; green, 2; blue, 27 }  ,fill opacity=0.16 ] (318.37,170.62) -- (309.15,188.17) -- (309.14,154.03) -- (312.08,144.28) -- (313.5,134.53) -- (312.08,123.48) -- (309.12,117.95) -- (309.11,97.54) -- (318.73,80.34) -- (318.55,125.48) -- (318.32,143.3) -- cycle ;
\draw  [draw opacity=0][fill={rgb, 255:red, 208; green, 2; blue, 27 }  ,fill opacity=0.16 ] (309.14,154.03) -- (312.08,144.28) -- (313.5,134.53) -- (312.65,126.08) -- (309.12,117.95) -- (332.79,116.98) -- (357.74,116.98) -- (378.17,117.3) -- (399.72,116.33) -- (412.48,115.03) -- (420.42,111.12) -- (428.08,101.37) -- (434.04,80.57) -- (451.62,91.62) -- (465.52,103.65) -- (476.3,112.1) -- (486.97,117.32) -- (479.42,121.53) -- (460.98,125.75) -- (436.59,127.7) -- (404.54,128.03) -- (379.87,127.7) -- (355.76,129.98) -- (336.19,134.85) -- (323.99,139.4) -- (314.92,145.9) -- cycle ;
\draw  [draw opacity=0][fill={rgb, 255:red, 208; green, 2; blue, 27 }  ,fill opacity=0.16 ] (127.78,171.11) -- (118.56,188.66) -- (118.52,98.03) -- (128.13,80.83) -- cycle ;

\end{tikzpicture}

    \caption{On the left, we see the chart $x$ where the red portion represents the intersection of the core of $c$ with the domain of $x$ and the blue arc the intersection of $g_\alpha.$ On the right, we see the result of the explained modification. The red portion in the right picture represents the embedding $j_0$ and the blue arc the embedding $j_1.$ Notice that this move simulaneously modifies the embedding $e$ as in \cref{adding dual picture 1} and takes an embedded connected sum of the blue arc with the "geometric dual" of the red portion in $V_1.$ This is a $0$-dimensional version of the move depicted in \cref{dual move picture}.}
    \label{grw move 1}
\end{figure}
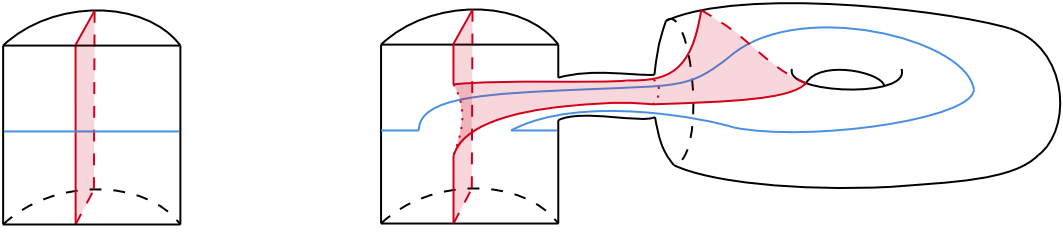

    We start by finding an embedded path in $\hb W_\infty$ from a point in $x(\partial (\partial_0 D^{n+1}_+\times D^n))$ to an embedded copy of $(V_1,W_{1,1})$ disjoint from $\vb W_\infty$ and the remaining of images under $F'$. This is possible because $n\geq 2$. We can thicken this path to define an embedded boundary connected sum \[x':((D^{n+1}_+\times D^n)\natural V_1,(\partial_0D^{n+1}_+\times D^n)\natural W_{1,1})\hookrightarrow (W_\infty,\hb W_\infty).\]There are two prefered disjoint embeddings 
    \begin{align*}
         j_0: (D^{n+1}_+\times \frac{1}{2} D^n,\partial_0D^{n+1}_+\times \frac{1}{2}D^n) & \hookrightarrow ((D^{n+1}_+\times D^n)\natural V_1,(\partial_0D^{n+1}_+\times D^n)\natural W_{1,1})\\
         j_1: \{0\}\times D^n & \hookrightarrow (\partial_0D^{n+1}_+\times D^n)\natural W_{1,1}
    \end{align*}which agree with the canonical inclusion near the (vertical) boundary: identify $((D^{n+1}_+\times D^n)\natural V_1,(\partial_0D^{n+1}_+\times D^n)\natural W_{1,1})\cong (S^n\times D^{n+1},S^n\times S^n\backslash (D^n_-,D^n_-)$ and take the submanifolds $(\frac{1}{2}D^n_+\times D^{n+1}_+,\frac{1}{2}D^n_+\times \partial_1D^{n+1}_+) $ and $D^n_+\times \{0\}$ respectively. We define the embedding $c'$ to agree with $c_\epsilon$ away from $x$ and to be the embedding $j_0$ on the remaining domain. Similarly, we define $g_\alpha'$ to agree with $g_\alpha$ away from $x$ and the map $j_1$ inside the domain of $x.$ It is clear to see that $c'$ and $g_\alpha'$ intersect transversely where $c$ and $g_\alpha$ do except the point $p.$ (See \cref{grw move 1} for a depiction of this move.) Observe that the $c'$ is the result of the construction before \cref{claim 5.6}, thus we can define a triple $(t+\epsilon, c',\L')$ in $X_0$. Let $F'_v$ be the simplicial map that agrees with $F'$ for all vertices different from $v$ and that $F'_v(v)=(t+\epsilon,c',\L').$ This map is homotopic to $F'$ by the same reason as before. By \cref{claim 5.6}, we still have that $K_\sigma(F'_v)=0$ and that $I_\sigma(F')\subseteq I_\sigma(F'_v).$ After repeating this process for every vertex of $\sigma$, we obtain a map $F''_\alpha$ and a map $g_\alpha''.$ However, the map $g_\alpha''$ is obtained by taking multiple embedded connected sums with sphere with standard (and in particular nullhomotopic) $\Theta$-structure. Thus the class $[l_{\hb W_\infty}\circ g_\alpha'']$ agrees with $[l_{\hb W_\infty}\circ g_\alpha]$ in $\pi_n(B^\partial).$ But, since $g_\alpha''$ is disjoint from the cores of $F''(\sigma),$ it lies in $I_\sigma(F'').$ By finite generation and applying this process to every simplex, we produce a map $F''$ homotopic to $F$ relative to $\partial I^k$ where $K_\sigma(F'')=0$ and $I_\sigma(F'')=\pi_n(B^\partial)$ for every simplex $\sigma.$ This finishes the proof of this statement.
\end{proof}

\begin{prop}\label{no discretize}
    The maps
    \(\de(W_\infty,\phi,\ell)_\bullet\to \deb(W_\infty,\phi,\ell)_\bullet\leftarrow \debd(W_\infty,\phi,\ell)_\bullet\) induce an equivalence on geometric realizations.
\end{prop}
\begin{proof}
    One proceeds exactly as \cite[Lemma 5.3]{GRWII} for the left equivalence and \cite[Lemma 5.7]{GRWII} for the right equivalence. 
\end{proof}

\begin{proof}[Proof of \cref{final contractability of middle left}.]
    Combine \cref{quasi fibration}, \cref{no discretize}, \cref{from no (v) to (v)}, \cref{from immersed to embedded cores} and \cref{contractibility easy}.
\end{proof}

\subsubsection{Contractability of right core complexes of middle dimension.}

In this section, we prove the analogous result for the right type, establishing \ref{step 3} from the introduction in this case.

\begin{teo}\label{contractability middle right}
    Let $P\in \smash{\CobbtL}$ and $K$ a $\Theta$-end. Let $\phi$ and $\ell:(0,\infty)\to \Thetaspace(T(D^n\times D^{n+1}_+),\Theta^*\gamma_{2n+1})$ be as in \cref{right resolution}. Then the map
    \(||\de_{\phi,\ell}(P,K|_\infty)_\bullet||\to \F(P,K|_\infty)\) is an equivalence, provided $(B,B^\partial)$ is $1$-connected.
\end{teo}

We proceed similarly as before by defining variants of $\de(W,\phi, l)$, deduce their contractability from \cite[Section 5]{GRWII} and compare them to $\de(W,\phi,\ell)$. Let $(s,W),\phi$ and $\ell$ be as in \cref{right resolution}.

    \begin{defn}\label{right resolution immersion}
    Define $\deb(W,\phi,\ell)_0$ to be the space of triples $(t,c,\L)$ where $t\in (0,\infty),$ 
    \[c:D^{n}\times D^{n}\hookrightarrow \hb W\]an immersion and $\L:[0,1]\to \Thetaspacev(D^{n}\times D^{n},(\theta^\partial)^*\gamma_{2n})$ satisfying the following properties:
    \begin{enumerate}[label=(\textit{\roman*})]
        \item $c$ is an embedding when restricted to $D^{n}\times \{0\}$ and for some $\delta>0,$ we have
        \(c(x,v)=\phi(\frac{x}{|x|},v+t\cdot e_1)+(1-|x|)\cdot e_0\) where $e_0$ is the first coordinate vector of $\RR\times \RR^{\infty}.$
        \item The image $C=c(D^{n}\times \{0\})$ is disjoint from $([0,s]\times \vb L)\cup (\{s\}\times Q)$ and $c^{-1}(\hb P)=\partial D^{n}\times D^{n}.$
        \item $\L(0)=c^*\ell_{\hb W}$ and $\L(1)=\ell_t|_{D^n\times \partial_0D^{n+1}_+}.$
        \item The map $\ell_W|_{W\backslash C}:W\backslash C\to B$ is strongly $n$-connected.
    \end{enumerate}Define $\deb(W,\phi,\ell)_p$ to be the subspace of $\deb(W,\phi,\ell)_0^{\times (p+1)}$ constisting of those tuples $(t_i;c_i;\L_i)_{i=0,\cdots,p}$ such that:
    \begin{enumerate}[label=(\textit{\roman*})]
        \item $t_0<t_1<\cdots<t_p.$
        \item The embeddings $c_i|_{D^{n}\times \{0\}}$ are pairwise disjoint.
        \item The map $\ell_{W}|_{W\backslash C}:W\backslash C\to B$ is strongly $n$-connected, for $C$ the union of the images of all $c_i|_{D^{n}\times \{0\}}.$
    \end{enumerate}This defines a semi-simplicial space $\deb(W,\phi,\ell)_\bullet.$ Denote by $\debd(W,\phi,\ell)_\bullet$ the semi-simplicial space where $\debd(W,\phi,\ell)_p$ is the set $\deb(W,\phi,\ell)_\bullet$ with the discrete topology. The identity induces a map of semi-simplicial spaces
    \[\debd(W,\phi,\ell)_\bullet\to \deb(W,\phi,\ell)_\bullet.\]
    
\end{defn}

\begin{defn}\label{right resolution boundary}
     Define $\dec(W,\phi,\ell)_0$ to be the space of triples $(t,c,\L)$ as in \cref{right resolution immersion} satisfying $(i)-(iii)$ along with 
    \begin{enumerate}[label=(\textit{\roman*'})]
        \setcounter{enumi}{3}
        \item The map $\ell_{\hb W}|_{\hb W\backslash \hb C}:\hb W\backslash C\to B^\partial $ is $n$-connected.
    \end{enumerate}Define $\dec(W,\phi,\ell)_p$ to be the subspace of $\dec(W,\phi,\ell)^{\times (p+1)}_0$ consisting of those tuples $(t_i;c_i;\L_i)_{i=0,\cdots,p}$ satisfying $(i)-(ii)$ and 
    \begin{enumerate}[label=(\textit{\roman*'})]
        \setcounter{enumi}{2}
        \item The map $\ell_{\hb W}|_{\hb W\backslash \hb C}:\hb W\backslash B^\partial $ is $n$-connected, where $C$ is the union of all images of $c_i.$
    \end{enumerate}This defines a semi-simplicial space $\dec(W,\phi,\ell)_\bullet$ with a map to $\debd(W,\phi,\ell)_\bullet.$
\end{defn}

In this case, we are able to extract contractability directly from the work of \cite{GRWII}, since our embedding is completely contained in the boundary. Similarly to last subsubsection, we denote by $\de(W_\infty,\phi,\ell)$ the strict colimit of the maps $K|_{[i,i+1]}\circ (-),$ and similarly for all variations of $\de.$

\begin{lemma}\label{input from grw}
    The space $||\dec(W_\infty,\phi,\ell)_\bullet||$ is contractible.
\end{lemma}
\begin{proof}
    Recall the semi-simplicial space $\overline{Y}^\delta(\hb W)_\bullet$ of \cite[Defn. 5.1]{GRWII} and notice that this is precisely $\dec(W,\phi,\ell)_\bullet$ for any $W.$ The claim follows by \cite[Prop. 5.4/5]{GRWII}.
\end{proof}

To compare the geometric realizations of the discrete semi-simplicial spaces defined above, we will employ a similar (but dual) strategy to \cref{from no (v) to (v)} once again using the infinite supply of $V_1$'s at hand.

\begin{prop}\label{from (v) boundary to complete}
    The map 
    \[\debd(W_\infty,\phi,\ell)_\bullet\to \dec(W_\infty,\phi,\ell)_\bullet\]induces an equivalence on geometric realizations, provided $(B,B^\partial)$ is $1$-connected.
\end{prop}
\begin{proof}
    Denote the source and target of this map by $X_\bullet$ and $Y_\bullet$. Let $k\geq 0$ and $F:(I^k,\partial I^k)\to (||Y_\bullet||,||X_\bullet||)$ be a map, which we can assume to be simplicial with respect to some triangulation of $I^k$, and $\sigma$ a simplex of $I^k$ for which $\ell_{W\backslash C_\sigma}$ is not (necessarily) strongly $n$-connected, where $C_\sigma$ is the union of all cores of $F(\sigma).$ We shall replace $F$ up to homotopy so that $F(\sigma)\in X_p.$ The inclusion $(W\backslash C_\sigma,\hb W\backslash \hb C_\sigma)\hookrightarrow (W,\hb W)$ is strongly $(n-1)$-connected. Since $n\geq 2$ and $\hb W\backslash \hb C_\sigma\to B^\partial$ is $n$-connected, we see $\pi_1(\hb W\backslash \hb C_\sigma)\cong \pi_1(B^\partial).$ Since the map of spaces $W\backslash C_\sigma\to W$ is an equivalence, it follows $\pi_1(W\backslash C_\sigma)\cong \pi_1(B)$ and thus $(W\backslash C_\sigma,\hb W\backslash \hb C_\sigma)$ is $1$-connected. By \cref{1 conn: old same as new}, we can detect strong $n$-connectivity on maps of relative homotopy groups. Set $K_\sigma(F)$ to be the kernel of the map \[\pi_{n}(W\backslash C_\sigma,\hb W\backslash \hb C_\sigma)\to \pi_n(B,B^\partial) \]and $I_\sigma(F)$ to be the image of the same map in degree $n+1.$ We see that $F(\sigma)$ lies in $X_p$ if and only if $K_\sigma(F)$ and $\pi_{n+1}(B,B^\partial)/I_\sigma(F)$ vanish, since we know $\ell_{\hb W\backslash \hb C}$ is $n$-connected. For $n\geq 2$, the maps above are of $\ZZ[\pi_1(B^\partial)]$-modules. We start by proving that $K_\sigma(F)$ is finitely generated as such. Note that the triad homotopy group (see \cite[Section 3.1]{BP}) $\pi_{n+1}(B,B^\partial,W\backslash C_\sigma)$ surjects into $K_\sigma(F).$ Moreover, by comparing the long exact sequences for the maps of pairs $W\backslash C_\sigma\to W\to B$ and the four-lemma, we have a surjection $\pi_{n+1}(W,\hb W,W\backslash C_\sigma)\to \pi_{n+1}(B,B^\partial,W\backslash C_\sigma).$ Since $W\backslash C_\sigma\to W$ is an equivalence (as it is equivalent to the attachment of a left handle), we have an isomorphism $\pi_{n+1}(W,\hb W,W\backslash C_\sigma)\cong \pi_{n}(\hb W,\hb W\backslash \hb C_\sigma).$ The latter group is a finitely generated $\ZZ[\pi_1(B^\partial)]$-module generated by the horizontal boundaries of belt (or meridian) spheres of $C_\sigma$ by the Hurewicz theorem, if $n\geq 3$, and normally finitely generated, if $n= 2$ (see \cite[166]{GRWII} or the proof of \cref{from no (v) to (v)}). In other words, the $\ZZ[\pi_1(B^\partial)]$-module $K_\sigma(F)$ is generated by the belt spheres $\{0\}\times \partial_1 D^{n+1}_+\subset D^n\times D^{n+1}_+$ composed with the embeddings of $F(\sigma).$ Similar to before, if $K_\sigma(F)=0,$ then we have a long exact sequence
    \[\cdots \to \pi_{n+1}(W\backslash C_\sigma,\hb W\backslash\hb C_\sigma)\to \pi_{n+1}(B,B^\partial)\to \pi_{n+1}(B,B^\partial,W\backslash C_\sigma)\to 0. \]As mentioned above, the $\ZZ[\pi_1(B^\partial)]$-module $\pi_{n+1}(B,B^\partial,W\backslash C_\sigma)$ is finitely generated and thus, $\pi_{n+1}(B,B^\partial)$ is generated by $I_\sigma(F)$ along with finitely many elements coming from lifts of the leftmost map.
    \\\\
    \textit{Killing $K_\sigma(F).$} We describe a modification on the embeddings $c_i$ and describe its main properties. This is analogous to \cite[166]{GRWII} and similar to the modification of $c$ to satisfy (D) in the proof of \cref{from immersed to embedded cores}. Let $v$ be a vertex of $\sigma$ and $F(v)=(t,c,\L)$ and let $j\geq 0$ be such that $\im(c)\subset W_j.$ Recall the embedding $\bar{e}:S^n\times \partial_0D^{n+1}_+\hookrightarrow W_{1,1}. $ Denote also by $\bar{e}$ the composition of this embedding with the inclusion of a fixed copy of $V_1$ in $K|_{[j,j+1]}$ seen now in $W_{j+1}.$ Choose an embedded path $\gamma$ in $\hb W_{j+1}$ from $c_\epsilon(x,0)$ for some point $x\in D^n_+$ and $\bar{e}(x,0)$ disjoint from $\hvb W_j$ and $C_\sigma$, and from $\im (c_\epsilon)$ and $\im(\bar{e})$ except at times $t=\{0,1\}.$ Using a thickening of this path, we can construct an embedded connected sum of the cores of $c_\epsilon$ and $\bar{e}$ and further thicken it to an embedding $c':D^n\times D^n \hookrightarrow \hb W_{j+1}$ which agrees with $c_\epsilon$ outside of a neighborhood of $x.$ As in \cref{from immersed to embedded cores} and \cref{from no (v) to (v)}, since the $\Theta$-structure $\bar{e}^*\ell_{\hb W_{j+1}}$ is standard, the structures $(c')^*\ell_{\hb W_{j+1}}$ and $c_\epsilon^*\ell_{\hb W_{j+1}}$ are homotopic. By choosing such homotopy and composing it with the path $\L_\epsilon$, we construct a path $\L'$ from $(c')^*\ell_{\hb W_{j+1}}$ to $\ell_{t+\epsilon}|_{D^n\times \partial_0D^{n+1}_+}.$ We remark that, once restricted to the horizontal boundary, this construction is exactly the one explained in \cite[166]{GRWII}. We conclude by \cite[Claim 5.6]{GRWII} that it does not affect the fact that $\hb W\backslash C_\sigma\to B^\partial$ is already $n$-connected, so we stay in $Y_\bullet.$ Denote $F_v:(I^k,\partial I^k)\to (||X_\bullet||,||Y_\bullet||)$ the map definedd on vertices by $F_v(w)=F(w)$ for $w\neq v$ and $F_v(v)=(t+\epsilon,c',\L')$. For $\epsilon$ small enough, we see that this map is well-defined and homotopic to $F$ relative to $\partial I^k.$ The following claim establishes properties of this construction.

    \begin{claim}\label{claim 5.6 but dual}
        The map $F_v$ is homotopic to $F$ relative to $\partial I^k$ and satisfies the following properties: 
     \begin{enumerate}[label=(\roman*)]
         \item for every simplex $\sigma=(v,\omega_1,\cdots,\omega_p)\in I^k,$ there exists a surjection $K_\sigma(F)\to K_\sigma(F_v),$ sending the class of the meridian of each simplex of $g$ to the meridian of each simplex of $F_v$ except $v$ and whose kernel contains the class of the meridian of $v.$ 
         \item for every simplex $\sigma\in I^k,$ we have $I_\sigma(F)\subseteq I_\sigma(F_v).$
     \end{enumerate}
    \end{claim}
    \begin{proof}[Proof of \cref{claim 5.6 but dual}.]
        This is essentially analogous to \cite[Claim 5.6]{GRWII}. We start with $(i).$ Let $X$ to be the union of $\im(c)=C$, the thickened path $\gamma$ and the $V_1$ used for the move above, in such a way that the image of modified embedding $c'$ lies in $X$ and $\hb X=\hb W\cap X. $ Let also $C_\omega$ be the image of the cores of $\omega_i.$ Consider the following maps
        \[K_\sigma(F')\leftarrow \ker(\pi_n(W\backslash (X\cup C_\omega),\hb W\backslash (\hb X\cup \hb C_\omega))\to \pi_n(B,B^\partial))\to K_\sigma(F).  \]We start by proving that the right map is an isomorphism. It suffices to prove that the map $\pi_n(W\backslash (X\cup C_\omega),\hb W\backslash (\hb X\cup \hb C_\omega))\to \pi_n(W\backslash (C\cup C_\omega),\hb W\backslash (\hb C\cup \hb C_\omega)) $ is an isomorphism. Observe that the following square is a pushout square of pairs
        \[\begin{tikzcd}
            (D^{2n}_+,\partial_0D^{2n}_+)\arrow[d]\arrow[r] & (W\backslash (X\cup C_\omega),\hb W\backslash (\hb X\cup \hb C_\omega))\arrow[d] \\
            (V_1,W_{1,1})\arrow[r] & (W\backslash (C\cup C_\omega),\hb W\backslash (\hb C\cup \hb C_\omega))
        \end{tikzcd}\]since the right map is obtained by taking a connected sum of pairs with $(V_1,W_{1,1})$. Given a pair of path connected spaces $(Y,Y')$, one finds an isomorphism \[\pi_i(Y\vee V_1,Y'\vee W_{1,1})\cong \pi_i(Y,Y')\oplus (\pi_i(V_1,W_{1,1})\otimes \coker(\ZZ[\pi_1Y']\to \ZZ[\pi_1Y])),\]for $i\leq n$ (for example, using the relative Hurewicz theorem or as a corollary of the Hilton-Milnor theorem). By taking $(Y,Y')$ to be the pair in the upper right corner in the square above, one sees that the right vertical map induces an isomorphism in relative homotopy groups in degrees at most $n.$ We proceed to prove that the left horizontal map above is surjective and its kernel contains the class of the meridian of $v$, which finishes $(i)$. Recall that $f$ is isotopic to the inclusion of $S^n\times \{*\}\subset S^n\times D^{n+1}$ for some $*\in \partial D^{n+1}.$ Thus, $V_1$ is obtained from a neighborhood of $f$ by attaching a left $(n+1)$-handle at the meridian of $f.$ We conclude that we have a (homotopy) pushout square of pairs 
        \[\begin{tikzcd}
            (\partial_1 D^{n+1}_+,\partial_{01}D^{n+1}_+)\arrow[d]\arrow[r] & (W\backslash (X\cup C_\omega),\hb W\backslash (\hb X\cup \hb C_\omega))\arrow[d] \\
            (D^{n+1}_+,\partial_0D^{n+1}_+)\arrow[r] & (W\backslash (C'\cup C_\omega),\hb W\backslash (\hb C'\cup \hb C_\omega))
        \end{tikzcd}\]where $C'$ denotes the image of $c'$ and the top map is the inclusion of the meridian sphere of $f$. By \ref{effect of left handles} in \cref{strong geo section}, the right map is strongly $(n-1)$-connected and moreover, the class of the meridian sphere of $f$ is in the kernel of the map on relative homotopy groups. We conclude that the map
        \[\pi_{n}(W\backslash (X\cup C_\omega),\hb W\backslash (\hb X\cup \hb C_\omega))\to \pi_n(W\backslash (C'\cup C_\omega),\hb W\backslash (\hb C'\cup \hb C_\omega))\]is surjective and the class of the meridian of $g'(v)$ (which is simply the meridian of $f$) is in the kernel. This implies $(i)$ by compositing to $\pi_n(B,B^\partial).$

        We finish by proving $(ii).$ Let $\psi:(D^{n+1},\partial D^{n+1})\to (W\backslash (C\cup C_\omega),\hb W\backslash (\hb C\cup \hb C_\omega))$ and suppose its image intersects the image of $c',$ otherwise, it is automatically in $I_\sigma(F')$. It must do so inside $W_{1,1}$ (in a finite set of points) as the path $\lambda$ can be isotoped away from it, since $n\geq 2,$ and the image of $c'$ is a collar of an embedding in $W_{1,1}.$ We can construct a new map $\psi'$ which agrees with $\psi$ outside this finite set of points given by the boundary connected sum with translates of the embedding $\bar{f}$ to remove intersection points exactly as in the proof of \cref{from no (v) to (v)} and explained in \cref{grw move 1}. The only difference is the fact that the chart $x$ is already in $V_1.$ Moreover, for this dual sphere $\bar{f}$, since the tangential structure on $V_1$ is assumed to be standard, we have $[l_W\circ \psi]=[l_W\circ \psi']\in \pi_{n+1}(B,B^\partial). $ Hence, the class induced by $\psi$ in $I_\sigma(F)$ also lies in $I_\sigma(F').$ This proves $(ii).$
    \end{proof}
    
    Apply the construction above and \cref{claim 5.6 but dual} for every interior vertex $v\in I^k.$ This produces a homotopic map $F'$ where the meridian of all vertices are trivial in $K_\sigma(F').$ Since this kernel is generated by such meridians, we have $K_\sigma(F')=0.$ By induction on the dimension of the simplices, we can obtain $F'$ such that $K_\sigma(F')=0$ for every simplex $\sigma$ of $I^k.$
    \\\\
    \textit{Adding elements to $I_\sigma(F').$} We apply exactly the same strategy as the analogous step in \cref{from no (v) to (v)} by choosing elements $\{\psi_i\}_i$ of $\pi_{n+1}(B,B^\partial)$ that generate this module along with $I_\sigma(F').$ To replace them by other representatives which are disjoint from $C_\sigma$, apply the same move as in \cref{from no (v) to (v)} with the small difference that the embedding $j_0$ is restricted to $D^{n+1}_+\times \{0\}$ and $j_1$ is extended to $\frac{1}{2}D^{n+1}_+\times D^n$ and the charts are modified in the analogous way. This move will replace both $C_\sigma$ and $\psi_i$ by connected sums with the embeddings $\bar{e}$ and $\bar{f}$ respectively, which have standard structures. Proceed exactly as in \cref{from no (v) to (v)} using \cref{claim 5.6 but dual} instead of \cref{claim 5.6}. This finishes the proof.
\end{proof}

\begin{prop}\label{no discretize in right}
    The maps
    \(\de(W_\infty,\phi,\ell)_\bullet\to \deb(W_\infty,\phi,\ell)_\bullet\leftarrow \debd(W_\infty,\phi,\ell)_\bullet\) induce an equivalence on geometric realizations.
\end{prop}
\begin{proof}
    One proceeds exactly as \cite[Lemma 5.7]{GRWII} for the right equivalence. We focus on the leftmost map for $p=0$. Observe that, by contracting the $t$-coordinate, $\de(W,\phi,\ell)_0$ is equivalent to (a point-set model for) the homotopy fiber of the map
    \[\Emb_{\vb}(D^{n}\times D^{n+1}_+,W)\to \Buncp(D^n\times D^{n+1}_+,\Theta^*\gamma)\]at the point $\ell_0$, where the source is the space of embeddings of pairs $c$ extending $\phi,$ taking an embedding $e$ to $e^*\ell_W.$ On the other hand, $\deb(W,\phi,\ell)_0$ is equivalent to the analogous space where $c$ is an embedding instead of only an immersion with embedded core (see \cite[Lemma 5.3]{GRWII}). This space is equivalent to (a point-set model for) the homotopy fiber of the analogous map
    \[\Emb_{\partial}(D^{n}\times \partial_0 D^{n+1}_+,\hb W)\to \Bun_\partial(D^n\times \partial_0 D^{n+1}_+,(\theta^\partial)^*\gamma_{2n})\]at $\ell_0.$ Moreover, the map above in $p=0$ is simply the induced map on the homotopy fiber induced by the restriction maps between the embedding and bundle spaces. We claim that both maps are equivalences: between embedding spaces one uses the parameterized isotopy extension theorem to deduce that it is a Serre fibration and its fiber identifies with the space of embeddings of $D^n\times D^{n+1}_+$ fixed in $D^n\times \partial_0D^{n+1}_+,$ which is contractible; and between bundle spaces one uses \cref{bundle map fibration} to deduce similarly that its fiber is contractible. Hence the map on homotopy fibers is an equivalence, which finishes the proof for $p=0.$ The higher cases are completely analogous.
\end{proof}

\begin{proof}[Proof of \cref{contractability middle right}.]
    Combine \cref{input from grw}, \cref{from (v) boundary to complete} and \cref{no discretize in right}.
\end{proof}

\subsubsection{Resolutions are degree-wise cuts.}\label{cuts section}

We establish now \ref{step 2} from the introduction, which involves specifying the data of $\phi$ and $\ell$ from the results above. Let $W:M\leadsto N$ be a morphism in $\smash{\CobbtL}$, $\sigma$ an attaching map and $\ell_t$ be a family of $\Theta$-structure on $\tr(\phi)$ indexed by $t\in (0,\infty)$ satisfying the hypothesis of \cref{main closure property}. For the remaining of this subsection, we fix $\phi\coloneqq \sigma|_{\partial D^n\times (0,\infty)\times \RR^{n}_+}$ when $\sigma$ is of right type and $\phi\coloneqq \sigma|_{\partial_1D^{n+1}_+\times (0,\infty)\times \RR^{n-1}}$ when $\sigma$ is of left type. By assumption, the image of $\phi$ lands away from the support of $W$, thus we see this map both in $N$ and in $M$ (as $N\backslash \supp(W)=M\backslash \supp(W)$). We define $\ell$ to be the restriction of the family $\ell_t$ to $D^n\times D^{n+1}_+$ for the right type and to $D^{n+1}_+\times D^n$ for the left type. Our first goal is to identify the space of $p$-simplices of $\de(W,\phi,\ell)$ in terms of the iterated surgeries $\chi^{p+1}(W,\sigma,\ell)$ (recall the definition from \cref{the multiple surgery}).

 Recall that we have the subspaces $\core_\sigma$ and $\cocore_\sigma$ of $\tr(\sigma),$ whose parameterization we fix. Given the multiple surgery $\chi^p(W,\sigma,\ell)$, one denotes by $\core_{\sigma,i}$ to be the union of $\core_{\sigma_i}\subset \tr(\sigma_i)$ with $[0,i-1]\times \im(\sigma_i)\subset \cup_{\!j< i}\tr(\sigma_{\!j})=\tr(\sigma^{i-1}).$ This is possible since $\sigma_i$ is disjoint from the traces of all $\tr(\sigma_{\!j})$ for $j<i.$ 

\begin{nota}
    Let $\L_i$ be the constant path in $\ell_{3i}.$ Given $(s,X)\in \F(\chi^p(N,\sigma,\ell),K|_i),$ we claim that the tuple
    \[(s+p+1,X\circ \tr(\sigma^{p+1});(3i;\core_{\sigma,i},\L_i)_{i=1,\cdots,p+1})\in \de_{\phi,\ell}(N,K|_i)_p.\]where $\tr(\sigma^{p+1})$ is the union of all $\tr(\sigma_i)$ for $i\leq p+1.$ Since $\tr(\sigma^p)$ is isotopy equivalent to the union of all cores, it follows that the complement of such cores is isotopy equivalent to $X.$ Thus, by assumption, it has strongly $n$-connected $\Theta$-structure. This establishes property $(iv)$ in the definition above. The remaining properties hold by construction. This defines a map
    \[(-)\circ \sigma^{p+1}:\F(\chi^{p+1}(N,\sigma,\ell),K|_i)\to \de_{\phi,\ell}(N,K|_i)_p.\]Similarly, we have the analogous map replacing $N$ by $M$ and using the cores of $\L_W(\tr(\sigma^{p+1}))$ coming from the cores of $\tr(\sigma^{p+1})$. We denote this map also by $(-)\circ \sigma^{p+1}. $\end{nota}

    \begin{lemma}\label{naturality of the gluing map}
        The following square
        \[\begin{tikzcd}
	{\F(\chi^{p+1}(N,\sigma,\ell),K|_i)} & {\de_{\phi,\ell}(N,K|_i)_p} \\
	{\F(\chi^{p+1}(M,\sigma,\ell),K|_i)} & {\de_{\phi,\ell}(M,K|_i)_p}
	\arrow["{(-)\circ \sigma^{p+1}}", from=1-1, to=1-2]
	\arrow["{(-)\circ \chi^{p+1}(W,\sigma,\ell)}", from=1-1, to=2-1]
	\arrow["{(-)\circ W}", from=1-2, to=2-2]
	\arrow["{(-)\circ \sigma^{p+1}}", from=2-1, to=2-2]
\end{tikzcd}\]is homotopy commutative.
    \end{lemma}
    \begin{proof}
        If we post-compose the horizontal maps with the projection maps $\de_{\phi,\ell}(-,Q)\to \F(-,Q),$ then we have the square
        \[\begin{tikzcd}[column sep=huge]
	{\F(\chi^p(N,\sigma,\ell),K|_i)} & {\F(N,K|_i)_p} \\
	{\F(\chi^p(M,\sigma,\ell),K|_i)} & {\F(M,K|_i)_p}
	\arrow["{(-)\circ \tr(\sigma^{p+1})}", from=1-1, to=1-2]
	\arrow["{(-)\circ \chi^{p+1}(W,\sigma,\ell)}", from=1-1, to=2-1]
	\arrow["{(-)\circ W}", from=1-2, to=2-2]
	\arrow["{(-)\circ \L_W(\tr(\sigma^{p+1}))}", from=2-1, to=2-2]
\end{tikzcd}\]which is homotopy commutative, since we have a diffeomorphism $W\circ \tr(\sigma^{p+1})\cong \L_{W}(\tr(\sigma^{p+1}))\circ \chi^{p+1}(W,\sigma,\ell),$ by \cref{path of rearrangement}. To verify that the original square commutes, it suffices to observe that the diffeomorphism can be chosen so it takes the cores of the bottom map to the ones defined using the top map. This is done exactly as \cite[157]{GRWII}.
    \end{proof}

    \begin{prop}\label{de are cuts}
        The map \[(-)\circ \sigma^{p+1}:\F(\chi^{p+1}(N,\sigma,\ell),K|_i)\to \de_{\phi,\ell}(N,K|_i)_{p}\]is an equivalence, for every $p\geq 0.$ The analogous result for $M$ instead of $N$ holds too.
    \end{prop}
    \begin{proof}
        We follow the strategy of \cite[156]{GRWII}. Let $E$ be the subspace of \[\F(N,K|_i)\times \text{Emb}_N(\tr(\sigma^{p+1}),[0,\infty)\times \RR^{\infty}_+)\times \Thetaspacev(\tr(\sigma^{p+1}),\Theta^*\gamma)^{[0,1]}\]where the second factor is the space of embeddings of pairs extending the inclusion of $N,$ consisting of those triples $((s,W'),e,\L)$ where the image of $e$ lies in $W'$, the complement has strongly $n$-connected $\Theta$-structure, $\L(0)=e^*\ell_{W'}$ and $\L(1)=l_{\tr(\sigma^{p+1})}.$ Taking such a triple to $((s,W'),(3i, e|_{\core_{\sigma,i}},\L|_{\core_{\sigma,i}})_{i=0,\cdots, p})$ provides a map $E\to \de_{\phi,\ell}(N,K|_i)_p$. This map is an equivalence since $\tr(\sigma^{p+1})$ is isotopy equivalent to the union of all its cores and (a collar on) $N.$ Consider the projection
        \[E \to \text{Emb}_N(\tr(\sigma^{p+1}),[0,\infty)\times \RR^{\infty}_+).\]We claim that this projection is a fibration by the parameterized isotopy extension theorem (see \cite[Appendix A]{steimle}). We sketch a proof here for completeness: if we restrict to the path components of $E$ of those triples where the first coordinate is a fixed cobordism $W':N\leadsto K|_i.$ Ignoring the third coordinate, we describe $E$ alternatively by the space of pairs $(e_0:\tr(\sigma^{p+1})\hookrightarrow W',e_1:W'\hookrightarrow [0,\infty)\times \RR^{\infty}_+)$ such that both $e_0$ and $e_1$ extend (in a collared way) the embedding of $N$ in $\{0\}\times \RR^{\infty}_+$ and additionally $e_1$ extends a (translated) embedding of $K|_i$ in $\{0\}\times \RR^{\infty}_+$ (to guarantee that it is disjoint from $K|_i$) balanced by the action of $\Diff_{\vb}(W')$. The projection above can now be seen as composition of $e_0$ and $e_1.$ This map is equivariant with respect to the action of $\Diff_{\{0\}\times \RR^\infty_+}([0,\infty)\times \RR^{\infty}_+).$ Moreover, by \cite[Appendix A]{steimle}, the target admits local cross sections and thus, by \cite[Thm. A]{Palais1960}, any equivariant map to it is a (locally trivial) Serre fibration.

        The base of this fibration is contractible by \cite[Thm. 2.7]{genauer}. The fiber over the standard inclusion of $\tr(\sigma^{p+1})$ is the equivalent to the subspace of $(t,W)\in\F(N,K|_i)$ such that $\tr(\sigma^{p+1})\subset W',$ which is homeomorphic to $\F(\chi^p(N,\sigma,\ell),K|_i).$ The map $\F(\chi^p(N,\sigma,\ell),K|_i)\to E$ is an equivalence and therefore, we have factored the map $(-)\circ \sigma^{p+1}$ into two equivalences. This finishes the proof.  
    \end{proof}

From \cref{naturality of the gluing map} and \cref{de are cuts}, we see that if $\chi^{p+1}(W,\sigma,\ell)\in \W$, then the induced map on resolutions is an abelian homology equivalence on $p$-simplices. This establishes \ref{step 2} of our strategy.

\subsubsection{Proof of the closure property.}\label{proof main section}

This section is dedicated to the proof of \cref{main closure property}. Let $W:M\leadsto N$ be a morphism in $\smash{\CobbtL}$, $\sigma$ an attaching map satisfying the hypothesis and $\ell_t$ be a family of $\Theta$-structure on $\tr(\sigma)$ indexed by $t\in (0,\infty).$ We assume that $\chi^p(W,\sigma,\ell)\in \W$ for every $p\geq 1.$ We wish to prove that $W\in \W.$ Recall the definition of $\phi$ from the previous subsubsection. Consider the square
\[\begin{tikzcd}
	{||\de_{\phi,\ell}(N,K|_\infty)_\bullet||} & {||\de_{\phi,\ell}(M,K|_\infty)_\bullet||} \\
	{\F(N,K|_\infty)} & {\F(M,K|_\infty)}
	\arrow[from=1-1, to=1-2]
	\arrow[from=1-1, to=2-1]
	\arrow[from=1-2, to=2-2]
	\arrow[from=2-1, to=2-2]
\end{tikzcd}\]obtained by \cref{pre functoriality}. By \cref{final contractability of middle left} and \cref{contractability middle right}, the vertical maps are weak equivalences, so it suffices to prove that the top map is an abelian homology equivalence. For $p\geq 0$, consider now the square
\[\begin{tikzcd}
	{\F(\chi^{p+1}(N,\sigma,\ell),K|_\infty )} & {\de_{\phi,\ell}(N,K|_\infty)_p} \\
	{\F(\chi^{p+1}(M,\sigma,\ell),K|_\infty)} & {\de_{\phi,\ell}(M,K|_\infty)_p}
	\arrow["{(-)\circ \sigma^{p+1}}", from=1-1, to=1-2]
	\arrow["{(-)\circ \chi^{p+1}(W,\sigma,\ell)}", from=1-1, to=2-1]
	\arrow["{(-)\circ W}", from=1-2, to=2-2]
	\arrow["{(-)\circ \sigma^{p+1}}", from=2-1, to=2-2]
\end{tikzcd}\]which commutes by \cref{naturality of the gluing map} and whose horizontal maps are equivalences by \cref{de are cuts}. The leftmost vertical map is an abelian homology equivalence by hypothesis, thus so is the right map. Thus the top map in the first square is a degree-wise abelian homology equivalence of semi-simplicial spaces. Thus, so is the map after taking geometric realizations (see \cite[158]{GRWII} for an argument). This finishes the proof of \cref{main closure property}.

\subsection{The higher dimensional case.}
The main goal of this subsection is to prove \cref{main closure property high}. We proceed similarly to the previous section. Recall the definition of $\mu_t$ from \cref{higher section}. We start by defining augmented semi-simplicial spaces over $\F(P,Q)$ for right, left and interior type.

\begin{defn}[Right resolution] \label{right resolution high}
    Let $n<k<2n$ be an integer, $P,Q\in \CobbtL$ and $(s,W)\in \F(P,Q).$ Let \[\chi:(\partial D^k\times \RR^{2n+1-k}_+\backslash D^{2n+1-k}_+,\partial D^k\times \partial\RR^{2n+1-k}_+\backslash \partial_0 D^{2n+1-k}_+)\hookrightarrow (P,\hb P) \]be an embedding and $\ell$ a $\Theta$-structure on $[0,2]\times \partial D^k\times \RR^{2n+1-k}_+$ such that
    \(\smash{\ell|_{\{0\}\times\partial D^k\times \RR^{2n+1-k}_+\backslash D^{2n+1-k}_+}=\chi^*\ell_P.}\) Let $\phi_t^\chi=\chi\circ \mu_t$ for every $t\in (2,\infty).$ Define $\de(W,\chi,\ell)_0$ to be the space of triples $(t,c,\L)$ where $t\in (2,\infty),$ 
    \(c:(D^{2n+1-k}_+\times D^k, \partial_0D^{2n+1-k}_+\times D^k)\hookrightarrow (W,\hb W)\) an embedding of pairs and $\L:[0,1]\to \Thetaspacev(T(D^{2n+1-k}_+\times D^k),\Theta^*\gamma_{2n+1})$ is a path of structures fixed as $\mu_t^*\ell$ on $\partial_1D^{2n+1-k}_+\times D^k$ satisfying the following properties:
    \begin{enumerate}[label=(\textit{\roman*})]
        \item \textit{(Collared expansion of meridian of $\chi$ by $t$ near $P$)} For some $\delta>0,$ we have
        \(c(x,v)=\phi_t^\chi(\frac{x}{|x|},v)+(1-|x|)\cdot e_0\) where $e_0$ is the first coordinate vector of $\RR\times \RR^{\infty}_+$ for all $x$ such that $1-|x|<\delta.$ 
        \item The image $(C,\hb C)=c(D^{2n+1-k}_+\times D^k,\partial_0D^{2n+1-k}_+\times D^k)$ is disjoint from $([0,s]\times \vb L)\cup (\{s\}\times Q)$ and $c^{-1}(P)=\partial_1D^{2n+1-k}_+\times D^k.$
        \item $\L(0)=c^*\ell_W$ and $\L(1)=\mu_t^*\ell.$
        
    \end{enumerate}This space is topologized in the same way as the resolutions from the previous subsection. Define $\de(W,\chi,\ell)_p$ to be the subspace of $\smash{\de(W,\chi,\ell)_0^{\times (p+1)}}$ consisting of those tuples $(t_i;c_i;\L_i)_{i=0,\cdots,p}$ such that:
    \begin{enumerate}[label=(\roman*)]
        \item $t_0<t_1<\cdots<t_p.$
        \item The embeddings $c_i$ are pairwise disjoint.
    \end{enumerate}This defines a semi-simplicial space $\de(W,\chi,\ell)_\bullet$ in the same way as in the last subsection.
    
\end{defn}

\begin{defn}[Left resolution]\label{left resolution high}
    Let $n+1<k<2n+1$ be an integer, $P,Q\in \CobbtL$ and $(s,W)\in \F(P,Q).$ Let \[\chi:(\partial_1 D^k_+\times \RR^{2n+1-k}\backslash D^{2n+1-k},\partial_{01} D^k_+\times \RR^{2n+1-k}\backslash D^{2n+1-k})\hookrightarrow (P,\hb P) \]be an embedding and $\ell$ a $\Theta$-structure on $[0,2]\times \partial_1D^{k}_+\times \RR^{2n+1-k}$ such that
    \(\smash{\ell|_{\{0\}\times\partial_1 D^k_+\times \RR^{2n+1-k}\backslash D^{2n+1-k}}=\chi^*\ell_P.}\) Let $\phi_t^\chi=\chi\circ \mu_t$ for every $t\in (2,\infty).$ Define $\de(W,\chi,\ell)_0$ to be the space of triples $(t,c,\L)$ where $t\in (2,\infty),$ 
    \(c:(D^{2n+1-k}\times D^k_+, D^{2n+1-k}\times \partial_0D^k_+)\hookrightarrow (W,\hb W)\) an embedding of pairs and $\L:[0,1]\to \Thetaspacev(T(D^{2n+1-k}\times D^k_+),\Theta^*\gamma_{2n+1})$ is a path of structures fixed as $\mu_t^*\ell$ on $\partial D^{2n+1-k}\times D^k_+$ satisfying the following properties:
    \begin{enumerate}[label=(\textit{\roman*})]
        \item \textit{(Collared expansion of meridian of $\chi$ by $t$ near $P$)} For some $\delta>0,$ we have
        \(c(x,v)=\phi_t^\chi(\frac{x}{|x|},v)+(1-|x|)\cdot e_0\) where $e_0$ is the first coordinate vector of $\RR\times \RR^{\infty}_+$ for all $x$ such that $1-|x|<\delta.$ 
        \item The image $(C,\hb C)=c(D^{2n+1-k}\times D^k_+,D^{2n+1-k}\times \partial_0D^k_+ )$ is disjoint from $([0,s]\times \vb L)\cup (\{s\}\times Q)$ and $c^{-1}(P)=\partial D^{2n+1-k}\times D^k.$
        \item $\L(0)=c^*\ell_W$ and $\L(1)=\mu_t^*\ell.$
    \end{enumerate}This space is topologized in the same way as the resolutions from the previous subsection. Define $\de(W,\chi,\ell)_p$ to be the subspace of $\de(W,\chi,\ell)_0^{\times (p+1)}$ consisting of those tuples $(t_i;c_i;\L_i)_{i=0,\cdots,p}$ such that:
    \begin{enumerate}[label=(\textit{\roman*})]
        \item $t_0<t_1<\cdots<t_p.$
        \item The embeddings $c_i$ are pairwise disjoint.
    \end{enumerate}This defines a semi-simplicial space $\de(W,\chi,\ell)_\bullet$.
    
\end{defn}

\begin{defn}[Interior resolution]\label{interior resolution}
    Let $n+1<k<2n+1$ be an integer, $P,Q\in \CobbtL$ and $(s,W)\in \F(P,Q).$ Let \[\chi:\partial D^k\times \RR^{2n+1-k}\backslash D^{2n+1-k}\hookrightarrow P\backslash\hb P \]be an embedding and $\ell$ is a $\Theta$-structure on $[0,2]\times \partial D^k\times \RR^{2n+1-k}$ such that
    \(\smash{\ell|_{\{0\}\times\partial D^k\times \RR\backslash D^{2n+1-k}}=\chi^*\ell_P.}\) Let $\phi_t^\chi=\chi\circ \mu_t$ for every $t\in (2,\infty).$ Define $\de(W,\chi,\ell)_0$ to be the space of triples $(t,c,\L)$ where $t\in (2,\infty),$ 
    \(c:D^{2n+1-k}\times D^{k}\hookrightarrow W\backslash\hb W\) an embedding of pairs and $\L:[0,1]\to \Thetaspacev(T(D^{2n+1-k}\times D^{k}),\Theta^*\gamma_{2n+1})$ is a path of structures fixed as $\mu_t^*\ell$ on $\partial D^{2n+1-k}\times D^k$, where $\Bun_{\partial}(T(D^{2n+1-k}\times D^{k}),\theta^*\gamma_{2n+1})$ satisfying the following properties:
    \begin{enumerate}[label=(\textit{\roman*})]
        \item \textit{(Collared translation of the meridian of $\chi$ by $t$ near $P$)} For some $\delta>0,$ we have
        \(c(x,v)=\phi_t^\chi(\frac{x}{|x|},v)+(1-|x|)\cdot e_0\) where $e_0$ is the first coordinate vector of $\RR\times \RR^{\infty}_+$ for all $x$ such that $1-|x|<\delta.$ 
        \item The image $C=c(D^{2n+1-k}\times D^k)$ is disjoint from $([0,s]\times \vb L)\cup (\{s\}\times Q)$ and $c^{-1}(P)=\partial D^{2n+1-k}\times D^{k}.$
        \item $\L(0)=c^*\ell_W$ and $\L(1)=\ell_t.$
    \end{enumerate}This space is topologized in the same way as the resolutions from the previous subsection. Define $\de(W,\chi,\ell)_p$ to be the subspace of $\de(W,\chi,\ell)_0^{\times (p+1)}$ consisting of those tuples $(t_i;c_i;\L_i)_{i=0,\cdots,p}$ such that:
    \begin{enumerate}[label=(\textit{\roman*})]
        \item $t_0<t_1<\cdots<t_p.$
        \item The embeddings $c_i$ are pairwise disjoint.
    \end{enumerate}This defines a semi-simplicial space $\de(W,\chi,\ell)_\bullet$.
    
\end{defn}

\begin{defn}\label{defn of resolutions}
    Given $P,Q\in \CobbtL$ and $\chi$ and $\ell$ as in one of the above definitions. Define $\de_{\chi,\ell}(P,Q)_p$ to be the space of pairs $(W,x)$ where $W\in \F(P,Q)$ and $x\in \de(W,\chi,\ell)_p.$ This space maps to $\F(P,Q)$ by forgetting $x$ for every $p\geq 0.$
\end{defn}

\subsubsection{Functoriality of the resolutions.}
We continue by remarking "functoriality" properties of these resolutions with respect to pre- and post-composition. This will play the same role as the analogous step for the middle dimensional case.

\begin{lemma}\label{post functoriality high}
    Let $W':Q\leadsto S$ be a morphism in $\CobbtL$ which is strongly $(n-1)$-connected relative to $Q$, then post-composition defines a map of augmented semi-simplicial spaces
    \[W'\circ (-):\de_{\chi,\ell}(P,Q)_\bullet\to \de_{\chi,\ell}(P,S)_\bullet.\] 
\end{lemma}
\begin{proof}
    By the same argument as in \cref{subfunctor}, post-composition defines a map $\F(P,Q)\to \F(P,S).$ Since the embeddings in $\de_{\chi,\ell}(P,Q)$ are disjoint from $Q$, post composition with inclusion in $W'$ induces data in $\de_{\chi,\ell}(P,S).$ The verification of the conditions above follows by unravelling definitions.
\end{proof}

\begin{cons}\label{reparameterization drama 2}
    We follow exactly the strategy of \cref{reparameterization drama}. Fix a diffeomorphism $\epsilon_s:(D^k\times D^{2n+1-k}_+,D^k\times \partial D^{2n+1-k}_+)\cong ((s\cdot e_1+D^k)\times D^{2n+1-k}_+,(s\cdot e_1+D^k)\times \partial_0 D^{2n+1-k}_+)\cup [0,s]\times (\partial D^k\times D^{2n+1-k}_+,\partial D^k\times \partial_0 D^{2n+1-k}_+)\eqqcolon B_s\cup A_s=D_s$. Given $(t,c,\L)\in \de(W,\chi,\ell)$ for $\chi$ of left type and $(s,W'):P'\leadsto P$ such that $\supp(W)\cap \chi=\emptyset,$ define $c_s:D_s\hookrightarrow W'\cup W$ by taking $B_s$ to $c$ and $A_s$ to $[0,s]\times c|_{\partial D^k\times D^{2n+1-k}_+}.$ The triple $\smash{(t,\epsilon^{-1}_s(c_s),\epsilon^{-1}_s(\L_s))}$, where $\L_s$ is given by extending the path of $\Theta$-structure $\L$ along the collar $A_t,$ satisfies properties $(i)-(iii).$ When $\chi$ is of right type, one proceeds in a similar way by adding now a collar $[0,s]\times (\partial_1D^k_+\times D^{2n+1-k},\partial_{01}D^k_+\times D^{2n+1-k})$ to $(D^k_+\times D^{2n+1-k},\partial_{0}D^k_+\times D^{2n+1-k}).$ For interior type, one proceed analogously. (For more detail, see \cite[Defn. 4.13]{GRWII}.)
\end{cons}

\begin{lemma}\label{pre functoriality high}
    Let $(s,W'):P'\leadsto P$ be a morphism in $\CobbtLs$ such that $\supp(W)\cap \im(\chi)=\emptyset$, then pre-composition defined above 
    \[(-)\circ W':\de_{\chi,\ell}(P,Q)_\bullet\to \de_{\chi,\ell}(P',Q)_\bullet\]is a map of augmented semi-simplicial spaces. Moreover, for $W'':Q\leadsto S,$ then the pre and post-composition maps commute strictly in the natural way. 
\end{lemma}
\begin{proof}
    This follows exactly as \cref{post functoriality high}.
\end{proof}

\subsubsection{Contractability of core complexes for higher handles.}
We establish the analog of \cref{final contractability of middle left} and \cref{contractability middle right} for higher handles. Notice that the statement has no mention of a $\Theta$-end, since the augmented semi-simplicial spaces are already resolutions before statibilizing with the $\Theta$-end. The analogous phenomenon was observed in \cite[Theorem 6.6]{GRWII}. 

\begin{prop}\label{contractability of higher}
        Let $P,Q\in \CobbtL,$ $\chi$ and $\ell$ be as in \cref{defn of resolutions}. Then \(||\de_{\chi,\ell}(P,Q)_\bullet||\to \F(P,Q)\) is an equivalence, provided $(B,B^\partial)$ is $0$-connected.
\end{prop}
\begin{proof}
    By the same argument as in \cref{quasi fibration}, it suffices to prove that $\de(W,\chi,\ell)_\bullet$ has contractible realisation for every $W\in \F(P,Q).$ We can replace $\de(W,\chi,\ell)_\bullet$ by $\deb(W,\chi,\ell)$ of $0$-simplices $(t,c,\L)$ where $c$ is only assumed to be an immersion whose restriction to its core is an embedding and higher simplices $(t_i,c_i,\L_i)_i$ if the cores of $c_i$ are pairwise disjoint. This is analogous to \cref{left resolution immersion}. By the same argument of \cref{no discretize}, it follows that the map $||\de(W,\chi,\ell)_\bullet||\to ||\deb(W,\chi,\ell)_\bullet||$ is an equivalence. One verifies that $\deb(W,\chi,\ell)$ produces a topological flag complex (augmented over a point) as in \cref{contractability of surgery data}. (See the proof of \cite[Thm. 6.6]{GRWII} for a similar argument.) Thus, it suffices to verify the conditions $(i)-(iii)$ of \cite[Thm. 6.2]{GRWStableModuli} for this augmented flag complex. Condition $(i)$ is vacuous because $X_{-1}$ is a point. Condition $(ii)$ is equivalent to $\de(W,\chi,\ell)_0$ being non-empty. We separate by cases and start with the case of $\chi$ being of right type. Consider the following diagram
    \[\begin{tikzcd}
	{(\partial_1 D^{2n+1-k}_+\times D^{k},\partial_{01} D^{2n+1-k}_+\times D^{k})} & {(P,\hb P)} & {(W, \hb W)} \\
	{(D^{2n+1-k}_+\times D^{k},\partial_{0} D^{2n+1-k}_+\times D^{k})} && {(B,B^\partial)}
	\arrow["\phi_3^\chi", from=1-1, to=1-2]
	\arrow[hook, from=1-1, to=2-1]
	\arrow[hook, from=1-2, to=1-3]
	\arrow["{\ell_W}", from=1-3, to=2-3]
	\arrow["\hat{c}", dashed, from=2-1, to=1-3]
	\arrow["{\mu_3^*\ell}", from=2-1, to=2-3]
\end{tikzcd}\]which admits a lift making the top triangle commute strictly and the lower one up to homotopy, since $\ell_W$ is strongly $n$-connected and $k> n$ by \cref{strong implies lifting}, since $(B,B^\partial)$ and thus $(W,\hb W)$ are $0$-connected. By Smale-Hirsch theory as in \cref{contractibility easy}, we may assume that $\hat{c}$ is an immersion of pairs, as $\mu_3^*\ell$ is covered by a bundle map. Moreover, since $2(2n+1-k)<2n+1$, this immersion can be assumed to be an embedding once restricted to $D^{2n+1-k}_+\times \{0\}.$ By scaling the tubular neighborhood, we can homotope this immersion to be an embedding $c$. Choose $\L$ to be the path of $\Theta$-structures given by the homotopy of the bottom triangle above, and thus $(3,c,\L)\in \de(W,\chi,\ell)_0.$ The right and interior types follow exactly the same strategy, where the appropriate lifting problem will have a solution by the connectivity of $\ell_W$ and the restriction on $k.$ We finish with condition $(iii).$ Given a $p$-simplex $v=\{(t_i,c_i,\L_i)\}_{i=0,\cdots, p}$ we have to find $(t,c,\L)$ such that $\{(t,c,\ell)\}\cup v$ is a $(p+1)$-simplex. Take $t=t_0-\epsilon$ for all $i,$ consider $c$ the embedding $c_0$ composed with an expansion on its vertical boundary (in $P$) so it lies in $\phi_t^\chi.$ This is done in the same way as $c_\epsilon$ in the proof of \cref{from immersed to embedded cores}. We can isotope $c$ relative to its vertical boundary such that it is disjoint to all $c_i$ by transversality \cite[Lemma 9.2]{genauer}, since the handle dimension of the cores is strictly smaller than half of the total dimension. Taking $\L$ to be $\L_0$ composed with the expansion in $P.$ We see that $\{(t,c,\ell)\}\cup v$ is a $(p+1)$-simplex.
\end{proof}

\subsubsection{Resolutions are degree-wise cuts.}

As in the middle dimensional case, we proceed now by specifying our choice of $\chi$ and $\ell$. Let $W:M\leadsto N$ be a morphism in $\smash{\CobbtLs}$, $\sigma$ an attaching map of index $k$ and $\ell$ be a $\Theta$-structure satisfying the hypothesis of \cref{main closure property high}. For the remaining of this subsection, we set $\chi\coloneqq \sigma|_{\partial D^k\times \RR^{2n+1-k}_+\backslash D^{2n+1-k}_+}$ if $\sigma$ is of right type, $\chi\coloneqq \sigma|_{\partial_1 D^k_+\times \RR^{2n+1-k}\backslash D^{2n+1-k}}$ for the left type and $\chi\coloneqq \sigma|_{\partial D^k\times \RR^{2n+1-k}\backslash D^{2n+1-k}}$ for the interior type. Once again, we see $\chi$ both in $N$ and in $M$.

Recall the definition of $\phi^\sigma$ from \cref{definition of phi}. We proceed exactly as in \cref{cuts section} to define a morphism $(p+1,\tr(\phi^{\sigma,p+1})):P\leadsto \chi^{p+1}(P,\phi^\sigma,\ell)$ given by the union of the traces $\tr(\phi^\sigma_{3i})$ for $1\leq i\leq p+1,$ whose $\Theta$-structure is induced by $\ell.$ Inside this cobordism, we have embeddings $\core_{\phi^\sigma,i}$ for $i=1,\cdots, p+1.$ We define a map 
\[(-)\circ \phi^{\sigma,p+1}:\F(\chi^{p+1}(N,\sigma,\ell),K|_i)\to \de_{\chi,\ell}(N,K|_i)_p\]by taking $(s,X)$ to $(s+p+1, X\circ \tr(\phi^{\sigma, p+1});(3i;\core_{\phi^\sigma,i}; \L_i)_{i=1,\cdots, p+1})$ where $\L_i$ is the constant path at $\mu_{3i}^*\ell.$ Similarly, we have the analogous map replacing $N$ by $M$ and using the cores of $\L_W(\tr(\phi^{\sigma, p+1}))$ coming from the cores of $\tr(\phi^{\sigma, p+1})$. We denote this map also by $(-)\circ \phi^{\sigma,p+1}. $

    \begin{lemma}\label{naturality of the gluing map high}
        The following square
        \[\begin{tikzcd}
	{\F(\chi^{p+1}(N,\sigma,\ell),K|_i)} & {\de_{\chi,\ell}(N,K|_i)_p} \\
	{\F(\chi^{p+1}(M,\sigma,\ell),K|_i)} & {\de_{\chi,\ell}(M,K|_i)_p}
	\arrow["{(-)\circ \phi^{\sigma,{p+1}}}", from=1-1, to=1-2]
	\arrow["{(-)\circ \chi^{p+1}(W,\sigma,\ell)}", from=1-1, to=2-1]
	\arrow["{(-)\circ W}", from=1-2, to=2-2]
	\arrow["{(-)\circ \phi^{\sigma,{p+1}}}", from=2-1, to=2-2]
\end{tikzcd}\]is homotopy commutative.
    \end{lemma}
    \begin{proof}
        This follows verbatim with the strategy of \cref{naturality of the gluing map}.
    \end{proof}

    \begin{prop}\label{de are cuts high}
        The map \[(-)\circ \phi^{\sigma,p+1}:\F(\chi^{p+1}(N,\sigma,\ell),K|_i)\to \de_{\chi,\ell}(N,K|_i)_p\]is an equivalence, for every $p\geq 0.$ The analogous claim for $M$ also holds.
    \end{prop}
    \begin{proof}
        This follow verbatim with the strategy of \cref{de are cuts}.
    \end{proof}

\subsubsection{Proof of the closure property.} We finish by proving \cref{main closure property high}. Let $W:M\leadsto N$ be a morphism in $\smash{\CobbtLs}$, $\sigma$ an attaching map of index $k$ and $\ell$ be a $\Theta$-structure satisfying the hypothesis of \cref{main closure property high}. Assume $\chi^{p}(W,\phi^\sigma,\ell)\in \W$ for every $p\geq 1.$ Recall the definition of $\chi$ from the previous subsubsection. Consider the square
\[\begin{tikzcd}
	{||\de_{\chi,\ell}(N,K|_\infty)_\bullet||} & {||\de_{\chi,\ell}(M,K|_\infty)_\bullet||} \\
	{\F(N,K|_\infty)} & {\F(M,K|_\infty)}
	\arrow[from=1-1, to=1-2]
	\arrow[from=1-1, to=2-1]
	\arrow[from=1-2, to=2-2]
	\arrow[from=2-1, to=2-2]
\end{tikzcd}\]obtained by \cref{pre functoriality}. By \cref{contractability of higher}, the vertical maps are equivalences. We proceed exactly as in \cref{proof main section} to deduce that the upper map is an abelian homology equivalence by using \cref{de are cuts high} instead of \cref{de are cuts}. This finishes the proof. The proofs for the left and interior cases are analogous.

\section{Stable homology and group completion.}\label{group comp section}


In this section, we use the main results of the last sections, namely \cref{surgery on morphisms,stable stability}, to prove a generalization of \cref{final with null in intro}. This will be crucial to deduce \cref{main no tang} and more generally, a version with tangential structures in the next section. We start by stating the main result of this section. To do so, consider the following definitions. Throughout this section, we fix again the data in \cref{the assumptions for stable stability}. Although our main result of this section will not feature the $L$, its proof requires us to choose a specific $L$ and thus establish preliminary results where a general $L$ is present. We specialize to our desired $L$ by the end of this section. In \cref{theta end in l}, we introduced the notion of a $\Theta$-end in $\CobbtL$, but we now consider a variation which does not use the $L.$

\begin{defn}\label{theta end without l}
    A \textit{$\Theta$-end} $K$ in $\Cobbt$ is a sequence of composable morphisms
    \[\{K|_{[i,i+1]}:K|_i\leadsto K|_{i+1}\}_{i\geq 0}\]in $\Cobbt$ such that:
    \begin{enumerate}[label=(\roman*)]
        \item for every $i\geq 0$, the inclusions $(K|_i,\partial K|_i)\hookrightarrow (K|_{[i,i+1]},\hb K|_{[i,i+1]})$ and $(K|_{i+1},\partial K|_{i+1})\hookrightarrow (K|_{[i,i+1]},\hb K|_{[i,i+1]})$ are strongly $(n-1)$-connected,
        \item for every $i\geq 0,$ there exists an embedding $\omega:(V_1,W_{1,1})\hookrightarrow (K|_{[i,i+1]},\hb K|_{[i,i+1]})$ such that $\omega^*\ell_{K_{[i,i+1]}}$ is a standard $\Theta$-structure in the sense of \cref{standard structure on V1}. 
    \end{enumerate}
\end{defn}

For a pair of objects $P,Q\in \Cobbt$, consider the map $\Cobbt(P,Q)\to \Omega_{[P,Q]}\B\Cobbt$ that takes a morphism to the path in $\B\Cobbt$ represented by this morphism. For an object $P\in \Cobbt$, let $\Nt(P)$ denote $\Cobbt(\emptyset, P)$ and let $\Ntn(P)$ be the subspace of $\Nt(P)$ consisting of those morphisms $(W,\ell_W)$ such that $\ell_W:(W,\hb W)\to (B,B^\partial)$ is strongly $n$-connected. Given a $\Theta$-end $K$ in $\Cobbt$, denote by $\Ntn(K|_\infty)$ the homotopy colimit over $i$ of the maps $K|_{[i,i+1]}\circ (-):\Ntn(K|_i)\to \Ntn(K|_{i+1})$ (which are well defined by the first assumption in \cref{theta end without l}, see \cref{subfunctor}.). 

We have a map \(\Ntn(K|_\infty)\to \Omega_{[\emptyset,K|_\infty]}\B\Cobbt\) where the target denotes the homotopy colimit of the maps $ \Omega_{[\emptyset,K|_i]}\B\Cobbt\to  \Omega_{[\emptyset,K|_{i+1}]}\B\Cobbt$ given by concatenating loops with the path represented by the morphism $K|_{[i,i+1]}.$ Recall that a map of spaces $f:X\to Y$ is \textit{acyclic} if for all local systems $\L$ on $Y$, the induced map $\H_*(X;f^*\L)\to \H_*(Y,\L)$ is an isomorphism in all degrees. We are now ready to state the main result of this section. Recall from \cref{the assumptions for stable stability} that we assume $2n+1\geq 7$ and $B^\partial$ is path-connected.

\begin{teo}\label{final with null}
    Assume $(B,B^\partial)$ is $1$-connected. Let $K$ be a $\Theta$-end in $\Cobbt$ such that $\partial K|_0\neq \emptyset$ and $\Ntn(K|_0)\neq \emptyset.$ Then the map 
    \[\Ntn(K|_\infty)\to \Omega_{[\emptyset,K|_\infty]}\B\Cobbt\]is acyclic.
\end{teo}

The majority of the remainder of this section is devoted to proving this statement. We will proceed similarly to \cite[Section 7]{GRWII}, so we attempt to make our exposition as close to this reference as possible. To prove this statement, we will need two preliminary versions of \cref{final with null}, which are stated and proved in the next two subsections. Finally, the last subsection is devoted to the proof of \cref{stable class thm}.

\subsection{The group completion argument.}

The main goal of this subsection is to prove the first preliminary version of the statement above using \cref{stable stability}. Recall that we have fixed an $L$ as in \cref{the assumptions for stable stability}. Given two objects $P,Q\in \smash{\CobbtL}$, denote by $\smash{\Fo(P,Q)}$ the subspace of $\smash{\CobbtLs(P,Q)}$ (see \cref{cobordism categories definition L cut}) of those $(W,\ell_W)$ such that $\ell_W:(W,\hb W)\to (B,B^\partial)$ is strongly $n$-connected. Given a $\Theta$-end $K$ in $\smash{\CobbtL}$ (see \cref{theta end in l}), we denote by $\smash{\Fo(P,K|_\infty)}$ the homotopy colimit of the post-composition maps by $K|_{[i,i+1]}$, similar to the definition above. Before stating the main result of this subsection, we introduce the following notion. For $k\geq 2$, we say that the pair $(B,B^\partial)$ is $(E_k)$ if there exists a strongly $k$-connected map $(X,X')\to (B,B^\partial)$ from a CW complex pair $(X,X')$ such that $X'$ has finite $k$-skeleton and the pair $(X,X')$ has finite $(k+1)$-skeleton (i.e. the $(k+1)$-skeleton of the relative CW filtration of $X$ relative to $X'$ is finite). We now state the first preliminary version of \cref{final with null}.

\begin{prop}\label{highly connected theta ends}
    Assume $(B,B^\partial)$ is $1$-connected and $(E_n)$ . Let $K$ be a $\Theta$-end in $\smash{\CobbtL}$ such that $\ell|_{K|_0}$ is strongly $(n-1)$-connected. If $P\in \smash{\CobbtL},$ then the map
    \[\Fo(P,K|_\infty)\to \Omega_{[P,K|_\infty]}\B\CobbtLs\]is acyclic.
\end{prop}

\begin{rmk}\label{theta staisfies En}
    Observe that if $\Theta$ satisfies the assumption of \cref{final with null}, there exists a nullbordism $W$ of $K|_0$ with a strongly $n$-connected map to $(B,B^\partial)$. Since $W$ and $\hb W$ are compact manifolds, they admit finite CW structures. Thus $(B,B^\partial)$ satisfies the condition $(E_n).$
\end{rmk}

\subsubsection{An easier case.} We start the proof of \cref{highly connected theta ends} by establishing an intermediate version, where the hypothesis on $K$ and $\Theta$ are stronger. 

\begin{prop}\label{stably highly connected theta ends}
    Assume $(B,B^\partial)$ is $1$-connected and $(E_n)$. Let $K$ be a $\Theta$-end in $\CobbtL$ such that $\ell|_{K|_0}$ is strongly $(n-1)$-connected and $\ell|_{K|_{[i,+\infty)}}$ is strongly $n$-connected for every $i$. If $P\in \smash{\CobbtL},$ then the map
    \[\Fo(P,K|_\infty)\to \Omega_{[P,K|_\infty]}\B\CobbtLs\]is acyclic.
\end{prop}

The proof of this result will use \cref{stable stability} in conjuction with the generalized group completion theorem for topological categories, as stated in \cite[Thm. A.1]{GRWII}. To use these results, we need to prove the following lemmas. The overall goal of these lemmas is to relate the morphism spaces of the category $\CobbtLs$ and the subfunctor $\F(-,K|_\infty)$, which is the object of study of \cref{stable stability}.

\begin{lemma}\label{take out fo}
    In the context of \cref{stably highly connected theta ends}, the inclusion \(\Fo(P,K|_i)\subseteq \F(P,K|_i) \) is an equality, for every $i\geq 0$.
\end{lemma}
\begin{proof}
    It suffices to prove that for any morphism $W:P\leadsto K|_i$ such that $\ell|_W$ is strongly $n$-connected, then $(K|_i,\hb K|_i)\to (W,\hb W)$ is strongly $(n-1)$-connected. Start by observing that the map $(K|_i,\hb K|_i)\to (K|_{[i,+\infty)},\hb K|_{[i,+\infty)}) $ is strongly $(n-1)$-connected. To do this, we show that $(K|_i,\hb K|_i)\to (K|_{[i,i+j]},\hb K|_{[i,i+j]})$ is strongly $(n-1)$-connected for all $j\geq 1.$ For $j= 1,$ this follows from the definition of a $\Theta$-end. By induction, it suffices to prove that $(K|_{[i,i+j]} ,\hb K|_{[i,i+j]})\to (K|_{[i,i+j+1])},\hb K|_{[i,i+j+1])}) $ is strongly $(n-1)$-connected. Observe the following homotopy pushout square
    \[\begin{tikzcd}
	(K|_{i+j} ,\hb K|_{i+j}) & (K|_{[i,i+j]},\hb K|_{[i,i+j]}) \\
	   (K|_{[i+j,i+j+1]}, \hb K|_{[i+j,i+j+1]}) & (K|_{[i,i+j+1]},\hb K|_{[i,i+j+1]})
	\arrow[from=1-1, to=1-2]
	\arrow[from=1-1, to=2-1]
	\arrow[from=1-2, to=2-2]
	\arrow[from=2-1, to=2-2]
\end{tikzcd}\]where the left vertical map is strongly $(n-1)$-connected by definition. Applying \cref{strong conn and pushouts}, we have that the right vertical map is strongly $(n-1)$-connected, as claimed. Now the following square
    \[\begin{tikzcd}
	(K|_{i},\hb K|_{i}) & (W,\hb W)\\
	   (K|_{[i,\infty)},\hb K|_{[i,\infty)}) & (B,B^\partial)
	\arrow[from=1-1, to=1-2]
	\arrow[from=1-1, to=2-1]
	\arrow[from=1-2, to=2-2]
	\arrow[from=2-1, to=2-2]
\end{tikzcd}\]implies that the composite $(K|_i,\hb K|_i)\to (B,B^\partial)$ is strongly $(n-1)$-connected. By \cref{closure properties strongly}, we deduce that $(K|_i,\hb K|_i)\to (W,\hb W)$ is strongly $(n-1)$-connected, since $n\geq 3$. This finishes the proof.
\end{proof}

Before we prove \cref{stably highly connected theta ends}, we need the following lemma about the condition $(E_n)$. For a map of pairs $f:(X,X')\to (Y,Y')$, recall the definition of the triad homotopy group $\pi_*(Y,Y',X)$ from \cref{triad homotopy}.

\begin{lemma}\label{some fn stuff}
    Let $k\geq 2$. Let $(B,B')$ be a $1$-connected pair satisfying $(E_k)$ such that $B'$ is path-connected. Then the following claims hold:
    \begin{enumerate}[label=(\textit{\roman*})]
        \item $B'$ satisfies Wall's finiteness condition $(F_k)$ (in the sense of \cite[57]{Wallfin})
        \item Let $(W,W')$ be a \text{CW} pair such that both $W$ and $W'$ have finite $k$-skeleton relative to $\emptyset$. Let $f:(W,W')\to (B,B')$ be a strongly $(k-1)$-connected map such that $f|_{W'}:W'\to B'$ is $k$-connected, then the triad homotopy group $\pi_{k+1}(B,B',W)$ is finitely generated as a $\ZZ[\pi_1 W']$-module.
    \end{enumerate}
\end{lemma}
\begin{proof}
    Let $\alpha:(X,X')\to (B,B')$ be strongly $k$-connected map from a CW pair $(X,X')$ such that $X'$ has finite $k$-skeleton and the pair $(X,X')$ has finite $(k+1)$-skeleton. For $(i)$, we observe that $\alpha|_{X'}:X'\to B'$ can be made into a weak equivalence by attaching cells of dimension at least $k+1$. Let $Y$ be the CW complex given by the result of that procedure. Then $Y$ satisfies Wall's $(F_k)$ condition by \cite[Thm. A]{Wallfin}. Since $B'$ is weakly equivalent to $Y$, we conclude that $B'$ satisfies Wall's $(F_k)$ condition.   

    We focus now on $(ii).$ From the hypothesis, one obtains that the spaces $W,W',B$ and $B'$ are path-connected, hence we will drop the basepoint from the notation of the relative and triad homotopy groups that appear in this proof. Additionally, we obtain that $f:W\to B$ is $k$-connected, since it factors as $W\to W\cup_{W'} B'\to B$, both of which are $k$-connected maps. We denote $\pi'\coloneqq \pi_1(W')\cong \pi_1(B')$ and $\pi\coloneqq \pi_1(W)\cong \pi_1(B)$. Consider the long exact sequence of the map of pairs $(W,B)\to (W',B')$ induced by $f$
    \[\cdots \to \pi_{k+1}(B',W')\overset{f_*}{\to} \pi_{k+1}(B,W)\to \pi_{k+1}(B,B',W)\to \pi_{k}(B',W')\to \cdots\]of $\ZZ[\pi']$-modules. Since $f|_{W'}$ is $k$-connected, the rightmost module vanishes, and hence $\pi_{k+1}(B,B',W)$ is isomorphic to the cokernel of $f_*$. Moreover, since $f|_W$ and $f|_{W'}$ are $k$-connected, the Hurewicz homomorphisms $\pi_{k+1}(B',W')\to \H_{k+1}(B',W';\ZZ[\pi'])$ and $\pi_{k+1}(B,W)\to \H_{k+1}(B,W;\ZZ[\pi])$ are isomorphisms. Under these isomorphisms, the map $f_*$ corresponds to the composite 
    \[\H_{k+1}(B',W';\ZZ[\pi'])\to \H_{k+1}(B',W';f^*\ZZ[\pi])\overset{f_*}{\to} \H_{k+1}(B,W;\ZZ[\pi]).\]Denote the kernel of $\ZZ[\pi']\to \ZZ[\pi]$ by $K$. Since $\pi'\to \pi$ is surjective, the left map above fits in the long exact sequence
    \[\cdots \to \H_{k+1}(B',W';\ZZ[\pi'])\to \H_{k+1}(B',W';f^*\ZZ[\pi])\to \H_{k}(B',W';K)\to \cdots.\]Since $f|_{W'}$ is $k$-connected, the $\ZZ[\pi']$-module $\H_{k}(B',W';K)$ vanishes. We conclude that $\pi_{k+1}(B,B',W)$ is isomorphic to the cokernel of the map of $\ZZ[\pi']$-modules $f_*:\H_{k+1}(B',W';f^*\ZZ[\pi])\to \H_{k+1}(B,W;\ZZ[\pi])$. Once again by $k$-connectivity of $f|_{W'}$, this cokernel is isomorphic to $\H_{k+1}(B,B'\cup_{W'}W;\ZZ[\pi]).$ Thus, it suffices to prove that the latter module is finitely generated as a $\ZZ[\pi]$-module, as this implies that it is finitely generated as a $\ZZ[\pi']$-module, as $\pi'\to \pi$ is surjective. We will do so by showing that the chain complex of $\ZZ[\pi]$-modules $C_*(B,B'\cup_{W'}W;\ZZ[\pi]) $ of relative chains with coefficients in $\ZZ[\pi]$ is quasi-isomorphic to a chain complex $C_*$ such that $C_*=0$ for $*\leq k$ and $C_{k+1}$ is finitely generated. We will use the following general facts about chain complexes of $R$-modules for a ring $R$. We say that a chain complex of $R$-modules has \textit{finite $k$-skeleton} if it is quasi-isomorphic to a chain complex which is degreewise finitely generated projective in degrees at most $k$.
    \begin{enumerate}[label=(\textit{\roman*})]
        \item\label{cofibs} \textit{Let $i\in \ZZ$, given a cofiber sequence of chain complexes of $R$-modules $A\to B\to C$ such that $A$ and $B$ have finite $i-1$ and $i$-skeletons, respectively, then $C$ has finite $i$-skeleton}: this follows since $C$ is quasi-isomorphic to the cone of $A\to B$, which in degree $d\in \ZZ$ is given by $A_{d-1}\oplus B_d$.
        \item\label{finits} \textit{Let $i\geq 0$ and let $C$ be a chain complex of projective $R$-modules supported in non-negative degrees. If $\H_*(C)=0$ for $*\leq i$, then $C$ is quasi-isomorphic to a chain complex $D$ supported in degrees at least $i+1$. If, additionally, $C_*$ is finitely generated projective for $*\leq i+l$ for a fixed $l\geq 0$, then $D_*$ can be assumed to be finitely generated projective for $*\leq i+l:$} This is essentially an algebraic analogue of Whitehead's cell trading lemma \cite[Lemma 15]{Whitehead1950} (see also \cite[Lemma 6.1]{pedersenFerry}). By induction on $i$, we can assume that $C_*=0$ for $*\leq i-1$. Since $\H_i(C)=0$, it follows that $d:C_{i+1}\to C_i$ is surjective. Choose a section $s:C_i\to C_{i+1}$ of $d$, which is possible since $C_i$ is projective. Let $D(C_i)$ be the chain complex concentrated in degrees $i$ and $i+1$ with value $C_i$ in both entries, and whose only possible non-trivial differential is the identity. Let $C'$ be the cone of the chain map $s:D(C_i)\to C$ given by $\id_{C_i}$ in degree $i$ and $s$ in degree $i+1$. In particular, $C'_*=C_i$ for $*\neq i+1,i+2$ and $C'_{*}=C_*\oplus C_i$ for $*=i+1,i+2$. Since $D(C_i)$ is acyclic, it follows that $C\to C'$ is a quasi-isomorphism. Moreover, we have a degreewise split injective map $(0,\id):D(C_i)\to C'$ which is the identity in degree $i$ and is $0\oplus \id_{C_i}:C_i\to C_{i+1}\oplus C_i$ in degree $i+1$. Let $D$ be the degreewise cokernel of $(0,\id)$, which is quasi-isomorphic to $C$, as it is also the cofiber of $(0,\id)$, as this map is degreewise injective. Moreover, we see that $D_*=0$ for $*\leq i$, $D_{i+1}=C_{i+1}$, $D_{i+2}$ is the direct sum of projective modules, and $D_*=C_*$ for $*\geq i+3$. Finally, one traces that the additional finite generation hypothesis on $C$ implies the desired finite generation on $D$.
    \end{enumerate}
    Start by observing that we have a cofiber sequence of $\ZZ[\pi]$-modules 
    \[C_*(W,W';\ZZ[\pi])\to C_*(B,B';\ZZ[\pi])\to C_*(B,B'\cup_{W'} W;\ZZ[\pi]).\]The chain complex $C_*(W,W';\ZZ[\pi])$ is the cofiber of $C_*(W';\ZZ[\pi])\to C_*(W;\ZZ[\pi])$, which is map of chain complexes with finite $k$-skeleta, by taking cellular chains (which are degreewise free modules) and using the assumptions on the skeleta of $W'$ and $W$. From \ref{cofibs} above, we deduce that $C_*(W,W';\ZZ[\pi])$ has finite $k$-skeleton. On the other hand, the chain complex $C_*(B,B';\ZZ[\pi])$ is the cofiber of the connecting map $C_{*+1}(B,B'\cup_{X'}X;\ZZ[\pi])\to C_*(X,X';\ZZ[\pi])$, where $(X,X')$ is the source of the map $\alpha$ defined in the first line of this proof. By taking cellular chains and the assumption on $(X,X')$, we obtain that $C_*(X,X';\ZZ[\pi])$ has finite $(k+1)$-skeleton. Since $B'\cup_{X'} X\to B$ is $(k+1)$-connected, we can apply \ref{finits} to obtain that $C_{*+1}(B,B'\cup_{X'}X;\ZZ[\pi])$ is quasi-isomorphic to a chain complex concentrated in degrees at least $k+1$. In particular, this former chain complex has finite $k$-skeleton. From \ref{cofibs}, we conclude that $C_*(B,B'\cup_{W'} W;\ZZ[\pi])$ has finite $(k+1)$-skeleton. Now applying \ref{finits} to this chain complex, we obtain that $C_*(B,B'\cup_{W'} W;\ZZ[\pi])$ is quasi-isomorphic to a chain complex concentrated in degrees at least $k+1$ and whose entry in degree $k+1$ is finitely generated, as $B'\cup_{W'} W\to B$ is $k$-connected. This finishes this proof.
\end{proof}

\begin{lemma}\label{take out n stably}
    In the context of \cref{stably highly connected theta ends}, the map induced by inclusion
    \[\Fo(P,K|_\infty)\to \CobbtLs(P,K|_\infty) \]is an equivalence.
\end{lemma}
\begin{proof}
    Start by noticing that $\Fo(P,K|_i)$ is a union of path components of $\CobbtLs(P,K|_i)$ for every $i\geq 0.$ So it suffices to prove that the induced map on $\pi_0$ is surjective after taking the colimit over $i.$ Moreover, it suffices to prove that given a morphism $W:P\leadsto K|_i$, there exists a $j\geq i$ such that $(W\cup K|_{[i,j]},\hb (W\cup K|_{[i,j]})\to (B,B^\partial)$ is strongly $n$-connected. 
    
    Since $(K|_i,\hb K|_i)\to (W,\hb W)$ and $(K|_i,\hb K|_i)\to (B,B^\partial)$ (see \cref{take out fo}) are strongly $(n-1)$-connected, we have that $(W,\hb W)\to (B,B^\partial)$ is strongly $(n-1)$-connected by \cref{closure properties strongly}. We show that there exists $j,$ such that the map $(W\cup K|_{[i,j]},\hb (W\cup K|_{[i,j]})) \to (B,B^\partial)$ is strongly $n$-connected. This is inspired by \cite[Lemma 7.7]{GRWII}. We first show that there exists $j,$ such that the map $\hb (W\cup K|_{[i,j]})\to B^\partial$ is $n$-connected. By \cref{some fn stuff}, $B^\partial$ is $(F_n)$ and so $\pi_n(B^\partial,\hb W)$ is a finitely generated $\ZZ[\pi_1(B^\partial)]$-module, since $n\geq 3$. This implies that $\ker(\pi_{n-1}(\hb W)\to \pi_{n-1}(B^\partial))$ is a finitely generated $\ZZ[\pi_1(B^\partial)]$-module. Choose a finite basis $\{a_p\}_{p}$ for this module. It follows from the hypothesis, there exists an integer $l\geq 0 $ such that the images of $a_p$ in $\pi_{n-1}(\hb (W\cup K|_{[i,i+l]}))$ vanish for all $p$. Similarly, we need to add a finite number of elements of $\pi_n(B^\partial)$ to the image of $\pi_n(\hb (W\cup K|_{[i,i+l]}))$ to generate $\pi_n(B^\partial)$ as $\ZZ[\pi_1(B^\partial)]$-modules. Thus there exists an $l'\geq l$ such that, additionally, the map $\pi_n(\hb (W\cup K|_{[i,i+l']}))\to \pi_n(B^\partial)$ is surjective. Therefore, the map $\hb (W\cup K|_{[i,i+l']})\to B^\partial$ is $n$-connected.
    
    Let $W'\coloneqq W\cup K|_{[i,i+l']}$ and consider now the group $\pi_{n+1}(B, B^\partial, W').$ By \cref{some fn stuff}, this group is finitely generated as a $\ZZ[\pi_1(\hb W')]$-module. This implies that $\ker(\pi_{n}(W',\hb W')\to \pi_{n}(B,B^\partial))$ is finitely generated as a $\ZZ[\pi_1(\hb W')]$-module. Choose a finite basis $\{b_q\}_q$ for this module. By hypothesis, there exists an integer $m\geq 0 $ such that the images of $b_p$ in $\pi_{n}(W'\cup K|_{[l',l'+m]}, \hb (W'\cup K|_{[l',l'+m]}))$ vanish for all $q.$ Similarly, we can generate $\pi_{n+1}(B,B^\partial)$ as a module over $\pi_1(B^\partial)\cong \pi_1(\hb W')$ by the image of $\pi_{n+1}(W'\cup K|_{[l',l'+m]}, \hb (W'\cup K|_{[l',l'+m]}))$ and finitely many elements. Thus, there exists an $m'\geq m$ such that, additionally, the map $\pi_{n+1}(W'\cup K|_{[l',l'+m']}, \hb (W'\cup K|_{[l',l'+m']}))\to \pi_{n+1}(B,B^\partial)$ is surjective. Therefore, the map $(W\cup K|_{[i,i+l'+m']},\hb(W\cup K|_{[i,i+l'+m']}) \to (B,B^\partial)$ is strongly $n$-connected. This finishes the proof.
\end{proof}

\begin{proof}[Proof of \cref{stably highly connected theta ends}]
    We apply \cite[Thm. A.14]{GRWII} to the category $\smash{\CobbtLs}.$ Condition $(i)$ follows as in the proof of \cite[Prop. 7.5]{GRWII}. To verify $(ii),$ we observe that the functor $\smash{\CobbtLs(-,K|_\infty)}$ is equivalent to $\smash{\F(-,K|_\infty)}$ by \cref{take out n stably} and \cref{take out fo}. The latter takes all morphisms to abelian homology equivalences by \cref{stable stability}, which is precisely $(ii)$. Thus, we obtain that the map
    \[\CobbtLs(P,K|_{\infty})\to \Omega_{[P,K|_\infty]}\B\CobbtLs\]is acyclic. By pre-composing with the map in \cref{take out n stably}, we finish the proof.
\end{proof}

\subsubsection{\texorpdfstring{$\Theta$-surgery.}{Theta-surgery.}}

The goal of this subsubsection is to provide a general construction which will allow us to produce a $\Theta$-end $K'$ satisfying the hypothesis of \cref{stably highly connected theta ends} from a $\Theta$-end $K$ satisfying only the hypothesis of \cref{highly connected theta ends}. In the next subsubsection, we will relate these $\Theta$-ends to give a proof of the latter statement. This construction takes the form of the following result. 

\begin{prop}\label{Krecks surgery}
    Assume $(B,B^\partial)$ is $(E_{n})$ and let $R$ be an object in $\CobbtL$ such that $\ell_R:(R,\hb R)\to (B,B^\partial)$ is strongly $(n-1)$-connected. If $ (R,\hb R)$ is $1$-connected, then there exists a morphism $W_R:R\leadsto R'$ in $\CobbtL$ satisfying the following properties:
    \begin{enumerate}[label=(\textit{\alph*})]
        \item\label{connectivity between kreck} The inclusions $(R,\hb R)\hookrightarrow (W_R,\hb W_R)\hookleftarrow (R',\hb R')$ are strongly $(n-1)$-connected;
        \item\label{connectivity to B kreck} $\ell:(W_R,\hb W_R)\to (B,B^\partial)$ is strongly $n$-connected. 
    \end{enumerate}
\end{prop}

 One can see this result as an analog of \cite[Prop. 4]{Kreck1999SurgeryAD} for triads. In fact, the strategy is analogous to the proof of loc.cit. (see also \cite[Lemma 7.6]{GRWII}). We start with two constructions that will be essential to the proof of \cref{Krecks surgery}. Let $(M,\ell_M)\in \smash{\CobbtL}$ and $k\geq 0$ be an integer. Assume we are given an embedding $e:S^k\to \hb M\backslash \hvb M$ and a map $\ell:D^{k+1}\to B^\partial$ extending $e^*\ell_M$ on $\partial D^{k+1}.$ We say that $e$ and $\ell$ \textit{occur as the trace of a right $\Theta$-surgery} if there exists a morphism $(W,\ell_W):M\leadsto M'$ which admits a  handle decomposition with a single right $(k+1)$-handle (recall \cref{defn of handles}) relative to $M$, whose attaching map is a thickening of $e$ such that there exists a homotopy commutive diagram 
\[\begin{tikzcd}
    (M\cup_{e} D^{k+1},\hb M\cup_{e} D^{k+1}) \arrow[r, "\iota"]\arrow[rd, "{\ell_M}\cup {\ell}"'] & (W,\hb W) \arrow[d, "\ell_W"] \\
    & (B,B^\partial)
\end{tikzcd}\]where $\iota$ is the inclusion map induced by the core of the single right $(k+1)$-handle. The following lemma provides a sufficient condition for this to be the case. This will be useful to construct morphisms in $\CobbtL$ with prescribed $\Theta$-structure.  

\begin{lemma}\label{right theta surgery}
    Let $e:S^k\to \partial M$ be an embedding and $\ell:D^{k+1}\to B^\partial$ be an extension of the $\Theta$-structure $e^*\ell_M$. Then $e$ and $\ell$ occur as the trace of a right $\Theta$-surgery, provided $k\leq n-1$.
\end{lemma}
\begin{proof}
    This is essentially \cite[Section 4.1]{GRWStableModuli} by applying it to $\partial M$ and extending the trace and $\theta^\partial$-structure to a collar.
\end{proof}

The second construction is analogous to the first but for left surgery. Let again $(M,\ell_M)\in \smash{\CobbtL}$ and $k\geq 0$ be an integer. Let \(e:(\partial_1D^{k+1}_+,\partial_{01}D^{k+1}_+)\to (M,\partial M)\) be an embedding along with a map $\ell:(D^{k+1}_+,\partial_0D^{k+1})\to (B,B^\partial)$ that extends $e^*\ell_{M}$ in $(\partial_1D^{k+1}_+,\partial_{01}D^{k+1}_+).$ We say that $e$ and $\ell$ \textit{occur as the trace of a left $\Theta$-surgery} if there exists a morphism $(W,\ell_W):M\leadsto M'$ which admits a handle decomposition with a single left $(k+1)$-handle (recall \cref{defn of handles}) whose attaching map is a thickening of $e$ such that there exists a homotopy commutive diagram 
\[\begin{tikzcd}
    (M\cup_{e} D^{k+1}_+,\hb M\cup_{e} \partial_{1}D^{k+1}_+) \arrow[r, "\iota"]\arrow[rd, "{\ell_M}\cup {\ell}"'] & (W,\hb W) \arrow[d, "\ell_W"] \\
    & (B,B^\partial)
\end{tikzcd}\]where $\iota$ is the inclusion map induced by the core of the single left $(k+1)$-handle. The next lemma is analogous to \cref{right theta surgery}.

\begin{lemma}\label{left theta surgery}
    Let $e:(\partial_1D^{l+1}_+,\partial_{01}D^{l+1}_+)\to (M,\partial M)$ be an embedding and $\ell_D:(D^{l+1}_+,\partial_0D^{l+1})\to (B,B^\partial)$ be an extension of the $\Theta$-structure $e^*\ell_M$. Then $\ell_D$ occurs as the trace of a left $\Theta$-surgery, provided $l\leq 2n-1$.
\end{lemma}
\begin{proof}
    We follow the strategy of \cite[Section 4.1]{GRWStableModuli}. The embedding $e$ induces an embedding \[e_1=[0,1]\times e:[0,1]\times (\partial_1D^{l+1}_+,\partial_{01}D^{l+1}_+)\to [0,1]\times (M,\partial M).\]Let $A\coloneqq(D^{l+1}_+\backslash \text{int} (\frac{1}{2} D^{l+1}_+),\partial_0D^{l+1}_+\backslash \text{int} (\frac{1}{2} \partial_0D^{l+1}_+))$ and identify it with $[0,1]\times (\partial_1D^{l+1}_+,\partial_{01}D^{l+1}_+).$ As in \cite[Section 4.1]{GRWStableModuli}, $\ell_M$ induces a collared bundle map of pairs $T(D^{l+1}_+,\partial_0D^{l+1}_+)|_A\to (\theta^*\gamma_{2n+1},(\theta^\partial)^*\gamma_{2n}).$ We start by proving that such a map admits an extension to a collared bundle map \(\smash{T(D^{l+1}_+,\partial_0D^{l+1}_+)\to (\theta^*\gamma_{2n+1},(\theta^\partial)^*\gamma_{2n})}\) which covers the map $\ell_D.$ Note that it suffices to find a collared bundle map \(\smash{T(D^{l+1}_+,\partial_0D^{l+1}_+)\to \ell_D^*(\theta^*\gamma_{2n+1},(\theta^\partial)^*\gamma_{2n})}\)over the pair $(D^{l+1}_+,\partial_0D^{l+1}_+)$ extending the given map on the subspace $A.$ Since both collared bundles are trivial, as they are bundles over a contractible pair, this problem is equivalent to finding a lift of the following lifting problem
    \[\begin{tikzcd}
        (\partial_1D^{l+1}_+,\partial_{01}D^{l+1}_+)\arrow[r]\arrow[d] & (V_l(2n+1),V_{l-1}(2n))\\
        (D^{l+1}_+,\partial_1D^{l+1}_+)\arrow[ur,dashed]
    \end{tikzcd}\]where $V_i(j)$ is the Stiefel manifold of $i$-frames on $\RR^j$ and the top map is induced by the trivialization over the subspace $A$. The upper right pair is $(2n-1)$-connected, thus a lift always exists if $l\leq 2n-1.$ Choose now such a collared bundle map of pairs $T(D^{l+1}_+,\partial_0D^{l+1}_+)\to (\theta^*\gamma_{2n+1},(\theta^\partial)^*\gamma_{2n}),$ which induces a collared vector bundle pair $(V,V')\to (D^{l+1}_+,\partial_0D^{l+1}_+)$ of dimension $2n-l$ and collar bundle map
    \[T(D^{l+1}_+,\partial_0D^{l+1}_+)\oplus (V,V')\to (\theta^*\gamma_{2n+1},(\theta^\partial)^*\gamma_{2n}),\]which is a fiberwise isomorphism and extends the given map. The disc bundle of $(V,V')$ is then identified with the normal bundle of $e_1$ once restricted to $A.$ Since the bundle $(V,V')$ has contractible base, then its disc bundle is diffeomorphic to $(D^{l+1}_+\times D^{2n-l},\partial_0D^{l+1}_+\times D^{2n-l})$ Consider now the pushout of pairs 
    \[\begin{tikzcd}
        (D(V),D(V'))|_A\arrow[d]\arrow[r] & (D(V),D(V'))\arrow[d] \\
        ([0,1]\times M,[0,1]\times \partial M) \arrow[r] & (W,\hb W)
    \end{tikzcd}\]which is a model for attachment of a left $(l+1)$-handle along a thickening of $e.$ This description comes with a preferred $\Theta$-structure given by the bundle maps defined above. In other words, $\ell_D$ occurs as a left trace of a surgery. This finishes the proof.
    \end{proof}

We remark that the discussion above did not use the fact that the dimension of the cobordisms is odd in any essential way. We are now ready to start the proof of \cref{Krecks surgery}. This will be done by constructing $W_R$ step by step increasing the connectivity of its $\Theta$-structure. We start by deducing the existence of an intermediate morphism, where \ref{connectivity to B kreck} is satisfied in the boundary.

\begin{lemma}\label{killing boundary kernel}
    Let $R$ be an object of $\CobbtL$ satisfying the conditions of \cref{Krecks surgery}. Then there exists a morphism $(W_0,\ell_0):R\leadsto R_0$ in $\CobbtL$ satisfying properties \ref{connectivity between kreck} and where the map $\ell_0:\hb W_0\to B^\partial$ is $n$-connected.
\end{lemma}
\begin{proof}
    This follows by applying the construction of the proof of \cite[Lemma 7.6]{GRWII} to $P=\hb R$ and extending it to a morphism in $\CobbtL$ using right handles (see \cref{extending handle decompositions}), which is possible by \cref{right theta surgery}. The property \ref{connectivity between kreck} follows from \ref{effect of right handles} from \cref{strong geo section}. The map $\ell_0$ is $n$-connected by property $(iv)$ in the proof of \cite[Lemma 7.6]{GRWII}. 
\end{proof}

\begin{lemma}\label{killing relative kernel}
    Let $R$ be an object of $\CobbtL$ satisfying the conditions of \cref{Krecks surgery}. Then there exists a morphism $(W_1,\ell_1):R\leadsto R_1$ in $\CobbtL$ satisfying properties \ref{connectivity between kreck} and where the map $\ell_1|_{\hb W_1}:\hb W_1\to B^\partial$ is $n$-connected and the map $(\ell_1)_*:\pi_{i}(W_1,\hb W_1)\to \pi_i(B,B^\partial)$ is an isomorphism for $i\leq n.$ 
\end{lemma}
\begin{proof}
    Let $(W_0,\ell_0):R\leadsto R_0$ be a morphism satisfying the properties of \cref{killing boundary kernel}. Observe that $\pi_1(\hb W_0)\cong \pi_1(B^\partial)$, since $n\geq 2.$ We start by observing that the triad homotopy group (recall from \cref{triad homotopy}) $\pi_{n+1}(B,B^\partial,W_0)$ is a finitely generated $\ZZ[\pi_1(B^\partial)]$-module, by \cref{some fn stuff}. 
    We construct now the surgery data in order to construct the cobordism $W_1$. The $\ZZ[\pi_1(B^\partial)]$-module
    \[\ker\left((\ell_0)_*:\pi_{n}(W_0,\hb W_0)\to \pi_{n}(B,B^\partial)\right)\]is finitely generated: consider the long exact sequence on homotopy groups 
    \[\cdots \to \pi_{n+1}(B,B^\partial)\to \pi_{n+1}(B,B^\partial, W_0)\to \pi_{n}(W_0,\hb W_0)\to \pi_{n}(B,B^\partial)\to \cdots.\]This is an exact sequence of $\ZZ[\pi_1(B^\partial)]$-modules, since $\hb W_0\to B^\partial$ is an isomorphism on $\pi_1$. However, the module $\pi_{n+1}(B,B^\partial, W_0)$ surjects to the kernel considered above and thus the latter is also finitely generated. Pick generators $\beta_1,\cdots,\beta_q$ of this kernel as a $\ZZ[\pi_1(B^\partial)]$-module. Since the map $\pi_{n}(R_0,\hb R_0)\to \pi_{n}(W_0,\hb W_0)$ is surjective, we can represent these generators by maps 
    \[\beta_i:(D^{n},\partial D^{n})\to (R_0,\hb R_0).\]To represent these elements by embeddings of pairs to the interior of $(R_0,\hb R_0)$, we use a result by Hudson \cite[Thm. 1]{Hudson_1972}, best stated for our purposes in \cite[Thm. 5.13]{BP}, which requires $n\geq 3$. In the notation of \cite{BP}, we use this result for $(P,\partial_0P,\partial_1P)\coloneqq(D^{n},\partial D^{n},\emptyset)$ and $(W,\partial_0W,\partial_1W)\coloneqq(R_0,\hb R_0,\vb R_0)$, since the inclusion $\hb R\to R$ is $1$-connected. Once such embeddings are found, we can make them disjoint by \cite[Thm. C.1]{BP}, as $n\geq 3$. Once again, use \cref{left theta surgery} along with nullhomotopies of $\ell$ restricted to $\beta_i$ to define the morphism $W_1:R\leadsto R_1$ given by the simultaneous attachment of left $(n+1)$-handles to $W_0$ along thickenings of the embeddings $\beta_i.$ We finish by verifying the properties claimed in the statement. Once again by \ref{effect of left handles} from \cref{strong geo section}, $(W_1,\hb W_1)$ is strongly $(n-1)$-connected relative to $(W_0,\hb W_0)$ and thus relative to $(R,\hb R).$ By \cref{dual handle} and \ref{effect of right handles} from \cref{strong geo section}, the inclusion of $(R_1,\hb R_1)$ to $(W_1,\hb W_1)$ is $(n-1)$-connected. The cobordism $\hb W_1 $ is obtained from $\hb W_0$ by attaching trivial $n$-cells along $\beta_i|_{\partial D^n}$ (by \cref{killing boundary kernel}), which implies that $\hb W_0\to \hb W_1$ admits a section collapsing these cells. Thus, the map $\pi_i(\hb W_0)\to \pi_i(\hb W_1)$ is an isomorphism for $i\leq n-1$ and injective for $i=n.$ Thus, the map $\hb W_1\to B^\partial$ is $n$-connected, as $\hb W_0\to B^\partial$ so is. Finally, the group $\pi_i(W_1,\hb W_1)$ is isomorphic to $\pi_i(W_0,\hb W_0)\cong \pi_i(R,\hb R)$ for $i< n$, and thus also to $\pi_i(B,B^\partial).$ For $i=n,$ the group is a quotient of $\pi_n(W_0,\hb W_0)$ by the submodule generated by the classes $\beta_i,$ and isomorphic to $\pi_n(B,B^\partial)$. This finishes the proof.
\end{proof}

To finish the proof of \cref{Krecks surgery}, recall the definition of triad connected sum from \cref{triad connect sum}. 

\begin{rmk}\label{tangential structure of connect sum}
   Given two $\Theta$-triads $W$ and $W'$ such that $\ell_W$ and $\ell_{W'}$ hit the same path component of $B.$ Then $W\natural W'$ admits a $\Theta$-structure extending $\ell_W$ and $\ell_{W'}$ using the following argument: Pick a path from $\ell_W(e(0))$ and $\ell_{W'}(e(0))$ in $B^\partial$ and define $\Theta$-structure on $(W\natural W',\hb (W\natural W'))\simeq (W\cup_{e(0)}[0,1]\cup_{e'(0)}W',\hb W\cup_{e(0)}[0,1]\cup_{e'(0)}\hb W')$ given by this path on $[0,1].$ In most cases, we can take embeddings such that $\ell_W(e(0))$ and $\ell_{W'}(e(0))$ coincide. In this case, we take the constant path.
\end{rmk}

\begin{proof}[Proof of \cref{Krecks surgery}]
    Let $W_1:R\leadsto R_1$ be a morphism satisfying the properties of \cref{killing relative kernel}. We can pick a finite collection of $\zeta_0, \zeta_1,\cdots,\zeta_m\in \pi_{n+1}(B,B^\partial)$ that generate $\pi_{n+1}(B,B^\partial)$ as a $\ZZ[\pi_1(B^\partial)]$-module along with the image of $\pi_{n+1}(W_2,\hb W_2).$ This is possible since $\pi_{n+1}(B,B^\partial,W_2)$ is a finitely generated $\ZZ[\pi_1(B^\partial)]$-module, by the same argument as in loc.cit. Let $(F_i,\hb F_i)$ be the $n$-disc bundle given by a lift of the form
    \[\begin{tikzcd}
        & &(\BO(n),\BO(n))\arrow[d]\\
        (D^{n+1},\partial D^{n+1})\arrow[r]\arrow[urr,dashed] &(B,B^\partial)\arrow[r] &(\BO(2n+1),\BO(2n)). 
    \end{tikzcd}\]which exists and is unique by obstruction theory since $\pi_i(\BO(2n+1),\BO(2n),\BO(n))\cong \pi_i(\BO(2n+1),\BO(2n))$ vanishes for $i<2n+1$. Once again, we endow such pair with the $\Theta$-structure given by the composite $(F_i,\hb F_i)\to (D^{n+1},\partial D^{n+1})\to (B,B^\partial). $ Let $W_P$ be the cobordism given by a triad connected sum of $W_2$ with the triad determined by $(F_i,\hb F_i).$ Similar to before, we can obtain $F_i$ from $D^{2n+1}_+$ by attaching a left $(n+1)$-handle so \ref{connectivity between kreck} is automatic by \ref{effect of right handles} and \ref{effect of left handles} from \cref{strong geo section}, and \cref{dual handle}. By construction, the map $\pi_{n+1}(W_R,\hb W_R)\to \pi_{n+1}(B,B^\partial)$ is surjective. By \cref{1 conn: old same as new}, we obtain \ref{connectivity to B kreck}. By the same argument as in \cref{killing relative kernel}, we do not lose the established properties of $W_1$. This completes the proof since $(W_R,\hb W_R)$ is $1$-connected, as $n\geq 2$.
\end{proof}

\subsubsection{\texorpdfstring{Proof of \cref{highly connected theta ends}.}{Proof of second preliminary version.}}
In this subsubsection, we use the general construction from the previous subsubsection and \cref{stably highly connected theta ends} to prove the main result of this subsection. We start with the following lemma which uses this construction.

\begin{lemma}\label{theta end kreck}
    Assume $(B,B^\partial)$ is $1$-connected and $(E_n)$. If $P\in \CobbtL$ such that $\ell|_{P}$ is strongly $(n-1)$-connected, then there exists a $\Theta$-end $K$ such that $K|_0=P$ and $\ell|_{K|_{[i,+\infty)}}$ is strongly $n$-connected for every $i\geq 0$. 
\end{lemma}
\begin{proof}
    Given any object $P\in \CobbtL$ such that $\ell|_P$ is strongly $(n-1)$-connected, it suffices to construct a morphism $W_P:P\leadsto P'$ in $\smash{\CobbtLo}$ such that:
    \begin{enumerate}[label=(\textit{\roman*})]
        \item $W_P$ is strongly $(n-1)$-connected relative to both ends;
        \item $W_P$ contains an embedded copy of $V_1$ with standard $\Theta$-structure (see \cref{standard structure on V1})
        \item $\ell|_{W_P}$ is strongly $n$-connected.
    \end{enumerate}Given such construction, $\ell|_{P'}$ is strongly $(n-1)$-connected, so we can iterate this contruction to produce the $\Theta$-end $W_{P}\cup W_{P'}\cup W_{(P')'}\cup \cdots.$ This is constructed by applying \cref{Krecks surgery} (producing a morphism satisfying $(i)$ and $(iii)$ of above) and doing a triad connected sum (recall \cref{triad connect sum}) with $V_1$ with a standard $\Theta$-structure (to satisfy $(ii)$). This does not affect $(i)$ or $(iii)$ since $W\to W\natural V_1$ is strongly $n$-connected. (This is possible since $(P,\hb P)$ is $1$-connected, given the assumption on $(B,B^\partial)$ and $n\geq 3.$) 
\end{proof}

We are now ready to prove the first preliminary version of our main result of this section.

\begin{proof}[Proof of \cref{highly connected theta ends}]
   We start by fixing a $\Theta$-end $\{P|_{[i,i+1]}:P|_i\leadsto P|_{i+1}\}$ with $P|_0=P$. For example, take $P|_{[i,i+1]}=\prescript{}{P|_i}{H}$ (as defined before \cref{stability for V_1}). For any two $\Theta$-ends $K$ and $K'$ such that $K|_0=K'|_0$, we construct a zig-zag of abelian homology equivalences 
   \[\Fo(P,K|_\infty)\to \cdots \leftarrow \Fo(P,K'|_\infty)\]which commutes with an analogous zig-zag of equivalences
   \[\Omega_{[P,K|_\infty]}\B\CobbtLs\to \cdots \leftarrow \Omega_{[P,K'|_\infty]}\B\CobbtLs.\]The proof follows by letting $K'$ be the $\Theta$-end constructed in \cref{theta end kreck} and applying \cref{stably highly connected theta ends} and \cref{theta staisfies En}. This implies that the map in the claim is an abelian homology equivalence. However, the target of this map has abelian fundamental group, and thus the map is actually acyclic, since all local systems are abelian (see \cite[Appendix A]{GRWII}). 

   We construct the aforementioned zig-zags in the following way
   \[\begin{tikzcd}
	{\Fo(P,K|_\infty)} & {\hocolim_j\hocolim_i \Fo(P|_{\!j},K|_i)} \\
	& {\Fo(P|_{\infty},K|_0)=\Fo(P|_{\infty},K'|_0)} \\
	{\Fo(P,K'|_\infty)} & {\hocolim_{\!j}\hocolim_i \Fo(P|_{j},K'|_i)}
	\arrow[from=1-1, to=1-2]
	\arrow[from=2-2, to=1-2]
	\arrow[from=2-2, to=3-2]
	\arrow[from=3-1, to=3-2]
\end{tikzcd}\]where the horizontal maps are induced by pre-composition by $P|_*$ and the vertical maps by post-composition by $K$. The horizontal maps are abelian homology equivalences by \cref{stable stability}, since they are the sequential colimit of abelian homology equivalences of the form 
\[\Fo(P|_i,K|_\infty)\to \Fo(P|_{i-1},K|_\infty)\]and such equivalences are closed under sequential colimits. Since $\ell_{K|_i}$ is strongly $(n-1)$-connected (since $\ell_{K|_0}$ is and $K|_{[i,i+1]}$ is strongly $(n-1)$-connected relative to both ends), we see that $\Fo(P|_j,K|_i)= \F(P|_j,K|_i)$ since $n\geq 2.$ We can thus apply \cref{stable stability}. The vertical maps are abelian homology equivalences by the analogous statement of \cref{stable stability} where stabilization is done on the left. More precisely, the top vertical map is a colimit of maps of the form
\[\Fo(P|_\infty,K|_i)\to \Fo(P|_\infty,K|_{i+1}).\]By the same argument as above, both spaces can be replaced by the analogous ones for $\F(-,-).$  Applying \cref{stable stability} to the $\Theta$-end $\overline{P}|_{-*}$ given by the reflection (see \cite[140]{GRWII}) and composing it with the reflection equivalence
\[\F(P,Q)\simeq \Cob^{\partial}_{\Theta,\vb \overline{L},n}(\overline{Q},\overline{P})\]where the $\Theta$-structure on $\overline{Q}$ is the pre-composition of $\ell_Q$ with the equivalence of vector bundle pairs $(-1)\oplus \id:\varepsilon^1\oplus TQ\simeq \varepsilon^1\oplus TQ$ (see \cite[Defn. 5.1]{GRWStableModuli} for more details on the reflection automorphism).
It is clear that the diagram above is compatible with the analogous diagram for the different (colimits of) loop spaces of $\B\CobbtL.$ All the maps in that diagram are induced by concatenation with paths, and hence are weak equivalences. This finishes the proof. 
\end{proof}

\subsection{\texorpdfstring{Weak maps of $\Theta$-ends.}{Weak maps of Theta-ends.}} The goal of this subsection is to prove a second preliminary version of \cref{final with null}, which builds on \cref{highly connected theta ends}. More precisely, this version takes the following form. Recall the definition of triad handle decomposition from \cref{defn of handles}.

\begin{prop}\label{final without null}
    Assume $(B,B^\partial)$ is $1$-connected and $(L,\hb L)$ admits a triad handle decomposition with handles of any type and of index at most $n-1$ relative to $\emptyset$, and the maps $\vb L\to L$ and $\hvb L\to \hb L$ are fundamental groupoid isomorphisms. Let $K$ be a $\Theta$-end and $P\in \CobbtLo$ such that $\F(P,K|_0)$ is non-empty. Then
    \[\F(P,K|_\infty)\to \Omega_{[P,K|_\infty]}\B\CobbtL\]is acyclic.
\end{prop}

\begin{rmk}
    If $n\geq 3$, then the condition that $\vb L\to L$ and $\hvb L\to \hb L$ are fundamental groupoid isomorphisms, follows from the condition that $(L,\hb L)$ admits a triad handle decomposition with handles of any type and of index at most $n-1$ relative to $\emptyset$, by \ref{effect of right handles} to \ref{effect of interior handles} from \cref{strong geo section}.
\end{rmk}

Once again, we proceed by constructing a $\Theta$-end satisfying the conditions of \cref{highly connected theta ends} from one which only satisfies the hypothesis from above. In this case, these $\Theta$-ends are related in a different way to the one from the last subsection. To make this precise, it is convenient to introduce the following definition. 

\begin{defn}[Weak maps of $\Theta$-ends.]
    Let $K'$ and $K$ be $\Theta$-ends in $\smash{\CobbtL}$. A \textit{weak map of $\Theta$-ends} $V:K'\leadsto K$ is a collection of morphisms $V_i:K'|_i\leadsto K|_i$ in $\smash{\CobbtL}$ for every $i\geq 0$ such that the composition $K'|_i\leadsto K'|_{i+1}\leadsto K|_{i+1}$ is diffeomorphic to the composition $K'|_i\leadsto K|_i\leadsto K|_{i+1}$ for every $i\geq 0$ and for every $i\geq 0,$ the inclusion $(K'|_i,\hb K'|_i)\to (V_i,\hb V_i)$ is strongly $(n-1)$-connected.    
\end{defn}

By fixing diffeomorphisms issuing the condition above, one produces a map
\[V_*:\F(P,K'|_\infty)\to \F(P,K|_\infty)\]by taking homotopy colimits of the homotopy coherent diagram given by the maps $V_i\circ (-):\F(P,K'|_i)\to \F(P,K|_i)$. The reason why, in the definition above, we do not fix diffeomorphisms is that the property that this map is an abelian homology equivalence does not depend on the choice of diffeomorphisms. More precisely, if such a map, given a certain choice, is an abelian homology equivalence, then any other map induced by another choice also is one. This follows from \cite[Lemma A.10]{GRWII} and was pointed out in \cite[183]{GRWII}. Therefore, the property that $V_*$ is an abelian homology equivalence is a well-defined property of the weak map of $\Theta$-ends $V$. 

\begin{lemma}\label{theta ends and maps}
    Assume the hypothesis of \cref{final without null}. Let $V:K'\leadsto K$ be a weak map of $\Theta$-ends, then the map 
    \[V_*:\F(P,K'|_\infty)\to \F(P,K|_\infty)\]is an abelian homology equivalence.
\end{lemma}
\begin{proof}
    As in \cref{highly connected theta ends}, we start by fixing a $\Theta$-end $\{P|_{[i,i+1]}:P|_i\leadsto P|_{i+1}\}$ with $P|_0=P$. Consider the following square
    \[\begin{tikzcd}
	{\F(P,K'|_\infty)} & {\F(P,K|_{\infty})} \\
	{\F(P|_\infty,K'|_\infty)} & {\F(P|_\infty,K|_\infty)}
	\arrow["{V_*}", from=1-1, to=1-2]
	\arrow[from=1-1, to=2-1]
	\arrow[from=1-2, to=2-2]
	\arrow["{V_*}", from=2-1, to=2-2]
\end{tikzcd}\]where the vertical maps are the maps induced by $V_*$ after taking homotopy colimits with respect to stabilization on the left by the $\Theta$-end. The vertical left map is the homotopy colimit of the abelian homology equivalences of the form
\(\smash{\F(P|_i,K'|_\infty)\to \F(P|_{i-1},K'|_\infty)}\), by \cref{stable stability}. Hence, the map to the homotopy colimit is also an abelian homology equivalence. The same applies for the right vertical map. On the other hand, the lower horizontal map is the homotopy colimit over $i$ of maps \(\smash{\F(P|_\infty,K'|_i) \to \F(P|_\infty, K|_i)} \) which are abelian homology equivalences by the dual of \cref{stable stability} (see the proof of \cref{highly connected theta ends}). This finishes the proof by $2$-out-of-$3$ property for abelian homology equivalences (see \cref{two out of three}).
\end{proof}


Following the result above, our goal is to find a weak map of $\Theta$-ends $K'\leadsto K$ such that $\ell|_{K'|_i}$ is strongly $(n-1)$-connected, for any $\Theta$-end $K$ satisfying the hypothesis of \cref{final without null}. The following lemma achieves precisely that.

\begin{lemma}\label{truncation theta end}
    Let $K$ be a $\Theta$-end and $P\in \CobbtLo$ such that $\F(P,K|_0)\neq \emptyset$, then there exists a weak map  of $\Theta$-ends $K'\leadsto K$ such that $\ell|_{K'|_i}$ is strongly $(n-1)$-connected, for every $i\geq 0$.
\end{lemma}
\begin{proof}
    Let $W\in \F(P,K|_0)$ and fix a triad handle decomposition of $(W,\hb W)$ relative to $(P,\hb P)$ (see \cref{defn of handles}). Define $V_0:K'|_0\leadsto K|_0$ to be the union of all right handles of index at least $n$, and interior and left handles of index at least $n+1$. Thus, by \ref{effect of right handles} to \ref{effect of interior handles} from \cref{strong geo section}, the inclusion $(K'|_0,\hb K'|_0)\hookrightarrow (V_0,\hb V_0)$ is strongly $(n-1)$-connected. The inclusion of $K'|_0\hookrightarrow W$ is strongly $(n-1)$-connected, since it is obtained by attaching right handles of index $n$, interior handles of index at least $n+1$ and left handles of index at least $n+2.$ Thus, the map $(K'|_0,\hb K'|_0)\to (B,B^\partial)$ is strongly $(n-1)$-connected. Consider the composition $K|_{[0,1]}\circ V_0:K'|_0\leadsto K|_1.$ By assumption, $K|_{[0,1]}\cong M\circ H_{K|_0}$ for some morphism $M$. The handle $H_{K|_0}$ is attached trivially to $K|_0$ by definition, so it can be made disjoint from the belts of the handles in $V$ since $n\geq 1$. Hence, by \cref{rearragement,factorizations and handle structures}, the composition $K|_{[0,1]}\circ V_0:K'|_0\leadsto K|_1$ factors as $V_1\circ H_{K'|_0}$ for some $V_1:K'|_1\leadsto K|_1.$ Clearly, the inclusions $(K'|_0,\hb K'|_0)\hookrightarrow (H_{K'|_0},\hb H_{K'|_0})\hookleftarrow (K'|_1,\hb K'|_1)$ are strongly $(n-1)$-connected. Thus, so is the map $(K'|_1,\hb K'|_1)\to (B,B^\partial).$ The inclusion $(K'|_1,\hb K'|_1)\hookrightarrow (V_1,\hb V_1)$ is strongly $(n-1)$-connected by the following argument: Consider the following (homotopy) pushout square of pairs
    \[\begin{tikzcd}
	{(K'|_1,\hb K'|_1)} & {(H_{K'|_0}, \hb H_{K'|_0})} \\
	{(V_1,\hb V_1)} & {(V_1\circ H_{K'|_0},\hb (V_1\circ H_{K'|_0})) }
	\arrow[from=1-1, to=1-2]
	\arrow[from=1-1, to=2-1]
	\arrow[from=1-2, to=2-2]
	\arrow[from=2-1, to=2-2]
\end{tikzcd}\]where both horizontal maps are strongly $(n-1)$-connected and the right vertical map is strongly $(n-1)$-connected since the composition $(K'|_0,\hb K'|_0)\to (H_{K'|_0},\hb H_{K'|_0}) \to (V_1\circ H_{K'|_0},\hb (V_1\circ H_{K'|_0})) \cong (K|_{[0,1]}\circ V_0,\hb (K|_{[0,1]}\circ V_0))$ is strongly $(n-1)$-connected. Thus, by \cref{strong conn and pushouts} and $n\geq 2,$ the leftmost vertical map is strongly $(n-1)$-connected. We can iterate this construction to obtain a weak map between $K'=\{H_{K'|_i}:K'|_i\leadsto K'|_{i+1}\}_{i\geq 0}$ and $K$ with the desired property.
\end{proof}
Before we prove \cref{final without null}, we return briefly to the difference between the categories $\smash{\CobbtLs}$ and $\smash{\Cobbtlo}$. Recall that the functor $(M\mapsto M^\circ):\smash{\Cobbtl}\to \smash{\CobbtL}$ is an isomorphism of topological categories, but under this map $\smash{\Cobbtls}$ does not map in general into $\smash{\CobbtLs}$. However, the next lemma gives a condition for that to happen. This will be necessary to use both \cref{surgery on morphisms} and \cref{stable stability}, which concern these two different categories. 

\begin{lemma}\label{cut or not cut is the same}
    The image of functor $(M\mapsto M^\circ):\smash{\Cobbtls}\to \smash{\CobbtL}$ is $\smash{\CobbtLs}$, provided the maps in $\vb L\to L$ and  $\hvb L\to \hb L$ induce an isomorphism on fundamental groupoids. In this case, the functor $(M\mapsto M^\circ):\smash{\Cobbtls}\to \smash{\CobbtLs}$ is an isomorphism of categories.
\end{lemma}
\begin{proof}
    Let $W:N\leadsto M$ be a morphism in $\Cobbtls$. Start by noticing that the assumption implies that the map $(M^\circ,\hb M^\circ)\to (M,\partial M) $ induces an isomorphism on fundamental groupoids: This follows from the left homotopy pushout square of pairs in the following diagram
    \[\begin{tikzcd}
	(\vb L,\hvb L) & (M^\circ,\hb M^\circ) \\
	    (L,\hb L) & (M,\partial M)
	\arrow[from=1-1, to=1-2]
	\arrow[from=1-1, to=2-1]
	\arrow[from=1-2, to=2-2]
	\arrow[from=2-1, to=2-2]
\end{tikzcd} \quad \quad \quad \begin{tikzcd}
	(M^\circ,\hb M^\circ) & (M,\partial M) \\
	    (W^\circ,\hb W^\circ) & (W,\hb W)
	\arrow[from=1-1, to=1-2]
	\arrow[from=1-1, to=2-1]
	\arrow[from=1-2, to=2-2]
	\arrow[from=2-1, to=2-2]
\end{tikzcd}\]and by Seifert-Van Kampen. Consider now the right homotopy pushout square of pairs above. Since the top map induces an isomorphism on fundamental groupoids, then by \cref{strong conn and pushouts} the right vertical map is strongly $(n-1)$-connected if and only if the leftmost is.
\end{proof}

We are now ready to finish the proof of the second preliminary version of \cref{final with null}.

\begin{proof}[Proof of \cref{final without null}.]
    Let $K'\leadsto K$ be a weak map such that $\ell|_{K'|_i}$ is strongly $(n-1)$-connected, which exists by \cref{truncation theta end}. Given $W\in \F(P,K'|_i),$ we have that $K'|_i\to W$ is strongly $(n-1)$-connected since $\ell|_{K'|_i}$ and $\ell|_W$ are strongly $(n-1)$ and $n$-connected respectively and $n\geq 2$. Thus, $\Fo(P,K'|_i)=\F(P,K'|_i)$ for every $i\geq 0.$ Consider the following commutative diagram
    \[\begin{tikzcd}
	{\F(P,K'|_\infty)} & & {\F(P,K|_\infty)}  \\
	{\Fo(P,K'|_{\infty})} \\
	{\Omega_{[P,K'|_\infty]}\B\CobbtLo} & {\Omega_{[P,K|_\infty]}\B\CobbtLo} & {\Omega_{[P,K|_\infty]}\B\CobbtL}.
	\arrow[from=1-1, to=1-3]
	\arrow[Rightarrow, no head, from=1-1, to=2-1]
	\arrow[from=1-3, to=3-3]
	\arrow[from=2-1, to=3-1]
	\arrow[from=3-1, to=3-2]
        \arrow[from=3-2, to=3-3]
\end{tikzcd}\]By \cref{highly connected theta ends}, the left vertical map is an acyclic and by \cref{theta ends and maps} the top map is an abelian homology equivalence. The bottom right map is an equivalence by \cref{surgery on morphisms} and \cref{cut or not cut is the same}, since $(L,\hb L)$ is obtained by handles of index at most $n-1$ from $\emptyset$ and the maps $\vb L\to L$ and $\hvb L\to \hb L$ are fundamental groupoid isomorphisms. Since the bottom map is an equivalence, we see that the right vertical map is an abelian homology equivalence and thus, acyclic since the fundamental group of the target is abelian. 
\end{proof}

\subsection{\texorpdfstring{Proof of \cref{final with null}.}{Proof of main statement.}} In this subsection, we prove the main result of this section using the preliminary versions established above. Until now, we have proved all our statements for general classes of submanifolds $L$. We will now make use of this additional generality to choose a specific $L$, which depends on the $\Theta$-end $K$ in the hypothesis of \cref{final with null}. We start by constructing this $L.$

\begin{cons}\label{the L}
    Let $K$ be a $\Theta$-end in $\Cobbt$ such that $\partial K|_0\neq \emptyset$ and pick a triad handle decomposition of $(K|_0,\partial K|_0,\emptyset)$ relative to $\emptyset$ with only left and interior handles, which is possible by \cref{existence of handle decomposition}. Let $L$ be the union of all handles of index at most $n-1$. This is a $\Theta$-triad by restricting. By isotoping the embedding $(K|_0,\partial K|_0)\hookrightarrow (\RR^{\infty}_+,\partial \RR^{\infty}_+),$ we can assume that $L=K|_0\cap ([0,+\infty)\times (-\infty,0]\times \RR^{\infty-2}).$ Therefore, we can assume $K|_0\in \Cobbtl.$ Note that, $(L,\hb L)$ is built from $(\vb L,\hvb L)$ by attaching right and interior handles of index at least $n+1$. Since $n\geq 2$, we see from \ref{effect of right handles} and \ref{effect of interior handles} from \cref{strong geo section} that $\hvb L\to \hb L$ and $\vb L\to L$ are fundamental groupoid isomorphisms. Hence, $(L,\hb L)$ satisfies the hypothesis of \cref{final without null}, since $\hvb L\neq \emptyset$ as $\partial K|_0\neq \emptyset$.
\end{cons}

The following lemma relates the notion of a $\Theta$-end in $\Cobbt$ to a $\Theta$-end in $\CobbtL$. This will allow us to use \cref{final without null} for our specific choice of $L$ to deduce our main result.

\begin{lemma}
    The embedding $K\hookrightarrow [0,\infty)\times [0,1)\times (-1,1)^{\infty-2}  $ may be isotoped, along with a bundle homotopy of its $\Theta$-structure, such that \(K\cap ([0,\infty)\times [0,1)\times (-\infty,0]\times \RR^{\infty-2})=[0,\infty)\times L   \) relative to $K|_0$ as a $\Theta$-manifold. Moreover, $K^\circ$ is a $\Theta$-end in $\smash{\CobbtLs}.$
\end{lemma}
\begin{proof}
    We start by showing that any embedding $\{i\}\times L\hookrightarrow K|_i$ can be extended to an embedding $[i,i+1]\times L\hookrightarrow K|_{[i,i+1]}.$ By \cref{strong geometrical connectivity}, we can find a triad handle decomposition of $K|_{[i,i+1]}$ relative to $K|_{i+1}$ with right handles of index at least $n$ and interior handles of index at least $n+1$. This induces a triad handle decomposition of $K|_{[i,i+1]}$ relative to $K|_i$ with left handles of index at most $n+1$ and interior handles of index at most $n.$ One can extend the embedding of $L$ if its image can be made disjoint from the attaching maps of such handles. It suffices to make the cores of $L$ disjoint from such attaching maps (since $L$ is isotopy equivalent to the union of its cores). The attaching maps of the left and interior handles are submanifolds of dimension at most $n$ and the cores of all handles in $\hb L$ have dimension at most $n-1.$ Since $n+(n-1)<2n,$ we can isotope the attaching maps of the left and interior handles in $K|_i$ to be disjoint from all cores in $\hb L$ by transversality. By iterating this procedure, we have an embedding of $e:[0,\infty)\times L\hookrightarrow K$ relative to $K|_0.$ It suffices to prove that we can isotope $K$ in $[0,\infty)\times [0,1)\times (-1,1)^{\infty-2}$ to have the claimed property. This follows from the isotopy extension theorem: We can isotope (as maps of pairs) the image of $[0,\infty)\times L\hookrightarrow K \hookrightarrow [0,\infty)\times [0,1)\times (-1,1)^{\infty-2}$ to be the inclusion (recall that $[0,\infty)\times L$ is a submanifold of the target). By the isotopy extension theorem, we can extend such an isotopy to an isotopy of $K.$ We can improve such an isotopy to be one of $\Theta$-triads by choosing a bundle homotopy of the $\Theta$-structure $e^*\ell|_K$ to $\ell|_{[0,\infty)\times L}$ (which exists since both are just extensions of $\ell|_L$). This is possible, since the embedding of $L$ is a Hurewicz cofibration. We finish by proving that $K^\circ =K\backslash [0,\infty)\times L$ is a $\Theta$-end in $\CobbtLo.$ The first condition is verified by considering the following pushout square
    \[\begin{tikzcd}
	{(M^\circ,\hb M^\circ)} & {(M,\partial M)} \\
	{(W^\circ,\hb W^\circ)} & (W,\hb W)
	\arrow[from=1-1, to=1-2]
	\arrow[from=1-1, to=2-1]
	\arrow[from=1-2, to=2-2]
	\arrow[from=2-1, to=2-2]
\end{tikzcd}\]for any morphism $W:M\leadsto M'$ in $\Cobbtlo.$ The top map is a pushout of the map $(\vb L,\hvb L)\to (L,\hb L)$ (see \cref{cut or not cut is the same}), which is an isomorphism on fundamental groupoids by assumption, and thus an isomorphism on fundamental groupoids, by the Seifert--Van Kampen theorem. Thus, by \cref{strong conn and pushouts}, if the right vertical map is strongly $(n-1)$-connected than so is the left one. Once again by transversality, the embedded copy of $V_1$ in $K|_{[i,i+1]}$ (which can be seen as trivially attached pair of left $(n+1)$ and right $n$-handle) can be made disjoint from $[i,i+1]\times L$ for each $i\geq 0.$ This verifies the second condition and finishes the proof.  
\end{proof}

We are now equipped for the proof of the main result of this section.

\begin{proof}[Proof of \cref{final with null}.]
    Pick $L$ as in \cref{the L}. An analogous contruction to \cite[Lemma 2.17]{GRWII} produces an equivalence $\F(\overline{L},K|_\infty^\circ)\to \Ntn(K|_\infty)$ (and thus the former is non-empty) and a commutative diagram
    \[\begin{tikzcd}
	{\F(\overline{L},K|_\infty^\circ)} & {\Omega_{[\overline{L},K|_\infty^\circ]}\B\CobbtL} \\
	{\Ntn(K|_\infty)} & {\Omega_{[\emptyset,K|_\infty]}\B\Cobbt}
	\arrow[from=1-1, to=1-2]
	\arrow[from=1-1, to=2-1]
	\arrow[from=1-2, to=2-2]
	\arrow[from=2-1, to=2-2]
	\arrow[from=2-1, to=2-2]
\end{tikzcd}.\]The proof of this statement follows analogously to loc.cit. Here $\overline{L}$ denotes the reflection of $L$ along the hyperplane $[0,1)\times \{0\}\times (-\infty,\infty)^{\infty-2}$ (recall the definition of $L$ above \cref{cob with L}). The vertical compositions are equivalences by once again the analog of \cite[Lemma 2.17]{GRWII} and \cref{have or not L in BCOb}. The top map is acyclic by \cref{final without null}, thus so is the bottom map.
\end{proof}

\subsection{Stable diffeomorphism classification.}\label{proof of kreck section}

In this subsection, we prove \cref{stable class thm} on the classification of manifold triads up to stable diffeomorphism using \cref{final with null}. This result is akin to Kreck's seminal work \cite{Kreck1999SurgeryAD} on stable diffeomorphism classification of even-dimensional manifolds with boundary. We now recall the setting of \cref{stable class thm}: Let $n\geq 3$ be an integer, $P$ be a compact $2n$-manifold with boundary, and $(N_i,\hb N_i,\vb N_i)$ be $(2n+1)$-dimensional manifold triads (see \cref{section of pairs of manifolds}) for $i=0,1$, together with an identification of $\vb N_i$ with $P$. Assume also that $\hb N_i$ is connected for $i=0,1.$

\begin{proof}[Proof of \cref{stable class thm}]
    We start with the "only if"-direction, namely, we assume that $N_0$ and $N_1$ are stably diffeomorphic relative to $P$ and fix $g\geq 0$ and a diffeomorphism $\smash{\hat{f}}$ between the triads $N_i\natural V_g$ for $i=0,1$ as in \ref{stably diffeo defn} from the introduction. The statement about relative Euler characteristics follows from the equality $\chi(N_i\natural V_g,\hb N_i\sharp W_g)= \chi(N_i,\hb N_i)+(-1)^{n+1}g$. We prove now the second statement. To define a stable normal $n$-type, take a Moore-Postnikov $n$-factorization of maps of pairs $\nu_0=\Theta^\perp\circ \ell_0:(N_0,\hb N_0)\to (B,B^\partial)\to (\BO,\BO)$ in the sense of the introduction (see also \cref{def moore} below). By definition, the map $\ell_0$ is strongly $n$-connected, so it suffices to produce a strongly $n$-connected map $\ell_1:(N_1,\hb N_1)\to (B,B^\partial)$ to conclude that $N_0$ and $N_1$ have the same stable normal $n$-type. Take a Moore-Postnikov factorization of pairs $\nu_{N_0\natural V_g}=\Theta_g^\perp\circ l^g:(N_0\natural V_g,\hb N_0\sharp W_g)\to (B_g,B^\partial_g)\to (\BO,\BO)$. One can check that there is an induced homotopy commutative diagram 
    \[\begin{tikzcd}
        (N_0,\hb N_0) \arrow[r, "\ell_0"]\arrow[d] & (B,B^\partial) \arrow[r, "\Theta^\perp"]\arrow[d, dashed, "b"]  & (\BO,\BO) \arrow[d, "\id"]\\
        (N_0\natural V_g,\hb N_0\sharp W_g) \arrow[r, "\ell^g"] & (B_g,B^\partial_g) \arrow[r, "\Theta_g^\perp"]  & (\BO,\BO)
    \end{tikzcd}.\]We claim that $b$ is a weak equivalence, which can be seen from the fact that the map $(*,*)\to (V_g,W_{g,1})$ is strongly $(n-1)$-connected, and that  the pair $(V_g,W_{g,1})$ is stably parallelizable and thus the map of pairs to $(\BO,\BO)$ classifying its stable normal bundle is nullhomotopic. Thus, without loss of generality, we can assume the map $\ell_0$ factors through $(N_0\natural V_g,\hb N_0\sharp W_g).$ We define now $\ell_1$ to be the composite of the inclusion $(N_1,\hb N_1)\hookrightarrow (N_1\natural V_g,\hb N_1\sharp W_g)$, the diffeomorphism $\smash{\hat{f}^{-1}}:(N_1\natural V_g,\hb N_1\sharp W_g)\cong (N_0\natural V_g,\hb N_0\sharp W_g)$ and the map $\ell^g$. It remains to check that $\ell_1$ is strongly $n$-connected. This can be seen again by naturality of Moore-Postnikov factorizations of pairs and the fact that the map $(*,*)\to (V_g,W_{g,1})$ is strongly $(n-1)$-connected and that $(V_g,W_{g,1})$ is stably parallelizable. We leave this check to the reader. We conclude that $N_0$ and $N_1$ have the same stable normal $n$-type $\Theta^\perp$. Finally, to see that $N_0$ and $N_1$ admit bordant $\Theta^\perp$-smoothings, we will construct a $\Theta^\perp$-nullbordism of $(N_0\cup_{P} N_1,\ell_0\cup (-\ell_1))$. This is inspired by \cite[Lemma 2.2]{Crowley2011}. Consider the $(2n+2)$-dimensional triad cobordism (recall \cref{relative cobordism}) $M_i$ obtained from $N_i\times [0,1]$ by attaching $g$-many trivial left $(n+1)$-handles (recall \cref{handles}), then its vertical boundary is diffeomorphic to the disjoint union of $N_i$ and $N_i\natural V_g:$ This can be seen by observing that $(V_1,W_{1,1})$ can be obtained from $(D^{2n+1}_+,\partial_0D^{2n+1}_+)$ by removing a neighborhood of a $(\partial_1D^{n+1}_+,\partial_{01} D^{n+1}_+)$ and gluing $(D^{n+1}_+\times S^n,\partial_0 D^{n+1}_+).$ Let $\smash{\hat{M}}$ be the triad bordism given by gluing $M_0$ and $M_1$ along $\smash{\hat{f}}:N_0\natural V_g\to N_1\natural V_g$. This is a triad cobordism from $N_0$ to $N_1$, which we can see as triad whose vertical boundary is $\smash{N_0\cup \partial_1\hat{M}\cup N_1}$, which is diffeomorphic to $N_0\cup_P N_1$. Thus, we can see $\smash{\hat{M}}$ as a nullbordism of $N_0\cup_PN_1.$ We proceed now to extend $\ell_0$ to a $\Theta^\perp$-structure $\smash{\hat{\ell}}$ on $\smash{\hat{M}}$ such that $\smash{\hat{\ell}|_{N_1}}$ is strongly $n$-connected. As the left handles of $M_0$ relative to $N_0$ are attached trivially, we can find a $\Theta^\perp$-structure $\ell_{M_0}$ on $M_0$ by extending $\ell_0$ trivially along the core of the handle, which is possible as the core is parallelizable. We must now extend this structure from $M_0$ to $\smash{\hat{M}}$. In other words, we must solve the following lifting problem of maps of pairs
    \[\begin{tikzcd}
        (M_0,\hb M_0) \arrow[d] \arrow[r, "\ell_{M_0}"] & (B,B^\partial) \arrow[d, "\Theta^\perp"]\\
        (\hat{M},\hb \hat{M}) \arrow[r, "\nu_{\hat{M}}"]\arrow[ru, dashed, "\hat{\ell}"] & (\BO,\BO)
    \end{tikzcd}.\]To do so, observe that $\smash{(\hat{M},\hb \hat{M})}$ is obtained from $(M_0,\hb M_0)$ by attaching $g$-many right $(n+1)$-handles, since $M_1$ is obtained from $N_1$ by attaching $g$-many left $(n+1)$-handles, by \cref{dual handle}. In particular, we have $\smash{\hat{M}}$ is equivalent to the union $\smash{M_0\cup_{\hb M_0}\hb \hat{M}}$ (see \ref{effect of right handles} from \cref{strong geo section}), and so it suffices to solve th lifting problem of maps of spaces given by restriction of the problem above of pairs to the sources of the pairs involved. By definition, the map $(\theta^\perp)^\partial:B^\partial\to\BO $ is $n$-coconnected. Since $\hb M_0\to \smash{\hb \hat{M}}$ is an $n$-connected map of CW complexes, we see that such a lift exists by obstruction theory. Moreover by \ref{effect of right handles} and \ref{effect of left handles} from \cref{strong geo section}, the map $(N_i,\hb N_i)\hookrightarrow(\smash{\hat{M},\hb \hat{M}})$ is strongly $(n-1)$-connected. We show now that $\smash{\hat{\ell}}|_{N_1}$ is strongly $n$-connected. Since $\ell_0$ is strongly $n$-connected, we deduce that $\smash{\hat{\ell}}$ is strongly $n$-connected by \cref{closure properties strongly}. We conclude that $\smash{\hat{\ell}}|_{M_1}$ is also strongly $n$-connected, since the map $(M_1,\hb M_1)\to (\smash{\hat{M},\hb \hat{M}})$ is strongly $n$-connected. On the other hand, we have $M_1\simeq N_1$ and $\hb M_1\simeq N_1\vee (S^{n})^{\vee g}$ by \ref{effect of left handles} in \cref{strong geo section}. Since the map $B^\partial\to \BO$ is $n$-coconnected and $\nu_1|_{(S^n)^{\vee g}}$ is nullhomotopic, we deduce $\smash{\hat{\ell}}|_{(S^n)^{\vee g}}$ is nullhomotopic. Thus, we conclude that $\partial N_1\to B^\partial$ is $n$-connected. It remains to show that $N_1\cup_{\partial N_1} B^\partial\to B$ is $(n+1)$-connected. Since $\smash{\hat{\ell}}|_{M_1}$ is also strongly $n$-connected, we see that $N_1\cup_{\partial N_1\vee (S^n)^{\vee g}} B^\partial\to B$ is $(n+1)$-connected. However, the target of the latter map is equivalent to $(N_1\cup_{\partial N_1} B^\partial)\vee (S^{n+1})^{\vee g}$ and the restriction of the map to $(S^{n+1})^{\vee g}$ is nullhomotopic since it is once post-composed with $\theta^\perp$, which is an $n$-coconnected map. We conclude that $N_1\cup_{\partial N_1} B^\partial\to B$ is $(n+1)$-connected, hence finishing the proof that $N_0$ and $N_1$ admit bordant $\Theta^\perp$-smoothings. This establishes the "only if"-direction.


    We move now to the "if"-direction. We start by proving the case $\hvb N_i\neq \emptyset$ and later deduce the general case from the special case. Start by choosing maps of pairs $\Theta^\perp:(B,B^\partial)\to (\BO,\BO)$ as in \smash{\ref{same normal type}} from the introduction, a strongly $n$-connected $\ell_i:(N_i,\hb N_i)\to (B,B^\partial)$ lift of $\nu_i$ as in \smash{\ref{bordant normal smoothings}} from the introduction such that $N_0\cup_P (-N_1)$ is $\Theta^\perp$-nullbordant. Fix an embedding of $(P,\partial P)$ into $(\RR^\infty_+,\partial \RR^\infty_+)$. Let $\Theta:(B',(B')^\partial)\to (\BO(2n+1),\BO(2n))$ be the map given by pulling back $\iota\circ \Theta^\perp$ along the stabilization map $(\BO(2n+1),\BO(2n))\to (\BO,\BO),$ where $\iota: \BO\to \BO$ is the map induced by taking a matrix to its inverse. Observe that $\iota\circ \nu_i$ is the classifying map for the stable tangent bundle of $(N_i,\hb N_i).$ Thus, $\ell_i$ lifts uniquely (up to homotopy) along $(B',(B')^\partial)\to (B,B^\partial)$ to a map $\ell'_i:(N_i,\hb N_i)\to (B',(B')^\partial),$ as the stable tangent classifier lifts to the unstable tangent classifier. Observe that $\ell'_i$ is also strongly $n$-connected since the map $(\BO(2n+1),\BO(2n))\to (\BO,\BO)$ is strongly $(2n-1)$-connected and $\iota$ is an equivalence. Therefore, by choosing embeddings of $N_0$ and $N_1$ into $(-\infty,0]\times \RR^\infty_+$ extending the embedding of $(P,\partial P)$ (using the fixed identifications), we can consider $(N_i,\ell'_i)\in \Ntn(P)$ as in \cref{final with null}. Consider a $\Theta$-end $K$ in the sense of \cref{theta end without l} given by taking connect sums with $V_1$ starting in $P\times [0,1]$. In particular, $K|_0=P$. Then, one can check that $N_0$ and $N_1$ are stably diffeomorphic if the classes $[N_0,\ell'_0]$ and $[N_1,\ell'_1]$ in $\pi_0(\Ntn(K|_\infty))$ agree. We proceed now to check the latter condition. By \cref{final with null}, this condition is equivalent to these classes agreeing in $\pi_0(\Omega_{[\emptyset, \vb N_0]}\B\Cobbt)$ after applying the map present in loc.cit. Observe that we can consider $(N_1,-\ell'_1)\in\Cobbt(\vb N_0,\emptyset)$ for the fixed map $\ell'_1$ as in \smash{\ref{bordant normal smoothings}} from the introduction, and concatenation by the path induced by $(N_1,-\ell'_1)$ in $\B\Cobbt$ induces an isomorphism $\pi_0(\Omega_{[\emptyset, P]}\B\Cobbt)\cong \pi_0(\Omega_{[\emptyset, \emptyset]}\B\Cobbt)$. By applying \cite[Main Cor. 4.6]{genauer} to $\pi_1$, there is an isomorphism \(\pi_0(\Omega_{[\emptyset, \emptyset]}\B\Cobbt)\cong \pi_0(\MTH)\) where $\MTH$ is the cofiber of the canonical map of Thom spectra $\Sigma^{-1} \MThb\to \MTc$ (recall the definition from the introduction). On the other hand, there is a stabilization map $s:\MTH\to \Sigma^{-2n-1}\MM\Theta^\perp$ (see proof of \cref{the gollinger lemma} below for more details). One can check that the composite
    \[\pi_0(\Ntn(K|_\infty))\overset{\cong}{\to} \pi_0(\Omega_{[\emptyset, \vb N_0]}\B\Cobbt)\overset{\cong}{\to} \pi_0(\Omega_{[\emptyset, \emptyset]}\B\Cobbt)\overset{\cong}{\to} \pi_0(\MTH)\overset{\pi_{0}(s)}{\to} \pi_{2n+1}(\MM\Theta^\perp)\cong \Omega_{2n+1}^{\Theta^\perp}   \]takes the class $[N_i,\ell_i']$ to the relative bordism class of $[N_i\cup_P N_1,\ell_i'\cup (-\ell'_1)].$ By hypothesis, the image of $N_0$ under this composite vanishes. On the other hand, the image of $N_1$ is $N_1\cup_P (-N_1)$, which vanishes by definition (by taking $N_1\times [0,1]$). Hence, the images of $N_0$ and $N_1$ under this composite agree. We conclude that the classes $[N_0,\ell'_0]$ and $[N_1,\ell'_1]$ agree on $\pi_0(\Ntn(K|_\infty))$ by the following fact (which we prove separately below): The map $\pi_{0}(s)\times \chi^{\text{rel}}:\pi_0(\MTH)\to \Omega_{2n+1}^{\Theta^\perp}\times \ZZ$ is injective, where $\chi^{\text{rel}}$ takes a $\Theta$-bordism class of manifolds with boundary to its relative Euler characteristic (one can check that the relative Euler characteristic is bordism invariant). This finishes the proof of the "if"-direction for the case $\hvb N_i\neq \emptyset$. 
    
    We finish by deducing the general case from this special case. Assume $\hvb N_i=\emptyset$ and let $N_i^\circ$ be the triad $(N_i, \partial N_i\backslash D^\circ_i, D_i)$ where $D_i$ is a codimension $0$ disc in $\partial N_i$. Observe that the hypothesis implies that $\chi(N_0,\hb N_0^\circ)=\chi(N_1,\hb N_0^\circ)$, since $\chi(N_i,\hb N_i^\circ)=\chi(N_i,\partial N_i)+1.$ Similarly, the pairs $(N_i,\partial N_i)$ and $(N_i,\hb N_i^\circ)$ have the same stable normal $n$-type $\Theta^\perp$, as the inclusion $\hb N_i^\circ \hookrightarrow \partial N_i$ is $(2n-1)$-connected. Moreover, we know that $N_0 $ and $N_1$ admit bordant $\Theta$-smoothings, so there exists $\ell_i$ and $W$ as in \smash{\ref{bordant normal smoothings}} from the introduction (which now has empty vertical boundary) such that $N_0$ and $N_1 $ agree in $\Omega^{\Theta^\perp}_{2n+1}$. Let $W$ be a $\Theta^\perp$-bordism between $N_0$ and $N_1$. We can assume that there is a path in $\hb W$ between a point in $D_0$ and a point in $D_1$ by the following argument: since $\ell_i$ is strongly $n$-connected and $\hb N_i$ is connected, we know that $\ell_i|_{\hb N_i}$ hit the same path component of $B^\partial.$ Take the triad $W$ and attach a right $1$-handle along the horizontal boundary in the path components hit by $D_0$ and $D_1$. Denote the result of this procedure by $W'$. Since $\hb N_0$ and $\hb N_1$ hit the same component of $B^\partial$, we can extend the $\Theta^\perp$-structure from $W$ to $W'$, making the latter into a $\Theta^\perp$-bordism between $N_0$ and $N_1$, which satisfies this extra assumption. One can consider a thickening of that path to define a $\Theta^\perp$-nullbordism of $N_0^\circ\cup_DN_1^\circ$. We leave this step to the reader. Thus, we see that $N_0^\circ\cup_DN_1^\circ$ vanishes in this bordism group. Hence, $N_0^\circ$ and $N_1^\circ$ admit bordant $\Theta^\perp$-smoothings. Thus, by the special case, we see that $N_0^\circ$ and $N_1^\circ$ are stably diffeomorphic. Since this diffeomorphism extends the chosen identification of their vertical boundaries, we see that it extends to a stable diffeomorphism between $N_0$ and $N_1.$ Thus, we have established the general case and hence finished the proof.
    \end{proof}

    \begin{lemma}\label{the gollinger lemma}
        The map $\pi_{0}(s)\times \chi^{\text{rel}}:\pi_0(\MTH)\to \Omega_{2n+1}^{\Theta^\perp}\times \ZZ$ is injective.
    \end{lemma}
    \begin{proof}
        This is inspired by \cite[Lemma 1.4]{KrannichExotic}. Let us start by defining the map $s$ in detail. Recall that we have the map $\Theta=(\theta,\theta^\partial):(B',(B')^\partial)\to (\BO(2n+1),\BO(2n))$ which is obtained by pulling back $\iota\circ \Theta^\perp$ along the stabilization map $(\BO(2n+1),\BO(2n))\to (\BO,\BO)$. Thus, the map $B^\partial\to \BO(2n+1)\to \BO$ is the classifying map of the stable vector bundle $-(\theta^\partial)^*\gamma$ of rank $-2n$. Applying the Thom spectrum functor to this classifying map produces a map $\MThb\to \Sigma^{-2n}\MM(\theta^\perp)^\partial$, where $(\theta^\perp)^\partial\coloneqq \Theta^\perp|_{B^\partial}.$ Similarly, we have a map $\MTc\to \Sigma^{-2n-1}\MM\theta^\perp$ given by the Thomified classifying map of the stable vector bundle $-\theta^*\gamma$. The map $s:\MTH\to \Sigma^{-2n-1}\MM\Theta^\perp$ is defined as the induced map on cofibers between the maps $\Sigma^{-1}\MThb\to \MTc$ and $\Sigma^{-2n-1}\MM(\theta^\perp)^\partial\to \Sigma^{-2n-1}\MM\theta^\perp.$ Let $\Theta_{2n+2}=(\theta_{2n+2},\theta^\partial_{2n+2}):(B'',(B'')^\partial)\to (\BO(2n+2),\BO(2n+1))$ be the map obtained by pulling back $\iota\circ \Theta^\perp$ along the stabilization map $(\BO(2n+2),\BO(2n+1))\to (\BO,\BO)$. Then, clearly the map $s$ factors as the composite $\MTH\to \Sigma\MTH_{2n+2}\to \Sigma^{-2n-1}\MM\Theta^\perp.$ The map $\Sigma\MTH_{2n+2}\to \Sigma^{-2n-1}\MM\Theta^\perp$ can be seen to be $1$-connected and thus an isomorphism on $\pi_0$ (this follows since the maps $\MThb_{2n+1}\to \Sigma^{-2n-1}\MM(\theta^\perp)^\partial$ and $\MTc_{2n+2}\to \Sigma^{-2n-2}\MM(\theta^\perp)^\partial$ are $0$-connected, see \cite[5]{KrannichExotic}). Thus, it suffices to prove that the map $\pi_{0}(s)\times \chi^{\text{rel}}:\pi_0(\MTH)\to \pi_{-1}(\MTH_{2n+2})\times \ZZ$ is injective. We shall now study the fiber  $\MTx\coloneqq \fib(\MTH\to \Sigma\MTH_{2n+2}).$ The fibers of the maps $\MTc\to \Sigma\MTc_{2n+2}$ and $\MThb\to \Sigma\MThb_{2n+1}$ are equivalent to $\Sigma^\infty_+B''$ and $\Sigma^{\infty}_+(B'')^\partial$, respectively (see \cite[Prop. 3.0.7]{gollinger}). Thus, we have a cofiber sequence of spectra $\Sigma^\infty_+B''\to \MTx\to \Sigma^\infty_+(B'')^\partial.$ We conclude that there is an isomorphism $\pi_0(\MTx)\cong \pi_0(\Sigma^\infty_+ (B'')^\partial)\oplus \im(\pi_0(\Sigma^\infty_+B'')\to \pi_0(\MTx))$. One checks that $B''$ and $(B'')^\partial$ are path-connected and thus $\pi_0(\Sigma^\infty_+B'')\cong \ZZ\cong \pi_0(\Sigma^\infty_+(B'')^\partial)$. We proceed by showing the following claims:
        \begin{enumerate}[label=(\textit{\alph*})]
            \item\label{gollinger input} \textit{The composite $\pi_0(\Sigma^\infty_+B'')\to \pi_0(\MTx)\to \pi_0(\MTH)$ is zero}: We start by producing a splitting of the cofiber sequence $\Sigma^\infty_+B''\to \MTx\to \Sigma^\infty_+(B'')^\partial$ at the level of $\pi_0$. To do so, we start by observing that the group $\pi_0(\MTx)$ has the following geometric description: It is isomorphic to the bordism group of triples $[M,E,\phi]$ where:
            \begin{enumerate}[label=\arabic*.]
                \item $M=(M,\hb M,\vb M)$ is a $(2n+2)$-dimensional triad (in the sense of \cref{section of pairs of manifolds});
                \item $E=(E,E^\partial)$ is a $(2n+2)$-dimensional collared vector bundle pair over $(M,\hb M)$ together with a collared bundle map $(E,E^\partial)\to ((\theta_{2n+2})^*\gamma_{2n+2},(\theta_{2n+2}^\partial)^*\gamma_{2n+1})$  (in the sense of \cref{collared pairs}):
                \item $\phi$ is a stable isomorphism of stable vector bundle pairs $(TM,T\hb M)\cong (E,E^\partial)$.
            \end{enumerate}Such a triple is \textit{nullbordant} if there exists a $(2n+3)$-dimensional $4$-ad $(W,\partial_0 W,\partial_1 W,\partial_2 W)$ (recall \cref{relative cobordism}) such that $(\partial_2 W,\partial_{02} W,\partial_{12} W)=(M,\hb M,\vb M), $ a $(2n+2)$-dimensional collared vector bundle pair $(F,F^\partial)$ over $(W,\partial_0 W)$ which restricts to $(E,E^\partial)$ on $(M,\hb M),$ and a stable isomorphism $(TW,T\partial_0 W)\cong (F,F^\partial)\oplus \varepsilon^1$ compatible with $\phi.$ One sees this by combining \cite[Thm. 3.1.5]{laures} with the observation that $\MTx$ is equivalent to the total cofiber of the following commutative square of spectra
            \[\begin{tikzcd}
                \Sigma^{-2}\MThb \arrow[d] \arrow[r] & \Sigma^{-1}\MThb_{2n+2} \arrow[d]\\
                \Sigma^{-1}\MTc \arrow[r] & \MTc_{2n+2}
            \end{tikzcd}\]On the other hand, as $\Sigma_+^\infty B''$ is equivalent to the cofiber of $\Sigma^{-1} \MTc\to \MTc_{2n+2}$, we see that $\pi_0(\Sigma^\infty_+B'')$ is isomorphic to the bordism group of triples $[M',E',\phi']$ where $M'$ is a $(2n+2)$-dimensional compact manifold with boundary, $E'$ is a $(2n+2)$-dimensional vector bundle over $M$ together with a bundle map $E'\to (\theta_{2n+2})^*\gamma_{2n+2},$ and $\phi':TM\cong E'$ is a stable isomorphism of bundles. Under these geometric desciptions, the map $\pi_0(\Sigma^\infty_+B'')\to \pi_0(\MTx)$ takes $[M',E',\phi']$ to the triple $[(M,\emptyset, \partial M), (E',\emptyset),\phi']$. This map admits a splitting, namely the map which takes a triple $[M,E,\phi]\in \pi_0(\MTx)$ to the triple $[M, E|_M,\phi|_M],$ where one forgets the triad structure on $M$. We conclude that we have an isomorphism $\pi_0(\MTx)\cong \pi_0(\Sigma^\infty_+B'')\oplus \pi_0(\Sigma^\infty_+(B'')^\partial).$ We now return to proving the initial claim. For that, we show that a generator $a\in \pi_0(\Sigma^\infty_+B'')\subset \pi_0(\MTx)$ is hit by the map $\pi_0(\MTH_{2n+2})\to \pi_0(\MTx)$, and hence maps to zero in $\pi_0(\MTH)$ by exactness. Let $x\in \pi_0(\MTH_{2n+2})$ be a class representing the triple $[D^{2n+2},TD^{2n+2},\phi]$ where $(TD^{2n+2},T\partial D^{2n+2})$ is seen with some $\Theta_{2n+2}$-structure. First, we show that $x$ maps to $a+b\in \pi_0(\MTx)$ where $a\in \pi_0(\Sigma^\infty_+ B'')$ is a generator and $b\in \pi_0(\Sigma^\infty_+(B'')^\partial):$ Consider the map $\chi^{\text{rel}}:\pi_0(\Sigma^\infty_+B'')\to \ZZ$ which takes $[M',E',\phi']$ to $\chi(M',\partial M')$ and observe that it is surjective and thus an isomorphism, as $B''$ is path-connected. However, the image $a$ of $x$ in $\pi_0(\Sigma^\infty_+B'')$ hits $1$ along $\chi^{\text{rel}}.$ Hence, $a$ is a generator of $\pi_0(\Sigma^\infty_+B'').$ Second, we prove that the image $b\in \pi_0(\Sigma^\infty_+(B'')^\partial)$ vanishes: Observe that the map $\pi_0(\Sigma^\infty_+(B'')^\partial)\to \pi_0(\MThb)$ is injective by the proof of \cite[Lemma 1.4]{KrannichExotic} (see also below in the proof of $(b)$). However, the class $b$ maps to $0$ along the latter map as it is in the image of $\pi_0(\MThb_{2n+1})$, by definition. Thus, $b$ vanishes. We conclude that $x\in \pi_0(\MTH_{2n+2})$ hits a generator of the summand $\pi_0(\Sigma^\infty_+B'')$ in $\pi_0(\MTx)$ and hence the map $\pi_0(\Sigma^\infty_+B'')\to \pi_0(\MTx)\to \pi_0(\MTH)$ is zero.   
            \item\label{krannich input} \emph{The composite $\ZZ\cong \pi_0(\Sigma^\infty_+(B'')^\partial)\hookrightarrow \pi_0(\MTx)\to \pi_0(\MTH)\overset{\chi^{\mathrm{rel}}}{\to} \ZZ $ is an isomorphism}: Consider the following diagram
            \[\begin{tikzcd}
                \pi_0(\MTx)\arrow[d]\arrow[r] & \pi_0(\MTH)\arrow[d] \arrow[r, "\chi^{\text{rel}}"] & \ZZ \arrow[d, "2\cdot"] \\
                \ZZ\cong \pi_0(\Sigma^\infty_+(B'')^\partial) \arrow[r] & \pi_0(\MThb) \arrow[r, "\chi"] & \ZZ
            \end{tikzcd}\]which commutes, since given a $(2n+1)$-manifold with boundary $M$, then $\chi(\partial M)=2\chi(M,\partial M).$ By the proof of \cite[Lemma 1.4]{KrannichExotic}, we can observe that the bottom composition is the map $2\cdot (-):\ZZ\to \ZZ$. Thus, we conclude that the generator of inclusion of the summand $\pi_0(\Sigma^\infty_+(B'')^\partial)$ hits $1\in \ZZ$ along the top composite, since its image along the right vertical map is $2$. This finishes the claim.
        \end{enumerate}Using \ref{gollinger input} and \ref{krannich input}, we finish the proof of this claim. Recall that it suffices to prove that $\pi_{0}(s)\times \chi^{\text{rel}}:\pi_0(\MTH)\to \pi_{-1}(\MTH_{2n+2})\times \ZZ$ is injective. By exactness, the kernel of $\pi_0(s)$ is the image of the map $\pi_0(\MTx)\to \pi_0(\MTH)$. By \ref{gollinger input}, this image equals the image of the summand $\pi_0(\Sigma^\infty_+(B'')^\partial)$. By \ref{krannich input}, this image is detected by $\chi^{\text{rel}}$. Hence, we conclude that $\pi_0(s)\times \chi^{\text{rel}}$ is injective.
    \end{proof}

    \begin{rmk}\label{grw and kreck}
        \cref{stable class thm} should be seen as an analog of \cite[Thm. C]{Kreck1999SurgeryAD} for odd-dimensional triads. Although in loc.cit. only the "if"-direction is stated, the converse direction holds as well (see \cite[Lemma 2.2]{Crowley2011}). A similar proof as above recovers \cite[Thm. C]{Kreck1999SurgeryAD} and its converse from the work of Galatius and Randal-Williams \cite[Thm. 1.5]{GRWII} (see \cite[Lemma 1.4]{KrannichExotic} for the analog of \cref{the gollinger lemma}).
    \end{rmk}

\section{General tangential structures.}\label{general}

In this section, we state and prove our main theorem. We start by defining the main object of study. Let $n\geq 3$ and $(N,\hb N,\vb N)$ be a $(2n+1)$-dimensional manifold triad (see \cref{section of pairs of manifolds}) and $\Theta=(\theta,\theta^\partial):(B,B^\partial)\to (\BO(2n+1),\BO(2n))$ be a map of pairs. Recall that we fix an inwards-pointing vector field on $\hb N$ inducing a collar on the vector bundle pair $(TN,T\hb N).$ On the other hand, the canonical collar on $(\gamma_{2n+1},\gamma_{2n})$ induces a collar on $(\theta^*\gamma_{2n+1},(\theta^\partial)^*\gamma_{2n})$ (see discussion below \cref{bun col is pullback}). Given a collared bundle map $\ell^v:(TN|_{\vb N},T\hb N|_{\hvb N})\to (\theta^*\gamma_{2n+1},(\theta^\partial)^*\gamma_{2n})$ and a smooth embedding of pairs $e^v:(\vb N,\hvb N)\hookrightarrow (\{0\}\times \RR^\infty_+,\{0\}\times \partial \RR^\infty_+)$ (and a collared vector bundle map $(T\vb N\oplus \varepsilon^1,T\hb N\oplus \varepsilon^1)\to (TN|_{\vb N},T\hb N|_{\hvb N})$ over the identity), we can see $(\vb N,\ell^v)$ as an object in $\Cobbt$.

\begin{defn}[Moduli space of $\Theta$-structures]\label{moduli space definition}
    Given the definitions above, we denote by $\Mt(N;\ell^v)$ the path components of $\Nt(\vb N)\coloneqq \Cobbt(\emptyset, \vb N)$ of all morphisms whose underlying submanifold pair is abstractly diffeomorphic to $(N,\hb N)$ relative to $\vb N$.
\end{defn}

\begin{rmk}\label{recovering bdiff}
One can check that the space $\Mt(N;\ell^v)$ is weakly equivalent to the balanced product
\[\Emb_{\vb}(N,(-\infty,0]\times \RR^\infty_+)\times_{\Diff_{\vb}(N)} \Buncp(TN,\Theta^*\gamma_{2n+1})\]of the actions of $\Diff_{\vb}(N)$ on the space of embeddings of pairs $e:(N,\hb N)\hookrightarrow ((-\infty,0]\times \RR^\infty_+,(-\infty,0]\times \partial \RR^\infty_+)$ extending $e^v$ (see \cref{section of pairs of manifolds}) and on the space $\smash{\Buncp(TN,\Theta^*\gamma_{2n+1})}$ of those collared bundle maps $(TN,T\hb N)\to (\theta^*\gamma_{2n+1},(\theta^\partial)^*\gamma_{2n})$ extending $\ell^v$ (see \cref{collared pairs}). The diffeomorphism group $\Diff_{\vb}(N)$ acts on the latter space by pre-composition with the differential of a diffeomorphism (which is a collared bundle map of pairs, see \cref{section of pairs of manifolds}). We remark that the action of $\Diff_{\vb}(N)$ on $\Emb_{\vb}(N,(-\infty,0]\times \RR^\infty_+)$ is free and admits local cross-sections (by an analogous argument to \cite[Thm. B]{Palais1960}, see also \cite[Theorem A.1]{steimle}). Moreover, it is contractible (see e.g. \cite[Thm 2.7]{genauer} for a similar argument) and thus $\Mt(N;\ell^v)$ is a model for the homotopy orbits \(\smash{(\Buncp(TN,\Theta^*\gamma_{2n+1}))_{\Diff_{\vb}(N)} }\). In particular, for $\Theta=\id_{\BO(2n+1)}$, then this is a model for the classifying space $\BDiff_{\vb}(N)$, by \cref{gamma is universal}. 
\end{rmk} 

By picking path components in $\hvb N$ (and assuming $\hvb N\neq \emptyset$), we have a morphism $(\prescript{}{\vb N}{H},\ell^{\std}):\vb N\leadsto (\vb N)'$ (see the definition above \cref{stability for V_1}) in $\Cobbt$ given by taking a triad connected sum of $\vb N\times [0,1]$ with $V_1$ at $\hb N\times \{1\}$, where $\ell^{\std}$ is chosen by requiring that the restriction to $V_1$ is standard (see \cref{standard structure on V1} for the definition). To make sense of the concept of standard $\Theta$-structure, we need a basepoint $\Theta$-structure (see \cref{the assumptions for stable stability}). In this context, we fix such a basepoint $\Theta$-structure by requiring that the restriction of $\ell^v$ to a disc in the chosen component of $\hvb N$ is homotopic to it. With this choice, a $\Theta$-structure $\ell^{\std}$ with the above requirement exists by \cref{the space of standard is connected}. By post-composing with this morphism, we have a map \[\smash{\Mt(N;\ell^v)\to \Mt(N\cup_{\vb N}(\prescript{}{\vb N}{H});\ell^v_1)},\]where $\ell^v_1\coloneqq \ell^{\std}|_{(\vb N)'} $. By letting $N_0\coloneqq N$ and $N_g$ be defined by the gluing $\smash{N_{g-1}\cup_{\vb N_{g-1}}(\prescript{}{\vb N_{g-1}}{H})}$ and the corresponding $\Theta$-structure on $\vb N_g$ be $\ell^v_g$, we can define the following space
\[\smash{\Mt(N_\infty;\ell^v):=\hocolim\left(\cdots \to \Mt(N_g;\ell^v_g)\to \Mt(N_{g+1};\ell^v_{g+1}) \to \cdots \right)} . \]Our main result describes the homology of the space above, assuming the pair $(N,\hb N)$ is $1$-connected. Our first goal is to state our main result. To do so, we start by introducing the necessary concepts. After that, we proceed by assembling the results of the previous sections to prove this result. We finish this section by giving some examples of computations.

\subsection{The main theorem.}

We will need the following definitions before we can state our main result. It will be convenient to fix a model structure on the category of pairs of spaces. Consider the \textit{projective model structure}\footnote{Here, we mean the Reedy model structure on the category of functors from the poset $\mathcal{C}\coloneqq \{0<1\}$ to the category of spaces, where we take $0$ to have degree $0$ and $1$ to have degree $1$ and $\smash{\overrightarrow{\mathcal{C}}=\mathcal{C}}$ and $\smash{\overleftarrow{\mathcal{C}}}$ to be the identity morphisms. See \cite[Thm. 15.3.4]{Hirschhorn2003-du}.} on the category of pairs of compactly generated spaces generated by the Quillen model structure on the category of spaces. Weak equivalences and fibrations in the category of pairs are taken objectwise and cofibrations are those maps $(A,A')\to (B,B')$ such that $A'\to B'$ and $A\cup_{A'}B'\to B$ are cofibrations, where here $A\cup_{A'}B'$ means the strict pushout. 

\begin{defn}\label{coconneted and moore postnikov}\label{def moore}
    Let $g:(Y,Y')\to (Z,Z')$ be a map of pairs and $k\geq 0$ an integer, we say $g$ is \textit{strongly $k$-coconnected} if $g|_{Y'}$ is $k$-coconnected (\textit{i.e.} all its homotopy fibers are $(k-1)$-truncated) and if all its total homotopy fibers $\tohofib(f)\coloneqq \hofib(Y'\to Y\times_Z Z')$ are $(k-1)$-truncated. A factorization $(X,X')\to (Y,Y')\to (Z,Z')$ of $g\circ f$ is called a \textit{Moore-Postnikov $k$-factorization} if $f$ is a strongly $k$-connected cofibration and $g$ is a strongly $k$-coconnected fibration.  
\end{defn} 

\begin{lemma}\label{mp exist and are unique}
    Let $f:(X,X')\to (Z,Z')$ be a map of pairs and $k\geq 0$ an integer, then there exists a Moore-Postnikov $k$-factorization $u\circ l:(X,X')\to (Y,Y')\to (Z,Z')$.
\end{lemma}
\begin{proof}
    Consider a Moore-Postnikov $k$-factorization $u'\circ \ell':X'\to Y'\to Z'$ of $f|_{X'}.$ We can assume $\ell'$ is a cofibration and $u'$ a fibration. We now define $Y$ to fit in a Moore-Postnikov $(k+1)$-factorization of $X\cup_{X'} Y'\to Y\to Z.$ Once again, we assume that the leftmost map is a cofibration and the right map a fibration. This produces a factorization $u\circ l:(X,X')\to (Y,Y')\to (Z,Z').$ By construction, $\ell$ is both a cofibration of pairs and strongly $k$-connected. One checks that $u$ is strongly $k$-coconnected since $Y'\to Z'$ and $Y\to Z$ are $k$ and $(k+1)$-coconnected, respectively. On the other hand, they are also fibrations. This finishes the construction.
    \end{proof}

One can check that these factorizations are unique up to homotopy (either directly or as consequence of \cref{type of pairs of maps inj}). Given a cofibration $\iota:(X,X')\to (Y,Y')$ and a fibration $u:(Y,Y')\to (Z,Z')$, we denote by $\Aut_X(u)$ the space of self weak equivalences of $(Y,Y')$ under $\iota$ and over $u$. Composition makes this space into a group-like topological monoid\footnote{This monoid is group-like since any such weak equivalence is a homotopy equivalence. This follows from the general fact that weak equivalences between fibrant-cofibrant objects in a model category are homotopy equivalences (see \cite[Thm. 7.5.10]{Hirschhorn2003-du}), and the fact that $(Y,Y')$ is fibrant-cofibrant in the category of pairs under $(X,X')$ and over $(Z,Z')$ (see \cite[Thm. 7.6.5.(3)]{Hirschhorn2003-du})}. 

Let $n\geq 3$ and $(N,\hb N,\vb N)$ be a $(2n+1)$-dimensional manifold triad and $\Theta:(B,B^\partial)\to (\BO(2n+1),\BO(2n))$ be a map of pairs. Given a collared bundle map $\ell:(TN,T\hb N)\to (\theta^*\gamma_{2n+1},(\theta^\partial)^*\gamma_{2n})$ and letting $\ell=u\circ \ell':(N,\hb N)\to (B',B'^\partial)\to (B,B^\partial)$ be a Moore-Postnikov $n$-factorization of the underlying map of spaces, then the monoid $\Aut_{\vb N}(u)$ acts on the collared bundle $(\theta'^*\gamma_{2n+1},(\theta'^\partial)^*\gamma_{2n})$ via collared bundle maps, where $\Theta'=(\theta',\theta'^\partial)$ is the composite $\Theta\circ u.$ 

Recall from the introduction that $\MTH'$ denotes the cofiber of the canonical map of Thom spectra $\Sigma^{-1}\MT\theta'^\partial\simeq \MT(\theta\oplus \varepsilon^1)\to \MT\theta'$ induced by the collar. Then, by the discussion above, the monoid $\Aut_{\vb N}(u)$ acts on the spectrum $\MTH'$. Finally, denote by $\Mt(N;\ell^v)_{\ell}$ the path component of $\ell$ in $\Mt(N;\ell^v)$, where $\ell^v\coloneqq \ell|_{\vb N}$. If $\hvb N\neq \emptyset$ and after picking a path component of $\hvb N$, denote by $\Mt(N_\infty;\ell^v)_{\ell}$ the homotopy colimit of the maps $\smash{\Mt(N_g;\ell^v_g)_{\ell_g}\to \Mt(N_{g+1};\ell^v_{g+1})_{\ell_{g+1}}}$, where $\ell_g$ is the gluing of the map $\ell$ with the maps $\ell^{\std}_{g'}$ on $\prescript{}{\vb N_{g'}}{H}$ for all $g'\leq g$. We are now ready to state our main theorem.

\begin{mainteoteo}\label{main}
    Let $n\geq 3$. Let $\Theta:(B,B^{\partial})\to (\BO(2n+1),\BO(2n))$ be a map of pairs, $(N,\hb N,\vb N)$ be a manifold triad such that $N$ is connected, $(N,\hb N)$ is $1$-connected and $\hvb N$ is non-empty, and $\ell:(TN,T\hb N)\to (\theta^*\gamma_{2n+1},(\theta^\partial)^*\gamma_{2n})$ be a $\Theta$-structure. Let $\ell=u\circ \ell':(N,\hb N)\to (B',(B')^\partial)\to (B,B^\partial)$ be a Moore-Postnikov $n$-factorization, then there exists a map
    \[\Mt(N_\infty;\ell^v)_{\ell} \to \left(\Omega^{\infty}\MTH'\right)_{\Aut_{\vb N}(u)} \]which is acyclic onto the path component it hits, where $\Theta'=\Theta\circ u$.
\end{mainteoteo}

By specializing this result to $\Theta=\id_{\BO(2n+1)}$ and by \cref{recovering bdiff}, we conclude \cref{main no tang}. The rest of this section is dedicated to the proof of this statement and presenting some examples of potential interest.

\begin{rmk}
    In the statement above, the map $\smash{\Mt(N_\infty;\ell^v)_{\ell} \to \left(\Omega^{\infty}\MTH'\right)_{\Aut_{\vb N}(u)}}$ is more precisely a morphism in the category $\Top[w.e.^{-1}]$, that is, the localization of the category of spaces at weak equivalences. The property of a map being acyclic is invariant under weak equivalence, so it is a well-defined property for morphisms in $\Top[w.e.^{-1}]$. More concretely, we construct a span of spaces $\smash{\Mt(N_\infty;\ell^v)_{\ell} \leftarrow X\to \left(\Omega^{\infty}\MTH'\right)_{\Aut_{\vb N}(u)}}$ where the leftmost map is a weak equivalence and the rightmost map is acyclic.
\end{rmk}

\subsection{\texorpdfstring{Proof of \cref{main}.}{Proof of the main theorem.}}

In this subsection, we assemble all the results from previous section into the proof of \cref{main}. This result will be a consequence of \cref{final with null} by specializing to a chosen $\Theta'$-end and by removing the condition of strong $n$-connectivity for the $\Theta$-structure. Before we proceed, we fix the standing assumptions hidden behind \cref{main}. 

\begin{ass}\label{ass for main}
    As mentioned before, the data above \cref{main} allows us to see $(\vb N,\ell^v)$ as an object of $\Cobbt$. Moreover, the collection of morphisms $\{((\prescript{}{\vb N}{H})_g,\ell^{\std}_g)\}_g$ is a $\Theta$-end in the sense of \cref{theta end without l}. We fix now further choices necessary to define the map in \cref{main}:
    
    \begin{enumerate}[label=$(\mathrm{\Roman*})$]
        \item A \textit{basepoint} $\Theta'$-structure on $(\RR^{2n+1}_+,\partial \RR^{2n+1}_+)$ such that after composing with $u$ it is homotopic to the basepoint $\Theta$-structure; this is possible by taking the pullback of $\ell'$ along inclusion of the fixed disc in the chosen component of $\hvb N$;
        \item\label{theta' struc on h} Choices of $\Theta'$-structures $\smash{\ell'^{\std}_g}$ on $\smash{(\prescript{}{\vb N}{H})_g}$ such that once restricted to $V_1$, they are standard (with respect to the basepoint $\smash{\Theta'}$-structure); this implies that $\smash{u\circ \ell'^{\std}_g}$ is a standard $\Theta$-structure once restricted to $V_1$ and thus by \cref{the space of standard is connected}, it is homotopic to $\ell^{\std}_g$ on $V_1$ and thus on $\smash{(\prescript{}{\vb N}{H})_g}.$ 
    \end{enumerate}As before, the data above allows us to see $(\vb N;\ell'^v)$ as an object of $\Cob_{\Theta'}^\partial$. The collection $\{((\prescript{}{\vb N}{H})_g,\ell'^{\std}_g)\}_g$ is then a $\Theta'$-end. 
\end{ass}

\subsubsection{Bundle maps.}

In \cref{group comp section}, we restricted ourselves to the study of the moduli spaces of $\Theta$-manifolds where the underlying map of the $\Theta$-structure is strongly $n$-connected. We start by describing the space where this condition is not assumed and relate it to a version of the former. Consider the fixed Moore-Postnikov factorization of maps of pairs \[\smash{\Theta:(B',B'^\partial)\overset{(u,u^\partial)}{\to} (B,B^\partial)\overset{\Theta'}{\to} (\BO(2n+1),\BO(2n))} \]from \cref{main}. Recall that $u$ is a fibration, by \cref{def moore}. We have a collared map of vector bundle pairs $u:(\theta'^*\gamma_{2n+1},(\theta'^\partial)^*\gamma_{2n})\to (\theta^*\gamma_{2n+1},\theta^*\gamma_{2n}).$ Post-composition with $u$ induces a map \[ \Buncpn(TN,\Theta'^*\gamma)\to \Buncp(TN,\Theta^*\gamma),\]where the source is the subspace of $\Bunc(TN,\Theta'^*\gamma)$ of those maps extending $\ell'^v\coloneqq \ell'|_{\vb N}$ whose underlying map of pairs is strongly $n$-connected. Moreover, this map factors through the homotopy orbits of the source by the action of $\Aut_{\vb}(u).$ See \cite[190]{GRWII} for the analogous definition in the closed case. The following claims that this map is a weak equivalence, under some hypothesis on $u.$

\begin{lemma}\label{galois theory for pairs}
    The induced map
    \(\smash{(\Buncpn(TN,\Theta'^*\gamma))_{\Aut_{\vb}(u)}\to \Buncp(TN,\Theta^*\gamma)}\) is a weak equivalence onto the path components which it hits. The same claim holds after replacing $N$ by $N_g$ with the $\Theta'$-structure from \ref{theta' struc on h} for any $g\geq 0$.
\end{lemma}
\begin{proof}
     We start by noticing that this claim is equivalent to proving that the homotopy fibers of this map are either empty or contractible. We start by studying the the map \(\smash{\Buncpn(TN,\Theta'^*\gamma)\to \Buncp(TN,\Theta^*\gamma)}\) and its homotopy fibers. By \cref{bundle map post comp}, this map is a Serre fibration, since $u$ is a fibration and $(\vb N,\hvb N)\to (N,\hb N)$ is a cofibration. Moreover, the (homotopy) fiber over a point $g\in \Buncp(TN,\Theta^*\gamma)$ identifies with the union of those path components of $\smash{\Map_{(\vb N,\hvb N)}^{(B,B^\partial)}((N,\hb N),(B',B'^\partial))}$ of strongly $n$-connected maps, by the same result. If this space is non-empty, we pick a point and call it $f:(N,\hb N)\to (B',B'^\partial)$. Consider the map 
     \begin{equation}\label{the thing we want to show is we}
         \smash{\Map_{(\vb N,\hvb N)}^{(B,B^\partial)}((B',B'^\partial),(B',B'^\partial))}\to \smash{\Map_{(\vb N,\hvb N)}^{(B,B^\partial)}((N,\hb N),(B',B'^\partial))}
     \end{equation}given by pre-composition with $f.$ We prove that this map is a weak equivalence onto the components hit. Restricting maps of pairs to the sources of the pairs induces a map \[\Map_{\hvb N}^{B^\partial}(B'^\partial,B'^\partial)\to \Map_{\hvb N}^{B^\partial}(\hb N,B'^\partial),\]which is a weak equivalence by the proof of \cite[Lemma 9.2]{GRWII}. Here, we use that $\hvb N\to \hb N$ and $\hvb N\to B'^\partial$ are cofibrations, $u^\partial:B'^\partial\to B^\partial$ is an $n$-coconnected fibration. By \cref{pairs to spaces in mapping spaces}, these restriction maps to the sources of the pairs are Serre fibrations. Here, we use that $(\vb N,\hvb N)\to (N,\hb N)$ and $(\vb N,\hvb N)\to (B',B'^\partial)$ are cofibrations and $u$ is a fibration. Now, to show that \eqref{the thing we want to show is we} is a weak equivalence, it suffices to show that the induced map on homotopy fibers of the restriction maps is a weak equivalence. For a map $f':B'^\partial\to B'^\partial$, the induced map on fibers agrees with the pre-composition map
     \[\Map_{\vb N\cup_{\hvb N}B'^\partial}^{B}(B',B')\to \Map_{\vb N\cup_{\hvb N}B'^\partial}^{B}(N\cup_{\hb N}B'^\partial,B')\]induced by $f:N\to B'$ (see \cref{redescription of fiber in mapping spaces}). This map is a weak equivalence, once again by the proof of \cite[Lemma 9.2]{GRWII} (and observing that it does not depend on the fact that the pair $(W,\partial W)$ is a manifold and its boundary). To apply this result, we use that the maps $\vb N\cup_{\hvb N}B'^\partial\to B'$ and $\vb N\cup_{\hvb N}B'^\partial\to N\cup_{\hb N} B'^\partial$ are cofibrations (which is a consequence of $\vb N\cup_{\hvb N} \hb N\to N $ being a cofibration), that $f:N\cup_{\hb N}B'^\partial\to B'$ is $(n+1)$-connected and $B'\to B$ is an $n$-coconnected fibration (which follows from the \cref{def moore}). We conclude that the map \eqref{the thing we want to show is we} is a weak equivalence. Observe now that, since $f$ is strongly $n$-connected, the subspace of the source of \eqref{the thing we want to show is we} consisting of weak equivalences hits the path components of strongly $n$-connected maps. On the other hand, by the uniqueness of Moore-Postnikov factorizations (see the discussion below \cref{mp exist and are unique}), we see that every component of the subspace of those strongly $n$-connected maps is hit by the subspace consisting of weak equivalences. By definition, the latter subspace is $\Aut_{\vb N}(u).$ Thus, we see that \eqref{the thing we want to show is we} restricts to a weak equivalence between $\Aut_{\vb N}(u)$ and the fiber over $g$ of the map \(\smash{\Buncpn(TN,\Theta'^*\gamma)\to \Buncp(TN,\Theta^*\gamma)}.\) By the same reasoning as in \cite[Lemma 9.2]{GRWII}, we conclude that after taking homotopy orbits on the source of the latter map, its homotopy fibers are either empty or contractible. This finishes the proof. The case of $N_g$ follows verbatim by observing that the proof above only used that $(\vb N,\hvb N)\to (N,\hb N)$ is a cofibration.  
\end{proof}

We now apply this result above to the moduli spaces $\Mt(-)$. Note that $\Aut_{\vb N}(u)$ acts on the category $\smash{\Cob^\partial_{\Theta'}}$ via continuous endofunctors by changing the $\Theta'$-structures by post-composition. Since the elements of $\Aut_{\vb N}(u)$ are weak equivalences $(B',B'^\partial)\to (B',B'^\partial)$ under the map $\ell'^v:(\vb N,\hvb N)\to (B',B'^\partial)$, the object $(\vb N,\ell'^v)$ is fixed under this action. This induces an action on $\N^\partial_{\Theta'}(\vb N)$ and thus on $\M^\partial_{\Theta',n}(N;\ell'^v),$ where $\smash{\M^\partial_{\Theta',n}(N;\ell'^v)}$ denotes the path components of $\smash{\Mt(N;\ell'^v)}$ of those $(W,\ell_W)$ where $\ell_W$ is strongly $n$-connected. Similarly, we have a map $\smash{(\M^\partial_{\Theta',n}(N;\ell'^v))_{\Aut_{\vb}(u)}\to \Mt(N;\ell^v)}$ given by post-composition with $u$. We deduce from \cref{galois theory for pairs} that this map is a weak equivalence onto the path components hit.

\begin{cor}\label{the moore postnikov application}
    Then the induced map
    \(\smash{(\M^\partial_{\Theta',n}(N;\ell'^v))_{\Aut_{\vb}(u)}\to \Mt(N;\ell^v)}\) is a weak equivalence onto the path components hit. The same claim holds after replacing $N$ by $N_g$ for any $g\geq 0$.
\end{cor}
\begin{proof}
    This follows from \cref{galois theory for pairs}, the equivalence \(\smash{\Mt(N;l)\simeq (\Buncp(TN,\Theta^*\gamma))_{\Diff_{\vb}(N)} }\) from \cref{recovering bdiff}, and the fact that the actions of $\Diff_{\vb}(N)$ and $\Aut_{\vb}(u)$ commute and thus, one can form the homotopy orbits in any order.
\end{proof}

    Note that the equivalence of \cref{the moore postnikov application} is natural with respect to the stabilization map $\smash{\M^\partial_{\Theta',n}(N_g;\ell'^v_g)}\to \smash{\M^\partial_{\Theta',n}(N_{g+1};\allowdisplaybreaks \ell'^v_{g+1})}$ on the target and the analogous maps after taking homotopy orbits by $\Aut_{\vb N_g}(u) $ and $\Aut_{\vb N_{g+1}}(u),$ respectively. These two monoids are equivalent for every $g\geq 0$, since $\vb N_{g}\cong \vb N_{g+1}$ and under such identification $\ell'^v_g|_{\vb N_g}$ and $\smash{\ell'^v_{g+1}|_{\vb N_{g+1}}}$ of \ref{theta' struc on h} are homotopic, since $\smash{\ell'^{\std}_g}$ and $\smash{\ell'^{\std}_{g+1}}$ are standard (and thus extend to a contractible space). Thus, we have an equivalence after taking homotopy colimits
    \[\left(\M^\partial_{\Theta',n}(N_\infty;\ell'^v) \right)_{\Aut_{\vb}(u)}\to \Mt(N_\infty;\ell^v)\]onto the path components which are hit. In the context of the $\Aut_{\vb N}(u)$-action defined above \cref{the moore postnikov application}, the maps
    \[\M^\partial_{\Theta',n}(N;\ell'^v)\to \N^\partial_{\Theta'}(\vb N)\to \Omega_{[\emptyset,\vb N]}\B\Cob^\partial_{\Theta'} \]are equivariant. We are now ready to prove our main result.

\begin{proof}[Proof of \cref{main}]
    Consider the notation and choices fixed in \cref{the assumptions for stable stability,ass for main}. We can apply \cref{final with null} to the $\Theta'$-end given by $K|_0\coloneqq \vb N$ and $K|_{[i,i+1]}\coloneqq (\prescript{}{\vb N}{H})_i$. Note that by assumption $\partial K|_0=\hvb N\neq \emptyset$, $(B',B'^\partial)$ is $1$-connected and $B'^\partial$ is path-connected, since $\ell'$ is strongly $n$-connected and $n\geq 3$. By restricing to the component of $N$, we obtain that the map $\smash{\M^\partial_{\Theta',n}(N_\infty;\ell'^v)\to \Omega_{[\emptyset,\vb N|_\infty]}\B\Cob^\partial_{\Theta'}}$ is acyclic. Consider now the span
    \[\Mt(N_\infty;\ell^v)\leftarrow  \left(\M^\partial_{\Theta',n}(N_\infty;\ell'^v)\right)_{\Aut_{\vb}(u)}\to \left(\Omega_{[\emptyset,(\vb N,\ell'^v)]}\B\Cob^\partial_{\Theta'}\right)_{\Aut_{\vb}(u)},\]where the leftmost map is an equivalence onto the hit path components by the discussion below \cref{the moore postnikov application}. Since sequential homotopy colimits and homotopy orbits commute and acyclicity is preserved after taking homotopy orbits, the right map is acyclic onto the path components hit. By \cite[Main Thm. 4.5]{genauer}, there exists an equivalence $\B\Cob^\partial_{\Theta'}\simeq \Omega^{\infty-1}\MTH'$, which is a composition of equivariant maps with respect to the actions of $\Aut_{\vb}(u)$ (see the zig-zag in p.538 of loc.cit). This finishes the proof, by taking the component of $\ell$ in the leftmost space and picking a $\Theta'$-nullbordism of $(\vb N,\ell'^v)$ (e.g. $(N,\hb N)$ itself) to identify $\Omega_{[\emptyset,(\vb N,\ell'^v)]}\Omega^{\infty-1}\MTH'$ with $\Omega^\infty\MTH'$.
\end{proof}

\subsection{Some examples.}\label{examples section}

In this subsection, we expand on a few examples and prove the claims made in the introduction. Throughout the entire subsection, we let $(N,\hb N,\vb N)$ be a connected $(2n+1)$-manifold triad such that $(N,\hb N)$ is $1$-connected for $n\geq 3$ and $\hvb N\neq \emptyset$. We start by describing the Thom spectra and automorphism monoid present in \cref{main} in specific cases.

\subsubsection*{Induced tangential structures.}
Let $\theta:B\to \BO(2n+1)$ be a map of spaces. We can define the map of pairs $\iota_*\theta:(B,B\times_{\BO(2n+1)} \BO(2n))\to (\BO(2n+1),\BO(2n)). $ One has that $\MT\iota_*\theta\simeq \Sigma^\infty_+B$. (See \cite[Prop. 3.1]{GMTW} for a proof when $\theta=\id_{\BO(2n+1)}$ and \cite[Prop. 3.0.7]{gollinger} in general.) Given a map of pairs $\Theta:(B,B^\partial)\to (\BO(2n+1),\BO(2n)),$ we denote by $\iota^*\Theta$ the underlying map $B\to \BO(2n+1).$ We are interested in knowing when the structure $\Theta_N:(B,B^\partial)\to (\BO(2n+1),\BO(2n))$ given by the strong Moore-Postnikov $n$-factorization (see \cref{coconneted and moore postnikov}) of $\tau_N:(N,\hb N)\to (\BO(2n+1),\BO(2n))$ is equivalent to $\iota_*\iota^*\Theta_N$.

\begin{lemma}\label{when is mp induced}
    In the context of the above, the map $\Theta_N$ is equivalent to $\iota_*\iota^*\Theta_N$ if and only if $(N,\hb N)$ is $n$-connected. More generally, let $\theta:B\to \BO(2n+1)$ be a map and $l$ a $\iota_*\theta$-structure on $(N,\hb N).$ Let $l=u\circ l':(N,\hb N)\to (B',(B')^\partial)\to (B, B\times_{\BO(2n+1)} \BO(2n))$ be a Moore Postnikov $n$-factorization of pairs, then $\Theta_l\coloneqq \Theta\circ u$ is equivalent to $\iota_*\iota^*\Theta_l$ if and only if $(W,\hb W)$ is $n$-connected.
\end{lemma}
\begin{proof}
    We prove the second statement, as the first one follows by applying the second to $\theta=\id_{\BO(2n+1)}.$ We start with the "only if"-direction. Suppose $\Theta_l:(B',(B')^\partial)\to (\BO(2n+1),\BO(2n))$ is equivalent to $\iota_*\iota^*\Theta_l,$ thus it follows that $(B',(B')^\partial)$ is $2n$-connected, as it is a pullback of a $2n$-connected map $\BO(2n)\to \BO(2n+1).$ Thus, we see that $\pi_i(W,\hb W)\cong \pi_i(B',(B')^\partial)$ for $i\leq n$ by \cref{1 conn: old same as new}, since $(W,\hb W)\to (B',(B')^\partial)$ is strongly $n$-connected and $(W,\hb W)$ is $1$-connected. For the "if"-direction, consider the factorization $(B')^\partial\to B'\times_{\BO(2n+1)} \BO(2n)\to B'$. The first map is $n$-coconnected by definition, and the second map is $2n$-connected, since it is the pullback of an $2n$-connected map. On the other hand, the pair $(B',(B')^\partial)$ is $n$-connected, since $(W,\hb W)$ is. Since $n\geq 1,$ we have that $(B')^\partial\to B'\times_{\BO(2n+1)}\BO(2n)$ is also $n$-connected, again by \cref{1 conn: old same as new}. Thus, it is an equivalence. This finishes the proof.
\end{proof}

\subsubsection*{Contractible automorphisms.} In certain cases, the grouplike monoid $\Aut_{\vb N}(u)$ is contractible. This eliminates substantial complexity in computing the (co)homology of the stable moduli space $\Mt(N_\infty)$ via \cref{main}. The following is an application of \cref{type of pairs of maps inj}.  

\begin{lemma}\label{when aut is contractible}
    Let $\Theta:(B,B^\partial)\to (\BO(2n+1),\BO(2n))$ and $l_N$ a $\Theta$-structure on $(N,\hb N)$. Denote a Moore-Postnikov $n$-factorization of pairs $u\circ l'_N:(N,\hb N)\to (B',(B')^\partial)\to (B,B^\partial) $ of $l_N$. Assume that $(N,\vb N)$ and $(\hb N, \hvb N)$ are $(n-1)$-connected, then $\Aut_{\vb N}(u)$ is contractible. Moreover, if $(B')^{\partial}\to B'\times_{B} B^{\partial}$ is an equivalence, the condition on $(\hb N, \hvb N)$ can be removed.
\end{lemma}

\subsubsection*{Moduli spaces of $h$-cobordisms stabilized.} We call a triad $(N,\hb N,\vb N)$ an \textit{h-cobordism} if the inclusions $\hb N\to N$ and $\vb N\to N$ are equivalences. Let $\theta: B\to \BO(2n+1)$ be a map and $l_N$ be a $\iota_*\theta$-structure on $N.$

\begin{cor}\label{no tang h cobs}
    For an $h$-cobordism $(N,\hb N,\vb N)$, there exists a map
    \[\mathscr{M}^\partial_{\iota_*\theta}(N_\infty;l^v_N)_{l_N} \to \Omega^\infty \Sigma^\infty_+ B'\]which is acyclic onto the path component it hits, where $N\to B'\to B$ is the classical Moore-Postnikov $n$-factorization of the map $l_N:N\to B$.
\end{cor}
\begin{proof}
    This follows by combining \cref{main}, \cref{when is mp induced}, \cref{when aut is contractible} and the fact that the underlying map of the strong Moore-Postnikov factorization of this map of pairs in this case is a classical Moore-Postnikov factorization of $\tau_N:N\to B$. We leave this check to the reader.
\end{proof}

\begin{rmk}
    By taking $B$ $n$-connected and $W$ to be the triad $(D^{n+1}_+,\partial_0D^{n+1}_+,\partial_1D^{n+1}_+) $ (as defined in \cref{triad homotopy}), this result recovers \cite[Thm. A*]{BP} and extends it to dimension $7$.
\end{rmk}

\subsubsection*{A rational computation.} Let $A$ be $2n$-dimensional manifold with non-empty boundary. Consider the triad $(A\times [0,1], A\times \{0\},\partial A \times [0,1]\cup A \times \{1\}) $, and observe that it satisfies the conditions of \cref{main no tang} (or \cref{main} for $\Theta=\id_{\BO(2n+1)})$. The space $\BDiff_{\vb}(A\times [0,1]\natural V_\infty)$ can be viewed as a space of $V_g$-stabilized concordances of the manifold $A$. We compute its rational cohomology as follows.

\begin{cor}\label{rational computation}
    Assume $A$ is aspherical and parallelizable. There exists an isomorphism
    \[\H^*(\BDiff_{\vb}(A\times [0,1]\natural V_\infty);\QQ)\cong \QQ[\H^*(A;\QQ)\otimes \QQ[p_{[(n+1)/4]},\cdots, p_n]] \]as graded commutative $\QQ$-algebras.
\end{cor}
\begin{proof}
    This follows by applying \cref{no tang h cobs} to the trivial $h$-cobordism $A\times [0,1]$ along with the following classical facts: the rational cohomology $\QQ$-algebra $\H^*(\Omega^\infty_0\Sigma^\infty_+X)$ is isomorphic to $\QQ[\widetilde{\H}^*(X;\QQ)]$ of a space $X$ with finite dimensional $\QQ$-homology in each degree; since $A$ is parallelizable and aspherical, the composition $A\to A\times \tau_{>n}\BO(2n+1)\to \BO(2n+1)$ is a Moore-Postnikov $n$-factorization, where $\tau_{>n}\BO(2n+1)$ is the $n$-connected cover of $\BO(2n+1)$; finally the $\QQ$-algebra $\H^*(\tau_{>n}\BO(2n+1);\QQ)$ is the free graded commutative $\QQ$-algebra on the Pontrjagin classes $p_{[(n+1)/4]},\cdots, p_n,$ where $p_i$ lives in degree $4i.$
\end{proof}

\subsubsection*{The initial structure.}

We finish by considering the \textit{initial} tangential structure on $N.$ We consider its tangent classifies $\tau_N:(N,\hb N)\to (\BO(2n+1),\BO(2n))$. The identity map of $N$ induces a $\tau_N$-structure on $N.$ In this case, taking the Moore-Postnikov factorization of $\id: (N,\hb N)\to (N,\hb N)$ gives the identity of $N.$ Thus, $\MTH$ is the cofiber of the map $\Sigma^{-1}\Th(-T\hb N)\to \Th(-TN),$ which is equivalent to the Spanier-Whitehead dual $\text{D}(N/\vb N)$ by Atiyah duality. By \cref{type of pairs of maps inj} for $k=l=0,$ we see that $\Aut_{\vb N}(\id)$ is contractible. Thus, we obtain the following result.

\begin{cor}
    There exists a map
    \[\mathscr{M}^\partial_{\tau_N}(N_\infty;\mathrm{id}|_{\vb N})_{\id} \to \Omega^\infty \textD(N/\vb N)\]which is acyclic onto the path component it hits.
\end{cor}

\appendix

\section{Appendix: Mapping spaces of pairs.}\label{mappin appendix}

\renewcommand{\theteo}{\thesection.\arabic{teo}}
\setcounter{teo}{0}
\makeatletter
\@addtoreset{teo}{section}
\makeatother

\renewcommand{\theprop}{\theteo}
\renewcommand{\thelemma}{\theteo}
\renewcommand{\thecor}{\theteo}

Recall from the discussion above \cref{def moore} that we consider the category of pairs of spaces $\Top^{[1]}$ as a model category with the projective model structure induced by the Quillen model structure on $\Top$. Given pairs $(A,A')$ and $(B,B')$, one can define the category of pairs of spaces under $(A,A')$ and over $(B,B').$ We consider this category with the under-over model structure (i.e. the model structure where maps are fibrations, cofibrations, or weak equivalences if and only if the underlying map in $\Top^{[1]}$ is, see \cite[Thm 7.6.5]{Hirschhorn2003-du}). This category is enriched in spaces by taking the space of maps \[\Map_{(A,A')}^{(B,B')}((X,X'),(Y,Y'))\] with the compact-open topology. We can also see this category as enriched in simplicial sets by taking $\text{Sing}_\bullet.$

\begin{prop}\label{pairs to spaces in mapping spaces}
    Let $(A,A')\to (X,X')$ be a cofibration of pairs, $(Y,Y')\to (B,B')$ a fibration of pairs, and $f:(A,A')\to (Y,Y')$ and $g:(X,X')\to (Y,Y')$ be maps of pairs. Then the restriction map
    \[\Map_{(A,A')}^{(B,B')}((X,X'),(Y,Y'))\to \Map_{A'}^{B'}(X',Y')\]is a fibration. Moreover, the homotopy fiber over a map $\alpha:X'\to Y'$ is equivalent to $\Map_{A\cup_{A'}X'}^B(X,Y),$ where $Y$ is seen as under $A\cup_{A'} X'$ using $g\circ f$ and $\alpha$.
\end{prop}
\begin{proof}
    Let $I\to I'$ be a trivial cofibration in $\Top$ and consider a diagram of the following form
    \[\begin{tikzcd}
        I'\arrow[r]\arrow[d] & \Map_{(A,A')}^{(B,B')}((X,X'),(Y,Y'))\arrow[d] \\
        I \arrow[r]\arrow[ur, dashed] & \Map_{A'}^{B'}(X',Y').
    \end{tikzcd}\]The data of such a lift is equivalent to the data of a map $I\times (X,X')\to (Y,Y')$ over $(B,B')$ such that its restriction to $I\times X'\to Y'$ is the bottom horizontal map, the restriction to $I'\times (X,X')\to (Y,Y')$ is the top horizontal map and the restriction to $I\times (A,A')\to (Y,Y')$ is constant at $f.$ In other words, we have to find a lift of the following diagram
    \[\begin{tikzcd}
        (I'\times X)\cup_{I'\times (A\cup_{A'} X')} I\times (A\cup_{A'} X') \arrow[d]\arrow[r] & Y\arrow[d] \\
        I\times X \arrow[r] \arrow[ur, dashed] & B   
    \end{tikzcd}\]which is possible if the left vertical map is a trivial cofibration. This map is the pushout product of $I'\to I$ and $A\cup_{A'}X'\to X$ (see \cite[Defn.4.2.1]{hovey2007model}). Since the category of spaces is a monoidal model category, by \cite[Defn. 4.2.6, Prop. 4.2.11]{hovey2007model}\footnote{Note that in this reference, our definition of spaces agrees with the definition of $k$-spaces, while compactly generated spaces are assumed to weakly Hausdorff.}, the pushout product of two cofibrations is a cofibration and trivial if one of them is. By assumption, $I'\to I$ is a trivial cofibration and $A\cup_{A'}X'\to X$ is a cofibration, thus the left vertical map above is a trivial cofibration. This finishes the proof.
\end{proof}

\begin{rmk}\label{redescription of fiber in mapping spaces} 
By \cref{pairs to spaces in mapping spaces}, we have a fiber sequence of the form
\[\Map_{A\cup_{A'}X'}^B(X,Y)\to \Map_{(A,A')}^{(B,B')}((X,X'),(Y,Y'))\to \Map_{A'}^{B'}(X',Y').\]It is convenient to have a variant of this fiber sequence where the middle space is the subspace of those maps $(X,X')\to (Y,Y')$ which are strongly $k$-connected (recall from \cref{strong conn defn}). To do so, observe that the map \(\smash{\Map_{A\cup_{A'}X'}^B(X,Y) \to \Map_{A\cup_{A'}Y'}^B(X\cup_{X'}Y',Y)}\) is a homeomorphism. Thus, we have fiber sequence analogous to the one above, where the base is the subspace of those maps $X'\to Y'$ which are $k$-connected, the total space is the subspace of those maps $(X,X')\to (Y,Y')$ which are strongly $k$-connected and the fiber is the subspace of $\smash{\Map_{A\cup_{A'}Y'}^B(X\cup_{X'}Y',Y)}$ of those maps which are $(k+1)$-connected.

\end{rmk}

We can describe the homotopy type of these mapping spaces given connectivity and coconnectivity of the pairs involved. Recall the definition of strong connectivity and coconnectivity from \cref{strong conn defn,def moore}.

\begin{prop}\label{type of pairs of maps inj}
    Let $k,l\geq 0$ be integers. Let $(A,A')\to (X,X')$ be a cofibration of pairs and $(Y,Y')\to (B,B')$ be a strongly $l$-coconnected fibration of pairs. If both $A\to X$ and $A'\to X'$ are $k$-connected, then the space 
    \[\Map_{(A,A')}^{(B,B')}((X,X'),(Y,Y'))\]is a $(l-k-2)$-type. In particular, it is contractible if $k\geq l-1,$ provided it is non-empty. Furthermore, if $Y'\to Y\times_B B'$ is an equivalence, then we can remove the $k$-connectivity assumption of $A'\to X'.$
\end{prop}
\begin{proof}
    We start by noticing that, under our assumptions, the space $\smash{\Map_{(A,A')}^{(B,B')}((X,X'),(Y,Y'))}$ is a model for the derived mapping space in the category of pairs under $(A,A')$ and over $(B,B')$, since $(X,X')$ is cofibrant and $(Y,Y')$ is fibrant in this category. However, the derived mapping space of a model category only depends on the class of weak equivalences (see \cite[Corollary 4.7]{Dwyerkan} together with the fact that this category is a closed simplicial model category in the sense of \cite[Defn. 2.2.2]{Quillen1967}). We define the following model structure with the same class of weak equivalences on $\Top^{[1]}.$ The \textit{injective model structure} on $\Top^{[1]}$ is model structure where the weak equivalences and cofibrations are objectwise and the fibrations are those maps $(X,X')\to (Y,Y')$ such that $X\to Y$ and $X'\to X\times_YY'$ are fibrations of spaces. Similarly to before, this induces a model structure on the over-under category. Since the injective model structure has the same weak equivalences, our space of interest is equivalent to the derived mapping space in this category. The latter is equivalent to the same mapping space but replacing $(Y,Y')\to (B,B')$ by an injective fibration, as this gives a fibrant replacement of the target. For simplicity of notation, we denote this replacement also by $(Y,Y').$ Let $Z$ be the mapping space $\smash{\Map_{(A,A')}^{(B,B')}((X,X'),(Y,Y'))}$. Let $i\geq 0$ and given a map $\alpha:S^i\to Z$ be a map, then $\alpha$ can be extended to $D^{i+1}$ if and only if the following lifting problem can be solved
    \[\begin{tikzcd}
        (D^{i+1}\times A \cup S^i\times X,D^{i+1}\times A'\cup S^i\times X')\arrow[d]\arrow[r] & (Y,Y')\arrow[d] \\
        (D^{i+1}\times X,D^{i+1}\times X')\arrow[r]\arrow[ur, dashed] &(B,B').
    \end{tikzcd}\]From the assumptions, we conclude the following facts:
    \begin{enumerate}[label=(\roman*)]
        \item the maps $D^{i+1}\times A \cup S^i\times X\to D^{i+1}\times X$ and $D^{i+1}\times A'\cup S^i\times X'\to D^{i+1}\times X'$ are $(k+i+1)$-connected cofibrations. This can be seen by induction on cells on the $k$-connected pairs $(X,A)$ and $(X',A')$.
        \item the maps $Y\to B$ and $Y'\to Y\times_B B'$ are $l$-coconnected fibrations. 
    \end{enumerate}By obstruction theory, we have no obstructions to solving the leftmost lifting problem
    \[\begin{tikzcd}
        D^{i+1}\times A \cup S^i\times X \arrow[d]\arrow[r] & Y\arrow[d] \\
        D^{i+1}\times X \arrow[r]\arrow[ur, dashed] & B
    \end{tikzcd} \quad \begin{tikzcd}
        D^{i+1}\times A'\cup S^i\times X'\arrow[d]\arrow[r] & Y'\arrow[d] \\
        D^{i+1}\times X'\arrow[r]\arrow[ur, dashed] & Y\times_B B'
    \end{tikzcd}\]if $k+i+1\geq l,$ that is, when $i\geq l-k-1.$ Provided $i\geq l-k-1$ and fixing a lift, we can solve the initial lifting problem if we can solve the right problem above where the bottom horizontal map is the unique map to the pullback induced by the map $D^{i+1}\times X'\to B'$ and the chosen lift $D^{i+1}\times X'\to D^{i+1}\times X\to Y'$. The same obstruction theoretic argument implies that such a problem can be solved provided $i\geq l-k-1.$ This implies that $Z$ is a $(l-k-2)$-type. If the map $Y'\to Y\times_B B'$ is an equivalence, then no conditions are needed on the left vertical map of the third lifting problem to solve it. Thus, no condition on the cofibration $A'\to X'$ is necessary to deduce that $Z$ is an $(l-k-2)$-type in this case.
\end{proof}

\begin{prop}\label{post composition by fibration}
    Let $(A,A')\to (X,X')$ be a cofibration of pairs, $(Y,Y')\to (Z,Z')\to (B,B')$ two composable fibrations of pairs, and $(A,A')\to (Y,Y')$ and $(X,X')\to (Y,Y')$ be maps of pairs. Then the post-composition map
    \[\Map_{(A,A')}^{(B,B')}((X,X'),(Y,Y'))\to \Map_{(A,A')}^{(B,B')}((X,X'),(Z,Z'))\]is a fibration.
\end{prop}
\begin{proof}
    The claim follows if we show that the simplicially enriched category of pairs of spaces under $(A,A')$ and over $(B,B')$ is a simplicial model category (in the sense of \cite[Defn. 9.1.6]{Hirschhorn2003-du}) when taken with the projective model structure, as this claim is a special case of axiom M7. This follows by combining the general fact that the category of diagrams indexed in a Reedy category of a simplicial model category is a simplicial model category, when taken with the Reedy model structure (see \cite[Thm. 15.3.4.(3)]{Hirschhorn2003-du}), and the fact that over-under categories of a simplicial model categories are simplicial model categories.
\end{proof}

We establish now similar properties for spaces of collared bundle maps.  

\begin{prop}\label{bundle map fibration}
    Let $\iota:(A,A')\to (X,X')$ be a cofibration of pairs, $(\xi,\xi')$ a collared vector bundle over $(X,X')$, $(\eta,\eta')$ a collared vector bundle over $(Y,Y')$ and $f: (\iota^*\xi,\iota^*\xi')\to (\eta,\eta')$ be a collared bundle map. The restriction map 
    \[\smash{\Bunca(\xi,\eta)\to \Bun_{A'}(\xi',\eta')}\]is a fibration. Moreover, the homotopy fiber over $\alpha:\xi'\to \eta'$ is equivalent to $\Bun_{A\cup_{A'} X'}(\xi, \eta).$
\end{prop}
\begin{proof}
    This statement will follow from the following claim, whose proof we leave to the reader: the restriction map $\Bun(\xi,\eta)\to \Bun(\xi|_{X'},\eta)$ is a fibration, if $X'\to X$ is a cofibration. We show now our desired statement using this claim. From the pullback decomposition of \cref{bun col is pullback}, we observe that it suffices to show that $\Bun_A(\xi,\eta)\to \Bun_{A'}(\xi|_{X'},\eta)$ is a fibration, as $\xi|_{X'}\cong \xi'\oplus \varepsilon^1$. For that, consider the following commutative diagram
    \[\begin{tikzcd}
        \Bun_A(\xi,\eta)\arrow[d]\arrow[r] & \Bun(\xi,\eta)\arrow[d] \arrow[r] & \Bun(\xi|_A,\eta)\arrow[d] \\
        \Bun_{A'}(\xi|_{X'},\eta)\arrow[r] & \Bun(\xi|_{X'},\eta)\arrow[r] & \Bun(\xi|_{A'},\eta)
    \end{tikzcd}\]where the rows are fibration sequences, by the claim above as $A\to X$ and $A'\to X'$ are cofibrations. One can see that the left vertical map is a fibration if the map 
    \[\Bun(\xi,\eta)\to \Bun(\xi|_{X'},\eta)\times_{\Bun(\xi|_{A'},\eta)}\Bun(\xi|_A,\eta)\]induced by the commutativity of the left square, where $\times $ denotes the strict pullback, is a fibration. However, this map is homeomorphic to the restriction map $\Bun(\xi,\eta)\to \Bun(\xi|_{X'\cup_{A'} A},\eta)$, which is fibration by the claim above using that $X'\cup_{A'} A\to X$ is a cofibration. This finishes the proof of the first claim. The claim about the homotopy fiber follows by observing that $\Bun_{A\cup_{A'} X'}(\xi, \eta)$ is the strict fiber of the restriction map $\Bun_A(\xi,\eta)\to \Bun_{A'}(\xi|_{X'},\eta).$
\end{proof}

\begin{rmk} \label{redescription fiber bundle maps}

Assume that the dimensions of $\xi$ and $\eta$ agree. Similarly to \ref{redescription of fiber in mapping spaces}, the (homotopy) fiber of the map above $\Bun_{A\cup_{A'}X'}(\xi,\eta)$ is homeomorphic to $\Bun_{A\cup_{A'}Y'}(\xi\cup \eta'\oplus_{\xi'\oplus \varepsilon^1} \varepsilon^1, \eta).$ On the other hand, by forgetting the bundle map, one obtains a map of fiber sequences between the one describes in \cref{bundle map fibration} and the one in \cref{pairs to spaces in mapping spaces}.\end{rmk}

\begin{prop}\label{bundle map post comp}
    Let $\iota:(A,A')\to (X,X')$ be a cofibration of pairs, $(\xi,\xi')$ a collared vector bundle over $(X,X')$, $(\eta_Y,\eta'_Y)$ a collared vector bundle over $(Y,Y')$, $(\eta_Z,\eta'_Z)$ a collared vector bundle over $(Z,Z')$ and $f: (\iota^*\xi,\iota^*\xi')\to (\eta_Y,\eta'_Y)$ be a collared bundle map. Given a collared bundle map $u:(\eta_Y,\eta_Y')\to (\eta_Z,\eta_Z')$ whose underlying map is a fibration, the map induced by post-composition with $u$ 
    \[\smash{\Bunca(\xi,\eta_Y)\to \Bunca(\xi,\eta_Z)}\]is a fibration. Moreover, the homotopy fiber over a map $g$ is equivalent to $\Map_{(A,A')}^{(Z,Z')}((X,X'),(Y,Y')),$ where $(X,X')$ is seen over $(Z,Z')$ using $g.$
    \end{prop}
\begin{proof}
    We start by showing that the map above is a Serre fibration. Let $I'\to I$ be a trivial cofibration in $\Top$ and consider a diagram of the form
    \[\begin{tikzcd}
        I'\arrow[r]\arrow[d] & \Bunca(\xi,\eta_Y)\arrow[d] \\
        I \arrow[r]\arrow[ur, dashed] & \Bunca(\xi,\eta_Z)
    \end{tikzcd}.\]Observe that the analogous lifting problem where collared bundle mapping spaces are replaced by mapping spaces of the underlying pairs can be solved by \cref{post composition by fibration}. In other words, we can find a map $l:I\times (X,X')\to (Y,Y') $ extending the underlying map of the top arrow above, such that after composing with the map $(Y,Y')\to (Z,Z')$, it agrees with the underlying map of the bottom map of the square above. The original lifting problem can be solved if we can find a collared bundle map $I\times(\xi,\xi')\to (l^*\eta_Y,l^*\eta_Y')$ over the identity of $(X,X')$ compatibly with the existing bundle maps given by the square above. To do that, observe that a collared bundle map $I\times (\xi,\xi')\to (l^*\eta_Y,l^*\eta_Y')$ over $I\times X$ is equivalent to a section of the map of pairs $(\Hom(\xi,l^*\eta_Y),\Hom(\xi',l^*\eta_Y'))\to I\times(X,X')$ which is induced by the collars of $\xi$ and $\eta$. Each individual map in this map of pairs is a fiber bundle and thus a fibration. We conclude that the solution of the lifting problem above is reduced to the following lifting problem
    \[\begin{tikzcd}
        I\times (A,A')\cup_{I'\times (A,A')} I'\times (X,X')\arrow[r] \arrow[d] &(\Hom(\xi,l^*\eta_Y),\Hom(\xi',l^*\eta_Y'))  \arrow[d] \\
        I\times (X,X')\arrow[r] \arrow[ru, dashed] &(\Hom(\xi,(u\circ l)^*\eta_Z),\Hom(\xi',(u\circ l)^*\eta_Z'))
    \end{tikzcd}\]where the bottom map is induced by the bottom map of the original lifting problem. Now, this can be solved uniquely since the right vertical map is a homeomorphism. This proves the first claim.

    The proof of the identification of the homotopy fibers follows the same idea of the proof of the first claim. One can identify the fiber of this map with the fiber of the map between the underlying mapping spaces. This fiber identifies precisely with $\Map_{(A,A')}^{(Z,Z')}((X,X'),(Y,Y')).$ We leave this verification to the reader.
\end{proof}

\printbibliography

\textsc{Department of Mathematics, Karlsruhe Institute of Technology, 76131 Karlsruhe, Germany}

\textit{Email address:} \url{joao.fernandes@kit.edu}

\end{document}